\documentclass[reqno]{amsart}
\usepackage[all]{xy}
\SelectTips{cm}{}

\usepackage{newtxtext}
\usepackage{newtxmath}
\usepackage{bm}

\usepackage{amscd}
\usepackage{xcolor}
\usepackage{cjhebrew}
\usepackage{pdfpages}
\usepackage{hhline}
\usepackage{multirow}
\usepackage{imakeidx}
\usepackage[mathscr]{eucal}
\usepackage{mathdots}
\usepackage{graphicx}

\numberwithin{equation}{section}

\newtheorem{Thm}{Theorem}[section]
\newtheorem{Lem}[Thm]{Lemma}
\newtheorem{Pro}[Thm]{Proposition}
\newtheorem{Cor}[Thm]{Corollary}

\theoremstyle{definition}
\newtheorem{Def}[Thm]{Definition}

\newtheorem{Rem}[Thm]{Remark}
\newtheorem{Table}[Thm]{Table}

\newtheorem*{Rem-intro}{Remark}

\DeclareMathOperator{\Tor}{Tor}
\DeclareMathOperator{\Hom}{Hom}
\DeclareMathOperator{\Image}{Im}

\DeclareMathOperator{\Ext}{Ext}

\DeclareMathOperator{\Prin}{Prin}

\newcommand{\io}{\iota}

\newcommand{\smat}[4]
	{ \left[ \begin{smallmatrix} {#1} & {#2} \\ {#3} & {#4} \end{smallmatrix} \right]  }

\newcommand{\crn}{\rm{ch}}
\newcommand{\pont}{\rm{Ph}}
\newcommand{\pt}{\rm{pt}}

\newcommand\crr{^{\scriptscriptstyle {\it CR}}}
\newcommand\crt{^{\scriptscriptstyle {\it CRT}}}
\newcommand\scrt{_{\scriptscriptstyle {\it CRT}}}

\newcommand{\Dirlim}{\varinjlim}
\newcommand{\Invlim}{\varprojlim}

\newcommand{\ind}{\operatorname{ind}}

\newcommand{\hil}{{\mathcal{H}}}
\newcommand{\Lin}{{\mathcal{L}}}

\newcommand{\Len}{{\widetilde{L}}}

\newcommand{\LL}{{\mathcal{L}}}

\newcommand{\ZZ}{{\mathbb{Z}}}
\newcommand{\QQ}{{\mathbb{Q}}} 

\newcommand{\CC}{{\mathbb{C}}}
\renewcommand\sc{_{\scriptstyle \CC}}
\newcommand\sr{_{\scriptstyle \RR}}
\newcommand\pr{^{\scriptstyle \RR}}
\newcommand\pc{^{\scriptstyle \CC}}
\newcommand\pf{^{\scriptstyle \FF}}
\newcommand\subf{_{\scriptstyle \FF}}
\newcommand\po{^{\scriptscriptstyle O}}
\newcommand\pu{^{\scriptscriptstyle U}}
\newcommand\ppt{^{\scriptscriptstyle T}}
\newcommand\so{_{\scriptscriptstyle O}}
\newcommand\su{_{\scriptscriptstyle U}}
\newcommand\sst{_{\scriptscriptstyle T}}

 \newcommand{\GL}{{\mathcal{G}}}

\newcommand{\NN}{{\mathbb{N}}}
\newcommand{\HH}{{\mathbb{H}}}

\newcommand{\RR}{{\mathbb{R}}}

\newcommand{ \Kalg}{K^{\rm{alg}}}

\newcommand{\Z}{{\mathbb{Z}}}
\newcommand{\Q}{{\mathbb{Q}}}

\newcommand{\KK}{{\mathcal{K}}}

\newcommand{\CCC}{{\mathcal{C}}}

\newcommand{\FF}{{\mathbb{F}}}
\newcommand{\cR}{{\mathcal{R}}}

\newcommand{\oo}{\mathbf O}
\newcommand{\uu}{\mathbf U}
\newcommand{\pp}{\mathbf {Sp}}

\newcommand{\bbu}{\mathbf{BU}}
\newcommand{\bbo}{\mathbf{BO}}
\newcommand{\bbsp}{\mathbf{BSp}}

 \newcommand{\real}{C\pr}

 \newcommand{\realo}{C_0\pr}

 \newcommand{\complex}{C\pc}

 \newcommand{\complexo}{C_0\pc}

  \newcommand{\iS}{C^{\RR}(iS^1)}
 \newcommand{{\BS}}{monograph}  
 
 \newcommand{\ev}{\mathrm{ev}}
 \newcommand\id{\mathrm{id}}
 \newcommand{\diag}{\mathrm{diag}}
  
 \newcommand{\topr}{{\bf{R}}}

\newtheorem{Exam}[Thm]{Example}
 \newcommand{\york}{Clifford\,}
 \newcommand{\lancaster}{Lie algebra\,}

\newcommand\TTT{\rule{0pt}{2.6ex}} 
\newcommand\BBB{\rule[-1.2ex]{0pt}{0pt}} 

\makeindex   
\makeindex[name=notation, title=Index of Mathematical Notation]

\begin{document}

\title[Real $K$-Theory for $C^*$-Algebras: Just the Facts \\  \today ]
{Real $K$-Theory for $C^*$-Algebras: just the facts\\   \today }

 \author{Jeffrey L. Boersema}
\address{Mathematics Department, Seattle University, boersema@seattleu.edu}
\author{Chaim (Claude) Schochet} 
\address{Mathematics Department, The Technion, clsmath@gmail.com}
 
 \begin{abstract} 
    
  This {{\BS}}     is intended to present    the basic properties of $KO$-theory for real $C^*$-algebras and to explain its relationship with complex $K$-theory and with 
   $KR$- theory.
     Whenever possible we will rely upon proofs in
  printed literature, particularly the work of Karoubi 
  \cite{K}, Wood \cite{Wood}, Schr\"oder \cite{Schr}, and more recent work of Boersema \cite{BoersemaKT}
   and J. M. Rosenberg \cite{R}.   In addition, we shall explain how $KO$-theory is related to the Ten-Fold Way in physics and point out how some deeper features 
   of $KO$-theory for operator algebras may provide powerful new tools there. Commutative real $C^*$-algebras {\emph{not}} of the form $\real (X)$ will play a special role.

   Unfortunately, there is 
   no single reference for $KO$-theory for operator algebras that begins to compare with Blackadar's wonderful exposition of complex $K$-theory \cite{Blackadar}.   This {\BS} is intended to provide a platform upon which mathematicians and mathematical physicists can rely 
   in order to use these new tools in their research.   As we are writing for a diverse audience of functional analysts, topologists, and physicists, we often present material well-known to one group of people and unfamiliar to another.  
   
   Part I is a general introduction to real $C^*$-algebras, presenting basic examples and constructions, 
   homotopy properties, and developing the general structure of the category of real $C^*$-algebras.  After describing general homology theories, we define $KO_*$, display its basic properties and lead up to Bott Periodicity.
   
   Part II consists of a deeper look at some of the more advanced tools of $KO$-theory. We discuss various relationships
   between $KO$ and $K$-theory,  display some concrete calculations of various $KO$-groups, 
   introduce $KO$-theory with coefficients, and demonstrate how elements of each $KO$-group may be 
   represented in terms of unitaries.

   Part III is devoted to $KO^*$ for topological spaces. For this we review basic topological constructions, fibrations, fibre bundles, and vector bundles.  We define $KO^*(X)$ and present a novel proof of Swan's theorem relating $KO^*(X)$ to $KO_*(C\pr (X))$  for compact spaces. 
   
   Part IV introduces Clifford algebras and their relationship with Bott Periodicity and moves on to symmetric
   spaces and quantum symmetry.   There is a section entitled ``Topology and Physics" written by our colleague Nadav Orion, which aims to make the path to physics clearer.    We summarize this path via 
   the Homotopy Road Map.

    \end{abstract}
\maketitle

\newpage
\tableofcontents

\newpage
   \section{\bf{Introduction}}
    \label{Section:Introduction}
 
 \vspace{.2in}
 \centerline{\bf{\emph{``Real $K$-theory behaves differently from complex $K$-theory." }}}
 \centerline {M.F. Atiyah \cite{AtiyahColVol2}, p. 4.  }
 
 \vspace{.2in}

 \vspace{.2in}
\centerline{\bf{A little history }}
 \vspace{.1in}

\index{K-theory@$K$-theory!history|(}

  First, a little history to put $KO$-theory into historical context.  Let $R$ be an arbitrary ring. The group 
  $K_0(R)$ was introduced by Alexander Grothendieck   in 1957 as part of his revolution in algebraic geometry, as a way to formulate what is now known as the Grothendieck-Riemann-Roch theorem. It was based upon coherent 
 sheaves on an algebraic variety.  It was extended using modules over  arbitrary rings and is   difficult to compute in general.\footnote{For everything you would want to know about the early history of 
 algebraic $K$-theory, see Weibel \cite{Whistory}.} The group  $K_1(R)$ was defined in 1962, $K_2(R)$        in 1967, and        $K_n(R)$
in 1970.

 At about the same time, Raoul Bott was studying the homotopy theory of the classical Lie groups and discovered Bott Periodicity - that homotopy groups
 of the classical unitary groups $\pi _n (\bigcup _k U_k) $ were periodic of period 2 and the homotopy groups of the classical orthogonal groups $\pi _n (\bigcup _k O_k) $ were periodic of 
 period 8.   This fit in, somehow, with Elie Cartan's earlier classification of symmetric spaces, though certainly that was not clear at the time, and  the connection to the Ten-Fold Way in physics was not understood.\footnote{ other than the fact that $2 + 8 = 10.$}
 
 Michael Atiyah and Friedrich Hirzebruch made the fundamental breakthrough in 1959 (cf. \cite{AH}), introducing $K$-theory for complex vector bundles and $KO$-theory for real vector bundles,\footnote{Atiyah used the notation $K(X)$ for complex bundles, $KO(X)$ for real bundles, and $KR(X)$ for Real bundles. Others used $KU(X)$ for the complex case.}
using Bott's periodicity result to show   
   that the associated cohomology theories on compact spaces were periodic (of period 2 and 8 respectively), and solving several major outstanding open questions in classical algebraic topology.  The Atiyah-Singer Index Theorem, perhaps the most important mathematics result in the 1960's, grew out of this milieu, and that is just one of many important results in algebraic topology established in the 1960's using $K$-theory for compact spaces.  
 
 Atiyah, Bott, and Arnold Shapiro (1964) \cite{ABS} (ABS) classified Clifford algebras, related them to   $K$-theory and $KO$-theory, and predicted the existence of a proof of the Periodicity theorem 
 based upon Clifford algebras. 
  
 As their paper was about to appear, Atiyah's student R. Wood produced (1966)\cite{Wood}  a   proof of the Periodicity Theorem along the lines suggested by ABS, 
 where Clifford algebras play an essential role and where  real Bott Periodicity is a consequence of eight homotopy equivalences.
 He took a fundamental 
 step forward, in that he stated the theorem in terms of the general linear groups of real and complex Banach algebras, and not just compact spaces, though his proof was rather sketchy.  (See also Milnor \cite{M3}).  This was followed by several detailed papers 
 by Max Karoubi (cf. \cite {K68}). Karoubi  generalized the Periodicity Theorem to Banach categories and Banach rings, proving periodicity in detail. Karoubi and Orlando Villamayor (\cite{KV}, \cite{KV71} ) explored     its relationship to algebraic $K$-theory. A few 
 years later Karoubi introduced Hermitian $K$-theory (\cite{Knumber14}, \cite{Knumber29}), where he proved what he called ``A Fundamental Theorem" which is 
 another road to Bott Periodicity (see the introduction to \cite{Knumber14}).  This theorem has important consequences nowadays in algebraic geometry 
 and in  quadratic forms over any field, not just the field of real numbers.
 Karoubi anticipated many of the results that we use in this {\BS}. In particular, the idea of 
 desuspension was used to produce the negative algebraic $K$-theory groups developed by Bass (\cite{Bass}, 1988)
 
 The first major application of complex $K$-theory for $C^*$-algebras was the   work of Brown-Douglas-Fillmore ( \cite{BDF1}(1973), \cite{BDF2}(1977) ), answering a question of Paul Halmos. Suppose that $T$ is an essentially normal operator.\footnote{meaning that $T$ is a bounded operator on complex Hilbert space and $T^*T - TT^*$ is a compact operator.} When is $T$ the sum of a normal operator and a compact operator? The answer turns out to depend 
 upon the complex $K$- homology theory of the essential spectrum of the operator and hence was best described in $C^*$-algebra terms. (Initially this involved their functor $Ext(X)$ for 
 the essential spectrum $X$, a compact subset 
 of the plane, 
 which was isomorphic to $K_1(X)$, but then their work was tied to the groups  $KK_*(A,B)$ for $C^*$-algebras $A$ and $B$  of Gennadi Kasparov, which definitely required complex 
 $K$-theory for $C^*$-algebras.
 

    This success led to a spectacular flowering of work in functional analysis in the 1970's that continues, using complex $K$-theory and Kasparov theory    as primary tools.  For example, the on-going classification 
 of simple, separable, nuclear  complex $C^*$-algebras relies heavily upon the Elliott invariant, which is built on complex $K$-theory.
 
 Real $K$-theory was not as popular. For example, there are several different modern texts on complex $K$-theory for $C^*$-algebras, such as 
 \cite{Blackadar}, \cite{Hatcher}, \cite{H}, \cite{K}, \cite {RLL},  but only \cite{Schr} for real $K$-theory.  However, there has been 
 a steady stream of results for $KO$ - both for spaces and for  real $C^*$-algebras during this period, and in the last ten years its usefulness in various areas of mathematics has been 
 established.
 For example, in Loring and Sorensen \cite{LS} (2016), $KO$-theory for real $C \sp *$-algebras is used to prove that almost commuting pairs of real $n \times n$ matrices can be approximated in a uniform way, independent of $n$, by exactly commuting pairs of matrices. This result extends to the real case what had been proven in the complex case much earlier by  Lin \cite{Lin}.
 
 There have been serious applications in algebraic topology as well. $KO_*$ plays a fundamental role in providing obstructions for the existence 
 of a metric of positive scalar curvature on a compact manifold.  It interacts with various versions of the Baum-Connes conjecture, and helps 
 bring light to orbifold string theories in physics. For all of these matters and much more, we refer to Rosenberg's deep and valuable survey article \cite{R}.
  
 There has been a recent resurgence of interest in various problems in physics related to the Cartan classification of symmetric spaces, though the physicists would better describe it as studying symmetries in the atomic structure of matter. This generally involves $KO$-theory (since 8 of the 10 classes are related to $KO$, as we shall explain).  Topological insulators and tilings are  especially 
 fertile areas for $K$ and for $KO$, as witness, for example, the work of Stone et al \cite{Stone}, Hastings and Loring \cite{HL},  and Akkermans, Don, Rosenberg, and Schochet \cite{ADRS}.
\index{K-theory@$K$-theory!history|)}
 
  \vspace{.2in}
\centerline{\bf{Section Summaries}}
 \vspace{.1in}

This {\BS} is designed to provide a basic   source for those in mathematics and physics who wish to use $KO$-theory in their work.  
There is no such source at present:
we mean to remedy that problem.

 We now briefly summarize our work.  
 \vspace{.2in}
 
 {\bf{Part 1: $KO_*$- Theory for Real $C^*$-algebras}}  
 
  \vspace{.1in}
  
{\bf{Section~\ref{Section:RealC*-Algebras}}} deals with the basics of real $C^*$-algebras, giving many examples. There is a  fundamental difference between the real and complex commutative settings. Every complex commutative $C^*$-algebra is isomorphic to an algebra of continuous complex-valued functions on a locally compact space, by the Gelfand theory.  The analogous statement is not true for real $C^*$-algebras and will lead to the 
possibility of ``desuspending" a real $C^*$-algebra. 

\vspace{.1in}

{\bf{Section~\ref{Section:NewC*-Algebras}}} focuses on 
  how to construct various kinds of $C^*$-algebras including sums, tensor products, and how to complexify.  We also introduce constructing a real $C^*$-algebra
  from a compact space with an involution and demonstrate how $\otimes$ and $\Dirlim$ interact.
  
  \vspace{.1in}
  
  {\bf{Section~\ref{Section:HomotopyConstructions}}}
 is all about homotopy constructions for real $C^*$-algebras. We recall the definitions of the cylinder, the cone, and the suspension and we introduce the desuspension of a real $C^*$-algebra 
 as well as  cofibrations.
 
 \vspace{.1in}
 
 {\bf{Section~\ref{Section:HomologyTheories}.}}  Here we introduce axioms for homology and cohomology theories on real and complex $C^*$-algebras and show how to use them.  Our immediate interest is in showing how much of the general theory of $KO_*$  goes through without invoking the depth of the Bott Periodicity theorem, as well as illustrating many of the similarities between real and complex $K$-theory.  We also introduce cofibre homology and cohomology theories.
 
 \vspace{.1in}

In {\bf{Section~\ref{Section:KODefinition}}}
 we introduce the homology theory $KO_*$ for real $C^*$-algebras, the main focus of this {\BS}. 
   There are several distinct approaches to $KO$. We use what we believe to be the simplest 
 approach, not relating much to algebraic $K$-theory, Hermitian $K$-theory, or other homotopy constructions.  We show that $KO_*$ is a homology theory and introduce some of its relationships to $K_*$. 
 Still no mention of Bott Periodicity.
 
 \vspace{.1in}

 {\bf{Section~\ref{Section:Products}.}}   Ordinary cohomology of a topological space is a graded ring, and this is also true for the $K$ and $KO$-theories when focused on spaces, but when we pass to non-commutative $C^*$-algebras the ring structure disappears.  
 Nevertheless, we show that some partial product structure is retained. 
 $KO_*(\RR)$ is a graded ring (with a complicated structure!), and each $KO_*(A)$ is a graded $KO_*(\RR)$-module. This module structure turns out to be very important when dealing with the real  K\"unneth 
 and Universal Coefficient Theorems.
 
  \vspace{.1in}
  
 In  {\bf{Section~\ref{Section:Bott}}} we (finally!) introduce Bott Periodicity, both in the complex and the real settings, and we explain how this allows us to extend $KO$-theory and $K$-theory to periodic homology theories defined 
 for all integers. 
 
 \vspace{.2in}

  {\bf{Part 2: A Deeper Look at the Structure of $KO$ and $K$}}
  
     \vspace{.1in}
 
 {\bf{Section~\ref{Section-United}}} 
  is devoted to CR and CRT united $K$-theory. These help us understand the complicated relationships between real and complex $K$-theory.  They are critical for 
  some important results in real $K$-theory. 
  
      \vspace{.1in}
      
  {\bf{Section~\ref{Section-Examples}.}} The section is entitled ``Examples and Calculations" which describes exactly what is contained in it.
     
     \vspace{.1in}
     
       {\bf{Section~\ref{Section-Kunneth}}}.  We present the K\"unneth Theorem which deals with computing in the complex setting $K_*(A\otimes B)$ in terms of $K_*(A)$ and $K_*(B)$. The real analog for computing $KO_*(A \otimes B)$ requires united $K$-theory for the full generality. We also present a vast generalization of the Atiyah-Hirzebruch spectral sequence, which deals with computing the $K$-theory of a filtered $C^*$-algebra, both in 
       the real and the complex setting. 

        \vspace{.1in}
     
     {\bf{Section~\ref{KKSection} }} is all about $KK$-theory for real and complex $C^*$-algebras. We sketch the relevant definitions and primary properties and state the Coefficient Coefficient Theorem, both in the real and complex settings.
 
        \vspace{.1in}
     
      {\bf{Section~\ref{Section-KwithCoefficients}.}} This  section contains our treatment of $KO$-theory with coefficients.  It is well-known that $H^*(X; G)$, the cohomology of a space with coefficients in a group $G$, is sometimes much more useful and certainly more tractable than when you take $G = \ZZ$.  We suspect that the same will 
  be true for $KO$ and hence we present how to define it.

 \vspace{.1in}   
     
  {\bf{Section~\ref{Section-KwithUnitaries}}} This section presents a novel approach to the groups $K_*(A)$ and  $KO_*(A)$, showing that each of the ten groups may be described concretely in 
  terms of unitary elements of the complexification of $A$ and matrix algebras over $A$.
    
 \vspace{.2in}

  {\bf{Part 3: $KO^*$-theory for Spaces}}

  \vspace{.1in}

  {\bf{Section~\ref{top spaces}}}  introduces basic constructions in topological spaces such as the product, the join, the smash product, and the suspension. It  establishes 
 the fundamentally important identification $F(SX, Y) \cong F(X, \Omega Y) $.   We will assume throughout the {\BS} that {\emph{all topological spaces are 
 compactly generated}}. We introduce CW-complexes briefly.
 
  \vspace{.1in}
 
  {\bf{Section~\ref{Section:Fibrations}.}} This is a classical homotopy section. We introduce fibrations, fibre bundles, principal bundles and briefly explain their important 
  properties.  We conclude by discussing the universal principal bundle $\oo \longrightarrow {\bf{EO}} \longrightarrow {\bf{BO}} $.
 
  \vspace{.1in}
 
  {\bf{Section~\ref{Section:BottRevisited}.}} In this short but very important section we revisit the Bott Periodicity theorem, state it in terms of the classifying space ${\bf{BO}}$ and
  present the first version of the Homotopy Road Map which connects Bott Periodicity via classical homotopy to Clifford algebras to symmetric spaces to physics.
 
  \vspace{.1in}
 
  {\bf{Section~\ref{Section:LowHomotopyGroups}.}} Bott Periodicity tells us to concentrate on the homotopy groups of spaces like $\uu /\oo $, and it is true that 
  \[
  \pi_k( \uu / \oo ) \cong \lim_{n \to \infty} \pi _k(U_n/O_n)  
   \]
   and so when $n$ is large  compared to $k$ then  
   \[
    \pi_k( \uu / \oo ) \cong   \pi _k(U_n/O_n)  
   \]
 but sometimes (and specifically, in some recent developments in physics) we need to know 
 $\pi _k(U_n/O_n) $ when $n$ is small. 
   So we pause to spend 
 a section computing such groups for $k = 0, 1, 2$.  

  \vspace{.1in}

  {\bf{Section~\ref{Section:RealVectorBundles}}}
   introduces real and complex vector bundles, mostly over compact spaces, and states the classification theorem with which isomorphism classes of vector bundles over $X$ are related to homotopy classes of maps from $X$ into the appropriate Grassmann manifold.

   \vspace{.1in}

  {\bf{Section~\ref{Section:KOforSpaces}}} introduces 
 $K$ and $KO$-theory for compact spaces, defined via complex and real vector bundles respectively. This material is very well-treated in various books and so 
 we simply state the highlights. We mention the Chern character and the Pontrjagin character which relate     $K$ and $KO$-theory for compact spaces to the cohomology of those spaces.
  
    \vspace{.1in}
 
  {\bf{Section~\ref{Section:RelatingKO-theory}}}  is devoted to establishing the relationship between $KO^*$-theory for spaces and $KO_*$-theory for real $C^*$-algebras and similarly for $K^*$-theory for spaces and $K_*$-theory  for  complex $C^*$-algebras. 
 We present a novel proof of Swan's classical result relating 
 the  $KO_*$-theory of a real $C^*$-algebra $C\pr (X)$ and the  topological $KO^*$-theory of the compact space $X$.  
 
   \vspace{.1in}
 
  {\bf{Section~\ref{Section:KRTheory}}} 
   is special.  Atiyah introduced $KR^*(X, \tau )$  (where $X$ is a compact space and $\tau $ is an involution of $X$) for specific applications in functional analysis, as we shall briefly explain.    We show that there is a natural isomorphism of the form
$$
 KR^*(X, \tau ) \cong KO_*(A),
 $$
  where   $A$ is a real commutative $C^*$-algebra constructed from $X$ and from the involution $\tau $.  

  \vspace{.2in}

  {\bf{Part 4: Clifford Algebras and Symmetric Spaces}}

     \vspace{.1in}

   {\bf{Section~\ref{Section:CliffordAlgebras}.}}  This section is a quick expos\'e of the   basic information about real Clifford algebras that we need going forward. We concentrate upon the 
  Lie groups of each Clifford algebra in the standard sequence and the homogeneous spaces constructed from them.  We emphasize the relationship between 
  this sequence of homogenous spaces and the sequence of homogeneous spaces obtained from the Bott Periodicity theorem, referring to  the Homotopy Road Map for reference.  We contrast the \york perspective and the \lancaster perspective for this list of spaces and explain their relationships. 
     
         \vspace{.1in}

 {\bf{Section~\ref{Section:CliffordAlgebrasBottPeriodicity}}}  gives a brief discussion of the Wood-Karoubi proof of Bott Periodicity and, specifically, its relationship to Clifford algebras. 
      
     \vspace{.1in}

 {\bf{Section~\ref{Section:CliffordAlgebrasComplex}.}}  In this short section we mention complex Clifford algebras and indicate how they too relate to the Homotopy Road Map via Bott Periodicity.

     \vspace{.1in}

 {\bf{Section~\ref{Section:SynmetricSpaces}}}  gives a basic introduction to symmetric spaces.  Every symmetric space arises as a homogeneous space of Lie groups and so 
  the discussion turns at once to associated Lie groups and Lie algebras. Compact symmetric spaces have been classified and there are ten different types. These ten 
  types correspond to the ten rows on the Homotopy Road Map and we present Step 2 of the Road Map, adding on this correspondence. We then present 
  a ``Timely Example", namely the symmetric space $SU_4/O_4$.  We present a basis for the Lie algebras that reflect their geometry and define a symmetry 
  $T$ that acts on the Lie algebra with very specific characteristics.  Appendix 20.1 is by Nadav Orion. In it he explains the basic connection between topology and physics in our context.
  
  \vspace{.1in}
 
 {\bf{Section~\ref{Section:QuantumSymmetry}}}  connects us directly with the Ten-Fold Way in physics.  We discuss the actions of Hamiltonians and the symmetries that 
  they may or may not respect. We introduce the $T,C, S$ symmetries and link them to the ten types of compact symmetric spaces, adding them in order to complete the Homotopy Road Map.
 
     \vspace{.1in}
     
     What's not in the {\BS}?        We have  not 
 discussed some of the applications of $KO$-theory in differential geometry, as this requires quite a bit of background to appreciate. 
 Nor have we discussed $KO$-homology theory for compact spaces. We have not discussed the Hopkins-Hovey isomorphism \cite{HH}
 \[
 MSpin_*(X) \otimes _{MSpin(pt)} KO_*(pt) \cong KO_*(X)
 \]
 relating $KO_*(X)$ and Spin bordism. 
We  have not discussed the generalization of $KO$-homology for spaces to 
  $KO^*(A)$,
 a cohomology theory for real $C^*$-algebras, defined most easily using the real version of Kasparov's $KK$-theory. There is some literature on each 
 topic, but these are substantially more technical  and, to this point, less central than the parts of the theory that we have covered.   Finally, we have not gone deeply into the physics, but that is truly for another time. 
 
This work arose in connection with the second author's collaboration with Erik Akkermans. We are most grateful to him and to his students, particularly to Nadav Orion, who have provided invaluable feedback and true insight in the development of this work.  We are also very grateful to Jean Bellissard, Joachim Cuntz,  
Yigal Kamel, 
John Klein, N. Christopher Phillips, Emil Prodan, Jonathan  Rosenberg, and Chuck Weibel, for providing us with various deep insights as this {\BS}  developed.

We especially want to acknowledge with great pleasure the extensive philosophical and detailed help that we have received from Max Karoubi, who, more than anyone, is responsible for the very existence of $KO$-theory for spaces and algebras as we know it.

 \newpage
 \part{$KO_*$-Theory for Real $C^*$-Algebras}

   \section{\bf{Real $C^*$-algebras}}
  \label{Section:RealC*-Algebras}
  
 Two points on  notation. First, we have been confused over the years by various papers using $C(X)$, $K_0(A)$, and the like, without specifying whether they had real or complex 
versions in mind.\footnote{and we each plead guilty to this sin.} Frequently authors deal with this by stating at the beginning that everything is complex or that everything is real.  In this {\BS} both possibilities will arise. We will be 
more precise. 

A second point. Usually the real numbers occur in this {\BS} as a real 
vector space or  with their real $C^*$-algebra structure, 
and when that happens we shall denote them as $\RR $.  Sometimes, though, we want to think of the real numbers
simply as a locally compact topological space, In that case 
we shall use the notation $\topr $.  We will label the points of $\topr $ as usual, so that we can write continuous functions, such as $x \mapsto x^2$ or   the   involution $x \mapsto -x$.
  
 \index[notation]{R@$\topr$}
 \index[notation]{R@$\RR$}
 
 For any \index{compact space} compact\footnote{We assume that locally compact spaces and compact spaces are always Hausdorff.}  space $X$ and normed algebra $A$, let  $C(X, A )$ denote the algebra of continuous functions from $X$ to $A$ with pointwise operations and the sup norm.   Similarly, for a locally compact space 
$X$, let $C_o(X, A)$ denote the algebra of continuous functions from $X$ to $A$ which vanish at infinity.  We then abbreviate as follows:

\index[notation]{C@$\real(X)$} 
\index[notation]{C@$\realo(X)$} 

\[
\real (X) = C(X, \RR ) \qquad\qquad \realo (X) = C_0(X, \RR )
\] 
\[
\complex (X) = C(X, \CC ) \qquad\qquad \complexo (X) = C_0(X, \CC)
 \]
  
 When developing $C^*$-algebras and their $K$-theory in the past, people generally used $C(X)$ for complex-valued functions on $X$ and 
  $K_*$ for complex $K$-theory, because almost everything involving $C^*$-algebras was over the complex numbers.  
  That is no longer the case, but, following tradition, we shall continue to use $K_*(A)$ for complex $K$-theory.
  We shall use $KO_*$ for real $K$-theory throughout\footnote{except in Section~\ref{KKSection}, as we shall explain there.}.   
  
 \begin{Def} \index{real $C^*$-algebra!definition} A {\emph{real $C^*$-algebra }}  is a real Banach $*$-algebra $A$ 
 \begin{enumerate}
 \item
   which is isometrically $*$-isomorphic  to a norm-closed
 $*$-algebra of bounded operators on real Hilbert space. 

 \item  or, equivalently, such that for all $a \in A$, 
\begin{itemize}
 \item $ \|a^*a \| = \| a \|^2$
 \item and the element $a^*a $ has spectrum contained in $[0, \infty )$. 
 \end{itemize}
 \end{enumerate}
 \end{Def}
 
 The reader should not expect it to be obvious to see that these two definitions are equivalent -- this is a deep result of 
 T. Palmer \cite{Palmer}. See the monograph by Goodearl \cite{Goodearl} for a development of theory to the point of proving that every abstractly defined real $C \sp *$-algebra can be isometrically represented on a real Hilbert space. We also note to a reader familiar with complex $C \sp *$-algebras: the presence of the second point of (2) is surprising and seemingly redundant. In fact, in the complex case the first point of (2) is sufficient in the definition and implies the second point. But in the real case, both points are necessary. 
 
 We let  $\mathscr{C}^*_{\rm{real}}$ denote the category of real $C^*$-algebras and 
 real $*$-homomorphisms (always assumed to be continuous; frequently we simply say ``maps").
 The category of complex $C^*$-algebras, with complex $*$-homomorphisms is denoted by $\mathscr{C}^*_{\rm{complex}}$.
 \index[notation]{C@$\mathscr{C}^*_{\rm{real}}$}
  \index[notation]{C@$\mathscr{C}^*_{\rm{complex}}$}

 We present a variety of examples.
 
 \begin{Exam} \label{Examples-RealC*-algebras}
 \begin{enumerate}
 \item
 \index{real $C^*$-algebra!finite dimensional}
 If $A$ is a finite-dimensional real $C^*$-algebra, then it is a direct sum of matrix algebras $M_n(\FF)$, where $\FF = \RR, \,\, \CC,\,\,$ or $\HH$, the quaternions. 
 
 \index[notation]{MNA@$M_n(A)$}
 \index[notation]{Hamilt@$\HH$}
  
 \item If $X$ is a locally compact space then $\realo (X)$ with the sup norm is a commutative real $C^*$-algebra. If $X$ is compact then $\real (X)$ is a unital commutative real $C^*$-algebra.

\item
If we take $\mathcal{H}\sr$ to be a real Hilbert space, then we have the algebra $\mathcal{B}\sr(\mathcal{H}\sr) $ of bounded linear operators on $\mathcal {H}\sr$ and the algebra of compact operators ${\mathcal{K}}\sr(\mathcal{H}\sr)$, which are both real $C \sp *$-algebras. In the case that $\mathcal{H}\sr$ is taken to be a separable infinite dimensional real Hilbert space, then we denote the corresponding operator algebras simply by $\mathcal{B}\sr$ and $\mathcal{K}\sr$ respectively.
 \index[notation]{Bounded@$\mathcal{B}\sr$}
 \index[notation]{Kompact@$\mathcal K \sr$}

\item If $\mathcal{S}$ is any finite or infinite subset of any real $C \sp *$-algebra $A$, then we let $C\sr ^*(\mathcal{S}) $ denote the smallest closed real $*$-subalgebra of $A$ which contains the set $\mathcal{S}$. 
In particular, for any bounded operator $T \in \mathcal{B}\sr$, we have $C\sr^*(T)$ as the closed real  *-subalgebra of $\mathcal{B}\sr(\mathcal{H}\sr)$ generated by $T$. It will be commutative if and only if  $T$ is normal (that is, $T T^* = T^* T$).


 \item
 \index{compact operators}
Any  closed and $*$-closed ideal of a real $C^*$-algebra is again a real $C^*$-algebra. 
For example, ${\mathcal{K}}\sr$ is a $*$-closed ideal in $\mathcal{B}\sr$.

\item
  Suppose that $A$ is a real $C^*$-algebra and $X$ is a locally compact space. Then, as mentioned above,  we may form $\realo (X, A)$, all continuous functions from $X$ to $A$ vanishing 
at infinity.  With pointwise operations and sup norm, this forms a real $C^*$-algebra.  Similarly, if $X$ is compact and $A$ is unital then $\real (X, A)$ is a unital real $C^*$-algebra. 

\index[notation]{CXA@$\real(X,A)$}

 \item Suppose that $X$ is a locally compact space and $\tau :X \to X$ is an involution of $X$; that is, a  homeomorphism of $X$ with $\tau ^2 = 1$.  Define
 \[
\realo (X, \tau ) =  \{ f \in \complexo (X) \mid f(\tau (x)) = \overline{f(x)} \text{ for all $x \in X$}  \} . 
 \]
 \index[notation]{CXt@$\realo (X, \tau )$}
 Then $\realo (X, \tau )$ is a real commutative $C^*$-algebra.    It is a nice exercise to show that
 \[
 \realo (X, \tau ) \cong \realo (X)
 \]
 if and only if $\tau  = \id $ is the identity homeomorphism.     See Theorem~\ref{Arens} for a complete discussion of real commutative $C^*$-algebras.

\item  The locally compact space $\topr $ has two obvious involutions, namely, the  identity map and
the homeomorphism $\hat\tau (x) = -x$.  This gives us two different real $C^*$-algebras ``based" on $\topr $, namely $\realo (\topr )$,  and 
$\realo (\topr, \hat\tau ). $
This latter $C^*$-algebra is going to be the core of ``desuspension" and so we state its definition explicitly, using Schr\"oder's notation:
\[
\realo(i\topr ) = \{ f \in \complexo (\topr ) \mid f(-x) = \overline{f(x)} \text{ for all $x \in \topr$ } \}.
\]
\index[notation]{Cr0R@$\real_0(i\topr)$}
 
\item  Similarly, the circle  $S^1 $ has three obvious  involutions. Regarding $S^1$ as the unit circle in the complex plane, they are
\[
z \overset{\id}\longmapsto  z
 \qquad\qquad z \overset{\hat{\tau}}\longmapsto \bar{z}      \qquad\qquad z \overset{\sigma}\longmapsto -z
\]
namely the identity map, complex conjugation, and the antipodal map.
   This gives us three different real commutative $C^*$-algebras arising from the circle $S^1$,
   \begin{align*}
   \real (S^1, \id) &= \{ f \in \complex(S^1) \mid f(z) = \overline{f(z)}, z \in S^1 \} = \real(S^1) \\
 \real (S^1, \hat\tau )  &=     \{ f \in \complex (S^1) \mid f( \overline{z} ) = \overline{f(z)}, z \in S^1  \} \\ 
  \real(S^1, \sigma) &=              \{ f \in \complex (S^1) \mid  f(-z)  = \overline{f(z)}, z \in S^1 \} \; . 
    \end{align*}
All three of these algebreas will be important in this work. Anticipating events, we let $\iS$ denote the second algebra above, so
$$\iS   =      \real (S^1, \hat\tau )  $$
\index[notation]{CR0S@$\real_0(i S^1)$}
\index[notation]{T@$T$}

\item   \index{Toeplitz algebra}  
The {\emph{real Toeplitz algebra}} ${\mathcal T}\sr $  is   defined as the
 real $C^*$-algebra\footnote{To be more precise,  take 
 $\mathcal{H}\sr = L^2\sr $
  to be real $L^2$-functions on the 
 circle and 
 $ \mathcal{H}\sr ^2 \subset L^2\sr $
  to be the real Hardy space, with canonical projection 
 \[
 P :  L^2\sr    \longrightarrow   \mathcal{H}\sr ^2 
\]
If  $\phi $ is a complex-valued continuous function on the circle then define 
 $T_\phi f = P(\phi f)$.   
}

\[
\mathcal{T}\sr  = \{ T_\phi  + K \mid \phi \in \iS , K \in \KK\sr \}.
\index[notation]{TR@$\mathcal{T} \sr$}
\]
with associated short exact sequence
\[
0 \longrightarrow \KK\sr \longrightarrow \mathcal{T}\sr   \overset{\sigma}\longrightarrow \iS \to 0 .
\]

\item  Let $X = \{ a, b \}$, the space with two discrete points.  It has two obvious involutions, namely the identity and the involution $\tau $ which interchanges the two points. 
In the first case, we have 
\[
\real (X, \id ) \cong \real (X) \cong \RR \oplus \RR 
\]
with pointwise operations. In the second case, we have
\[
\real (X, \tau ) \cong \CC
\]
the complex numbers, regarded as a real $C^*$-algebra. This is yet another example of a real, commutative unital $C^*$-algebra which is {\emph{not}} of the form $\real (X)$. 
   
\end{enumerate}
\end{Exam}

Recall that every complex unital commutative $C^*$-algebra $A$ is isometrically isomorphic to $\complex (X )$  for some compact 
space $X$, and a non-unital commutative $C^*$-algebra is isometrically isomorphic to $\complexo (X )$  for some locally compact 
space $X$.  The situation is profoundly different for real $C^*$-algebras, as witness the previous example.   Here is the full story.

 \begin{Thm}(Arens-Kaplansky \cite{AK}, Thm. 9.1)\label{Arens} 
 \index{real $C^*$-algebra!commutative} Suppose that $A$ is a commutative real 
 $C^*$-algebra. Then: 
 \begin{enumerate}
 \item
 If $A$ is unital, 
 then there exists a  compact space $X$ and an involution  
   $\tau : X \to X$   such that     $A \cong \real (X, \tau )$.  
  \item
If $A$ is not unital,  then there exists a locally  compact space $X$ and an involution 
  $\tau : X \to X$   such that    $A \cong \realo (X, \tau )$.
  \item There is an equivalence of categories between the category of commutative real unital $C^*$-algebras (with unital $*$-homomorphisms) and the category of compact spaces with involution (with continuous involution-respecting maps).
  \item Similarly, there is an equivalence of categories between the category of commutative real $C^*$-algebras 
  (with non-degenerate\footnote{A homomorphism $f \colon A \rightarrow B$ is non-degenerate if $\{f(a_j)\}$ is an approximate unit in $B$ whenever $\{a_j\}$ is an approximate unit in $A$.} $*$-homomorphisms)
  and the category of locally compact spaces with involution (with continuous involution-respecting maps).
\end{enumerate}  
  \end{Thm}
  
  \begin{Rem} The real commutative $C^*$-algebras 
  $\realo (i\RR )$ and
   $\real (iS^1)$ 
   will both be very important in our story.
   They should be regarded as the poster children for 
  commutative real $C^*$-algebras that seem to have no complex analog. Or stated better, there can be multiple real $C \sp *$-algebras that are real analogs of a given complex $C \sp *$-algebra. For example, in (9) above all three real $C \sp *$-algebras -- $\real (S^1)$, $\real (iS^1)$, and $T$ --  are plausible analogs of $\complex (S^1)$. The notion of {\it complexification} in the next section will clarify what we mean by this.
     \end{Rem}
  
 \begin{Def} \index{simple $C \sp *$-algebra} \label{Def-simpleC*-alg}
A real or complex $C \sp *$-algebra $A$ is {\it simple} if $A$ has no non-trivial proper closed two-sided ideals. In other words, $A$ is simple if any $I \subseteq A$ that satisfies $A I A \subseteq I, I^* = I,$ and $\overline{I} = I$, must satisfy $I = \{0\}$ or $I = A$.
\end{Def}

Examples of simple real $C \sp *$-algebras include $M_n(\RR)$ and $\KK\sr$. However, if $X$ has more than one point, then $\real(X)$ is not simple since an non-trivial idea is
$$I_x = \{ f \in \real(X) \mid f(x) = 0 \}$$
for any $x \in X$.

  
 \newpage
   \section{\bf{New $C^*$-algebras from Old Ones}}
  \label{Section:NewC*-Algebras}
  
 In this section we review several methods of constructing new   $C^*$-algebras from old ones.

\begin{Def} 
\index{complexification} 
\index[notation]{Ac@$A\sc$}
Let $\mathcal{H}$ be a real Hilbert space.  Given a real $C^*$-algebra $A$ of bounded operators on $\mathcal{H}$,
 its {\emph{complexification}} is defined as the complex algebra
 \[
 A\sc = A + iA
 \]
of bounded operators on $\mathcal{H}\sc$. Further,  
\begin{enumerate}
\item  $A\sc $ is         complete in the induced norm and hence is a complex $C^*$-algebra, independent of the representation of $A$.\footnote{We refer the reader to Goodearl \cite{Goodearl} pp. 66-67.}
\item  If $f: A \to B$ is a morphism of real $C^*$-algebras then it extends in an obvious way to $f\sc : A\sc \to B\sc $
by 
\[
f\sc (s + it) = f(s) + if(t).
\]
\item Complexification is a covariant functor from the category of real $C^*$-algebras to the category of complex $C^*$-algebras, denoted
\[
c : \mathscr{C}^*_{\rm{real}} \rightsquigarrow \mathscr{C}^*_{\rm{complex}}  .
\]
\end{enumerate}
\end{Def}

Note that any vector $\xi $ in $\mathcal{H}\sc $ can be written uniquely as 
$\xi = x + iy
$
where $x, y \in \mathcal{H}$, and then 
\[
\overline{\xi} =  \overline{x + iy} = x - iy .
\]

The complexification functor is not an equivalence of categories. 
Not every complex $C^*$-algebra arises as the complexification of a real $C^*$-algebra.  Further, two non-isomorphic real $C^*$ algebras can have isomorphic 
complexifications. For instance,  $M_2(\RR)$ is not isomorphic to  $ \HH $ as a real $C^*$-algebra, and yet 
\[
M_2(\RR)\sc \cong (\HH) \sc  \cong M_2(\CC).
\]
 Similarly, 
 the  three different examples of involutions on $S^1$ all lead to  isomorphic complexifications:
 \[
 \real (S^1)\sc \cong \real (iS^1)\sc \cong \real (S^1, \sigma )\sc \cong \complex (S^1) .
 \]
 
\begin{Pro}
Let $A$ be a real $C^*$-algebra. Then there is an involutive $*$-antiautomorphism\footnote{An {\it antiautomorphism} of a ring $R$ is a map $\sigma \colon R \rightarrow R$ which is additive and antimultiplicative, that is $\sigma (ab) = \sigma (b) \sigma(a)$ for all $a,b \in R$.}
$\sigma : A\sc  \to A\sc $ given by 
\[
\sigma (s + it) = s^* + it^* \; 
\]
which satisfies
$$A = \{ a \in A\sc \mid a^* = \sigma(a) \} \; .$$
\end{Pro}

Both of the statements of this Proposition are easy to verify directly. The fact that $A$ can be recovered from $\sigma$ means that there is an equivalence between the category of real $C^*$-algebras and the category of pairs $(B, \sigma)$ where $B$ is a complex $C^*$-algebra and $\sigma$ is a involutive *-antiautomorphism of $B$.

\begin{Exam}  \label{realforms}
Here are some examples of how this categorical equivalence works for particular real $C^*$-algebras.
\begin{itemize}
\item The real $C \sp *$-algebra $\RR$ corresponds to the pair $(\CC, \id)$.
\item The real $C^*$-algebra $M_n(\RR)$ corresponds to the pair $(M_n(\CC), \text{Tr})$ where $\text{Tr}$ is the matrix transpose operator.
\item 
\index[notation]{1@$\sharp$}
The $C^*$-algebra of quaternions $\HH$ corresponds to the pair $(M_2(\CC), \sharp)$ where
$$\sharp \colon \begin{bmatrix} a & b \\ c & d \end{bmatrix} \mapsto \begin{bmatrix} d & -b \\ -c & a \end{bmatrix} \; .$$
\item Let $\mathcal{H}$ be a real Hilbert space and let $\mathcal{H}\sc$ be the complexification. The real $C^*$-algebra of bounded operators $\mathcal{B}\sr(\mathcal{H})$ corresponds to the pair
$(\mathcal{B}(\mathcal{H}\sc), \text{Tr})$ where $\text{Tr}$ is the transpose operator defined by
$$\text{Tr}(T)(x) = (\overline{T^*(\overline{x}))} \; .$$
\item For a topological space $X$, the real $C^*$-algebra $\real (X)$ corresponds to the pair $(\complex ( X), \id)$.
\item For a topological space $X$ with continuous involution, the real $C^*$-algebra $\real (X, \tau)$ corresponds to the pair $(\complex ( X), \tau^*)$
where $\tau^*$ is the involution on $\complex(X)$ defined by $\tau^*(f)(x) = f(\tau(x))$.
\end{itemize}
\end{Exam}

\begin{Rem}
 We note that for a real $C^*$-algebra, we could alternatively focus on the involution
 $$\phi \colon s + it \mapsto s - it \; $$
 on $A\sc$ and recover $A$ as the set of fixed points, that is points which satisfy $\phi(a) = a$. Note that $\phi$ is multiplicative but is conjugate-linear, as it satisfies $\phi(\lambda a) = \overline{\lambda} \phi(a)$.
 Under this perspective, there is an equivalence between the category of real $C^*$-algebras and the category of pairs $(B, \phi)$, where $B$ is a complex $C^*$-algebra and $\phi$ is a conjugative-linear involutive $*$-automorphism of $B$.
 \end{Rem}
 
  \begin{Def}
 Suppose that $A$ is a real $C^*$-algebra, unital or not.  
 Define its {\emph{unitization}} \index{unitization} $A^+$ by
  $A^+ = A \oplus \RR $ additively, with multiplication given by 
  \[
  (a, r) \times (b,s) = (ab + rb + as, rs)
  \]
  and norm given by
  \[
  \| (a,r) \| = \sup \{ \| ab + rb \| \mid b \in A, \,\,\| b \| \leq 1 \}.
  \] 
 Then $A^+$  is a real unital $C^*$-algebra, $A$ is a closed ideal in $A^+$ and  
 $A^+/A \cong \RR$, so there is a natural short 
 exact sequence
 \[
 0 \longrightarrow A \longrightarrow A^+ \overset{\pi }\longrightarrow \RR \longrightarrow 0 .
 \]
 If $A$ already were unital, then $A^+ \cong A \oplus \RR$ as real $C^*$-algebras.
 \end{Def}

Here are some examples:
 \begin{enumerate}
 
 \item If $X$ is a locally compact space, then 
 \[
 \realo (X)^+ \cong \real (X^+) 
 \]
 where $X^+$ denotes the one-point compactification of $X$.   For example, $\realo (\topr )^+ \cong \real (S^1)$ since the one-point compactification of the real numbers is the circle.

 \item  Similarly,
   \[
 \realo (i\topr )^+ \cong \iS . 
 \]
   
 \end{enumerate}

 \begin{Def}  
 \index[notation]{MnA@$M_n(A)$} 
For any real $C \sp *$-algebra, we let $M_n(A)$ denote the real $*$-algebra of $n \times n$ matrices with entries in $A$, equipped with the usual addition, multiplication, and scalar multiplication. The $*$-action on $M_n(A)$ is given by
$$\begin{bmatrix} a_{11} & a_{12} & \cdots & a_{1n} \\ 
			a_{21} & a_{22} & \cdots & a_{2n} \\
			\vdots & \vdots & \ddots & \vdots \\
			a_{n1} & a_{n2} & \cdots & a_{nn} \end{bmatrix}^* 
			= 
\begin{bmatrix} a_{11}^* & a_{21}^* & \cdots & a_{n1}^* \\ 
			a_{12}^* & a_{22}^* & \cdots & a_{n2}^* \\
			\vdots & \vdots & \ddots & \vdots \\
			a_{1n}^* & a_{2n}^* & \cdots & a_{nn}^* \end{bmatrix}  \; .$$

To define the norm on $M_n(A)$, we follow Section~1.3 of \cite{RLL} in the complex case. Let $\phi \colon A \rightarrow \mathcal{B}\sr(\mathcal{H}\sr)$ be a a faithful representation of $A$ on a real Hilbert space. Then $\phi$ extends to a faithful *-algebra representation 
$$M_n(\phi) \colon M_n(A) \rightarrow \mathcal \mathcal{B}\sr(\mathcal{H}^n\sr) \; . $$
Then we give $M_n(A)$ the norm inherited from the norm on $\mathcal{B}\sr(\mathcal{H}^n\sr)$.
 
 \end{Def}
 
 \begin{Def} \index{direct sum of $C^*$-algebras} Given real $C^*$-algebras $A_1, A_2, \dots A_n$ we may define their {\emph{direct sum}} $A_1 \oplus \dots \oplus A_n$  just as in the complex setting, by taking 
 their algebraic direct sum, and the evident $ *$ and
 close it in the evident
  norm.   Similarly, we may take arbitrary sums $\oplus _j A_j $ by taking the algebraic direct sum and giving it the appropriate 
 norm,  and closing it just as in the complex setting.      It is not hard to show that the direct sum operation respects complexification. That is, 
 \[
 (A \oplus B)\sc   \cong A\sc \oplus B\sc 
 \]
 and similarly for arbitrary direct sums, as explained below.
 \end{Def}

 \begin{Def}
 \index{direct limit of $C^*$-algebras} 
 \index[notation]{limit@$\Dirlim _j A_j$}
 More generally,  given a sequence 
 \[
 A_1 \longrightarrow A_2 \longrightarrow \dots \longrightarrow A_n \longrightarrow \dots
 \]
 of real $C^*$-algebras, we may form the {\emph{direct limit}}\footnote{Sometimes called the {\emph{injective limit}} } $\,A = \lim_{j \to \infty} A_j$ with canonical maps 
 \[
 \iota _n : A_n \to A \;  .
 \]
  \index{direct limit of $C^*$-algebras}
 We refer to Blackadar (\cite{Blackadar} p.17) for details.    If each $A_j$ is finite-dimensional, then  $\lim_{j \to \infty}A_j$ is what is called an {\emph{AF algebra}}.  These have been studied 
 very deeply, and we refer to T. Giordano \cite{Giordano} for details.
 The direct limit operation allows us to define precisely the sum $\oplus _{j= 1}^\infty  A_j $, namely 
 \[
  \bigoplus _{j= 1}^\infty  A_j = \lim_{j \to \infty} \left[ A_1\oplus A_2 \oplus \dots \oplus A_j \right]. 
 \]
 The direct limit operation also respects complexification. That is, there is a natural isomorphism
 \[
 \lim_{n \to \infty} \left[ (A_n) \sc \right] \overset{\cong}\longrightarrow  \left[ \lim_{n \to \infty}A_n   \right] \sc
 \]
 \end{Def}

 The definition and basic properties of the tensor product of two complex $C^*$-algebras are well-exposed in several sources, notably Blackadar \cite{Blackadar}.  For 
 real $C^*$-algebras we rely upon Boersema \cite{BoersemaKT}.

 \begin{Def}  \index{tensor product of $C^*$-algebras} 
 \index[notation]{AB@$A \otimes\sr B$}
   Suppose that $A$ and $B$ are real $C^*$-algebras.  Then we may form their {\emph{tensor product}} $A \otimes B$ by taking their algebraic tensor product  over $\RR $ and then completing 
 that algebra in a suitable $C^*$-norm.   \footnote{We will use the symbol $\otimes$ when it is clear from context that we are working with real or with complex $C^*$-algebras.  If it 
 is not clear then we shall use $\otimes\sr $ and $\otimes \sc$.}

 The complication here is that in general, just as in the complex case, there is more than one ``suitable" norm on the algebraic tensor product - there is a minimal norm, a maximal norm, and possibly many in between. The minimal norm has the very nice property of being concrete.  The maximal norm has the very nice property of preserving exact sequences. We want both properties and hence {\bf{we make a blanket assumption: whenever we write $A \otimes B$     then at least one of the two $C^*$-algebras has the property that its complexification is \index{nuclear} nuclear.\footnote{For example, any commutative complex $C^*$-algebra is nuclear.  Please see Blackadar \cite{Blackadar} for definitions of nuclear $C^*$-algebras as well as for other operator algebra 
 terms that may seem mysterious.}
 }}
 This does the job, because, using the minimal tensor product norm, 
 \begin{enumerate}
 \item The complexification of a tensor product of two real $C^*$-algebras is isomorphic to the  tensor product of the complexification of the two real $C^*$-algebras, and
 \item A short exact sequence of real $C^*$-algebras is exact if and only if the complexification is exact. 
  \end{enumerate}
    We are particularly indebted to J. Cuntz for clarifying this matter for us. 
  \end{Def}

\noindent
\begin{Pro}
\begin{enumerate}
\item
Let $A_n$ be a directed system of complex $C \sp *$-algebras and let $B$ be a complex $C \sp *$-algebra. 
Then there is an isomorphism
$$ 
\lim_{n \to \infty} (A_n \otimes B) \cong \left( \lim_{n \to \infty}  A_n \right) \otimes B.
$$

\item Similarly, if each $A_n$ is unital and the directed system preserves the identity then there is a natural map 
\[
g: B \longrightarrow  \lim_{n \to \infty} (A_n \otimes B)
\]
defined as a limit of the sequence of maps $B \hookrightarrow A_n \otimes B$ given by $b \mapsto 1\otimes b$.

\item The same results hold in the category of real $C^*$-algebras.

\end{enumerate}

 \end{Pro}
 
\noindent
\begin{proof}
 From the statement for complex $C \sp *$-algebras, we prove the statement for real $C \sp *$-algebras by working in the complexification.  So consider first the complex case.
 
Assume first that $A_n$ is unital for all $n$ and that all of the connecting homomorphisms are unital. Assume also that $B$ is unital. The general case can be easily deduced from this. Let $A =  \lim_{n \to \infty}  A_n$ and let
$$\alpha_n \colon A_n \rightarrow A$$ 
be the natural homomorphism into the limit, which is also unital.

The sequence of maps 
$$\alpha_n \otimes 1 \colon A_n \otimes B \rightarrow A \otimes B$$
is coherent with respect to the inductive sequence $\{ A_n \otimes B\}_n$, and so there is a homomorphism in the limit
$$\theta \colon \lim_{n \to \infty} (A_n \otimes B) \rightarrow A \otimes B \; .$$

For the other direction, there are maps 
$$\beta_n \colon A_n \rightarrow A_n \otimes B$$
given by $\beta_n(a) = a \otimes 1$. 
If we compose $\beta_n$ with the limit map
$$A_n \otimes B \rightarrow \lim_{n \to \infty}  (A_n \otimes B) \; $$
then we have a system of coherent homomorphisms
$$ A_n \rightarrow\lim_{n \to \infty} (A_n \otimes B)  \; $$
which yields the map 
$$f \colon A \rightarrow \lim_{n \to \infty} (A_n \otimes B) $$
in the limit. There is also a map
$$g \colon B \rightarrow \lim_{n \to \infty}(A_n \otimes B) $$
defined by $b \mapsto 1 \otimes b$.
Check that the images of $f$ and $g$ commute in $ \lim_{n \to \infty} (A_n \otimes B) $, and so by the universal property of the maximal tensor product \cite{Mu}  we obtain a homomorphism
$$\psi \colon A \otimes B \rightarrow \lim_{n \to \infty} (A_n \otimes B)  \; $$
which restricts to $f$ and $g$ on $A \otimes 1$ and $1 \otimes B$ respectively.

Finally, check that $\phi$ and $\psi$ are inverses to each other.
\end{proof}

 \begin{Exam}
     Here are some   more examples. First establish the notation 
     $$A^{\otimes n} = A \otimes \dots \otimes A \quad \text{($n$ times.)} $$ 
     
     \begin{enumerate}
 \item 
 \[
 \realo (\topr )^{\otimes n} \,\cong\, \realo(\topr ^n) 
 	  \]
\item
	  \[
	   \iS ^{\otimes n}  \, \cong\,  \real (S^1 \times \dots \times S^1 , \hat\tau ^n )
	  \]
where 
$\hat\tau ^n : (S^1 \times \dots \times S^1) \longrightarrow (S^1 \times \dots \times S^1) $ is given by 
\[
\hat\tau ^n (x_1, \dots , x_n ) =   	(\bar{x}_1,\bar{x}_2, \dots , \bar{x}_n  )   .
 \]
 
\ \item
 If $A$ is a real $C^*$-algebra and the associated conjugate linear  involution on $A\sc $ is $\sigma $, then 
  \[
  \realo (\topr ) \otimes A \cong C(\topr, A)
  \cong \{ f \in C (\topr , A\sc ) \mid f(x) = \sigma (f(x))^* \}.
 \]

 \item
 If $A$ is a real $C^*$-algebra and the associated conjugate linear  involution on $A\sc $ is $\sigma $, then 
  \[
  \realo (i\topr ) \otimes A \cong \{ f \in \complexo (\topr , A\sc ) \mid f(-x) = \sigma (f(x))^* \}.
 \]
 \item  Finally, 
 \[
 \realo (X_1, \tau _1) \otimes  \realo (X_2, \tau _2) \cong \realo (X_1 \times X_2 , \tau _1 \times \tau _2 ). 
 \]

 \end{enumerate}
 \end{Exam}

 The section is titled ``New $C^*$-algebras from Old Ones" but our final example should be titled ``Old $C^*$-algebras from Old Ones." 
 This fact will lead eventually, as we shall see, to a strong connection between real and quaternionic $K$-theory. 
 
 \begin{Pro} \label{H-times-H}
 (cf. \cite{ABS} 4.1)
 There is a natural isomorphism 
 \[
 \HH \otimes\sr \HH \cong M_4(\RR ) .
 \]
 \end{Pro}
 
 \begin{proof}
 We start by considering $\HH$ in terms of the complex $C^*$-algebra with involution $(M_n(\CC), \sharp)$, as in Example~\ref{realforms}.
 Then $\HH \otimes_{\RR} \HH$ corresponds to  
 $$(M_2(\CC) \otimes_{\CC} M_2(\CC), \sharp \otimes \sharp) = (M_4(\CC), \sharp \otimes \sharp)$$
where $\sharp \otimes \sharp$ is the involution on $M_4(\CC)$ defined by
$$\sharp \colon \begin{bmatrix} A& B \\ C & D \end{bmatrix} \mapsto \begin{bmatrix} D^\sharp & -B^\sharp \\ -C^\sharp & A^\sharp \end{bmatrix} \; .$$

We claim that there is an isomorphism $ (M_4(\CC), \rm{Tr}) \cong (M_4(\CC), \sharp \otimes \sharp)$ as complex $C \sp *$-algebras with involution. It will follow from this that there is a real $C \sp *$-algebra isomorphism $M_4(\RR) \cong \HH \otimes_{\RR} \HH $.

To produce this isomorphism, let
$$Q = \tfrac{1}{\sqrt{2}} \begin{bmatrix} 1&0 & 0 & -i \\ 0 & 1 & i & 0 \\ 0 & i & 1 & 0 \\ -i & 0 & 0 & 1 \end{bmatrix} \; .$$
Check that $Q$ is a unitary and defines an isomorphism $\phi \colon M_4(\CC) \rightarrow M_4(\CC)$
by $$\phi(A) = Q A Q^* \; .$$
Then check that 
$$\phi(A^{\rm{Tr}}) = \phi(A)^{\sharp \otimes \sharp}$$
which therefore gives us the desired isomorphism $ (M_4(\CC), \rm{Tr}) \cong (M_4(\CC), \sharp \otimes \sharp)$.
\end{proof}

 \begin{Rem} One of the referees helpfully observed that one can also see this
  by viewing $\HH  $  
  as $\RR ^4$ and letting it act on itself by left and right multiplication, giving two unital representations with commuting images.
 \end{Rem}

 
 \newpage
   \section{\bf{Homotopy Constructions} } \label{Section:HomotopyConstructions}
 
 In this section we present some of the basic homotopy constructions that we shall need in subsequent sections.  This is important philosophically, 
 since it illustrates how much the homology theories $K_*$ and $KO_*$ have in common and opens the door for other possibilities in the future.\footnote{For 
 example, the homology theories $K_*(A; G) $ (complex $K$-theory with coefficients in an abelian group $G$ such as $\ZZ _2$ or $\QQ$ ) have been found to be quite useful in certain contexts. 
 It may turn out that similar constructions starting with $KO_*$ may also be of service. All of these theories are homology theories.
 See Section~\ref{Section-KwithCoefficients}.}
 
 To set up,  one needs first of all to specify a category which will be the domain of the theory.  For $K_*$ we sometimes use, for example, 
 \begin{enumerate}
 \item All complex $C^*$-algebras and their homomorphisms
 \item All separable complex $C^*$-algebras and their homomorphisms
 \item The bootstrap \index{bootstrap category} category\footnote{First defined by Schochet  in \cite{TopII}. This category is an important part of the hypotheses for the K\"unneth Theorem and the Universal Coefficient Theorem }   \index{universal coefficient theorem}
 (a certain subcategory of the previous example, which includes commutative $C^*$-algebras but is much larger), 
 \item All commutative complex $C^*$-algebras and their homomorphisms\footnote{which, by the Gelfand transform, corresponds to 
 all compact spaces and continuous functions}.
 \end{enumerate}
 
 For $KO_*$ we will use a similar list:
  \begin{enumerate}
 \item All real $C^*$-algebras and their homomorphisms
 \item All separable real $C^*$-algebras and their homomorphisms
 \item The real bootstrap \index{bootstrap category!real bootstrap category}category (a certain subcategory of the previous example, analogous to the complex bootstrap category above, which includes all commutative real $C^*$-algebras)
 \item  The pre-bootstrap \index{bootstrap category!pre-bootstrap category} category, which consists of all real $C^*$-algebras $A$ whose complexifications $A\sc$ lie in the complex bootstrap category. This includes the real bootstrap category, by Corollary~4.3 of \cite{BoersemaKT}.
 \item All commutative real $C^*$-algebras and their homomorphisms\footnote{which,  as we have noted, does {\bf{not}}  correspond to 
 all compact spaces and continuous functions. It is larger!}.
 \end{enumerate}

 The second author's paper \cite{TopIII} lays out almost all of the constructions that we need, but restricted to the complex categories specified above, as 
 the concern then was $K_*$.  These constructions carry over immediately to the real category. There is one very striking addition 
 that comes up - you cannot ``desuspend" in the complex world, but in the real world you can.\footnote{This is slightly exaggerated.  If $A$ is a real $C^*$-algebra then 
 we have real $C^*$-algebras $SA$ and $VA$ defined below, and they behave very differently.  In the complex world the analogous $C^*$-algebras are isomorphic as complex 
 $C^*$-algebras.}     Stay tuned! 
 
 We begin by defining some elementary constructions. Let $A$ be a real $C^*$-algebra.

 \begin{Def}  
 \index{cylinder} 
 \index[notation]{IA@$IA$}
 The {\emph{cylinder}} of $A$, denoted $IA$, is $C( [0,1], A)$ with canonical maps 
 $p_t : IA \to A $ given by $p_t(\xi ) = \xi (t)$.

 \end{Def}

 \begin{Def}
 \index{cone} 
 \index[notation]{CA@$CA$}
 The {\emph{cone}}   of $A$, denoted by $CA$, is defined by
 \[
 CA = \{ \xi \in IA \, | \, \xi (0) = 0 \} = \ker ( p_0 : IA \to A)
 \]
 with associated short exact sequence
 \[
 0 \longrightarrow CA \longrightarrow IA \overset{p_0}\longrightarrow A \longrightarrow 0.
 \]
\end{Def}

 \begin{Def} \label{Def-Sus}
 \index{suspension!of a $C^*$-algebra} 
 \index[notation]{SA@$SA$}
The  {\emph{suspension}} of $A$, denoted $SA$,  is defined by 
 
 \[
 SA = \{ \xi \in IA \,|\,  \xi (0) = \xi (1) = 0 \} \cong C_0(\topr, A)      \cong  \realo(\topr)\otimes A \cong S\RR \otimes A
 \] 
 with associated short exact sequence
 \[
 0 \longrightarrow SA \longrightarrow  CA \overset{p_1}\longrightarrow A \longrightarrow 0.
 \]

For each $n = 0,1,2,...$ we define the {\emph{n-th suspension}} of $A$, as 
 \[
 S^nA = C_o(\topr ^n, A) \cong \realo (\topr ^n) \otimes A ,
 \] 
 and we note at once that 
 \[
 S(S^nA) \cong S^{n+1}A .
 \]

  \end{Def}

 For example, if $A = \real (X)$, with $X$  compact, then 
 \[
 S^nA =  \realo (\topr ^n, C(X, \RR)) \cong \realo (\topr ^n \times X)  
  \]
 and 
 \[ 
 (S^nA)^+ \cong \real (S^n \times X) .
 \]
 
Taking $X = x_0$ a point, we have $\realo (x_0) = \RR$, and so 
\[
S^n\RR = \realo(\topr^n) 
\]
and 
\[
(S^n\RR)^+ = \real (S^n).
\]

Next we come finally to desuspension.   A construction similar to ours was first considered 
 by Karoubi (\cite{Knumber27}, \cite{Knumber29}) in the framework of Banach rings that is related to the 
 development of negative algebraic $K$-theory of Bass \cite{Bass} (cf. Weibel 
 \cite{Whistory}).  
 
 \begin{Def} {\label{Def-DeSus}} \index{desuspension}   
 \index[notation]{VA@$VA$}
 The {\emph{desuspension}} of $A$, denoted $VA$, is defined by
 \[
 VA = \realo (i\topr ) \otimes A
 \]
 and the {\emph{n'th desuspension}} of $A$, denoted $V^nA$, is defined by
 \[
 V^nA = V(V^{n-1}A) \,\cong\, \realo (i\topr )^{\otimes n}  \otimes A    \,\cong\, \realo (\topr ^n, \hat{\tau} ^{\otimes n} ) \otimes A.
 \]
 where $\hat\tau (x) = -x$ for $x \in \RR^n$ as in Example~\ref{Examples-RealC*-algebras}, Part (8).
 \end{Def}
 
 We will need a little machinery in order to justify the use of the word ``desuspension".   It will soon become clear.   
  
  \begin{Def} \label{Def-AllSus}
  \index[notation]{snA@$S^n A$}
  We combine the suspension and desuspension notations to write, for all integers $n$,
  \[ S^n A = \begin{cases} S^n A & n \geq 0 \\ V^{-n} A& n < 0 \end{cases}
\]
 \end{Def}
  
 \begin{Def} 
 \index{homotopic} Suppose that $f_0, f_1 : A \to B$ are homomorphisms of real $C^*$-algebras.  We say that $f_0$ is {\emph{homotopic}} 
 to $f_1$ (written $f_0 \simeq f_1 $) if there exists a homomorphism 
 \[
 H : A \longrightarrow  IB
 \]
 such that $p_0H = f_0$ and $p_1 H = f_1$.   We frequently write $f_t(a) = H(a)(t)$ and then think of $\{f_t \}$ as a continuous family of homomorphisms 
 from $f_0$ to $f_1$. The first definition encodes what we mean by ``continuous family."
 \end{Def}
 
This definition of homotopy generalizes the definition of homotopy between compact spaces.   If $f, g : X \to Y$ are basepoint preserving continuous maps, then $f \simeq g$ 
 if and only if the induced maps $f^*, \, g^* : \real (Y) \to \real (X)$ are homotopic 
 as real basepoint preserving maps\footnote{Suppose that \index{basepoint preserving}
 $x_0$ is the basepoint of $X$ and $y_0$ is the basepoint of $Y$. Then we say that a homomorphism $f : \real(Y) \to \real(X) $ is {\emph{basepoint preserving}} if 
 $f( \tau)(x_0) = \tau(y_0)$ for all $\tau \in \real(Y)$. Further, we say that two homomorphisms $f_0, f_1 \colon \real(Y) \to \real(X)$ are basepoint preserving homotopic if there is a homotopy $\{f_t\} $ of   $C^*$-algebra maps that is basepoint preserving for each $t$.        }  
.
 
 We note that homotopy is an equivalence relation, and we denote the set of homotopy classes of maps from $A$ to $B$ by $[A,B]$. 
 
 Recall that a topological space $X$ is {\emph{contractible}} if \index{contractible} the identity map $i: X \to X$ is homotopic to a constant map of the form $f(x) = x_o \,\,\, \forall x \in X$. 
 
 \begin{Def} 
 \index{contractible}
 A real $C^*$-algebra $A$ is {\emph{contractible}} if the identity map $i : A \to A$ is homotopic to the trivial map $A \to A$ sending every element of $A$ 
 to $0$. 
 \end{Def}
 
 The simplest examples of contractible real $C^*$-algebras are those of the form $A = \real_0(X \backslash x_0)$ where $X$ is contractible and $x_0 \in X$. 
 Here is another important example.
 
 \begin{Pro} 
 \index{cone}
 If $A$ is a real $C^*$-algebra then $CA$, the cone of $A$, is contractible. 
 \end{Pro}
 
 \begin{proof} Define $f_s : CA \to CA   $ by
 \[
 f_s(\xi ) (t) = \xi (st) .
 \]
  
 This is well-defined  since $f_s(\xi )(0) = \xi (0) = 0$ for all $s$.  Then it is easy to see that 
 $f_1 \colon CA \to CA$ is the identity map and that $f_0$ is the zero map.
The basepoint of $CA$ is the function $\phi  : [0,1] \to A$ with $\phi (t) = 0 $ for all $t$.   Then 
\[
f_s(\phi) (t) = \phi (st) = 0 \,\,\forall t
\]
and hence the homotopy preserves basepoints.

\end{proof}

 Next we introduce pullbacks.
 
 \begin{Def}
 \index{pullback}
 Suppose given real $C^*$-algebras $A, B, C$ and maps
 \[
 \begin{CD}
 @.     A        \\
 @.  @VfVV \\ 
B @>g>>     C.
\end{CD}
\]
The {\emph{pullback}} $D$ 
 \[
 \begin{CD}
D               @>{\hat g}>>            A        \\
 @V{\hat f}VV                 @VfVV \\ 
B                       @>g>>            C
\end{CD}
\]
is defined by
\[
D = \{ (a, b ) \in A \oplus B \mid f(a) = g(b) \}
\]
and
\[
 \hat f : D \to B  \quad and \quad  \hat g : D \to A 
 \]
 are  the evident projections.

\end{Def}

The notion of ``pullback" occurs in many different settings.\footnote {and different names. Sometimes it is called the fibre product or the Cartesian square, for instance.} Indeed, later in this {\BS} we will need to take the pullback 
of a vector bundle!  It is a categorical definition and enjoys the following universal property.  In the notation of the diagram 
above, suppose that one is given a  real $C^*$-algebra $L$ and maps 
\[
r: L \longrightarrow A \qquad and \qquad s: L \longrightarrow B
\]
such that $fr = gs: L \to C$. Then there exists a unique map $k: L \to D$ such that  
\[
r = \hat g k : L \longrightarrow D \longrightarrow A \qquad and \qquad s = \hat f k: L \longrightarrow D \longrightarrow B .
\]

\begin{Def} 
\index{mapping cone}
The {\emph{mapping cone}} $Cf$  of a map $f: A \to B$ is defined to be the pullback 
 \[
 \begin{CD}
Cf               @>{\pi (f) }>>            A        \\
 @VVV                 @VVV \\ 
CB                    @>{p_1}>>            B
\end{CD}
\]
with associated short exact sequence
\[
0 \longrightarrow SB \longrightarrow Cf \overset{\pi(f)}\longrightarrow A \longrightarrow 0.
\]
 
\end{Def}

\begin{Def} 
\index{mapping cylinder}
The {\emph{mapping cylinder}} $Mf$  of a map $f: A \to B$ is defined to be the pullback 
 \[
 \begin{CD}
Mf               @>{\hat g}>>            IB        \\
 @V{p}VV                 @V{p_1}VV \\ 
A                      @>f>>            B.
\end{CD}
\]
with associated short exact sequence

\[
0 \longrightarrow Cf \longrightarrow Mf \longrightarrow B \longrightarrow 0.
\]
\end{Def}

\begin{Pro} Suppose given the short exact sequence and map  of real $C^*$-algebras 
\[
\begin{CD}
@.           @.             @.            A       \\
@.       @.        @.         @VfVV    \\
0 @>>>  J      @>>>     B       @>g>>  B/J       @>>> 0 .
\end{CD}
\]
Then the pullback
\[
\begin{CD}
@.           @.             D    @>>>        A       \\
@.       @.        @VVV         @VfVV    \\
0 @>>>  J      @>>>     B       @>g>>  B/J       @>>> 0 .
\end{CD}
\]
extends, up to isomorphism, to a commuting diagram of short exact sequences
\[
\begin{CD}
0 @>>>          J  @>>>             D    @>>>        A     @>>> 0  \\
@.       @VV{\id}V        @VVV         @VVfV    \\
0 @>>>  J      @>>>     B       @>g>>  B/J       @>>> 0 .
\end{CD}
\]

\end{Pro}

Here is a very important example. 

\begin{Def} The {\emph{anonymous\footnote{What else do you say when one's identity is taken away?}  real Toeplitz algebra}} 
\index{Toeplitz algebra} 
\index[notation]{TR0@$\mathcal{T}_{\RR 0} $}
$\mathcal{T}_{\RR 0} $ is defined by the pullback diagram
\[
\begin{CD}
\mathcal{T}_{\RR 0}  @>>>    \realo (i\topr ) \\
@VVV        @VVV      \\
     \mathcal{T}_{\RR }   @>>>   \iS
     \end{CD}
\]
and fits into the commuting diagram of short exact sequences
\[
\begin{CD}
0 @>>>       \KK \sr      @>>>          \mathcal{T}_{\RR 0}  @>>>    \realo (i\topr)   @>>>     0 \\
@.        @VV{\id}V              @VVV        @VVV      \\
   0 @>>>       \KK \sr      @>>>      \mathcal{T}_{\RR }   @>>>   \iS      @>>>     0.
     \end{CD}
\]
\end{Def}

\begin{Def} 
\index{homotopy equivalence}
A {\emph{homotopy equivalence}} or {\emph{equivalence}} $f: A \to B$ of real $C^*$-algebras is a map $f$ for which there 
exists a map $g: B \to A $ (called a {\emph{homotopy inverse}}) such that $fg \simeq 1_B$ and $gf \simeq 1_A $.    If there is a homotopy equivalence from $A$ to $B$ then we say that 
the two real $C^*$-algebras are {\emph{homotopy equivalent}} and we write $A \simeq B$.  This is an equivalence relation on real $C^*$-algebras.
\end{Def}

This again is a straight-forward generalization of the notion for spaces.  For example, it is easy to see that the annulus
$$X = \{ (x,y) \in \RR^2 \mid 1 \leq x^2 + y^2 \leq 2 \}$$is homotopy equivalent to a circle.  As $C \sp *$-algebras,  $\real(X)$ and $\real(S^1)$ are homotopy equivalent.

At a much deeper level we have the following result. See Proposition~1.5.1 in Schr\"oder's book \cite{Schr} for a proof, which follows the idea of Cuntz in \cite{Cuntz1984} in the complex case.

 \begin{Thm}  Let $e_1 :    \iS      \longrightarrow  \RR $ be evaluation at $1$ and let
  \[
  \pi : \mathcal{T}\sr \to \RR 
  \]
   be defined as the composition
 \[
  \mathcal{T}\sr \longrightarrow   \iS  \overset{e_1}\longrightarrow  \RR .
 \]
 Then  $\pi : \mathcal{T}\sr \to \RR $ is a homotopy equivalence with homotopy inverse 
  $\iota  :    \RR \to \mathcal{T}\sr   $ the inclusion of the identity.
 Thus there is a split short exact sequence\footnote{A short exact sequence  \index{split exact sequence}
 \[
 0 \to A \to B \overset{p}\to C \to 0 
 \]
 is said to {\emph{split}} if there is a map $i: C \to B$   in the category such that $pi = \id$.  This implies that there is an isomorphism $B \cong A\oplus C$.  Sequences need not 
 split: for example, the short exact sequence of abelian groups 
 \[
 0 \to \ZZ \overset{\times 2}\longrightarrow \ZZ \to \ZZ _2 \to 0
 \]
 doesn't split: $\ZZ$ is NOT isomorphic to $\ZZ \oplus \ZZ _2 $. }
  
 \[
  0 \to     \mathcal{T}_{\RR 0}      \longrightarrow    \mathcal{T}\sr    \overset{\pi}\longrightarrow \RR   \to 0
 \]
 with $\pi $ a homotopy equivalence.  Further, for any real $C^*$-algebra $A$, there is a split short exact sequence
  \[
  0 \to     \mathcal{T}_{\RR 0} \otimes A     \longrightarrow    \mathcal{T}\sr  \otimes A   \overset{\pi \otimes 1}\longrightarrow A   \to 0
 \]
 with $\pi \otimes 1$ a homotopy equivalence. 
  \end{Thm}

  Next we move to cofibrations. For everything to do with this topic, please see the second author's paper \cite{TopIII}. 
  
  \begin{Def} A map $p: A \to B$ of real $C^*$-algebras is a {\emph{cofibration}} if whenever given a diagram 
    \index{cofibration}

  \[
\xymatrix{ & A\ar[d]^p   \\ 
D \ar[ur]^{H_0} \ar[r]^{h_t} & B}
\]

  \noindent with real $C^*$-algebra $D$, a homotopy of maps 
\[
  h_t : D  \longrightarrow B
  \]
  and a map 
  \[
  H_0 : D \longrightarrow A
 \]
  with $pH_0 = h_0$, then there exists a homotopy of maps 
  \[
  H_t : D \to A
  \]
   such that the diagram

   \[
\xymatrix{ & A\ar[d]^p   \\ 
D \ar[ur]^{H_t} \ar[r]^{h_t} & B}
\]
  
   \noindent commutes for $0 \leq t \leq 1$.

\end{Def}
  
  \begin{Rem}  
  The word  ``cofibration" appears  in classical algebraic topology. It has to do with a subspace $Y \subset X$ sitting nicely inside $X$. 
  For example, the circle $S^n$ sits nicely inside the $n+1$ disk $D^{n+1}$ but the set of numbers $\{ 1, 1/2, 1/3, 1/4, \dots \} $  does not sit nicely inside the closed unit interval.

    We introduced the word  ``cofibration" in \cite{TopIII} (1984)  in the context of $C^*$-algebras   
  so that the map 
  \[
  f: X \to Y
  \]
   is a cofibration if and only if the map 
   \[
   f^*:  C(Y) \to C(X) 
   \]
    is a cofibration.   However, several authors have chosen to quote \cite{TopIII} but call this a  ``fibration" (actually  ``Schochet fibration").\footnote{For example, Bunke \cite{Bunke} states the definition in much more elegant terms. He says that a map $A \to B$ is a 
  {\emph{Schochet fibration}} \index{fibration!Schochet fibration} if the map 
  \[
  f_* : \Hom (D, A) \longrightarrow \Hom (D, B)
  \]
  is a Serre fibration of topological spaces for every $C^*$-algebra $D$. If $f: A \to B$ is a Schochet fibration with $J = \ker(f) $ then the associated short exact 
  sequence 
  \[
  0 \longrightarrow J \longrightarrow A \longrightarrow B \longrightarrow 0
  \]
  is a {\emph{Schochet exact sequence.}} }
  We will stick with  ``cofibration."    The results in \cite{TopIII} are presented for complex $C^*$-algebras, but the parts having to do with cofibration apply to the real setting with 
  no serious modification needed, so we will state facts but offer no proofs. 
  \end{Rem}

  Here are some basic examples of cofibrations.  We mention the item and then insert the reference in \cite{TopIII} where the facts are verified.
  
  \begin{enumerate}
 \item 
  The map $p_t : IB \longrightarrow B$ is a cofibration for each $t \in [0,1]$.    (\cite{TopIII} Lemma 1.3)
 \item 
  In general, the map $A^+ \to A^+/A \cong \RR $ is not a cofibration.
 \item  
  If $p: A \to B$ is a cofibration and $f: D \to B$ is a map 
  \[
  \begin{CD}
  C @>>> A \\
  @V{f^*(p)}VV       @VpVV \\
  D @>f>>     B
  \end{CD}
  \]
  then the pullback $f^*(p) : C \to D$  is a cofibration. (\cite{TopIII} Prop. 1.5) 
 \item  Given $f: A \to B$ then the mapping cone
  \[
  Cf \overset{p}\longrightarrow A 
  \]
  is a cofibration.  (Definition 2.1)

  \end{enumerate}

 \newpage
   \section{\bf{Homology and Cohomology Theories on Real and Complex $C^*$-algebras}}
   \label{Section:HomologyTheories}
 
 In this section we give axioms for homology and cohomology theories on $C^*$-algebras.  We also point out some  important and powerful consequences of these axioms. 
 In truth, our main interest at the moment is $KO_*$.   However, we believe that using the axioms makes the properties clearer.  Certainly it makes the parallel with 
 $K_*$ for complex $C^*$-algebras clearer and it makes the discussion of $KO$-theory with coefficients much easier, as we shall see.

 In addition, we introduce axioms for  cofibre homology and cohomology theories for $C^*$-algebras. This is not immediately relevant for $KO_*$-theory but we believe that it belongs 
 in this work as we shall explain.     
 
 When needed, we will specify whether we are talking about a theory defined on a category of real $C^*$-algebras (call them ``real" theories) or on complex $C^*$-algebras (call them ``complex" theories.).  However, in the generality that we work this is seldom necessary.  The proofs for most of the results in this section are found in \cite{TopIII} which deals with complex theories and generalizes without effort to real theories.

 {\bf{Note: Nothing in this section   uses Bott Periodicity. }}
  
  \begin{Def}\label{D:homology}  
  \index{homology theory} A {\emph{homology theory}} is a sequence $\{ h_n \}$ of covariant functors 
  from an admissible category\footnote{Roughly speaking, a category is {\emph{admissible}} if it is closed under the various operations (direct sum, direct limits, 
  quotients, pullbacks, etc.) that we need to make our constructions. For instance, the natural numbers are closed under addition and multiplication but not under 
  subtraction.    For now, take the category to be all real $C^*$-algebras or else all complex $C^*$-algebras. One could also, for example, take the category of all real separable $C^*$-algebras.   Similarly, when working with homology theories on complex $C^*$-algebras we would use the category of complex 
  $C^*$-algebras or the category of complex separable $C^*$-algebras, or something smaller such as the bootstrap category. } 
  ${\mathcal {C}} $ of  real or complex  $C^*$-algebras  to abelian groups which satisfies the following axioms: \index{homology theory}\index{admissible category} 

  {\emph{Homotopy Axiom}}.\index{homotopy axiom}  Let $H : A \to IB$ be a homotopy from $f = p_0 H $ to $g = p_1 H$ 
 in  ${\mathcal {C}} $. Then 
 \[
 f_* = g_* : h_n(A) \longrightarrow h_n(B) \quad \text{for all n}.
 \]

 {\emph{Exactness Axiom.}} \index{exactness axiom}  Let
 \[
 0 \longrightarrow J \overset{i}\longrightarrow A \overset{j}\longrightarrow A/J  \longrightarrow 0
 \]
 be a short exact sequence in  ${\mathcal {C}} $.  Then there is a map 
 \[
 \partial : h_n(A/J) \longrightarrow h_{n-1}(J)
 \]
 and a long exact sequence 
 \[
\dots \to h_n(J) \overset{i_*}\longrightarrow h_n(A)  \overset{j_*}\longrightarrow h_n(A/J )  \overset{\partial }\longrightarrow h_{n-1}(J ) \to \cdots 
 \]
  The map $\partial $ is natural with respect to morphisms of short exact sequences.

  The axioms imply that if $A \cong \oplus _{i= 1}^k A_i$ in $\mathcal {C}$ with $k$ finite  then the natural map 
  \[ 
   \oplus _{i = 1}^k h_n(A_i)  \longrightarrow h_n(A)
   \]
   is an isomorphism. 
  A further axiom is sometimes assumed. If it holds then the homology theory is said to be {\emph{additive}}. \index{homology theory!additive homology theory}

  {\emph{Additivity Axiom}}\index{additivity axiom}  Let $A = \oplus _{i = 1}^\infty A_i$ in ${\mathcal{C}}$.  Then the natural maps $A_j \to A$ induce 
  an isomorphism 
  \[
  \oplus _i h_n(A_i) \overset{\cong}\longrightarrow h_n(A) .  \]
 for each $n$. 
 
 \end{Def}
 
 Henceforth, when we speak of a homology theory (or a cohomology theory) it is understood to be defined on an admissible category of real or of complex 
 $C^*$-algebras.
 
 \begin{Thm} (\cite{TopIII}, Theorem 5.1) \index{direct limit of $C^*$-algebras} \index{direct limit of $C^*$-algebras} Suppose that $h_*$ is an additive homology theory  and we are given a sequence of $C^*$-algebras 
 \[
 A_1 \longrightarrow A_2 \longrightarrow \dots \longrightarrow A_n \longrightarrow \dots
 \]
 with $\lim_{j \to \infty} A_j = A$.  Then the natural maps $h_*(A_j) \to h_*(A) $ induce an isomorphism
 \[
\lim_{j \to \infty}h_*(A_j) \overset{\cong}\longrightarrow h_*(A) . 
 \]
 \end{Thm}
 
 \begin{Cor}
 Suppose that $h_*$ is an additive homology theory and we are given  a sequence of  $C^*$-algebras 
 \[
 A_1 \longrightarrow A_2 \longrightarrow \dots \longrightarrow A_n \longrightarrow \dots
 \]
 with $\lim_{j \to \infty} A_j = A$.  In addition suppose given a  $C^*$-algebra $B$  with either $B$ nuclear or $\otimes$  the maximal tensor product. 
  Then the natural maps 
  \[
  h_*(A_j\otimes B) \to h_*(A\otimes B) 
  \]
   induce an isomorphism
 \[
\lim_{j \to \infty} h_*(A_j \otimes B ) \overset{\cong}\longrightarrow h_*(A\otimes B) . 
 \]
 \end{Cor}

 \begin{Def} 
 \index{h@$h$-contractible}
 Suppose that $A$ is a $C^*$ -algebra and $h_*$ is a homology theory for which $h_*(A)$ is defined. 
   We say that $A$ is {\emph{$h$-contractible}} if $h_n(A) = 0$ for all $n$. \index{h@$h$-contractible}
 \end{Def}
 
 For example, the cone $CA$ of a $C^*$-algebra   is $h$-contractible for any homology theory $h$.

 \begin{Pro} Suppose that $h_*$ is a homology theory on $C^*$-algebras.   Then 
  the boundary map    induces an isomorphism
 \[
 \partial : h_n(A)\overset{\cong}\longrightarrow h_{n-1}(SA)     
  \index{suspension!of a $C^*$-algebra} 
 \]
  for all $n$.
  \end{Pro}
  
  \begin{proof}
  There is a short exact sequence 
  \[
  0 \longrightarrow SA \longrightarrow CA  \longrightarrow A \longrightarrow 0
  \]
   with associated homology long exact sequence 
  \[
  \dots \rightarrow h_n(CA) \longrightarrow h_n(A)  \overset{\partial}\longrightarrow h_{n-1}(SA) \longrightarrow h_{n-1}(CA) \rightarrow\dots     .
  \]
  Since $h_*(CA) = 0$, the long exact sequence degenerates to   short sequences of the form
  \[
  \dots \rightarrow 0 \longrightarrow h_n(A)  \overset{\partial}\longrightarrow h_{n-1}(SA) \longrightarrow 0        \rightarrow\dots
  \]
  which implies the result       \end{proof}
  
  Recall:
\begin{itemize}
\item 
  $\mathcal{T} \sr $ is the  real Toeplitz algebra   \index{Toeplitz algebra}
  \item 
  $\mathcal{T} _{\RR 0} $ is the  real anonymous Toeplitz algebra  \index{Toeplitz algebra!anonymous Toepliz algebra}
  \item
  $VA = \realo(i\topr) \otimes A$  is the desuspension  of  a real $C^*$-algebra $A$ \index{desuspension}
 \end{itemize}

  \begin{Thm}   Suppose that $h_*$ is a  real homology theory  and let $A$ be an arbitrary real $C^*$-algebra. Then 
  \begin{enumerate}
  
  \item The homotopy equivalence $\pi : \mathcal{T}\sr  \xrightarrow{~\cong~} \RR$ extends to a homotopy equivalence 
  \[
  \pi \otimes 1 : \mathcal{T}\sr \otimes A    \xrightarrow{~\cong~}   A .
  \]

  \item There is a natural isomorphism
  \[
  (\pi \otimes 1)_*  :  h_*(\mathcal{T}\sr \otimes A ) \overset{\cong }\longrightarrow h_*( A).
  \]

  \item  The real $C^*$-algebra   $\mathcal{T}_{\RR 0} \otimes A $ is $h_*$-contractible.

\item  We have shown in the previous section that there is a short exact sequence of the form

\[
0 \longrightarrow  \KK\sr \longrightarrow   \mathcal{T}_{\RR 0 }  \longrightarrow  V\RR     \longrightarrow 0 . 
  \]
This produces, for any $A$, the short exact sequence 
\[
0 \longrightarrow A \otimes \KK\sr \longrightarrow   \mathcal{T}_{\RR 0 }\otimes A  \longrightarrow VA \longrightarrow 0.
  \]
 The long exact homology sequence associated to this short exact sequence
degenerates to isomorphisms
\[
\partial : h_n(VA) \overset{\cong}\longrightarrow h_{n-1}(A) .
  \]
  \end{enumerate}
  \end{Thm}
  
  \begin{proof}
  Part (1) is immediate from the definitions, and of course, a homotopy equivalence of $C^*$-algebras induces an isomorphism on homology, establishing Part (2). 
  
  For Part (3) we recall that there is an exact sequence 
\[  
  0 \to \mathcal{T}_{\RR 0}  \longrightarrow \mathcal{T}_{\RR }  \overset{\pi}\longrightarrow \RR \to 0  
  \]
  and when we tensor with $A$ the resulting sequence 
  \[  
  0 \to \mathcal{T}_{\RR 0} \otimes A \longrightarrow \mathcal{T}_{\RR } \otimes A  \overset{\pi\otimes 1}\longrightarrow A \to 0  
  \]
  remains exact.    
  The associated homology long exact sequence looks in part like   
  \[
 h_{n+1}( \mathcal{T}_{\RR } \otimes A )      \overset{\cong}\longrightarrow  h_{n+1}(A)  \overset{\partial}\longrightarrow   h_{n}(\mathcal{T}_{\RR 0} \otimes A)
 \to  h_{n}( \mathcal{T}_{\RR } \otimes A )       \overset{\cong}\longrightarrow      h_{n}(A)   
  \]
    and this implies\footnote{This is a matter of pure algebra: suppose given an exact sequence of abelian groups
    \[
    G_1 \overset{\alpha _1}\longrightarrow 
     G_2 \overset{\alpha _2}\longrightarrow 
 G_3 \overset{\alpha _3}\longrightarrow 
 G_4 \overset{\alpha _4}\longrightarrow  G_5
 \]
 with $\alpha _1 $ and $\alpha _4 $ isomorphisms.    Then 
 \[
 \ker (\alpha _2 ) \cong \Image( \alpha _1 ) = G_2
 \]
 since $\alpha _1$ is surjective, 
 which says that $\alpha _2 \equiv 0$.   Similarly, 
 \[
\Image(\alpha _3 ) = \ker(\alpha _4 ) = 0
 \]
  since $\alpha _4 $ is injective, which says that $\alpha _3 \equiv 0$.  Then $\alpha _2 = \alpha_ 3 \equiv 0$ implies that $G_3 = 0$.
  } 
  that   $h_{n}(\mathcal{T}_{\RR 0} \otimes A) = 0$ for all $n$, establishing Part (3). 
    
    For Part (4) we simply write out the associated homology long exact sequence and then observe that, by Part (3), every third term vanishes. We are thus left with 
    \[
0 \longrightarrow   h_n(VA) \overset{\partial}\longrightarrow h_{n-1}(A) \longrightarrow 0
  \]
     and this implies that $\partial $ is an isomorphism in each dimension.

  \end{proof}
  
  The following corollary describes both suspension, which also is available in the complex $K$-theory world, and desuspension, which is totally 
  hidden in the complex $K$-theory world, since 
  \[
  VA \otimes \CC \cong SA\otimes \CC = S(A\otimes \CC),
  \]
  a fact related to Bott Periodicity.
   
 Here's a simple way to emphasize the relationship between suspension and desuspension.
 
  \begin{Cor}\label{T: desuspension} 
Let  $h_*$ be a  homology theory. Then for each $C^*$-algebra $A$ in the category:
\begin{enumerate}
\item If $h_*$ is a {\bf{real or a complex}} homology theory, then there is a natural {\bf{suspension}} isomorphism 
  \[
  h_n(A) \cong h_{n-1}(SA)
  \]
  induced from a natural short exact sequence
  \[
  0 \to SA \longrightarrow B \longrightarrow A \to 0 \index{suspension!isomorphism}
  \]
  where $B$ is a contractible algebra.
  \item If $h_*$ is a {\bf{real}}  homology theory, then there is 
  a natural {\bf{desuspension}} isomorphism 
  \[
  h_n(A) \cong h_{n+1}(VA)
    \]
  induced from a natural short exact sequence 
  \[
  0 \to A \otimes \KK \sr \longrightarrow B \longrightarrow VA \to 0
  \]
   where again $B$ is a contractible algebra.
  \end{enumerate}
  \end{Cor}

  Next we turn our attention to the creation of new homology theories.  This technique has been used before, for example in \cite{TopIV} to build $K_*(A ; \ZZ _n)$, complex 
  $K$-theory with coefficients $\ZZ _n$ from $K_*(-)$.

  \begin{Pro}  \label{Proposition:homologytheory}
  \begin{enumerate}
  \item  Suppose that $h_*$ is a real homology theory and $N$ is a  real $C^*$-algebra such that $N\sc $ is nuclear.  Then $h_*( - \otimes N) $ is a real homology 
  theory as well.   If $h_*$ is additive then so is $h_*( - \otimes N) $.
  
  \item Suppose that $h_*$ is a complex homology theory and $N$ is a  complex nuclear $C^*$-algebra  Then $h_*( - \otimes N) $ is a complex homology 
  theory as well.   If $h_*$ is additive then so is $h_*( - \otimes N) $.
  \end{enumerate}
  
  \end{Pro} 
  
  \begin{proof}  We check the real case; the complex case is simpler.   We need to check the homotopy and exactness axioms. 
  
  Suppose first that $H : A \to IB $ is a homotopy from $f = p_0H $ to $g = p_1H$.  Then there is a natural continuous map 
  \[
  H_N : A\otimes N \longrightarrow I(B\otimes N)
  \]
   and this implies that $f\otimes I_N $ is homotopic to $g \otimes I_N$.  Hence 
   \[
   (f \otimes I_N)_* =  (g \otimes I_N)_*  : h_*(A \otimes N) \longrightarrow h_*(B\otimes N)
   \]
   as required. So the homotopy axiom does hold.
   
   Exactness is easy.  Suppose given a short exact sequence
   \[
   0 \to J \longrightarrow A \longrightarrow A/J \to 0 .
   \]
   Tensor the sequence with $N$ and the resulting sequence 
   \[
   0 \to J\otimes N \longrightarrow A\otimes N  \longrightarrow (A/J)\otimes N  \to 0 
   \]
   is still exact, since $N\sc $ is exact. Apply $h_*$ and obtain a long exact sequence.   
   
   Finally, suppose that $h_*$ is additive and given $A = \oplus _1^\infty A_j $.  Then the natural maps $A_j \to A$ yield maps $A_j\otimes N \to A\otimes N$ and 
   hence to 
   \[
   h_*(A_j \otimes N) \longrightarrow h_*(A \otimes N)
   \]
   so the additivity of $h_*$ implies that there is a natural isomorphism
   \[
   \oplus _j h_n(A_j\otimes N)  \overset{\cong}\longrightarrow h_n(A \otimes N).
   \]

   \end{proof}

  \begin{Cor}\label{C:Hhomology}
  If $h_*$ is a real homology theory and $\HH $ is the quaternions, then   $h_*(- \otimes \HH )$ is a real homology theory as well.
  \end{Cor}
  This will be of interest when we take $h_* = KO_*$.

  Next we discuss cofibre homology theories.\index{cofibre homology theory}  These are thoroughly discussed in the complex setting in \cite{TopIII} and so, as before, we simply give specific references
  to that paper. 
  
  \begin{Def}  (\cite{TopIII} p. 437) 
  \index{homology theory!cofibre homology theory}
   A {\emph{cofibre (real or complex) homology theory}}  is a sequence $\{h_n\}$ of covariant functors from an admissible category $\mathcal C$ of
    $C^*$-algebras
  to abelian groups which satisfies the following axioms:
  
  {\emph{Homotopy Axiom}}. \index{homotopy axiom} (same as for homology theories)
  
  {\emph{Suspension Axiom}}. \index{suspension!axiom} For all $A$, the natural suspension map 
  \[
  \sigma _A ; h_n(A) \longrightarrow h_{n-1}(SA)
  \]
  is an isomorphism. 
  
  {\emph{Cofibre Exactness Axiom}}. \index{exactness axiom!cofibre exactness axiom}  If $f: A \to B$ is a cofibration with mapping cone $Cf$ then for all $n$ the sequence
  \[
  h_n(Cf) \longrightarrow h_n(A) \longrightarrow h_n(B)
  \]
  is exact.
  
  \end{Def}

  The exactness and homotopy axioms of a homology theory imply the cofibration exactness axiom (\cite{TopIII}, Prop. 3.7), and hence any homology theory is a cofibre homology theory.       The converse is false: see \cite{TopIII}, Remark 8.6. 
  
  Here are two basic properties of cofibre homology theories. 
  
  \begin{Pro} (\cite{TopIII}, 3.7 and 3.8)    Suppose that $h_*$ is a cofibre homology theory on an admissible category. Then:
  \begin{enumerate}
  \item If $f: A \to B$ is a cofibration with $J = \ker(f)$ then there is a natural long exact sequence
  \[
  \cdots h_n(J) \longrightarrow h_n(A) \overset{f_*}\longrightarrow h_n(B) \longrightarrow h_{n-1}(J) \cdots
  \]
 \item if $A_1, \dots A_n$ are in the category then for each $k$ the natural map
 \[
 \oplus _{j=1}^n  h_k(A_j) \longrightarrow h_k(A_1 \oplus \dots \oplus A_n )
 \]
 is an isomorphism.
 \end{enumerate}
 \end{Pro}
 
 \begin{Rem}   (\cite{TopIII} 8.5)   Here is a family of examples of  cofibre homology theories which are generally {\emph{not}} homology theories. All $C^*$-algebras are in an admissible category.
  Fix some $C^*$-algebra $D$.  Define
  \[
  D_n(A) = \lim_{k \to \infty}   [S^{k-n}D, S^kA ].
  \]
 Then $D_*(A)$ is a cofibre homology theory. 
 \end{Rem}

  \begin{Rem}  Why discuss cofibre homology theories if they seem unrelated to $KO$-theory? The basic point to understand is that, in classical algebraic topology, 
  most homology and cohomology theories are in fact only cofibre theories. Topological $K$ and $KO$-theories are actually quite exceptional. We anticipate that in years to come 
  that there will arise other theories to help us understand real $C^*$-algebras. Given the precedent in classical algebraic topology, there is no reason to suppose that these 
  theories will satisfy the axioms that we have stated for homology theories, but they {\emph{will}} be cofibre theories.  It is already the case that cofibrations\footnote{known there as ``Schochet fibrations" as discussed in Section~\ref{Section:HomotopyConstructions} .}  have appeared in 
  several papers, for instance in Bunke \cite{Bunke}.  
    
  \end{Rem}

  Next we turn to cohomology theories.  These are common and well-known in classical algebraic topology but not so common (yet) in $C^*$-algebra theory.  
  
  Here are the axioms for a cohomology theory, taken from \cite{TopIII}.

  \begin{Def} 
  \index{cohomology theory} 
  A {\it{(real or complex) cohomology theory}}  is a sequence  $\{ h^n \}$  of 
  contravariant\footnote{A functor $F$ is {\it{contravariant}} if applying it to a map $f: A \to B$ 
  in the category produces a map $f^* : F(B) \to F(A).$ }  functors from an admissible category of (real or complex, respectively) $C^*$-algebras to abelian groups which satisfies the following
  axioms:
  
  {\it{Homotopy axiom}}. \index{homotopy axiom}  Let $h : A \to IB$ be a homotopy from $f_0 = p_0h$ to $f_1 = p_1h$.  Then for each $n$, 
  \[
  f_0^* = f_1^* : h^n(B) \longrightarrow h^n(A).
  \]
  
  {\it{Exactness axiom.}}. \index{suspension!axiom} Let
  \[
  0 \longrightarrow J \overset{i}\longrightarrow A \overset{j}\longrightarrow B \longrightarrow 0
  \]
  be a short exact sequence in the category.  Then for each $n$ there is 
  a map 
  \[
  \delta : h^n(J)  \to h^{n+1}(B) 
  \]
   and a long exact sequence
  \[
  \dots \longrightarrow h^n(B) \overset{j^*}\longrightarrow         h^n(A) \overset{i^*}\longrightarrow        h^n(J) \overset{\delta ^*}\longrightarrow h^{n+1}(B)  \longrightarrow \dots
  \]
  The map $\delta $ is natural with respect to morphisms of short exact sequences.
  \end{Def}

  The following axiom is sometimes also assumed. If so then the cohomology theory is said to be {\it{additive}}.\index{cohomology theory!additive cohomology theory}
  
  {\it{ Additivity axiom }} \index{additivity axiom}  Let $A \bigoplus _{i = 0} ^\infty A_i $. Then the natural maps 
  \[
  h^n(A) \longrightarrow h^n(A_i)
  \]
  induce an isomorphism 
  \[
  h^n(A)  \overset{\cong}\longrightarrow \prod _i h^n(A_i) .
  \]

  The example of cohomology theory that has been most studied is constructed as follows. Let $KK_*( -, -)$ be Kasparov's bivariant theory, which we will discuss briefly 
  in Section~\ref{KKSection}.
  \index{Kasparov} \index{KK-theory@$KK$-theory}
   Kasparov shows  under quite general circumstances,   that if  $A$ a fixed $C^*$-algebra then $KK(A, - )$ is a homology theory.  In particular, 
  $KK_*(\CC , B) \cong K_*(B)$.  
  Similarly, for fixed $B$, that $KK(-, B)$ is a cohomology 
  theory.  
  
  Just as we have both homology and cofibre homology theories, we have both cohomology and cofibre cohomology theories. 
  
  \begin{Def}
    \index{cohomology theory!cofibre cohomology theory}
  A {\it{(real or complex) cofibre cohomology theory}} is a sequence of contravariant functors from an admissible category of $C^*$-algebras to abelian groups which 
  satisfy the homotopy axiom and the following axioms: 
  
  {\it{Cofibre axiom.}} \index{cofibre axiom} If $p : A \to B$ is a cofibration then the sequence 
  \[
  h^n(B) \overset{p^*}\longrightarrow h^n(A) \overset{\pi (p)^* }\longrightarrow h^n(Cp)
  \]
  is exact for each $n$. 
  
  {\it{Suspension axiom}}. \index{suspension!axiom} There is a natural isomorphism 
  \[
  \sigma ^A :  h^n(SA)  \xrightarrow{~\cong~}  h^{n+1}(A) .  
  \]
  
  \end{Def}

  \begin{Pro}  \cite{TopIII}. 
  A homology theory {\it{is}} a cofibre homology theory and 
  a cohomology theory {\it{is}} a cofibre cohomology theory. 
  
  \end{Pro}

  There is rather a dramatic difference between homology and cohomology theories in their dealing with direct limits.  Here is the cohomology version.

 \begin{Thm} (\cite{TopIII}, Theorem 7.1) Suppose that $h^*$ is an additive cohomology theory  and we are given a sequence of $C^*$-algebras 
 \[
 A_1 \longrightarrow A_2 \longrightarrow \dots \longrightarrow A_n \longrightarrow \dots
 \]
 with $\lim_{j \to \infty}  A_j = A$.  Then the natural maps $h^n(A_j) \to h^n(A_j) $ induce a map 
 \[
 h^n(A) \longrightarrow \Invlim h^n(A_j)
 \]
   which fits 
 into a natural long exact sequence\footnote{ Recall that if 
 \[
 G_1 \overset{\alpha _1}\longleftarrow G_2 \overset{\alpha _2}\longleftarrow  G_3 \overset{\alpha _3}\longleftarrow \cdots
\]
is an inverse sequence of abelian groups and we define $\Psi : \prod G_i \to \prod G_i $ by 
\[
\psi (g_i) = (g_i - \alpha _i g_{i+1} )
\]
then  $\ker \Psi = \Invlim G_i   $ and $\rm{coker}\, \Psi = \Invlim^1 G_i $ \index{lim@$\lim^1$}.  To give you a feeling of what might happen, note that if each $G_i$ is finite or countable, then 
$\Invlim^1 G_i   $ is either zero or uncountable, a fact due to B. Gray.
For further information on $\Invlim^1$, cf.  May \cite{May}, p. 146 and Schochet \cite{Pext}.}
 \[
  0 \longrightarrow    {\Invlim}^1 h^{n-1}(A_j)         \longrightarrow           h^n(A) \longrightarrow \Invlim h^n(A_j)   \longrightarrow 0.
 \]
 \end{Thm}

  
  \newpage

    \section{\bf{$KO_*$-theory for Real $C^*$-Algebras: Definitions and Basic Properties}}
  \label{Section:KODefinition}
  
    {\bf{  Note: Nothing in this section   uses Bott Periodicity in any manner. }}

  Sometimes matters of taste intervene in mathematical research.  We are at that point.  $K$-theory can be developed in terms of projections or in terms of projective modules. Analysts tend to like projections and algebraists and topologists like projective modules\index{projective module}.  Here's how the connection goes. 
  
  Let $A$ be a unital ring (for example, a unital real $C^*$-algebra).  In the following, we choose to work with left $A$-modules throughout: one must choose left or right at the outset.
   A {\emph{finitely generated free $A$-module}} \index{finitely generated free module} 
   is a left $A$-module of the form $A \oplus \dots \oplus A$, with a finite number (say $k$) of copies of $A$.  Think of this as the identity function 
  \[
    A \oplus \dots \oplus A  \overset{1_k}\longrightarrow   A \oplus \dots \oplus A
    \]
    which in turn we may regard as the identity matrix $ 1_k \in M_k(A) $.    Conversely, if we start with $1_k \in M_k(A)$, regard it as a function and take its range, then we obtain $A \oplus \dots \oplus A$.   Along the way we have chosen a basis.

  \index{projection}  A {\emph{projection}} $p$ is an element of $M_k(A)$ with $p = p^* = p^2$.    It too induces a map 
  \[
    A^k \overset{p}\longrightarrow   A^k
    \]
  but in general the map will not be injective or surjective. Let us denote by $P$ the image of the projection.  Then $P$ is a finitely generated $A$-module. The orthogonal 
  complement of $p$ (call it $q$) has image that is orthogonal to $P$ and $P \oplus Q \cong A^k$.  This is exactly what is required to say that $P$ is a {\emph{finitely generated projective $A$-module. }} Again, the process may be reversed.
  
 The upshot is that projections in finite-dimensional matrix algebras over $A$ are exactly the same as finitely generated projective $A$-modules. When developing $K$-theory in terms of projective modules, \index{projective module} we regard two abstractly isomorphic finitely generated projective $A$-modules as equivalent. When developing $K$-theory in terms of projections (elements in $A$ or in $M_n(A)$) we also need to consider certain equivalence classes.
 
 \begin{Def}
   \index{projection!unitary equivalence}
     \index{projection!Murray-von Neumann equivalent}
  \index{projection!homotopic}

 Suppose that $A$ is a unital real $C^*$-algebra. Then we define three equivalence relations on the set of projections of $A$ as follows. Let $p$ and $q$ be projections 
 in $A$. 
 \begin{enumerate}
 \item $p$ and $q$ are {\emph{unitarily equivalent}}, written $p \equiv_u q$, if there exists a unitary\footnote{A unitary is an element which satisfies $uu^* = u^*u = 1$.} $u \in A$ such that $q = upu^*$.
\item $p$ and $q$ are {\emph{Murray-von Neumann equivalent}}, written $p \equiv_{m} q$, if there exists a partial isometry\footnote{A partial isometry is an element  $a$ for which $a ^*a$ is a projection.} $a \in A$ such that $p = a a^*$ and $a^* a = q$.
\item $p$ and $q$ are {\emph{homotopic}}, written $p \equiv_h q$, if there exists a continuous path $p_t$ of projections in $A$ such that $p_0 = p$ and $p_1 = q$.
\end{enumerate}
 \end{Def}
 
 \begin{Pro} \label{equivalence}
 \index{projection! semigroup of}
For any unital real $C^*$-algebra $A$, the relations $\equiv_u$, $\equiv_{m}$, and $\equiv_h$ are each equivalence relations. Furthermore, they are related by the following implications for any projections $p,q \in A$:
\begin{itemize}
\item if $p \equiv_u q$ then $p \equiv_{m} q $,
\item if $ p \equiv_{h} q$ then $p \equiv_{u} q$,
\item if $p \equiv_m q$ then $\begin{bmatrix} p & 0 \\ 0 & 0 \end{bmatrix} \equiv_u \begin{bmatrix} q & 0 \\ 0 & 0 \end{bmatrix}$ in $M_2(A)$,
\item if $p \equiv_u q$ then $\begin{bmatrix} p & 0 \\ 0 & 0 \end{bmatrix} \equiv_h \begin{bmatrix} q & 0 \\ 0 & 0 \end{bmatrix}$ in $M_2(A)$,
\end{itemize}
 \end{Pro}
See Section~4.2 of \cite{Blackadar} or Section 2.2 of \cite{RLL} for details. The definition of $KO_0(A)$ in what follows is based on unitary equivalence, but because of Proposition~\ref{equivalence} we could use either of the other equivalence relations with the same result.
 
 \begin{Def} Suppose that $A$ is a unital real $C^*$-algebra.  
 \index[notation]{projections1@$\mathcal{P}(A)$}
 \index[notation]{projections2@$\mathcal{P}_n(A)$}
 \index[notation]{projections3@$\mathcal{P}_\infty(A)$}
 Let $\mathcal{P}(A) $ denote the set of unitary 
 equivalence classes of projections in $A$. 
  Let $M_n(A) = A \otimes M_n(\RR) $ and 
 $\mathcal{P}_n(A) = \mathcal{P}(M_n(A))$. 
 The inclusion $M \to M\oplus 0$ of $M_n(A) $ into $M_{n+1}(A)$ induces inclusions 
 \[
\mathcal{P}_n(A)    \subset     \mathcal{P}_{n+1}(A).
 \]
Let  $ \mathcal{P}_\infty (A) $ be the union of 
  the $\mathcal{P}_n(A)$.  
  \end{Def}

   Note that if 
  $p \in M_r(A) $ 
  and $q \in M_s(A) $ then we can  
   form $p\oplus q \in M_{r+s}(A)$ by putting $p$ and $q$ into the upper left and lower right 
  corners of a square $r+s$ matrix respectively, with $0$'s in the upper right and lower left.  We do this even if $p = q$.  
  This gives  $ \mathcal{P}_\infty (A) $ the structure of an abelian semigroup.\footnote{An abelian semigroup \index{abelian semigroup}  is a set with a binary operation that we denote $+$ that is associative and commutative. There are bad classic jokes connected with this concept.  {\emph{What is purple and commutes? An abelian grape.  What is lavender and commutes? An 
  abelian semigrape.}}} 
  
  If $H$ is an abelian semigroup then there is an associated abelian group $\rm{Groth}(H)$, the \index{Grothendieck group} Grothendieck group of $H$. It may be defined as the quotient of $H \times H$ 
  under the equivalence relation 
  \[
  (h_1, h_2) \sim  (k_1, k_2) \quad \text{whenever  $\exists z \in H$ with} \quad h_1 + k_2 + z = h_2 + k_1 + z
  \]
  This can also be understood as taking formal differences of elements of $H$.  The classic example is the natural numbers $\NN$ with $\rm{Groth}(\NN) = \ZZ$.  In terms of formal 
  differences, we think of $-3$ as $[4] - [7] $. Of course we want that to be the same as $[5] - [8]$ and the equivalence relation does just that, 
  since $4 + 8 = 7 + 5$.
  
  The Grothendieck construction is natural and is a covariant functor from abelian semigroups to abelian groups.  If $H$ is an abelian semiring (all the ring axioms except 
  elements may not have additive inverses) then $\rm{Groth}(H)$ is a commutative ring. We will see this happen when we restrict attention to $ \mathcal{P}_\infty (C\pr (X))$. 
  
 \begin{Def} 
 \index{K-theory@$K$-theory!definition}
 \index[notation]{KOA@$KO_*(A)$}
 Suppose that $A$ is a real unital $C^*$-algebra. Define the real $K$-theory group of $A$ by
 \[
  KO_0(A)  = \rm{Groth} (\mathcal{P}_\infty (A) ).  
\]
\index[notation]{GrothP@$\rm{Groth} (\mathcal{P}_\infty (A) )$}
 Thus an element of $KO_0(A)$ 
 may be represented
  as a formal difference $[p] - [q]$, where 
  $p$ and $q$ are projections in matrix rings over $A$, 
 with the following relations:
 \begin{enumerate}
 \item  For any orthogonal matrix $u$, 
 \[
 [p] = [upu^{-1}] .
 \]
 \item For any projection $p$, 
 \[
 [p] + [0] = [p].
 \]
 \item
 \[
 [p] + [q] = [ p\oplus q ]. 
 \]
 \end{enumerate}

Suppose that $f: A \to B$ is a homomorphism of real $C^*$-algebras. If $p$ is a projection in $M_n(A)$ then $f(p)$ is a projection in $M_n(B)$.  The map 
$[p] \mapsto [f(p)] $ is well-defined on equivalences classes and respects the operators, so there is an induced
 homomorphism 
  \[
  f_* : KO_0(A) \longrightarrow KO_0 (B) 
 \]
  given on generators by
 \[
 f_*([p] ) = [f(p) ].
 \]
 It is easy to see that
 \begin{enumerate}
 \item
  If $i :A \to A$ is the identity map then $i_* : KO_0(A) \to KO_0(A)$ is the identity map.  
  \item
   Given real $C^*$-algebra maps
 \[
 A \overset{f}\longrightarrow B \overset{g} \longrightarrow C ,
 \]
then
  \[
  (gf)_* = g_*f_*  : KO_0(A) \longrightarrow KO_0(C) .
  \]
  \item If $f_0, f_1 : A \to B$  are homotopic, then 
  \[
  (f_0)_* = (f_1)_* : KO_0(A) \longrightarrow KO_0(B) .
  \]
   \end{enumerate}
 Thus $KO_0 $ is a homotopy-invariant covariant functor from unital $C^*$-algebras and maps to abelian groups.   \end{Def}

 \begin{Rem} It is easy to be tempted to think that one can use a similar procedure to produce a ring structure
  on $KO_0(A)$. The problem is the following: if $P$ and $Q$ are left modules over a non-commutative ring $A$, 
 then $P \otimes _A Q $ is an abelian group but not necessarily an $A$- module.  Thus 
 $\rm{Groth} (\mathcal{P}_\infty (A) )$ does not have a product and
 the construction grinds to a halt.   If $A$ is commutative then the product does make sense and we discuss this in 
 future sections.

 \end{Rem}
   
 We wish to make the relationship between $KO_0(A)$ and $\Kalg_0(A)$ clear from the beginning, and so we take a short detour. 
 For basic information on algebraic $K$-theory, see Weibel \cite{Wbook}.
 
 \begin{Def} 
\index{K-theory@$K$-theory!algebraic}
\index[notation]{KalgA@$\Kalg 0(A)$}
 Suppose that $A$ is a unital ring. Then we define the {\emph{algebraic $K$-theory group }} 
 $\Kalg_0(A) $ 
 by forming the 
 category 
$\mathcal{P}_\infty (A)$  
 of isomorphism classes of  finitely generated projective left $A$-modules, or, equivalently, the 
 category of equivalence classes of projections in matrix rings over $A$,  and then defining
 \[
 \Kalg_0(A)  = \rm{Groth} (\mathcal{P}_\infty (A) ).  
 \]
     \end{Def}

 Note that this definition makes sense for any unital ring.  There is a $\Kalg_1(A)$ and (by high-tech methods) $\Kalg_j(A)$ for all $j$, but these
 functors are NOT periodic in general.

   If $A$ is a real $C^*$-algebra then  
   \[
   KO_0(A) \cong \Kalg_0(A)
   \]
   simply by comparing definitions.
If $A$ is a  complex $C^*$-algebra then
\[
 K_0(A) \cong \Kalg_0(A).
 \]  
 This statement does not hold for $K_j$ or $KO_j$, for $j \geq 1$.  For example, for a complex $C \sp *$-algebra $A$, we have
$K_1(A)$ is a quotient of $\Kalg_1(A)$. 
 
The simplest case is, of course, $A = \RR $.   Then $M_n(\RR )$ is simply the ring of $n \times n$ real matrices, a projection $p$ is simply an $n \times n $ 
matrix with $p = p^2 = p^*$.   If $p \in M_n(\RR )$ is a projection then there is an orthogonal matrix $u$ such that 
\[
upu^{-1} = \diag(1,1,1,\dots 1, 0,0,0,\dots 0),
\] 
the diagonal matrix with $j   $ 1's down the main diagonal and $n - j$ 0's following.  This has rank $j$ with $0 \leq j \leq n$ and defines a dimension function 
\[
\dim : \mathcal{P}_n(\RR) \longrightarrow \{0,1,\dots , n\} .
  \]
  These dimension functions coalesce to form an isomorphism 
  \[
  \dim : KO_0(\RR ) \overset{\cong}\longrightarrow \ZZ .
  \]
  (Exactly the same argument gives $K_0(\CC ) \cong \ZZ.$ )

   So far we have defined $KO_0(A)$ only for unital $A$.  
   Note that if $A$ is a real $C^*$-algebra with \index{unitization} unitization $A^+$ then there is a canonical short exact sequence  
   \[
   0 \longrightarrow A \longrightarrow A^+ \overset{\pi}\longrightarrow  (A^+/A) \, \cong \RR \longrightarrow 0.
  \]
  
  \begin{Def} 
  \index{K-theory@$K$-theory!definition}
  For any real $C^*$-algebra $A$, we define   
  \[
  KO_0(A) = \ker \big[  KO_0(A^+) \overset{\pi} \longrightarrow KO_0(\RR ) \cong \ZZ \big]
   \]
   so there is a natural\footnote{In this context,``natural" means that if there is a map $f: A \to B$ of $C^*$-algebras, then the diagram
   \[
   \begin{CD}   
    0 @>>>  KO_0(A) @>>>  KO_0(A^+) @>>>  \ZZ @>>> 0   \\
    @.  @VVf_*V    @VVf_*V      @VV{\id}V @.      \\
       0 @>>>  KO_0(B) @>>>  KO_0(B^+) @>>>  \ZZ @>>> 0   
   \end{CD}
   \]  
   commutes.   The sequences split because $\ZZ $ is a free abelian group.
   However, the splitting maps $\ZZ \to   KO_0(A^+)$  and    $\ZZ \to   KO_0(B^+)$ are in general 
   not canonical, and hence the diagram
   \[
   \begin{CD}
   KO_0(A^+)  @>\cong >> KO_0(A) \oplus \ZZ       \\
   @VVf_* V       @VV{f_*\oplus \id }V    \\
      KO_0(B^+)  @>\cong >> KO_0(B) \oplus \ZZ  
   \end{CD}
   \]
   does {\emph{not}} in general commute.
   }        short exact sequence
   \[
   0 \longrightarrow KO_0(A) \longrightarrow KO_0(A^+) \longrightarrow \ZZ \longrightarrow 0 .
  \]
  which splits unnaturally, and hence there is an unnatural isomorphism 
  \[
  KO_0(A^+) \cong  KO_0(A) \oplus \ZZ   .
  \]
 
 Using exactly the same idea, for any complex $C^*$-algebra $A$, we define 
 \[
  K_0(A) = \ker \big[  K_0(A^+) \overset{\pi} \longrightarrow K_0(\CC ) \cong \ZZ \big].
   \]
     so there is a natural   split short exact sequence
   \[
   0 \longrightarrow K_0(A) \longrightarrow K_0(A^+) \longrightarrow \ZZ \longrightarrow 0  
  \]
  and an (unnatural) isomorphism
  \[
  K_0(A^+) \cong  K_0(A) \oplus \ZZ   .
  \]

  A homomorphism $f: A \to B$ induces a homomorphism $f_* : KO_0(A) \to KO_0(B) $ and $KO_0$ is a covariant functor from the category of real $C^*$-algebras and homomorphisms 
  to abelian groups.   The same idea carries through in the complex setting.\footnote{Karoubi conjectured in
  the   1970's that if $A$ is a complex $C^*$-algebra then the algebraic and $C^*$-algebraic $K$-theory groups of $A \otimes \KK$ were 
  isomorphic, and similarly for the real setting.  Brown and the second author proved this for $A = \CC$ in 1976 \cite {BS}.  The general case was established 
  by Suslin and Wodzicki (1990 and 1994)\cite{KW}.  To get a feeling for this, look at \cite{BS}.    
   It turns out that $K_1^{alg}(\KK ) \cong K_1(\KK ) = 0.$  It is easy to show that  $K_1(\KK ) = 0$ but quite difficult to show 
  that $K_1^{alg}(\KK ) =    0.$
  }
    
  \end{Def}

 \begin{Rem}  \label{ComplexificationMap}
  \index{complexification~homomorphism}
 For each real $C^*$-algebra $A$   there is an inclusion map $A \hookrightarrow A\sc$
    which is a homomorphism.  This induces a natural complexification map \index{complexification}
    \[
  c: KO_0(A) \longrightarrow K_0(A\sc)
  \]
  which produces a natural transformation of theories.  It has an easy description.  If $p$ is a projection in $A$ then it is also a projection in $A\sc $.  If $p$ and 
  $q$ are unitarily equivalent projections in $A$ then this also holds in $A\sc$.  Passing to $M_n(A)$ and $M_n(A\sc)$, the same holds. Therefore there is a homomorphism

  \[
  c:  KO_0(A) \longrightarrow K_0(A\sc )
  \]
 defined by $c([p]) = [p]$, which is well-defined and well-behaved for all real $C^*$-algebras.

   We can describe this map in a different way.  If we regard $A\sc $ as a real $C^*$-algebra this induces a natural transformation
  \[
c: KO_0(A) \longrightarrow KO_0(A\sc ).
  \]
 A moment's thought about projections or about finitely generated projective modules \index{projective module} will convince the reader that 
  \[
  KO_0(A\sc )    \cong  K_0(A\sc )
  \]
  and so the composition gives a natural homomorphism
   \[
  c:  KO_0(A) \longrightarrow K_0(A\sc ) .
  \]

    \end{Rem}

    Here's an important example: $A = \realo (\RR )$  so that $A^+  = \real (S^1)  $.  Then there is a natural short exact sequence
   \[
   0 \longrightarrow  \realo (\topr )\longrightarrow \real (S^1) \longrightarrow \RR \longrightarrow 0
   \]
   and we shall see later that this induces a short exact sequence 
   \[
   \begin{CD}
   0 @>>> KO_0( \realo (\topr)) @>>>    KO_0 ( \real (S^1)) @>>> KO_0(\RR ) @>>> 0  \\
   @.     @VV\cong V    @VV\cong V      @VV\cong V     @.   \\
   0 @>>>   \ZZ_2    @>>>    \ZZ \oplus \ZZ _2 @>>> \ZZ     @>>> 0
   \end{CD}
   \]
   
   For the corresponding analysis in the case of $C \sp *$-algebras, we have the short exact sequence
    \[
   \begin{CD}
   0 @>>> KO_0( \complexo (\topr)) @>>>    KO_0 ( \complex (S^1)) @>>> KO_0(\CC ) @>>> 0  \\
   @.     @VV\cong V    @VV\cong V      @VV\cong V     @.   \\
   0 @>>>  0    @>>>    \ZZ  @>>> \ZZ     @>>> 0
   \end{CD}
   \]
   From this we see that $c \colon KO_0(A) \rightarrow K_0(A\sc)$ is not generally an isomorphism.

  We move on now to the higher $KO$-groups.
  
  \begin{Def}
  \index{general linear group}
  If $A$ is a unital real $C^*$-algebra, then define
$
  GL_j(A) $ to be the topological group of invertible matrices in $M_j(A)$.  If $A$ is not unital,  define 
  \[
  GL_j(A) = \{ x \in GL_j(A^+) \mid x \equiv 1_j  \mod M_j(A) \} .
  \]
\index[notation]{GLn@$GL_n(A)$}
\index[notation]{GLzinf@$GL_\infty(A)$}

 Then $GL_j(A) $ is a closed normal subgroup of $GL_j(A^+)$.   We may embed 
 \[
 GL_j(A) \to GL_{j+1}(A) 
 \]
  by $x \mapsto \diag(x,1)$ and let 
  \[
  GL_\infty (A) = \lim_{j \to \infty} GL_j(A).
  \]
It is easy to see that $GL_j(-) $ is a covariant functor from real $C^*$-algebras to topological groups.  
  \end{Def}

  We note that it is straight-forward to prove\footnote{using CW-complexes, which we will discuss in \S~\ref{top spaces}.}  that the natural maps 
  \[
  GL_j(A) \longrightarrow GL_\infty (A)
  \]
  induce an isomorphism 
  \[
  \lim_{j \to \infty} \pi_k (GL_j(A)) \xrightarrow{~\cong~} \pi _k(GL_\infty (A)) \quad \text{for all $k$.}
  \]

  There are two ways to define $KO_n(A)$ for $n = 1, 2,  3,\dots $ and both ways have their pros and cons.

\begin{Def} {\bf{		FIRST DEFINITION }} \label{defn:realKOi}  
\index{K-theory@$K$-theory!definition}
For $n \geq 1$, define\footnote{We always take the identity of the topological group to be the basepoint for homotopy purposes, and we always use based spaces and basepoint preserving maps when dealing with compact spaces. In general $\pi _n(X)$ is an  abelian group for $n  >  1$ and only a group for 
$n= 1$.  
   However, $\pi _1$ of a topological group is always abelian. }  
  \[
  KO_n(A)  = \pi _{n-1} (GL_\infty (A) ).
  \]
  If $f : A \to B$ is a homomorphism of real $C^*$-algebras then it induces  a map
  \[
  f: GL_\infty (A) \to GL_\infty (B)
  \]
  and hence 
  \[
   f_* : KO_n(A) \to KO_n(B).
   \]
  Then each  $KO_n(-) $ is a covariant functor from real $C^*$-algebras to abelian groups.   We write $KO_*(A)$ for the graded abelian group 
  $\{ KO_n(A) \}_{n \geq 0}$ (we later extend this to $n \in \ZZ$ using Bott periodicity).  The same construction for complex $C^*$-algebras yields $K_n(A)$. 
       \end{Def}
     
     For many purposes it is convenient to use the fact that $O_n(A)$, the group of orthogonal $n \times n$ matrices, is 
     a strong deformation retraction of $GL_n(A)$.\footnote{For example, the set $O_1(\RR) = \{\pm 1\} $ is a strong deformation 
     retraction of $GL_1(\RR ) = \RR - \{0\}$, and $O_2  $ (rotations and reflections of the plane - topologically two circles)  is a strong deformation retraction of $GL_2(\RR )$.} 
     The compact Lie group $O_n(\RR ) $ is naturally contained in $ GL_n(\RR )$ and there is a strong deformation retraction\footnote{ also known as Gram-Schmidt orthogonalization}
     which deforms $GL_n(\RR)$ down to $O_n(\RR )$. 
  Thus the two spaces are homotopy equivalent. This implies that $O_\infty (\RR )$ has the same homotopy type as $GL_\infty (\RR)$, and an analogous argument 
  demonstrates that   $O_\infty (A )$ has the same homotopy type as $GL_\infty (A)$ for any real (or complex!)  $C^*$-algebra $A$.  Thus we may as well define 
   \[
  KO_n(A)  = \pi _{n-1} (O_\infty (A) )
  \]
  and this is the way it is done in many places in the literature.

  \begin{Def} {\bf{ SECOND DEFINITION}}   \label{defn:realKOi.v2}
  \index{K-theory@$K$-theory!definition}
   For $n \geq 1$, define 
   \[
   KO_n(A) = KO_0(S^nA).
   \]
     If $f : A \to B$ is a homomorphism of real $C^*$-algebras then it induces  
  \[
  f: S^nA \longrightarrow S^nB 
  \]
  and hence 
  \[
   f_* : KO_n(A) \to KO_n(B).
   \]
  Then each  $KO_n(-) $ is a covariant functor from real $C^*$-algebras to abelian groups. As before, we write $KO_*(A)$ for the graded abelian group 
  $\{ KO_n(A) \}_{n \geq 0}$.  
   \end{Def}
   
   \begin{Thm} 
   For all real unital $C^*$-algebras $A$, 
    \[
  \pi _{n-1} (GL_\infty (A) ) \cong KO_0(S^n A) 
    \]
and hence the two definitions of $KO_n(A)$ agree for all $A$.
   \end{Thm}
   
 This result has been ``well-known" for fifty years or so, but it is remarkably difficult to put together a rigorous proof.\footnote{Jon Rosenberg (private communication) suggested the following approach. First, define bigraded $K$-theory via Clifford algebras \index{Clifford algebra} and use the periodicity for Clifford algebras to 
   show that $KO_*(A)$ is periodic of period 8.  The resulting homology theory may be represented by an infinite loop space, and you then show that this infinite loop space may be taken 
   to be determined by $GL_\infty (A)$.  For information on Clifford algebras and periodicity please see future sections. The technical details of this approach are well beyond the scope 
   of this {\BS}.}  
 Herbert  Schr\"oder (\cite{Schr} Theorem 1.4.6)     briefly discusses the proof of this theorem and then refers to Karoubi's book \cite{K} III, 1.6 (it is actually 1.11).  The book does prove periodicity for $K_*$ but not for $KO_*$.  For a detailed proof of that you have to go to an earlier paper of Karoubi, namely \cite{K68}.

  \begin{Thm}  $KO_*$ is an additive homology theory on the category of real $C^*$-algebras.
  \end{Thm}
  
  \begin{proof} 
  Schr\"oder \cite{Schr} establishes the homotopy axiom (in somewhat greater generality) as Proposition~1.4.2 
  and the additivity axiom as Proposition~1.4.8.
  \end{proof}
  
  \begin{Rem} We previously defined the complexification natural transformation  \index{complexification~homomorphism}
  \[
  c:  KO_0(A) \longrightarrow K_0(A\sc )
  \]
  for all real $C^*$-algebras $A$.  Using the second definition (and its analog for complex $K$-theory),  it is easy to show that this natural transformation 
  extends to a natural transformation 
  \[
  c:  KO_n(A) \longrightarrow K_n(A\sc )
  \]
   of homology theories.
   \end{Rem}

  We note a very important and immediate consequence of this result. It justifies our use of ``desuspension" to describe $VA$. 
 
 \index{suspension!of a $C^*$-algebra} \index{desuspension}
 
  \begin{Cor} Suppose that $A$ is a real $C^*$-algebra. Then:
  \begin{enumerate}
  \item 
  the natural short exact sequence 
  \[
  0 \longrightarrow  SA \longrightarrow CA \longrightarrow A \longrightarrow 0
  \]
  induces a natural {\bf{suspension}} isomorphism
  \[
  KO_n(A) \cong KO_{n-1}(SA) \; .
  \]
   \item  \index{Toeplitz algebra!Toeplitz extension} There is a natural {\bf{desuspension}}  isomorphism\footnote{The boundary map in the $KO$-theory short 
exact sequence associated to the Toeplitz extension
 \[
 0 \longrightarrow A \otimes \KK\sr\longrightarrow   \mathcal{T}_{\RR 0 }\otimes A  \longrightarrow VA \longrightarrow 0  
  \] 
has the form 
  \[
  KO_{n+1}(VA) \overset{\cong}  \longrightarrow KO_n(A \otimes \KK\sr ) \cong KO_n(A)
  \]
 which we write right to left to make the comparison with suspension more dramatic.}
\[
  KO_n(A) \cong KO_{n+1}(VA) \; .
  \]

   \end{enumerate}
  \end{Cor}
  
  \begin{proof} This is immediate from Corollary \ref{T: desuspension}.
  \end{proof}
   
  \begin{Cor} For any real $C^*$-algebra $A$, there are natural isomorphisms
  \[
  KO_m(S^n A)\cong KO_{m+n}(A)    
   \]  
   where $S^n A$ is defined for positive and negative integers according to Definition~\ref{Def-AllSus}.
  
  \end{Cor}

   \begin{Cor} 
   \index{compact operators}
   For all real $C^*$-algebras $A$ and   integers $k \geq 1$ and $n \geq 0$, 
   \[
   KO_n(A) \cong KO_n( M_k(A) )   \cong KO_n(A \otimes \KK \sr  )
   \]
  \end{Cor}
  
  \begin{proof}  It suffices to consider the case $n = 0$, since by the previous corollary, $KO_n(A) \cong KO_0(S^nA)$.  The group $KO_0(A)$ is generated by 
  projections in matrix rings over $A$, and we use the simple observation that 
  \[
  M_i (M_j(A)) \cong M_{ij}(A)
  \]
  to show that, after passing to equivalence classes, we have 
  \[
 \mathcal{P}(A)_\infty = \mathcal{P}_\infty(M_k(A))       =       \mathcal{P}_\infty(A\otimes \mathcal{K}\sr ).
  \]
  \end{proof}

  \begin{Thm} (classical) \label{Thm-classical} Here are the $KO$-groups of the real $C^*$-algebras $\RR$, $\HH$ and $\CC$. 
  \index{quaternions}
\[
  \begin{tabular} { c c c c c c c c c}
  n & $0$ & $1$ & $2$ & $3$ & $4$   & $5$    &  $6$    & $7$ \\
  \hline
  $KO_n(\RR)$ &  $ \ZZ$ & $ \ZZ_2 $   & $\ZZ_2$      &  $0$    & $\ZZ$  & $0$    & $ 0$   & $0$     \\
   $KO_n(\HH)$ &  $ \ZZ$ & $0 $   & $0$      &  $0$    & $\ZZ$  & $\ZZ_2$    & $ \ZZ_2$   & $0$     \\ 
    $KO_n(\CC)$ &  $ \ZZ$ & $ 0 $   & $\ZZ $      &  $0$    & $\ZZ$  & $0$    & $ \ZZ$   & $0$     \\
    \end{tabular}
   \]

  \end{Thm}
  
  Let's look at some examples. 
    \index{general linear group} \index{special orthogonal group} \index{classical group}
    Steenrod \cite{St} \S 22 is an excellent reference for concrete calculations.  Of course he was 
  calculating homotopy groups of classical Lie groups, as $K$-theory had not been invented then.
  \begin{enumerate}
  \item 
  The group $KO_0(\RR) = \ZZ$.  We discussed this previously: a projection in $M_k(\RR)$ is determined up to equivalence by its rank. The group is generated by 
  formal differences of projections in matrix algebras over $\RR$, and so the dimension (rank) function gives the isomorphism.
  
  \item We have
  \[
   KO_1(\RR ) = \pi _0(GL_\infty (\RR ) ) \cong \pi _0 \left( \lim_{n \to \infty}  GL_n(\RR ) \right) = \lim_{n \to \infty} \pi _0(O_n).
  \]
  The topological group $O_n$ has two path components, namely the rotations $SO_n$ and the reflections, and the inclusions $O_n \subset O_{n+1} $ respect this, so  
  \[
  KO_1(\RR) = \lim_{n \to \infty}\pi _0(O_n) = \lim_{n \to \infty}  \ZZ_2 = \ZZ_2 .
  \]
  
  \index[notation]{On@$O_n$}
  \index[notation]{Son@$SO_n$}
  
Alternatively, we can describe the non-trivial element of $KO_1(\RR) = \ZZ_2$ in terms of a projection as follows. We have
\begin{align*} KO_1(\RR) = KO_0(S \RR) &= \ker(KO_0(S\RR^+) \rightarrow KO_0(\RR)) \\
	&= \ker((\ev_1)_* \colon KO_0(C\pr(S^1)) \rightarrow KO_0(\RR)) 
\end{align*}
The non-trivial element $\eta$ is represented by 
$$\eta = \left[ p  \right] - \left[ \smat 1000  \right] $$
where the projection $p \in M_2(C\pr(S^1))$ is defined by
$$p(x,y) = \tfrac{1}{2} \begin{bmatrix} 1+x & y \\ y & 1-x \end{bmatrix}$$ for $(x,y) \in S^1$.
We check that $\ev_1(p) =\smat 1000 $, so that $\eta \in \ker (e_1)_*$.
Of course, the projection $p$ corresponds to the M\"obius vector bundle over the circle. See Theorem 4.2 of Boersema \cite{Boersema2020}

\item We could also do the same thing for the non-trivial element of 
\[
KO_2(\RR) = KO_0(S^2 \RR).
\]
 This time the corresponding projection defined on $(x,y,z) \in S^2$ is
$$p(x,y,z) =\tfrac{1}{2} \begin{bmatrix} 1 + x & 0 & y & z \\ 0 & 1+ x &-z & y \\ y & -z & 1-x & 0 \\ z & y & 0 & 1-x \end{bmatrix}$$
See Theorem 4.4 of \cite{Boersema2020}.

  \item 
   For subsequent groups we need only look at the path component of the identity of $GL_\infty (\RR )  $ and hence really only at the special orthogonal group $SO_n$.  So 
  \[
  KO_2(\RR) \cong \pi _1(SO_n)   \qquad \text{$n$ large}
  \]
  Now $SO_2 \cong S^1$ and so $\pi _1(SO_2) \cong \ZZ$, but $SO_3 \cong RP^3$ as a topological space (where $RP^3$ is real projective space) and hence $\pi _1(SO_3) = \ZZ_2$. 
  For $n > 2$, any map $S^1 \to SO_n$ is homotopic\footnote{by the \index{cellular approximation theorem}
  Cellular Approximation Theorem, cf. G. Whitehead \cite{W} p. 77 and see \S 9 ahead.} to a map taking values in $SO_3 $ and hence $\pi _1(SO_n) = \ZZ_2$ for $n > 2$, confirming that 
  \[
  KO_2(\RR) \cong \ZZ_2 .  
  \]
  \item  To calculate
   \[
  KO_3(\RR) \cong \pi _2(SO_n)  \qquad \text{$n$ large}
  \]
we may use again the fact that $SO_3 \cong RP^3$ to deduce that $\pi _2(SO_3 ) = 0$.   Steenrod shows (p. 116) that $\pi _2(SO_4) = 0$. 
Cellular Approximation  tells us that $\pi _2(SO_n) = 0 $ for $n>4$ and hence  
   \[
  KO_3(\RR) \cong 0 .
  \]
  There is an alternate approach which is much simpler:  a 1936 result of \'E. Cartan tells us that if $G$ is any compact Lie group then $\pi _2(G) = 0.$
  
  \item Steenrod calculates  $\pi_j(SO_n) $ for $j = 3,4$ in \S\S 22-24, and these imply that 
  \[
  KO_4(\RR) \cong \ZZ   \qquad \text{and} \qquad     KO_5 (\RR) \cong 0.      
  \]
  \item The story of the computation of the group 
  \[
  KO_6(\RR) \cong \pi _5(SO_n) = 0
  \]
   is very interesting 
  historically.  Its computation depended upon calculating $\pi _5(S^3)$ and that group was the subject of controversy. Pontrjagin had announced that   $\pi _5(S^3) = 0$, which 
  conflicted with predictions stemming from Bott's work on the Periodicity Theorem  - the second author heard Bott mention once that he withheld publication of his results on periodicity for that reason. Eventually
  Pontrjagin and (independently) G. Whitehead determined that $\pi _{n+2}(S^n) \cong \ZZ_2$.  This all went on during the publication of Steenrod's book \cite{St}  (1951) - you see mention of it in \S 24.11 and the end of  \S 21. 
  
 \end{enumerate}
 

   \newpage
           \section{\bf{External and Internal Products }}
           \label{Section:Products}

   In this section we describe the external products on real and complex $K$-theory and how they fit together 
   coherently. We will see that this also implies the existence of
   internal products in the commutative setting, the existence of a ring structure for $KO_*(A)$ in the commutative case, and more generally a $KO_*(\RR)$-module structure for $KO_*(A)$ in all cases.  The product structure is also foundational to the K\"unneth Theorem which we discuss in Section~\ref{Section-Kunneth}.
      
    \vspace{.1in}
 
{\bf{External Products}}

    \vspace{.1in}

 \begin{Def}
 \index{product on $K$-theory!external product}
 An homology theory $h_*$ on real or complex $C^*$-algebras is said to have {\emph{external products}}
 if for any $C \sp *$-algebras $A$ and $B$ (or perhaps any nuclear $C \sp *$-algebras) in the category there is a family of natural transformations
 \[
 h_m(A) \otimes h_n(B) \longrightarrow h_{m+n}(A\otimes B)
 \]
 which are associative and graded commutative: if $x \in h_m(A)$ and $y \in h_n(B)$ then 
 \[
 yx = (-1)^{mn} xy .
 \]
 \end{Def}

   \begin{Pro} \label{Prop-external product}
   \begin{enumerate}
  \item
   For real $C ^*$-algebras $A$ and $B$ (with nuclearity assumed), there exists a natural external product map
 \[
  \mu\so \colon KO_m(A) \otimes KO_n(B) \rightarrow KO_{m+n}(A \otimes\sr  B) .
 \]
   \item Similarly, for complex 
    $C ^*$-algebras $A$ and $B$ (with nuclearity assumed), there exist a natural  external product map
 \[
  \mu\su \colon K_m(A) \otimes K_n(B) \rightarrow K_{m+n}(A \otimes\sc B).
 \]
  \item  These product maps respect complexification. That is, the diagram
  \[
  \begin{CD}
    KO_m(A) \otimes KO_n(B) @> \mu_O>>  KO_{m+n}( A \otimes\sr B) @>=>>  KO_{m+n}( A \otimes\sr B)  \\
 ^{c\otimes c}  @VVV^c      @VVV^c      @VVV    \\ 
  K_m(A\sc) \otimes K_n(B\sc) @>  \mu_U>>  K_{m+n}( A\sc \otimes\sc B\sc) @> \cong >>  K_{m+n}( (A \otimes\sr  B)\sc)  
  \end{CD}
  \]
   commutes, where the vertical maps are complexification.
   \end{enumerate}
   \end{Pro}

 \begin{proof} (Sketch) 
 If $p,q$ are projections in
 $M_j(A)$ and $M_k(B)$, respectively, then the Kronecker product\footnote{If $M$ and $N$ are matrices representing linear transformations $V_1 \to W_1$ and $V_2 \to W_2$ with respect to fixed choices of bases 
  $\{v_i^1\} $,   $\{v_j^2\} $,   $\{w_k^1\} $,   $\{w_l^2\} $,  
  then their Kronecker 
  product $M \otimes N$  is the matrix representing  the linear transformation $V_1 \otimes V_2 \to W_1 \otimes W_2 $ with respect to bases $\{v_i^1 \otimes v_j^2 \}$ and 
  $\{w_k^1 \otimes w_l^2 \}$.}
gives $p \otimes q$, a projection in $M_{j+k}(A \otimes B)$. This induces a map
 \[
 KO_0(A) \otimes KO_0(B) \rightarrow KO_0(A \otimes B)
 \]
  (skipping over details regarding how this provides a well-defined map vis a vis the Grothendieck construction). 
 Taking suspensions, we then obtain maps
 \[
 KO_0(S^mA) \otimes KO_0(S^n B) \rightarrow KO_0(S^{m+n}(A \otimes B))
 \]
  for all $m, n \geq 0$, which gives $\mu_O$.
 
 The analogous construction  gives the complex pairing. 
  
 \end{proof}

    \vspace{.1in}

{\bf{Internal Products}}
 
    \vspace{.1in}
 
 \index{product on $K$-theory!internal product}
 
 Internal products for $K$-theory do not exist in general for $C \sp *$-algtebras, but we will see that they do in the commutative setting, which is to say in the topological setting. 
  
 \begin{Def}
 An homology theory $h_*$ on a category $\mathcal A$ of real or complex $C^*$-algebras is said to have {\emph{internal products}}
 if for any (perhaps with nuclearity assumed) algebras $A$  in the category there is a family of  natural transformations
 \[
 h_m(A) \otimes h_n(A) \longrightarrow h_{m+n}(A)
 \]
 which is associative and graded commutative: if $x \in h_m(A)$ and $y \in h_n(A)$ then 
 \[
 yx = (-1)^{mn} xy .
 \]
 In other words, $h_* $ is a functor from $\mathcal A$ to the category of graded commutative rings.
 \end{Def}
 
\begin{Pro}  
\begin{enumerate}
\item
Restrict $K_*$ to the category of commutative complex $C^*$-algebras.  Then $K_*$ has an internal product.
\item
Restrict $KO_*$ to the category of commutative real $C^*$-algebras.  Then $KO_*$ has an internal product.
\end{enumerate}
\end{Pro}

\begin{proof} 
For (1), 
the natural map $m: A \otimes A \to A$ given by $x \otimes y  \to xy$ is a homomorphism of $C^*$-algebras since $A$ is commutative, and hence it induces
a map 
\[
m_* \colon K_{m+n}(A \otimes A) \longrightarrow K_{m+n}(A) \; .
\]
Then there is an internal product defined by composition
\[ m_* \circ \mu\pu  \colon K_m(A) \otimes K_n(A) \longrightarrow K_{m+n}(A)  \; . \]
 In the unital commutative case, we have $A = \complex(X)$ for some compact space $X$, and then the multiplication map $m$ corresponds to the diagonal 
map $\Delta \colon X \to X \times X $ which is used to obtain internal products for topological $K$-theory.

The proof in the real case is identical to the complex case in the setting of $C \sp *$-algebras. In the unital commutative case, the multiplication map $m$ on a commutative real $C \sp *$-algebra corresponds to the diagonal map 
\[
(X, \tau ) \longrightarrow (X, \tau ) \times (X, \tau ) \; 
\]
on topological spaces with involution.
\end{proof}

  \begin{Pro} \label{Prop-commutative ring}
  \begin{enumerate}
  \item
Complex $K$-theory $K_*$ is a functor from the category of commutative complex $C^*$-algebras to the category of 
 graded commutative rings. 
 \item
 Real $K$-theory  $KO_*$ is a functor from the category of commutative real $C^*$-algebras to the category of 
 graded commutative rings. 
 \end{enumerate}
  \end{Pro}
  
  That is, for any commutative real $C^*$-algebra $A$, there is a natural graded commutative ring structure on $KO_*(A)$  and 
a map $\alpha \colon A \to B$ induces a map of graded commutative rings $KO_*(A) \to KO_*(B)$.  The analogous statement holds in the complex world. 
 
\begin{Rem} If $A$ is a unital commutative real $C \sp *$-algebra, then the graded ring $KO_*(A)$ is actually a unital ring with identity element $[1]$. Similarly if $A$ is a unital commutative complex $C \sp *$-algebra then $K_*(A)$ is a ring with identity. For a non-unital $C \sp *$-algebra, the internal product on $KO_*(A)$ or $K_*(A)$ is not necessarily well behaved, and it could very well be that all products vanish.

For example, take $A = S\RR$. Then there is a non-trivial element $x \in KO_{-1}(\RR) = \ZZ$ which generates $KO_*(S\RR)$ as a module over $KO_*(\RR)$ (see the next section). We have $x^2 = 0$ because $KO_{-2}(\RR) = 0$. Furthermore, for any elements $y,z \in KO_*(S \RR)$ we have $y = \theta_1 \cdot x$ and $z = \theta_2 \cdot x$ for some elements $\theta_1$ and $\theta_2 \in KO_*(\RR)$. Then
$$y \cdot z = (\theta_1 \cdot x) \cdot (\theta_2 \cdot y) = (\theta_1 \cdot \theta_2) \cdot (x \cdot x) = 0 \; .$$

 \end{Rem}
 
  

  
 \vspace{.1in}
{\bf{The ring $KO_*(\RR)$ and modules over $KO_*(\RR)$.}}
 \vspace{.1in}
 
 Since $\RR$ is a commutative $C \sp *$-algebra it immediately follows that $KO_*(\RR)$ is a ring. The rest of this section is devoted to describing this ring and the structure that this provides to $KO_*(A)$ for any real $C \sp *$-algebra $A$.
 
 \begin{Pro} \label{prop-module}
$KO_*(\RR)$ is a graded commutative ring. Furthermore, for any real $C ^*$-algebra $A$,  $KO_*(A)$ is a graded module over $KO_*(\RR)$.
Similarly, 
$K_*(\CC)$ is a graded commutative ring; and for any complex $C ^*$-algebra $A$,  $K_*(A)$ is a graded module over $K_*(\CC)$.
 \end{Pro}
 
 \begin{proof}
 The first sentence follows from Proposition~\ref{Prop-commutative ring}.
 
 Then for any real $C ^*$-algebra $A$, we have $\RR \otimes A \cong A$ which give maps
 $$KO_m(\RR) \otimes KO_n(A) \rightarrow KO_{m+n}(\RR \otimes A) = KO_{m+n}(A) \; $$
 using the external product in Proposition~\ref{Prop-external product}.
 
 The complex case is essentially identical.
 \end{proof}

 We now discuss the ring structure of $KO_*(\RR)$. Recall from Theorem~\ref{Thm-classical} that $KO_*(\RR)$ has period 8 and that the groups in degrees $0$ through $7$ are as follows.  \[
  \begin{tabular} { c c c c c c c c c}
  n & $0$ & $1$ & $2$ & $3$ & $4$   & $5$    &  $6$    & $7$ \\
  \hline
  $KO_n(\RR)$ &  $ \ZZ$ & $ \ZZ_2 $   & $\ZZ_2$      &  $0$    & $\ZZ$  & $0$    & $ 0$   & $0$     \\
    \end{tabular}
   \]
   We specify generators for these groups as follows:
   \begin{enumerate}
   \item      The generator of $KO_0(\RR) = \ZZ$ is the class of the unit $[1]$, and it is easy to check that 
   \[
   [1] \cdot [p]= \mu([1] \otimes [p]) =  [p]  {\text{\qquad for any\quad }} p \in \mathcal{P}_k(A).
   \]
    In particular, $[1]$ is a unit in the ring $KO_*(\RR)$
    and serves as the identity under multiplication.
 
 \item
 Let $\eta \in KO_1(\RR)$ be the non-zero element. Then it turns out that $\eta^2$ is also the non-zero element of $KO_2(\RR)$. 
 
 \item We take $\xi$ to be a generator of $KO_4(\RR) = \ZZ$.
 
 \item The Bott element $\beta_O$ is a generator of $KO_8(\RR)= \ZZ$.    It turns out that $\xi^2 = 4 \beta_O$. 
  \end{enumerate}
 
 \begin{Thm}(Karoubi \cite{K}, III.5.19)  \label{thm:KO-relations}

  $KO_*(\RR)$ is the graded commutative unital ring with generators
  \begin{itemize}
  \item     ~$\eta$ in degree 1, 
  \item ~$\xi$ in degree 4, 
  \item ~$\beta\so$ in degree 8, 
  \end{itemize}
where $\beta\so$ has an inverse, and subject to the relations:
 \begin{itemize}
 \item ~$2 \eta = 0$
 \item ~$  \eta^3 = 0$
 \item ~$  \xi^2 = 4 \beta\so$
 \item ~$ \xi \eta = 0$
 \end{itemize}
 \end{Thm}
 
 \begin{Rem}
 This result can be rephrased as follows: the ring $KO_*(\RR )$ is isomorphic to the quotient of the graded polynomial ring over the integers generated by $\eta, \xi, \beta_O,$ 
 and $\beta_O ^{-1} $ modulo the ideal generated by the relations  specified above.
 \end{Rem}
 
 
  
 \begin{Rem} 
  \index{quaternions}
  From Theorem~6.18 we have
 \[ KO_n(\RR)\otimes \QQ  = \begin{cases} \QQ &  n = 0,\, 4 \\  0 & \text{otherwise.}  \end{cases}
 \]
 It turns out more generally that $KO_*(A)\otimes \QQ $ is periodic with period $4$ for any real $C \sp *$-algebra $A$. This is also true for  $KO_*(A) \otimes \Z[\tfrac{1}{2}]$.
 See Corollary~\ref{Cor:RationalPeriodicity} below.
  \end{Rem}

Next we introduce one of the other homology theories that may be constructed from $KO_*$.

\begin{Def}  (Karoubi) {\emph{Quaternionic   $K$-theory}}  \index{K-theory@$K$-theory!quaternionic}
is defined by 
\index[notation]{KHA@$KH_*(A)$}
\[
KH_*(A) = KO_*(A \otimes \HH )
\]
\end{Def}

\begin{Rem}(Karoubi) \label{Karoubi1}
We know from Corollary~\ref{C:Hhomology} that $KH_*$ is indeed an additive homology theory.   There is a natural exterior pairing
\[
KH_r(A) \otimes KH_s(A)  = KO_r(A \otimes \HH )\otimes  KO_s(A \otimes \HH ) \longrightarrow KO_{r+s}(A \otimes A \otimes \HH \otimes \HH ).
\]
Since $\HH\otimes \HH \cong M_4(\RR)$ as remarked in \S 1, this simplifies to 
\[
KH_r(A) \otimes KH_s(A) \longrightarrow KO_{r+s}(A \otimes A  )
\]
and if $A$ is real commutative we obtain a pairing
\[
KH_r(A) \otimes KH_s(A) \longrightarrow KO_{r+s}(A)
\]
which, as we shall see, is related to real Bott Periodicity.

\end{Rem}

\begin{Cor}\label{C:KH}
For any real $C \sp *$-algebra $A$, we have
\[
KH_*(A) \cong KO_{*+4}(A).
\]
\end{Cor}
  
  Thus $KH_*(A)$ is nothing but a shift of $KO_*(A)$. Since $KO_*(A)$ has period 8 (by Bott periodicity), the same is true for $KH_*(A)$.

  \begin{proof}
  \index[notation]{sz@$\Sigma^k M$}
For any graded abelian group $M$, let $(\Sigma ^kM)_i = M_{i + k}$ (see Definition~\ref{def-GroupSuspension}).  Then
\begin{align*}
KH_*(A) &\cong KO_*(A \otimes \HH) \\
	&\cong KO_*(A) \otimes_{KO_*(\RR)} KO_*(\HH) &&\text{by the K\" unneth Formula, Theorem~\ref{RealKF}} \\
	&\cong KO_*(A) \otimes_{KO_*(\RR)} \Sigma^4 KO_*(\RR) && \text{by Theorem~\ref{Thm-classical}} \\
	&\cong \Sigma^4 (KO_*(A) \otimes_{KO_*(\RR)}  KO_*(\RR)) \\
        &\cong \Sigma^4 KO_*(A). \\
\end{align*}

\end{proof}
For the relationship between $KH_*$, $KSp_*$ and other related theories, please see \cite {KW}, particularly Example 2.16.

 \newpage
   \section{\bf{Bott Periodicity} }
       \label{Section:Bott}
 
 In this section we show how Bott Periodicity provides real power to  $KO_*$ and $K_*$.    Complex $K$-theory is periodic 
 of period 2. In contrast, 
 
 \begin{Thm} (Bott Periodicity) \cite{Bott57}, \cite{K68}, \cite{Wood}. 
   \index{Bott Periodicity}
   For any real $C^*$-algebra $A$, there is a natural isomorphism
 \[
 KO_{j + 8}(A) \cong KO_j(A).
  \]

  \end{Thm} 
  
  This theorem was initially established by Bott \cite{Bott57} in 1957 in the context of the homotopy groups of the classical Lie groups, i.e.,   before there really was a $KO$-theory. 
  There are many, many other proofs of Bott's theorem, originating in several different 
  areas of math, illustrating its importance.  
  
   In the late 1960's,  Wood \cite{Wood}   and  Karoubi \cite{K68}  (see also Milnor \cite{M3})  proved it in the context of real (and also complex) Banach algebras and Banach categories.
 Here is a vague outline of Wood's proof. (We shall return to this in much greater detail in Part 3.)  Recall from Section~\ref{Section:KODefinition} that 

  \[
  KO_n(A) \cong \pi _{n-1} (GL_\infty (A) ).
     \]
  Wood proves that the connected component of the identity of $GL_\infty (A)$  is homotopy equivalent to the connected component of the identity of 
the 8th loop space   $\Omega ^8 GL_\infty (A) $ of  $GL_\infty (A)$.\footnote{See Section~\ref{top spaces} for the nice properties of  loop spaces $\Omega ^kX$.}
  This implies that, for each $n$, that there is an isomorphism
  \[
\pi _j( GL_\infty (A) )              \cong      \pi _j  (\Omega ^8 GL_\infty (A) )
  \]
  But
  \[
      \pi _j  (\Omega ^8 GL_\infty (A) )  \cong     \pi _{j+ 8}  (GL_\infty (A) )
  \]
 which translates to

   \[
 KO_{j + 8}(A) \cong KO_j(A) .
  \]
  
  Note that just as $GL_n(\RR) $ deformation retracts down to $O_n(\RR)$,  the group $GL_n(A) $ deformation retracts down to $O_n(A)$ and so we may use 
  orthogonal groups rather than the larger general linear groups. 
  
  \begin{Cor} Using periodicity, we may extend $KO_*(A)$ for all $* \in \ZZ .$
  \end{Cor}
  
  Recall that we let $\oo = \lim_{n \to \infty}  O_n $. 
  
    \index[notation]{oz@$\oo$}

  \begin{Cor}   For each $n > 0$,
  \[
  KO_n(\RR ) \cong \pi _{n-1}(\oo  ) \cong \pi _0 (\Omega ^{n-1}\oo).
  \]
    \end{Cor}
  
  This is a useful way to think about $KO_n(\RR )$, since 
 each of the spaces $\Omega ^k\oo  $ may be  constructed rather concretely from classical Lie groups, as we shall see.

  We now recap the principal properties of $KO_*$ for future use.   Recall that $\mathscr{C}^*_{\rm{real}}$ denotes the category of real $C^*$-algebras 
  and their homomorphisms.
  
  \begin{Thm} 
    \index{homology theory}
  Real $K$-theory $KO_*(A)$  is a  $\ZZ$-graded additive homology theory defined on the category $\mathscr{C}^*_{\rm{real}}$  
  and taking values in $\ZZ$-graded abelian groups. That is, it satisfies
  the following axioms:
  \begin{enumerate}
    \item {\underline{Homotopy Axiom}}. \index{homotopy axiom}  Let $H : A \to IB$ be a homotopy from $f = p_0 H $ to $g = p_1 H$ 
 in  $ \mathscr{C}^*_{\rm{real}}  $. Then 
 \[
 f_* = g_* : KO_n(A) \longrightarrow KO_n(B) \quad \text{for all n}.
 \]
 \vspace{.1in}
\item  {\underline{Exactness Axiom}}.  \index{exactness axiom}  Let
 \[
 0 \longrightarrow J \overset{i}\longrightarrow A \overset{j}\longrightarrow A/J  \longrightarrow 0
 \]
 be a short exact sequence in $ \mathscr{C}^*_{\rm{real}}$.  Then for each $n$  there is a map 
 \[
 \partial : KO_n(A/J) \longrightarrow KO_{n-1}(J)
 \]
 and a long\footnote{{\emph {very}} long! It extends infinitely in both directions.}
  exact sequence 
 \[
\dots \to KO_n(J) \overset{i_*}\longrightarrow KO_n(A)  \overset{j_*}\longrightarrow KO_n(A/J )  \overset{\partial }\longrightarrow KO_{n-1}(J ) \to \dots 
 \]
  The map $\partial $ is natural with respect to morphisms of short exact sequences. 
  \vspace{.1in}
 \item  {\underline{Additivity Axiom.}}  \index{additivity axiom} Let $A_i \in  \mathscr{C}^*_{\rm{real}} $ and let
 \[
 A = \oplus _{i = 1}^\infty A_i .   
 \]
   Then the natural maps $A_j \to A$ induce 
  an isomorphism 
  \[
  \oplus _i KO_n(A_i) \overset{\cong}\longrightarrow KO_n(A) 
  \]
 for each $n$.

 In addition, $KO_* $ satisfies the following properties:

 \vspace{.1in}
 \item{\underline{Stability}} \index{stability}  For all real $C^*$-algebras $A$, positive integers $k$ ,and integers $n$, 
   \[
   KO_n(A) \cong KO_n( M_k(A) )   \cong KO_n(A \otimes \KK \sr  ).
   \]

 \vspace{.1in}
\item {\underline{Periodicity:}} \index{Bott Periodicity} For any real $C^*$-algebra $A$ and $j \in \ZZ$, there is a natural Bott Periodicity isomorphism 
\[
KO_{j+8}(A) \cong KO_j(A).
 \]
 Stated in homotopy terms, Periodicity asserts that $\Omega ^8\oo \simeq \oo$. 
 
 \end{enumerate}
  \end{Thm}
  
  Here are the analogous properties for complex $C^*$-algebras and their $K$-theory.
    Recall that $\mathscr{C}^*_{\rm{complex}}$ denote the category of complex $C^*$-algebras 
  and their homomorphisms.

  \begin{Thm} \label{T:K_*}
    \index{homology theory}
 Complex  $K$-theory $K_*(A)$  is  a $\ZZ$-graded additive homology theory defined on the category 
 $\mathscr{C}^*_{\rm{complex}} $
  and taking values in $\ZZ$-graded abelian groups. That is, it satisfies
  the following axioms:
  \begin{enumerate}

  \item {\underline{Homotopy Axiom}}. \index{homotopy axiom}  Let $H : A \to IB$ be a homotopy from $f = p_0 H $ to $g = p_1 H$ 
 in  $\mathscr{C}^*_{\rm{complex}}$. Then 
 \[
 f_* = g_* : K_n(A) \longrightarrow K_n(B) \quad \text{for all n}.
 \]
  \vspace{.1in}
\item  {\underline{Exactness Axiom.}} \index{exactness axiom} Let
 \[
 0 \longrightarrow J \overset{i}\longrightarrow A \overset{j}\longrightarrow A/J  \longrightarrow 0
 \]
 be a short exact sequence in  $\mathscr{C}^*_{\rm{complex}} $.  Then for each $n$ there is a map 
 \[
 \partial : K_n(A/J) \longrightarrow K_{n-1}(J)
 \]
 and a long exact sequence 
 \[
\dots \to K_n(J) \overset{i_*}\longrightarrow K_n(A)  \overset{j_*}\longrightarrow K_n(A/J )  \overset{\partial }\longrightarrow K_{n-1}(J ) \to \cdots 
 \]
  The map $\partial $ is natural with respect to morphisms of short exact sequences. 
 
   \vspace{.1in}
 \item  {\underline{Additivity Axiom.}}  \index{additivity axiom} Let $A = \oplus _{i = 1}^\infty A_i$ in $ \mathscr{C}^*_{\rm{complex}}$.  Then the natural maps $A_j \to A$ induce 
  an isomorphism 
  \[
  \oplus _i K_n(A_i) \overset{\cong}\longrightarrow K_n(A) 
  \]
 for each $n$.

 In addition, $K_* $ satisfies the following properties:
 
 \vspace{.1in}
 \item{\underline{Matrix Stability:}} \index{stability}  For all complex $C^*$-algebras $A$ and positive integers $k$ and integers $j$, 
   \[
   K_j(A) \cong K_j( M_k(A) )   \cong K_j(A \otimes \KK \sc  )
   \]

\vglue ,1in 
\item {\underline{Periodicity:}}  \index{Bott Periodicity} For any complex $C^*$-algebra $A$, there is a natural Bott Periodicity isomorphism 
\[
K_{j+2}(A) \cong K_jA).
 \]
  Stated in homotopy terms, Periodicity asserts that $\Omega ^2\uu \simeq \uu$.

 \end{enumerate}
  \end{Thm}

\begin{Rem} (M. Karoubi)

\index{K-theory@$K$-theory!quaternionic}
We know from Corollary  \ref{C:Hhomology} and Corollary \ref{C:KH}  that $KH_*$ is indeed a homology theory.   It is periodic of period 4; the isomorphism is given by tensoring with the class of the 
identity in $KO_0(\HH )$. 
   In Remark \ref{Karoubi1} we noted that if $A$ is a real commutative $C^*$-algebra then 
   there is a natural map
   \[
KH_r(A) \otimes KH_s(A) \longrightarrow   KO_{r+s}(A).
\]
Take $A = \RR $ and $r = s = 4.$ Then we have
  \[
KH_4(A) \otimes KH_4(A) \longrightarrow   KO_{8}(A).
\]

  Now 
  \[
  KH_4(\RR ) \cong \pi_3 (GL_\infty (\HH) )
  \]
  by the first definition of $KO_*$.  By CW-approximation (see Section~\ref{top spaces}), 
  \[
  \pi_3 (GL_\infty (\HH)) \cong \pi_3 (GL_1 (\HH)) \cong \pi_3(S^3) \cong \ZZ
  \]
  The resulting map 
   \[
\ZZ \cong KH_4(A) \otimes KH_4(A) \longrightarrow   KO_{8}(A) \cong \ZZ.
\]
sends $1\otimes 1$ to the Bott generator of $KO_8(\RR )$.
 
 For further information on this matter, and in particular the relationship of $KH_*$ to symplectic $K$-theory, please see Karoubi and Weibel \cite{KL}.
  \end{Rem}

  \newpage

  \part{A Deeper look at the Structure of $KO$ and $K$}
  \vglue .3in

    \section{\bf{United $K$-theory -- Connecting Real and Complex $K$-theory}} \label{Section-United}
   
   \vglue .3in

Up to this point, we have been studying $K$-theory for complex $C \sp *$-algebras and $K$-theory for real $C \sp *$-algebras more or less in parallel, comparing and contrasting their behavior as we go. As it turns out, there is an important relationship between the $K$-theory for a real $C \sp * $-algebra $A$ and the $K$-theory of its complexification $A\sc$, which we will explore.
To be precise, we are considering the relationship between these two functors on the category of real $C \sp *$-algebras
\begin{align*}
A & \rightsquigarrow KO_*(A) \\
A & \rightsquigarrow K_*(A\sc) \; 
\end{align*}
and we consider both of these as $\ZZ$-graded abelian groups (with periodicity).
This relationship incluces natural transformations from one to the other, algebraic relations between these transformations, and a long exact sequence that involves both invariants. In this section we will describe this structure and then show how it can be packaged into a single invariant called ``united $K$-theory". United $K$-theory comes in two variations: the first is ``CR united $K$-theory'', which involves the two functors $KO_*(A)$ and $K_*(A\sc)$, and the second is ``CRT united $K$-theory" which also involves a third invariant called self-conjugate $K$-theory
$$A  \rightsquigarrow KT_*(A) \; .$$

For most of this information one should refer to Karoubi \cite{K},  Bousfield \cite{Bousfield}, Boersema \cite{BoersemaKT}.

We will use the language of  {\it natural transformations} of functors from category theory throughout this section, so we record the definition here. For this definition, it helps to have in mind functors from the category of real or complex $C \sp *$-algebras to the category of groups or rings. Accordingly, in this definition we will use the word ``homomorphism" for all the maps in the category $\mathcal{C}_2$.

\begin{Def}[cf. Section~A.3 of \cite{Weibel:HomAlg}] \label{Def:nat-trans} \index{natural transformations}
Let $\mathcal{F}$ and $\mathcal{G}$ be functors from a category {$\mathcal{C}_1$} to another category $\mathcal{C}_2$.
Then a natural transformation $\theta$ from $\mathcal{F}$ to $\mathcal{G}$ assigns to each object $A \in \mathcal{C}_1$ a homomorphism 
$$\theta_A \colon \mathcal{F}(A) \rightarrow \mathcal{G}(A) \; $$ 
in the category $\mathcal{C}_2$.
Furthermore, for each map
$\alpha \colon A \rightarrow B$
of objects in $\mathcal{C}_1$
the diagram 
  \[ \xymatrix{
& \mathcal{F}(A) \ar[r]^{\mathcal{F}(\alpha) } \ar[d]^{\theta_A} 
& \mathcal{F}(B)  \ar[d]^{\theta_B} \\
& \mathcal{G}(A) \ar[r]^{\mathcal{G}(\alpha) } 
& \mathcal{G}(B)
}  \]
must commute (the diagram consists of objects and homomorphisms in $\mathcal{C}_2$).
\end{Def}

We first remind the reader of the structure of real $K$-theory by itself. For any real $C \sp *$-algebra $A$, we know from Proposition~\ref{prop-module}
that $KO_*(A)$ is a $\ZZ$-graded module over the ring $KO_*(\RR)$, 
 \[
  \begin{tabular} { c c c c c c c c c}
  n & $0$ & $1$ & $2$ & $3$ & $4$   & $5$    &  $6$    & $7$ \\
  \hline
  $KO_n(\RR)$ &  $ \ZZ$ & $ \ZZ_2$ &  $ \ZZ_2$ & $0$ & $ \ZZ$ & $0$ & $0$ & $0$   \\
    \end{tabular}
   \]
so that multiplication by elements of $KO_*(\RR)$ define graded natural transformations on $KO_*(A)$ for any $A$. For example, for the non-trivial element $\eta \in KO_1(\RR)$ and for any element $x \in KO_n(A)$ we have $\eta \cdot x \in KO_{n+1}(A)$ using the product described in Proposition~\ref{prop-module}.
This gives the natural transformation $\eta_n$ displayed below. In the same way we obtain natural transformations corresponding to $\xi \in KO_4(\RR)$ and $\beta\so\ \in KO_8(\RR)$.

 In summary we have graded natural transformations
\begin{align*}
\eta_n &\colon KO_n(A) \rightarrow KO_{n+1}(A) \;  \\
\xi_n &\colon KO_n(A) \rightarrow KO_{n+4}(A) \;  \\
(\beta\so)_n &\colon KO_n(A) \rightarrow KO_{n+8}(A)  \; 
\end{align*}
where $\beta\so$ is an isomorphism
and these natural transformations are subject to the same relations listed in Theorem~\ref{thm:KO-relations}:
\begin{align*}
2 \eta & = 0 & \eta^3 & = 0 & \xi^2 &= 4 \beta \so &
\eta \xi &= 0
\end{align*}
 
For any complex $C \sp *$-algebra, $K_*(B)$ is a module over the ring $K_*(\CC)$. Thus in particular, $K_*(A\sc)$ is a module over $K_*(\CC)$. The structure of $K_*(\CC)$ is given by
  \[
  \begin{tabular} { c c c c c c c c c}
  n & $0$ & $1$ & $2$ & $3$ & $4$   & $5$    &  $6$    & $7$ \\
  \hline
  $K_n(\CC)$ &  $ \ZZ$ & $ 0 $   & $\ZZ$      &  $0$    & $\ZZ$  & $0$    & $ \ZZ$   & $0$     \\
    \end{tabular}
   \]
 with period 2. 
 The generator of $K_0(\CC)$ is $[1]$ and the invertible element $\beta \in K_2(\CC)$ implements the Bott Periodicity map on $K_*(A)$. 
 Therefore, there is another natural transformation
 $$\beta_n \colon K_n(A\sc) \rightarrow K_{n+2}(A\sc)$$
 for any real $C \sp *$-algebra, which is the usual Bott isomorphism.


We now define the three additional natural transformations which connect the two different versions of $K$-theory.

\begin{Pro} 
\index{complexification}
\index{realification}
Let $A$ be a real $C \sp *$-algebra. Then there exist natural transformations
\begin{align*}
c_n &\colon KO_n(A) \rightarrow K_n(A\sc ) \\
r_n &\colon K_n(A\sc ) \rightarrow KO_n(A) \\
\psi_n &\colon K_n(A\sc ) \rightarrow K_n(A\sc ) \\
\end{align*}
\end{Pro}

\begin{proof}
The complexification map $c \colon KO_*(A) \rightarrow K_*(A\sc) $ is defined in Remark~\ref{ComplexificationMap}, as the map induced on $K$-theory by the inclusion map $A \hookrightarrow A\sc$.
For a homomorphism of real $C \sp *$-algebras $\alpha \colon A \rightarrow B$, we have a commutative diagram
\[ \xymatrix{
A \ar[r]^\alpha \ar@{^{(}->}[d] &  B \ar@{^{(}->}[d]  \\
A\sc \ar[r]^{\alpha\sc} & B \sc
}\]
which induces the commutative diagram
\[ \xymatrix{
KO_n(A) \ar[r]^{\alpha_*} \ar[d]^{c_n} &  KO_n(B) \ar[d]^{c_n}  \\
K_n(A\sc) \ar[r]^{\alpha_*} & K_n(B \sc)
}\]
The naturality of $c_n$ follows from this.

There is also a standard inclusion $*$-homomorphism $\CC \rightarrow M_2(\RR)$ given by 
\[
x + iy \mapsto \begin{bmatrix} x & -y \\ y & x \end{bmatrix}  \; .
\]
 Tensoring by $A$, we obtain the homomorphism $A\sc \rightarrow M_2(A)$, which induces a map in $K$-theory so we have, for each $n$,
\[
 r_n \colon K_n(A\sc ) \rightarrow KO_n(M_2(A)) \cong KO_n(A).
 \]  
 The conjugate-linear map $A\sc \rightarrow A\sc$ defined by $a + ib \mapsto a - ib$ induces the map $\psi$.
 
 The naturality of $r_n$ and $\psi_n$ follow as for $c_n$.
 \end{proof}
  
  \begin{Rem}
  Even though $KO_n(A)$, and $K_n(\sc)$ are periodic, it it important to consider both to be $\ZZ$-graded modules here, so that both are defined for all $n \in \ZZ$. This facilitates being able to write the natural transformations above for all $n$. It would be awkward to try to express the natural transformations if $KO_n(A)$ was only defined for $0 \leq n \leq 7$ and if $K_n(A\sc)$ was only defined for $n = 0,1$.    
  \end{Rem}
  
 For example, if $A = \RR$ then we have 
$ KO_0(\RR) = \ZZ$ and $K_0(\RR\sc) = K_0(\CC) = \ZZ$. As both groups are generated by the class of the projection $1 \in \RR \subset \CC$, it follows that $c_0 \colon KO_0(\RR) \rightarrow K_0(\RR\sc)$ is an isomorphism. As the identity $1\sc \in \CC$ maps to the element $\diag(1,1) = 1_2 \in M_2(\RR)$ under the inclusion $\CC \subset M_2(\RR)$, the map $r_0 \colon K_0(\RR\sc) \rightarrow KO_0(\RR)$ is multiplication by 2.
 
 If $A = \HH$, we also have $KO_0(\HH) = K_0(\HH\sc ) = \ZZ$. \index{quaternions} However 
 \[
 c_0 \colon KO_0(\HH) \rightarrow K_0(\HH\sc)
 \]
  is multiplication by 2. Indeed, $KO_0(\HH)$ is generated by the class of the identity projection $1_\HH \in \HH$. But under the standard inclusion $\HH \subset \HH\sc = M_2(\CC)$, the identity of $\HH$ corresponds to the projection $1_2 \in M_2(\CC)$ and $[1_2]$ is twice the generator of $K_0(\CC)$. Then it follows from the following proposition that 
  \[
  r_0 \colon K_0(\HH\sc) \rightarrow KO_0(\HH)
  \]
   is an isomorphism.
   
   
 \begin{Pro} \label{CR-relations1}
 The natural transformations defined above satisfy the following relations:
\begin{align*}
rc &=  2 & cr &= 1 + \psi &  \psi^2 &= 1 & \psi c & = c & r \psi &= r  \\
\end{align*}
 \end{Pro}
  
 \begin{proof}
 \begin{enumerate}
 \item
 The composition 
$
 A \rightarrow A\sc \rightarrow M_2(A)
 $
  satisfies
   \[
  p \mapsto \begin{bmatrix} p & 0 \\ 0 & p \end{bmatrix}  \; ,
  \]
   and from this it follows that $rc([p]) = 2[p]$ for any projection in $M_n(A)$. 
   
   \item
   
   The composition $A\sc \rightarrow M_2(A) \rightarrow M_2(A\sc)$
 satisfies 
 \[
 a + ib \mapsto \begin{bmatrix} a & -b \\ b & a \end{bmatrix} 
 \]
  for any projection $p = a +ib$. Conjugating the latter by the unitary
 \[
 u = \tfrac{1}{\sqrt{2}} \begin{bmatrix} 1 & i \\ i & 1  \end{bmatrix} 
 \]
  yields 
  \[
   \begin{bmatrix} a + ib & 0 \\ 0 & a - ib \end{bmatrix} \; .
  \]
   Therefore 
   \[
   cr([p]) = (1 + \psi)([p])
   \]
    for any projection.
  \item
 The relation $\psi^2 = 1$ follows directly from the fact that $a + ib \mapsto a - ib$ is an involution (that is, it squares to the identity) on $A\sc$. 
   \item
 It is easy to check the relation $\psi c = c$ holds exactly on the level of $C^*$-algebra maps. 
   \item
 Finally to verify the relation $r \psi = r$, note that if $p = a + ib$ is a projection in $A\sc$ (or in $M_n(A\sc)$ for some $n$) then
 \[
 r \psi(p) =  \begin{bmatrix} a & b \\ -b & a \end{bmatrix} \; .
 \]
  But this projection is unitarily equivalent to
 $ \begin{bmatrix} a & -b \\ b & a \end{bmatrix} $, by conjugating by $u =  \begin{bmatrix} 1 & 0 \\ 0 & -1 \end{bmatrix}$. 
 It follows that $r \psi[p] = r[p]$ for any projection in $M_n(A\sc)$.
 \end{enumerate}
 \end{proof}

 With these additional operations, we have further relations.

 \begin{Pro} \label{CR-relations2}
 The natural transformations defined above satisfy the following relations:
\begin{align*}
c \beta\so &= \beta^4 c &  \beta\so r & = r \beta^4  &  \psi \beta &= - \beta \psi & \\
\eta r &= 0 & c \eta &= 0 & r \beta c &= \eta^2 & r \beta^{-1} c &= 0 \\
\end{align*}
\end{Pro}
 
 See \cite{Hewitt} or \cite{BoersemaKT}, but especially 
see \cite{Bousfield} for further discussion of these relations. 
We have now described all of the structure necessary for the definition of united $K$-theory.

\begin{Def} \index{united $K$-theory!$CR$ united $K$-theory}
\index[notation]{Kcrr@$K\crr(A)$}
Let $A$ be a real $C \sp *$-algebra. Then {\it CR united $K$-theory} refers to the invariant
$$
	K\crr(A) = \{ KO_*(A) , K_*(A\sc), \Theta\crr \}  
$$
\end{Def}

\begin{Rem}
\begin{enumerate}
\item In this definition, the $\Z$-graded group $KO_*(A)$ implicitly carries structure of a module over the ring $KO_*(\RR)$ which means that it has all the structure described in Proposition~\ref{prop-module} and Theorem~\ref{thm:KO-relations}, including the natural transformations $\eta, \xi, \beta\so$.
\item The $\ZZ$-graded group $K_*(A\sc)$ is a module over $K_*(\CC)$ which really only means that it is a $\ZZ$-graded abelian group with a periodicity isomorphism $\beta$ of degree 2.
\item $\Theta\crr = \{c, r, \psi\}$ represents the set of natural transformation among $KO_*(A)$ and $K_*(A\sc)$. These satisfy the relations 
stated here in Propositions~\ref{CR-relations1} and \ref{CR-relations2}.

\end{enumerate}
\end{Rem}

The long exact sequence given in the next proposition can be found in Section~1.11 of  Schr\"oder \cite{Schr} or Theorem~1.18 of \cite{BoersemaKT}, in the setting of real $C \sp *$-algebras. In the topological setting, it can be found in
 Theorem~5.18 of \cite{K} or Theorem~1.11 of \cite{Bousfield}. We find this long exact sequence indispensable for $K$-theory calculations.

\begin{Pro} \label{Pro-LES}
For any real $C^*$-algebra $A$, there is a long exact sequence
$$ \cdots \rightarrow KO_n(A) \xrightarrow{~\eta~} KO_{n+1}(A)  \xrightarrow{ ~c ~} K_{n+1}(A\sc ) 
	\xrightarrow {~ r \beta^{-1} } KO_{n-1}(A) \rightarrow \cdots $$
\end{Pro}
 
 \begin{proof}
 We start with the inclusion homomorphism $c \colon \RR \rightarrow \CC$. 
 Then the mapping cone construction from Section~\ref{Section:HomotopyConstructions} gives a fibration
 $$0 \rightarrow S\CC \hookrightarrow Cc \rightarrow \RR \rightarrow 0$$
  which gives rise to a long exact sequence 
$$
\cdots \longrightarrow  KO_{n+1}(A \otimes \CC) \rightarrow
KO_{n}(A \otimes Cc) \rightarrow 
KO_{n}(A) \xrightarrow{\partial}
KO_{n}(A\otimes \CC) \rightarrow \cdots \, . 
$$
We know that $\partial = c_*$
and it still remains to identify the other maps. For that we refer to the proof of Theorem~1.18 of  Boersema \cite{BoersemaKT}    for details. 
\end{proof}

The utility of this combined structure is illustrated by the following corollaries.
 
   \begin{Cor} \label{Cor-isomorphism} Suppose that $f: A \to B$ is a homomorphism of real $C^*$-algebras. Then the following are equivalent:
  \begin{enumerate}
  \item The map 
  \[
  f_* : KO_n(A) \longrightarrow KO_n(B)
  \]
  is an isomorphism for all $n$. 
  
  \item The map 
  \[
  f_* : K_n(A\sc) \longrightarrow K_n(B\sc)
  \]
  is an isomorphism for all $n$. 
 \end{enumerate}
\end{Cor}

\begin{proof}
Consider the map $f_*$ applied to the short exact sequence of Proposition~\ref{Pro-LES}.
\[ \xymatrix{
\cdots \ar[r]
& KO_{n-1}(A) \ar[r]^\eta \ar[d]^{f_*}
& KO_{n}(A) \ar[r]^c \ar[d]^{f_*}
& K_n(A\sc) \ar[r]^{r \beta^{-1}} \ar[d]^{f_*}
& KO_{n-2}(A) \ar[r] \ar[d]^{f_*}
&\cdots   \\
\cdots \ar[r]
& KO_{n-1}(B) \ar[r]^\eta 
& KO_{n}(B) \ar[r]^c 
& K_n(B\sc ) \ar[r]^{r \beta^{-1}}  
& KO_{n-2}(B) \ar[r]
&\cdots   
}\]

If $f_* \colon KO_n(A) \rightarrow KO_n(B)$ is an isomorphism for all $n$, then the five-lemma implies that $f_* \colon K_n(A\sc ) \rightarrow K_n(B\sc )$ is an isomorphism for all $n$.
For the converse, we refer to Section~2.3 of Bousfield \cite{Bousfield}.
\end{proof}

\begin{Cor}
Suppose that $A$ is a real $C^*$-algebra. Then the following are equivalent:
\begin{enumerate}
\item $KO_n(A) = 0 $ for all $n$. 

\item $K_n(A\sc) = 0 $ for all $n$.   
\end{enumerate}
\end{Cor}

\begin{proof}
If $KO_n(A) = 0$ for all $n$, then the long exact sequence of Proposition~\ref{Pro-LES} immediately implies that $K_n(A\sc ) = 0$ for all $n$.
Conversely, suppose that $K_n(A\sc ) = 0$ for all $n$. In this case, the long exact sequence implies that $\eta_n$ is an isomorphism for all $n$. However, the relation $\eta^3 = 0$ then implies that $KO_n(A) = 0$ for all $n$.
\end{proof}
 
We now introduce the new version of $K$-theory, known as ``self-conjugate $K$-theory", which will be necessary for CRT united $K$-theory.

\begin{Def} 
\index{K-theory@$K$-theory!self-conjugate $K$-theory}
\index[notation]{KTA@$KT_*(A)$}
Let $A$ be a real $C \sp *$-algebra. Then {\it self-conjugate $K$-theory} is the invariant
$$KT_*(A) = KO_*(A \otimes T) \; .$$
Here $T$ is the commutative $C \sp *$-algebra
$$T = \real(S^1, \sigma)$$ and $\sigma(x) = -x$ is the antipodal map on $S^1$ (see Example~\ref{Examples-RealC*-algebras}, Part (11)).
\end{Def}

Self-conjugate $K$-theory is a homology theory by Proposition~\ref{Proposition:homologytheory}. Furthermore, since $T$ is unital and commutative, $KO_*(T)$ is a unital ring and $KT_*(A)$ is a module over the ring $KO_*(T)$. The structure of this ring is given by the proposition below.
 
 \begin{Pro} \label{Prop-KO(T)} 
 The real $K$-theory of $T$ has period 4 and in degrees $i = 0,1,2,3$ it is given by

 \[
  \begin{tabular} { c c c c c c c c c}
  n & $0$ & $1$ & $2$ & $3$  \\
  \hline
  $KO_n(T)$ &  $ \ZZ$ & $ \ZZ_2$ & $0$  & $\ZZ$         \\
    \end{tabular}
   \]
 Since $T$ is commutative, $KO_*(T)$ is a ring. The generator $\beta_T \in KO_4(T) \cong \Z$ is an invertible element in the ring, implementing the periodicity.
 \end{Pro}
 
 \begin{proof}
 There is a unital surjective map $\ev \colon T \rightarrow \CC$ defined by point-evaluation
  and an associated short exact sequence of real $C \sp *$-algebras
$$0 \rightarrow S \CC \rightarrow T \xrightarrow{\ev} \CC \rightarrow 0\; .$$
 Therefore, there is an induced long exact sequence
 $$\cdots \rightarrow KO_n(S\CC) \rightarrow KO_n(T) \xrightarrow{\ev_*} KO_n(\CC) \rightarrow KO_{n-1}(S\CC) \rightarrow \cdots \; . $$
 
 Since $\ev$ is unital, it follows that 
 $$(\ev_*)_0  \colon KO_0(T) \rightarrow KO_0(\CC)$$
 is surjective. Also since $KO_0(S\CC) = 0$, the sequence implies that $(\ev_*)_0$ is an isomorphism. This is enough to get the calculation started. The rest of the calculation can be completed using this long exact sequence, together with the long exact sequence
 of Proposition~\ref{Pro-LES} taken with $A = T$. See Corollary~1.6 of \cite{BoersemaKT} for details.
 \end{proof}
 
 \begin{Rem}
 We denote the non-zero (additive) generators of $KT_*(\RR) \cong KO_*(T)$ by
 \begin{align*} 
 \eta\sst &\in KO_1(T) \cong \Z_2 \\
 \omega &\in KO_3(T) \cong \Z \\
 \beta\sst &\in KO_4(T) \cong \Z
 \end{align*}
 Then in the graded commutative ring $KT_*(\RR) = KO_*(T)$  the element $\beta\sst$ is invertible, and the elements $\eta\sst$ and $\omega$ satisfy $\eta\sst^2 = 0$, $\omega^2 = 0$,  and $\eta\sst \omega = 0$.
 These produce natural transformations
 \begin{align*}
 (\eta\sst)_i & \colon KT_i(A) \rightarrow KT_{i+1}(A) \\
 \omega_i & \colon KT_i(A) \rightarrow KT_{i+3}(A) \\
 (\beta\sst)_i & \colon KT_i(A) \rightarrow KT_{i+4}(A)
 \end{align*}
 where $\beta\sst$ is an isomorphism.
  \end{Rem}
 
 \begin{Rem}
Looking ahead, by Theorem~\ref{T:AtiyahR} there is an isomorphism $KO_*(T) \cong KR(S^1, \sigma)$ where $KR$ refers to Atiayah's ``Real $K$-theory". In particular, $KO_0(T)$ is isomorphic to the Grothendieck group of the semi-group of equivalence classes of pairs $(E, \tau)$ where $E$ is a vector bundle over $S^1$ and $\tau$ is a conjugate-linear involution on $E$ which intertwines with the involution $\sigma$ on $S^1$. See Section~\ref{Section:KRTheory} and for details.
\end{Rem}

\begin{Def} \index{united $K$-theory!$CRT$ united $K$-theory}
\index[notation]{Kcrt@$K\crt(A)$}
Let $A$ be a real $C \sp *$-algebra. Then  {\it united $K$-theory} or {\it CRT united $K$-theory} refers to the invariant
$$
	K\crt(A) = \{ KO_*(A) , K_*(A\sc), KO_*(A \otimes T), \Theta\crt \} 
$$
\end{Def}

\begin{Rem}
 $\Theta\crt$ represents the set of natural transformations among $KO_*(A)$, $K_*(A\sc)$, and $KT_*(A)$. This includes the transformations in $\Theta\crr$ previously discussed, as well as several others (namely: $\gamma, \zeta, \tau, \varepsilon$) and there is a long list of relations that these natural transformations satisfy, analogous to the list in Propositions~\ref{CR-relations1} and \ref{CR-relations2}. Formally, the algebra of natural transformations is parametrized by $KK\pr_*(B_1, B_2)$ where $B_i \in \{\RR, \CC, T\}$ and the complete set of relations is encoded by the product structure of maps
$$KK\pr_*(B_1, B_2) \otimes KK\pr_*(B_2, B_3) \rightarrow KK\pr_*(B_1, B_3) \; .$$
We refer the reader to Section~\ref{KKSection} for an introduction to $KK$-theory and this product in $KK$-theory.
These natural transformations and relations are fully described in \cite{Bousfield} and \cite{BoersemaKT}.

In fact, it turns out that the map $\ev$ in the proof of Proposition~\ref{Prop-KO(T)} induces the natural transformation $\zeta$, and the inclusion map $S\CC \hookrightarrow T$ induces the natural transformation $\gamma$. See \cite{BoersemaKT}.
\end{Rem}

\begin{Rem}
We wish to clarify the relationship between $CR$ united $K$-theory and $CRT$ united $K$-theory.
$CR$ united $K$-theory is clearly cleaner and more straightforward. In many situations $CR$ united $K$-theory is sufficient, in light of the crucial results in the dissertation of Beatrice Hewitt in \cite{Hewitt}. Those results imply that if $A$ and $B$ are real $C \sp *$-algebras, then 
$$K\crr(A) \cong K\crr(B) \quad \text{if and only if} \quad K\crt(A) \cong K\crt(B) \; .$$
Therefore, there is no information lost in restricting to $CR$ united $K$-theory. Furthermore, we have that $K\crr(A) \cong K\crr(B)$ holds if and only if $A$ and $B$ are $KK$-equivalent, see Definition~\ref{Def-KKequiv} and Corollary~\ref{Cor-KKequiv} below.

Though $CR$ united $K$-theory is sufficient for much of what we do, there are many situations where we still need to consider $CRT$ united $K$-theory. In particular, the most general form of the K\"unneth Theorem, in Section~\ref{Section-Kunneth}, is articulated in terms of $CRT$ united $K$-theory and cannot be reduced to or expressed in terms of $CR$ united $K$-theory. The realization theorem below for real $C \sp *$-algebras, Theorem~\ref{Thm-realization}, also requires the $CRT$ united $K$-theory.
\end{Rem}

The next Proposition, which we state here without proof, is found in Section~1.4 of \cite{BoersemaKT}. Then the Corollaries below follow from results in Section~2.3 of \cite{Bousfield}.

\begin{Pro} \label{Pro-LES2}
For any real $C^*$-algebra $A$, there are two long exact sequences
$$ \cdots \rightarrow KO_n(A) \xrightarrow{~\eta^2~} KO_{n+2}(A)  \xrightarrow{ ~\varepsilon ~} KT_{n+2}(A) 
	\xrightarrow {~ \tau \beta\sst^{-1} } KO_{n-1}(A) \rightarrow \cdots $$
and
$$  \cdots \rightarrow K_n(A\sc) \xrightarrow{~1 - \psi ~} K_{n}(A\sc)  \xrightarrow{ ~\gamma ~} KT_{n-1}(A) 
	\xrightarrow {~ \zeta } K_{n-1}(A\sc) \rightarrow \cdots $$
\end{Pro}

  \begin{Cor} \label{Cor-isomorphism} Suppose that $f: A \to B$ is a homomorphism of real $C^*$-algebras. Then the following are equivalent:
  \begin{enumerate}
  \item The map 
  \[
  f_* : KO_n(A) \longrightarrow KO_n(B)
  \]
  is an isomorphism for all $n$. 
  
  \item The map 
  \[
  f_* : K_n(A\sc) \longrightarrow K_n(B\sc)
  \]
  is an isomorphism for all $n$. 
  
   \item The map 
  \[
  f_* : KT_n(A) \longrightarrow KT_n(B)
  \]
  is an isomorphism for all $n$. 

 \end{enumerate}
\end{Cor}

\begin{Cor}
Suppose that $A$ is a real $C^*$-algebra. Then the following are equivalent:
\begin{enumerate}
\item $KO_n(A) = 0 $ for all $n$. 
\item $K_n(A\sc) = 0 $ for all $n$.   
\item $KT_n(A) = 0$ for all $n$.
\end{enumerate}
\end{Cor}

Finally, we finish this section with the description of the target category for united $K$-theory, as defined in \cite{Bousfield} and \cite{BoersemaKT}.
 The target category for $CR$ united $K$-theory is the category of abstract $CR$-modules and similarly the target category of $CRT$ united $K$-theory is the category of abstract $CRT$-modules. 
 
 \begin{Def} \label{Def-CRTmodule}

 \begin{enumerate}
 \index{CR-module@$CR$-module}
\item A $CR$-module is an object 
 $$M = \{ M\po, M\pu, c, r, \psi \}$$ where $M\po$ is a graded module over the ring $KO_*(\RR)$ and $M\pu$ is a graded module over the ring $K_*(\CC)$; and $r,c,\psi$ are homomorphisms
 \begin{align*}
 c & \colon M\po \rightarrow M\pu \\
r & \colon M\pu \rightarrow M\po \\
 \psi & \colon M \pu \rightarrow M\pu
 \end{align*}
 which satisfy all the relations given in Propositions~\ref{CR-relations1} and \ref{CR-relations2}, as well as many others. The complete list of relations is given in Section~8 of \cite{Hewitt}.
 
 Let $M = (M\po, M\pu, c,r, \psi)$ and $N = (N\po, N\pu, c,r, \psi)$ be $CR$-modules. Then a $CR$-module homomorphism $\alpha \colon M \rightarrow N$ is a pair of maps
 \begin{align*}
 \alpha\po \colon M\po \rightarrow N\po \\
 \alpha\pu \colon M\pu \rightarrow N\pu
 \end{align*} 
 that intertwine with the natural transformations in the obvious way. 
 
 The category $\mathcal{CR}$ is the category whose objects consist of all $CR$-modules and whose morphisms are all $CRT$-module homomorphisms.
  
 \item  \index{CRT-module@$CRT$-module}
 A $CRT$-module is an object 
 $$M = \{ M\po, M\pu, M\ppt, r, c, \psi, \gamma, \tau, \varepsilon, \zeta  \}$$ where $M\po$ is a graded module over the ring $KO_*(\RR)$, $M\pu$ is a graded module over the ring $K_*(\CC)$, and $M\ppt$ is a graded module over the ring $K_*(T)$; and there rest of the objects are homomorphisms among the three components of $K$-theory, which satisfy an appropriate list of relations, which include the relations in Propositions~\ref{CR-relations1} and \ref{CR-relations2}. The complete list of relations is given in Section~1.9 of \cite{Bousfield}.
 
 A $CRT$-module homomorphism is defined similarly to the definition of a $CR$-module homomorphism above, and the category of all $CRT$-modules is denoted by $\mathcal{CRT}$.
  \end{enumerate}
 \end{Def}
 
 \begin{Thm} (Theorem~1 of \cite{BoersemaSurj}) \label{Thm-realization}
 Let $M = \{ M\po, M\pu, M\ppt, r, c, \psi, \gamma, \tau, \varepsilon, \zeta  \}$ be a countable $CRT$-module such that three sequences are exact:
 \[ \cdots \rightarrow M\po_i \xrightarrow{\eta} M\po_{i+1} \xrightarrow{c} M_{i+1}\pu \xrightarrow{r \beta\su^{-1}} M_{i-1}\po \rightarrow \cdots \]
 \[ \cdots \rightarrow M\po_i \xrightarrow{\eta^2} M\po_{i+2} \xrightarrow{\varepsilon} M_{i+2}\ppt \xrightarrow{\tau \beta\sst^{-1}} M_{i-1}\po \rightarrow \cdots \] 
  \[ \cdots \rightarrow M\pu_i  \xrightarrow{~1 - \psi ~} M\pu_i  \xrightarrow{ ~\gamma ~} M\ppt_{i-1} 
	\xrightarrow {~ \zeta } M\pu_{i-1} \rightarrow \cdots  \]
 Then there exists a real $C \sp *$-algebra $A$ such that $K\crt(A) \cong M$.
 \end{Thm}
 
 Note if the sequences in Theorem~\ref{Thm-realization} are exact then the $CRT$-module $M$ is said to be {\it acyclic}. 
 	\index{CRT-module@$CRT$-module!acyclic $CRT$-module} 
 In fact, the $C \sp *$-algebra $A$ can be taken to be purely infinite and c-simple, and in the pre-bootstrap category (see Theorem~\ref{Thm-realization2} and the terminology defined leading up to that). 
\newpage
      \section{\bf{Examples and Calculations }} \label{Section-Examples}

 In this section, we present a number of different examples of the results of calculations of $K$-theory for real $C \sp *$-algebras and topological spaces with involution. Because of the structure discussed in the previous section, it is convenient to consider both the ``real" and ``complex" $K$-theory together. That is, for a real $C \sp *$-algebra $A$ we will often consider both $KO_*(A)$ and $K_*(A\sc)$. 
 In Sections~\ref{Section:KOforSpaces} and \ref{Section:RelatingKO-theory}, we will talk about $K$-theory for topological spaces. For now, we note that for a compact space $X$ we have
 \[
 KO^{-*}(X) \cong KO_*(\real(X)) \qquad  {\text{and}} \qquad  K^{-*}(X) \cong K_*(\complex(X))
 \]
 which will help provide some context for some of the calculations in this section.
 

\begin{Exam}
Let $A = \RR$. Then we also $KO_n(\RR) = KO^{-n}( \{ x \})$ where $\{x\}$ designates a one-point space.
Then expanding on Theorem~\ref{Thm-classical}, we have
\[
  \begin{tabular} { c c c c c c c c c}

  n & $0$ & $1$ & $2$ & $3$ & $4$   & $5$    &  $6$    & $7$ \\

  \hline

  $KO_n(\RR)$ &  $ \ZZ$ & $ \ZZ_2 $   & $\ZZ_2$      &  $0$    & $\ZZ$  & $0$    & $ 0$   & $0$     \\

  $K_n(\RR\sc)$ &  $ \ZZ$ & $0$   & $\ZZ$  & $0$  & $\ZZ$ & $0$   & $\ZZ$ & $0$    \\

  $c_n$ & 1 & 0 & 0 & 0 & 2 & 0 & 0 & 0 \\
  $r_n$ & 2 & 0 & 1 & 0 & 1 & 0 & 0 & 0 \\
  $\psi_n$ & 1 & 0 & -1 & 0 & 1 & 0 & -1 & 0 \\
  $\eta_n$ & 1 & 1 &  0 & 0 & 0 & 0 & 0 & 0 \\
    \end{tabular}
   \]
   This table shows the groups of real and complex $K$-theory, as well as the natural transformations.
   
 In the case at hand, each group is generated by a single element, so we can express the homomorphism by a single integer. For example, the table shows that 
   $$c_4 \colon KO_4(\RR) \rightarrow KO_4(\RR\sc)$$ is the map from $\ZZ$ to $\ZZ$ given by multiplication by 2. 
   In the examples below we use matrices to represent the homomorphisms when one of the groups has more than one generator.
 
  \end{Exam}

   \begin{Exam}

We can consider $\CC$ as a real $C \sp *$-algebra. When we do that we find there is an isomorphism

\[
\CC_{\CC} \cong \CC \otimes_{\RR} \CC \cong \CC \oplus \CC
\]
(this is a special case of Part (2) of Proposition~\ref{Pro-c-simple}).

Thus  $K_*(\CC\sc ) = K_*(\CC \oplus \CC)  $, so we have the following:

\[
  \begin{tabular} { c c c c c c c c c}

  n & $0$ & $1$ & $2$ & $3$ & $4$   & $5$    &  $6$    & $7$ \\

  \hline

  $KO_n(\CC)$ &  $ \ZZ$ & $0$   & $\ZZ$  & $0$  & $\ZZ$ & $0$   & $\ZZ$ & $0$    \\

  $K_n(\CC\sc )$ &  $ \ZZ^2$ & $0$   & $\ZZ^2$  & $0$  & $\ZZ^2$ & $0$   & $\ZZ^2$ & $0$    \\

   $c_n$ & $\left( \begin{smallmatrix} 1 \\ 1 \end{smallmatrix} \right)$ & 0 &  $\left( \begin{smallmatrix} -1 \\ 1 \end{smallmatrix} \right)$& 0

   & $\left( \begin{smallmatrix} 1 \\ 1 \end{smallmatrix} \right)$ & 0 & $\left( \begin{smallmatrix} -1 \\ 1 \end{smallmatrix} \right)$ & 0 \\

  $r_n$ & $\left( \begin{smallmatrix} 1 & 1 \end{smallmatrix} \right)$ & 0 & $\left( \begin{smallmatrix} -1 & 1 \end{smallmatrix} \right)$ & 0

  & $\left( \begin{smallmatrix} 1 & 1 \end{smallmatrix} \right)$ & 0 & $\left( \begin{smallmatrix} -1 & 1 \end{smallmatrix} \right)$ & 0 \\

$\psi_n$ & $\left( \begin{smallmatrix} 0 & 1 \\ 1 & 0 \end{smallmatrix} \right)$ & 0 & $\left( \begin{smallmatrix} 0 & -1 \\ -1 & 0 \end{smallmatrix} \right)$ & 0 
	& $\left( \begin{smallmatrix} 0 & 1 \\ 1 & 0 \end{smallmatrix} \right)$ & 0 & $\left( \begin{smallmatrix} 0 & -1 \\ -1 & 0 \end{smallmatrix} \right)$ & 0 \\

  $\eta_n$ & 0 & 0 &  0 & 0 & 0 & 0 & 0 & 0 \\

    \end{tabular}
   \]

  \end{Exam}

  \begin{Exam} \index{suspension!the suspension algebra}  \label{Ex-suspension}
  
  The suspension algebra is $S\RR = \real_0(\topr )$. Since 
  \[
  KO_n(S\RR) = KO_{n+1}(\RR) \,\,\,\, \text{and}\,\,\,\, K_n(S\CC) = K_{n+1}(\CC)
  \]
  we have
 \[
  \begin{tabular} { c c c c c c c c c}
  n & $0$ & $1$ & $2$ & $3$ & $4$   & $5$    &  $6$    & $7$ \\
  \hline
  $KO_n(S\RR)$ &   $ \ZZ_2 $   & $\ZZ_2$      &  $0$    & $\ZZ$  & $0$    & $ 0$   & $0$  & $\ZZ$    \\
  $K_*( (S\RR)\sc)$ &   $0$   & $\ZZ$  & $0$  & $\ZZ$ & $0$   & $\ZZ$ & $0$ & $\ZZ$    \\
    \end{tabular}
   \]
  \end{Exam}
  
\index{desuspension} \index{suspension!the desuspension algebra} 
    \begin{Exam} \label{Ex-desuspension}
    The desuspension algebra is $S^{-1} \RR = V\RR = \real_0(i\topr )$ (recall Definition~\ref{Def-AllSus}).  
    \index[notation]{Srr@ $S^{-1} \RR$}

     Since   
     \[
     KO_n(S^{-1} \RR ) = KO_{n-1}(\RR)\,\, \,\,  \text{but}\,\,\,\, (S^{-1} \RR) \sc \simeq S\RR \sc,
     \]
      we have 
 \[
  \begin{tabular} { c c c c c c c c c}
  n & $0$ & $1$ & $2$ & $3$ & $4$   & $5$    &  $6$    & $7$ \\
  \hline
  $KO_n(S^{-1} \RR)$ & $0$ & $\ZZ$ &   $ \ZZ_2 $   & $\ZZ_2$      &  $0$    & $\ZZ$  & $0$    & $ 0$    \\
  $K_n((S^{-1} \RR)\sc ) $ &   $0$   & $\ZZ$  & $0$  & $\ZZ$ & $0$   & $\ZZ$ & $0$ & $\ZZ$    \\
    \end{tabular}
   \]
   Compare the first row with the first row of the previous (suspension) example to see the difference in how the dimension shifts.
  \end{Exam}

  \begin{Def}  \label{def-GroupSuspension}
  \index{suspension!of a graded group}
  (See also the proof to Corollary~\ref{C:KH}.)
 Suppose that $M = (M_i)$ is a graded group (this includes the situation where $M$ is a $CR$-module or a $CRT$-module).  For $k \in \Z$ we define the {\it $k$th suspension} $\Sigma^k M$ of $M$ as follows:
  \[(\Sigma^k M)_i = M_{i+k}  \; .\]
  \end{Def}
     
With this notation, Examples~\ref{Ex-suspension} and \ref{Ex-desuspension} may be generalized as follows:

\begin{Exam}
 \[
  KO_*( \real_0(\topr ^n)) = KO_*(S^n \RR) =  \Sigma^n KO_*(\RR)
  \]
   and
  \[
  KO_*( \real_0(i \topr ^n))  =  KO_*(S^{-n} \RR) =  \Sigma^{-n} KO_*(\RR).
  \]
Now let $n = p + q$ and let  $\tau^{p,q}$ be the involution on $\topr ^n$   given by the identity on $p$ coordinates and by multiplication by $-1$ on the remaining $q$ coordinates. Then  
\[
KO_*( \real_0(\topr ^n, \tau^{p,q})) = \Sigma^{p-q} KO_*(\RR).
   \]
  \end{Exam}
   
   \begin{Exam} \label{example-circle}
   Consider the algebra $A = \real(S^1)$. Then 
   $$KO^{-n}(S^1) = KO_n(\real(S^1))$$ and
 \[
  \begin{tabular} { c c c c c c c c c}
  n & $0$ & $1$ & $2$ & $3$ & $4$   & $5$    &  $6$    & $7$ \\
  \hline
  $KO_n(\real (S^1))$ & $\ZZ \oplus \ZZ_2$ &   $ \ZZ_2^2 $   & $\ZZ_2$      &  $\ZZ$    & $\ZZ$  & $0$ & $0$ & $\ZZ$   \\
  $K_n(\complex (S^1))$ &   $\ZZ$    & $\ZZ$  & $\ZZ$   & $\ZZ$ &    $\ZZ$    & $\ZZ$  & $\ZZ$   & $\ZZ$   \\
    \end{tabular}
   \]
  \end{Exam}
  
  \begin{proof}
 Regard  $S^1$ as the usual subset of the complex plane and let 
  \[
  \ev_1 \colon \real(S^1) \rightarrow \RR
  \]
   be evaluation at $1 \in S^1$. This is a surjective real $C \sp *$-algebra homomorphism with kernel isomorphic to $\real_0( \topr )$, 
   and hence there is a  short exact sequence
  \[ 0 \rightarrow \real_0(\topr ) \rightarrow \real(S^1) \xrightarrow{\ev_1} \RR \rightarrow 0 \; . \]
   Furthermore, there is a splitting homomorphism $s \colon \RR \rightarrow \real(S^1)$ which is given by constant functions. This implies that the long exact sequence
  \[ \cdots \rightarrow KO_*( \real_0(\topr )) \rightarrow KO_*( \real (S^1)) \xrightarrow{(\ev_1)_*} 
  			KO_*( \RR) \rightarrow \cdots \]
splits, and we have an isomorphism
\[ KO_*( \real(S^1)) \cong KO_*(\RR) \oplus KO_*( \real_0(\topr )) 
	\cong KO_*(\RR) \oplus \Sigma KO_*(\RR)\; . \]
	
The complex groups  $K_*( \complex(S^1))$ are computed using the same split sequence.
\end{proof}

     \begin{Exam} \label{example-circle2}
    Recall that $iS^1$ denotes the circle with the involution $z \mapsto \overline{z}$
    and 
    $$\real(iS^1) = \{f \in \complex(S^1) \mid \overline{f(z)} = f(\overline{z}) \}$$ is the corresponding real $C^*$-algebra. Then we have:
 \[
  \begin{tabular} { c c c c c c c c c}
  n & $0$ & $1$ & $2$ & $3$ & $4$   & $5$    &  $6$    & $7$ \\
  \hline
  $KO_n(\real(iS^1))$ & $ \ZZ$ & $\ZZ \oplus \ZZ_2$ &   $ \ZZ_2^2 $   & $\ZZ_2$      &  $\ZZ$    & $\ZZ$  & $0$ & $0$    \\
  $K_n(\complex (S^1))$ &   $\ZZ$    & $\ZZ$  & $\ZZ$   & $\ZZ$ &    $\ZZ$    & $\ZZ$  & $\ZZ$   & $\ZZ$   \\
    \end{tabular}
   \]
  \end{Exam}
 
   \begin{proof}
   Since $1 \in S^1$ is a fixed point of the involution $z \mapsto \overline{z}$, the condition that $\overline{f(z)} = f(\overline{z})$ 
   implies that $f(1) \in \RR$. So as in the previous calculation, we have a map $\ev_1 \colon \real(i S^1) \rightarrow \RR$; this time the kernel is isomorphic to $S^{-1} \RR = \real_0( i\RR)$. So there is a split short exact sequence  
    \[ 0 \rightarrow \real_0(i \topr ) \rightarrow \real(iS^1) \xrightarrow{\ev_1} \RR \rightarrow 0 \; . \]
   Again, as in the previous calculation, the fact that the sequence splits implies that 
   \[ KO_*( \real(iS^1)) \cong KO_*(\RR) \oplus KO_*( S^{-1} \RR)
   	\cong KO_*(\RR) \oplus \Sigma^{-1} KO_*(\RR)\;  . \]
\end{proof}

\begin{Exam} More generally, we have
\begin{align*}
 KO_*( \real(S^n)) &= (1 + \Sigma^{n}) KO_*(\RR)  \;   \\
 \text{and~} KO_*( \real(iS^n)) &= (1 + \Sigma^{-n}) KO_*(\RR)  \; .
 \end{align*}
 This follows as in Examples~\ref{example-circle} and ~\ref{example-circle2} using the split short exact sequences
 \[ 0 \rightarrow \real_0(\topr^n ) \rightarrow \real(S^n) \xrightarrow{\ev_1} \RR \rightarrow 0 \;  \]
 and
 \[ 0 \rightarrow \real_0(i \topr^n ) \rightarrow \real(i S^n) \xrightarrow{\ev_1} \RR \rightarrow 0 \; . \]
 
  \end{Exam}
  
      \begin{Exam}
     For $T = C(S^1, \sigma)$ where $\sigma$ is the antipodal involution. Expanding on the result of Proposition~\ref{Prop-KO(T)} we have (from Corollary~1.6 of     \cite{BoersemaKT})
 \[
  \begin{tabular} { c c c c c c c c c}
  n & $0$ & $1$ & $2$ & $3$ & $4$   & $5$    &  $6$    & $7$ \\
  \hline
  $KO_n(\real (S^1, \sigma))$ & $ \ZZ$ & $ \ZZ_2$ &   $ 0 $   & $\ZZ$      &   $ \ZZ$ & $ \ZZ_2$ &   $ 0 $   & $\ZZ$      \\
  $K_n(\real (S^1, \sigma))$ &   $\ZZ$    & $\ZZ$  & $\ZZ$   & $\ZZ$ &    $\ZZ$    & $\ZZ$  & $\ZZ$   & $\ZZ$   \\
  $c_n$ & 1 & 0 & 0 & 2 & 1 & 0 & 0 & 2 \\
  $r_n$ & 2 & 1 & 0 & 1 & 2 & 1 & 0 & 1 \\
  $\psi_n$ & 1 & -1 & -1 & 1 & 1 & -1 & -1 & 1 \\
  $\eta_n$ & 1 & 0 &  0 & 0 & 1 & 0 & 0 & 0 \\
    \end{tabular}
   \]
  \end{Exam}

\begin{Exam}

Let $A = \real(S^1 \times S^1)$ be the algebra of functions on the torus. Then
\[KO_*(\real (S^1 \times S^1)) = KO_*(\RR) \oplus \Sigma KO_*(\RR)^2 \oplus \Sigma^2 KO_*(\RR)  \; .\]
Thus,
\[
  \begin{tabular} { c c c c c c c c c}
  n & $0$ & $1$ & $2$ & $3$ & $4$   & $5$    &  $6$    & $7$ \\
  \hline
  $KO_n (\real(S^1 \times S^1) ) $ & $\ZZ \oplus \ZZ_2^3$ &   $ \ZZ_2^3 $   & $\ZZ \oplus \ZZ_2$      &  $\ZZ^2$    & $\ZZ$  & $0$ & $\ZZ$ & $\ZZ^2 \oplus \ZZ_2$   \\ ~ \\ \vspace{-.7cm} \\
  $K_n( \complex(S^1 \times S^1))  $ &   $\ZZ^2$    & $\ZZ^2$  & $\ZZ^2$   & $\ZZ^2$ &    $\ZZ^2$    & $\ZZ^2$  & $\ZZ^2$   & $\ZZ^2$   \\
    \end{tabular}
   \]
  \end{Exam} 
 
 \begin{proof}
Since $KO_*( \real(S^1)) = KO_*(\RR) \oplus \Sigma KO_*( \RR)$, we can use the version of the K\"unneth Theorem expressed in Theorem~\ref{RealKF}, Part (1), in the next section.
\index{K\"unneth Theorem}
 We have
 \begin{align*}
  KO_*( \real(S^1 \times S^1)) &= K O_*( \real(S^1) \otimes_{\RR} \real(S^1) ) \\
 	&\cong KO_*( \real(S^1)) \otimes_{KO_*(\RR)} KO_*( \real(S^1))  \\
	&= (1 + \Sigma) KO_*(\RR) \otimes_{KO_*(\RR)} (1 + \Sigma) KO_*(\RR)   \\
	&= (1 + \Sigma)^2 KO_*(\RR) \; .
\end{align*}
 \end{proof}

\begin{Exam}
More generally, consider the $n$-torus $S^1 \times \overset{n}\cdots \times S^1$. Then
\[ 
KO_* (\real (S^1 \times \overset{n}\cdots \times S^1))  \cong
 \left(  \prod_{i = 1}^n (1 + \Sigma) \right) KO_*(\RR)\; . \]

Even more generally, consider a generalized torus where each factor is a sphere of arbitrary dimension: $X =\prod_{i = 1}^n S^{k_i}$. Then
\[ KO_*( \real(X)) \cong \left( \prod_{i = 1}^n (1 + \Sigma^{k_i}) \right) KO_*(\RR)\; . \] 

\end{Exam}
    
    
\newpage

\section{\bf{The K\"unneth Theorem and the Spectral Sequence}} \label{Section-Kunneth}
\index{K\"unneth Theorem|(}

In this section, we present two of the deeper tools used for computing $K$-theory for $C^*$-algebras: the K\"unneth Theorem for tensor products and the Atiyah-Hirzebruch spectral sequences for an algebra with a filtration. We will present these results without proof, but for reference we will include descriptions for both complex $C ^*$-algebras and real $C^*$-algebras. We jump ahead a bit and also state the corresponding results in the language of real and complex $K$-theory for topological spaces. See Sections~\ref{Section:KOforSpaces} and \ref{Section:RelatingKO-theory}.
 
K\"unneth theorems in general deal with the problem determining $h_*(A \otimes B) $ in terms of $h_*(A)$ and $h_*(B)$.   One has to put some constraints on the algebras 
 in order for the theorem to hold.  The K\"unneth Theorem for $K$-theory involves some serious homological algebra.  For $K_*(A)$ it was established by Schochet \cite{TopII} (1984) where the 
 bootstrap category was first introduced, and for $KO_*(A)$ it is due to Boersema \cite{BoersemaKT} (2002), building on the work of \cite{TopII} and utilizing united $K$-theory.   
    
  \begin{Thm}  {\bf{The K\"unneth Theorem (complex case)}} \label{ComplexKF}
 
 Suppose that $A$ and $B$ are complex $C^*$-algebras and $A$ is in the bootstrap category 
 (this includes all commutative $C \sp *$-algebras $A$).
   \begin{enumerate}
  \item If $K_*(B)$ is torsion-free, then 
  the natural external product map is 
  an isomorphism of graded groups 
%
  \[
  K_*(A) \otimes K_*(B) \overset{\cong}\longrightarrow K_*(A\otimes B) .
  \]
  \item
  In particular, if $X$ and $Y$ are compact spaces and  $K^*(Y)$ is torsion-free,
  then the natural external product map is
  an isomorphism of graded groups
  \[
  K^*(X) \otimes K^*(Y) \overset{\cong}\longrightarrow K^*(X \times Y) .
  \]
  \item In general, there is a natural long exact sequence of graded groups
   \[
  0 \longrightarrow K_*(A) \otimes K_*(B) \overset{\alpha}\longrightarrow K_*(A\otimes B) \overset{\beta}\longrightarrow \Tor(K_*(A), K_*(B)) \longrightarrow 0
  \]
 where $\alpha $ has degree $0$, $\beta$ has degree $-1$, and  $\Tor$ is the classical homological algebra functor on pairs of abelian groups.\footnote{$\Tor$ is the derived functor associated with tensor products. For details on $\Tor$  see p. 154 of \cite{MacLaneBook} or Chapter 3 of \cite{Wbook}. If $G$ or $H$ is torsion-free, then $\Tor(G,H) = 0$ which indicates how Parts (1) and (2) follow from Part (3).}
 

\end{enumerate}
\end{Thm} 
 
   \begin{Rem}
  Typically, complex $K$-theory is considered as just two groups $K_*(A) = (K_0(A), K_1(A))$. Then this is considered as graded groups over $\ZZ_2$ and the tensor products are tensor products of $\ZZ_2$-graded groups. Explicitly, if $M$ and $N$ are $\ZZ_2$-graded groups, then $M \otimes N$ is also $\ZZ_2$-graded and given by
  \begin{align*}
 (M \otimes N)_0 &= M_0 \otimes N_0 \oplus M_1 \otimes N_1 \\
 \text{and~}  (M \otimes N)_1 &= M_0 \otimes N_1 \oplus M_1 \otimes N_0 \; .  \\
  \end{align*}
  
  Similarly, $\Tor(M,N)$ is $\Z_2$-graded and 
   \begin{align*}
 \Tor(M,N)_0 &=  \Tor(M_0, N_0) \oplus \Tor(M_1, N_1) \\
 \text{and~}   \Tor(M,N)_1 &= \Tor(M_0, N_1) \oplus \Tor(M_1, N_0)  \; .  \\
  \end{align*}
  
  Alternatively, let us consider complex $K$-theory to be graded over $\ZZ$ with an invertible Bott map $K_*(A) \rightarrow K_*(A)$ of degree 2. This is equivalent to saying that
  $K_*(A)$ is a module over the ring $K_*(\CC)$ (which is also graded over $\ZZ$). If $M$ and $N$ are $\ZZ$-graded modules over $K_*(\CC)$ then so are $M\otimes_{K_*(\CC)} N$ and $\Tor_{K_*(\CC)} (M,N)$. See the Main Theorem of \cite{BoersemaKT}.
   
  This latter perspective matches what we do in the real case below, where we consider $KO_*(A)$ to be a $\ZZ$-graded module over the ring $KO_*(\RR)$ as in Proposition~\ref{prop-module}, and we take tensor products over $KO_*(\RR)$ in Part (1) of the following.
  \end{Rem}

  \begin{Thm}   {\bf{The K\"unneth Theorem (real case)}} \label{RealKF}
  
 Suppose that $A$ and $B$ are real  $C^*$-algebras and $A$ is in the pre-bootstrap category 
 (this includes all commutative real $C \sp *$-algebras $A$).
 \begin{enumerate}
 \item
 Suppose that $B$ has the property that
$K_*(B\sc)$ is torsion-free. Then the natural external product map is 
  an isomorphism of $KO_*(\RR)$-modules\  
  \[
  KO_*(A) \otimes_{KO_*(\RR)} KO_*(B) \overset{\cong}\longrightarrow KO_*(A\otimes B) .
  \]
  
  \item
  In particular, if $X$ and $Y$ are compact spaces and $K^*(Y)$ is torsion-free,
  then the natural external product map is
  an isomorphism of $KO_*(\RR)$-modules
    \[
  KO^*(X) \otimes_{KO^*(\pt)} KO^*(Y) \overset{\cong}\longrightarrow KO^*(X \times Y) .
  \] 
  \item
In general, there 
  is a natural short exact sequence of $CRT$-modules
\[
  0 \longrightarrow K\crt(A) \otimes\scrt K\crt(B) \overset{\alpha}\longrightarrow K\crt(A\otimes B) \overset{\beta} \longrightarrow  \Tor\scrt(K\crt(A), K\crt(B)) \longrightarrow 0
\]
   where $\alpha $ has degree $0$, $\beta$ has degree $-1$; and where $\otimes\scrt$ and $\Tor\scrt$ are versions of tensor product and $\Tor$ operators defined in the category  
   $\mathcal{CRT}$ of all $CRT$-modules\footnote{See \cite{BoersemaUCT} for the definitions.}.
\end{enumerate}
\end{Thm}

\begin{Rem}
\begin{enumerate}
\item Examples of real $C \sp *$-algebras $B$ that enjoy the property that $K_*(B\sc)$ is free include
$B = \RR, \CC, \HH, \iS, \realo (i\topr), \real (S^n)$ as well as any countable direct sum or  tensor product of these, or matrix algebra over any of these. Examples also include any real $C \sp *$-algebra of the form
$\real(X, \tau)$ where $X$ is a sphere or torus of any dimension and $\tau$ is any involution.

The condition that $K_*(B\sc)$ is free (as a group) is equivalent to the condition that $K\crt(B)$ is a free {$CRT$}-module, by \cite{Bousfield}.

\item 
We emphasize that the tensor product $KO_*(A) \otimes_{KO_*(\RR)} KO_*(B)$ in Part (1) (and also Part (2))
is a tensor product over the ring $KO_*(\RR)$. See Section~{IV. 5} in \cite{Hungerford} for the definition of tensor products over a ring. 
If $M$ and $N$ are both modules over the $KO_*(\RR)$, then
$$(M \otimes_{KO_*(\RR)} N)_n \cong  \left( \bigoplus_{i = 0}^7 M_i \otimes N_{n-i} \right)/K$$
where the group $K$  is such that in the quotient we have
$$( \eta \cdot m) \otimes  n = m \otimes (\eta \cdot n) \quad \text{and} \quad (\xi \cdot  m) \otimes n = m \otimes (\xi \cdot n)$$
for all $m \in M$ and $n \in N$.

In particular, it is not the case that
$$KO_n(A  \otimes B) \cong \bigoplus_{i = 0}^7 (KO_i(A) \otimes KO_{n-i}(B)) \; $$
(even under the hypothesis that $K_*(B\sc)$ is torsion-free)
but rather $KO_n(A \otimes B)$ is a quotient of the direct sum on the right side.
We can see that this is not an isomorphism just by considering $A = B = \RR$.

Indeed, $KO_*(\RR \otimes \RR) = KO_*(\RR) = (\Z, \Z_2, \Z_2, 0, \Z, 0, 0, 0)$.
On the other hand, by taking
$$M_n = \bigoplus_{i = 0}^7 (KO_i(\RR) \otimes KO_{n-i}(\RR))$$
we obtain 
$$M_n = (\Z^2, \Z_2^2, \Z_2^3, \Z_2^2, \Z^2 \oplus \Z_2, \Z_2^2, \Z_2^2, 0) \; .$$

\item The general K\"unneth Theorem for real $C^*$-algebras cannot be expressed without the language of united $K$-theory and 
$CRT$-modules. In particular, the isolated product map on real $K$-theory
\[
  KO_*(A) \otimes_{KO_*(\RR)} KO_*(B) \longrightarrow KO_*(A\otimes B) .
  \]
 always exists but is not necessarily injective or surjective apart from the situation of Part (1) of Theorem~\ref{RealKF}, as pointed out already by Atiyah in a footnote in \cite{AtiyahKT}.
 A counter-example here can be obtained by taking $A = B = \CC$.
\item 
As we saw in Section~\ref{Section-United}, $K\crt(A \otimes B)$ is a $CRT$-module with three components: a real part, a complex part, and a self-conjugate part. The real part is equal to $KO_*(A \otimes B)$. Therefore the real part of the sequence in Part (3) above can be used to calculate $KO_*(A \otimes B)$. However to do this we must compute the 
$CRT$-tensor product $K\crt(A) \otimes\scrt K\crt(B)$ which is a subtle object. In particular the real part of 
$K\crt(A) \otimes\scrt K\crt(B)$ is not just equal to $KO_*(A) \otimes_{KO_*(\RR)} KO_*(B)$, but must be calculated using the entire $CRT$-structure of both factors. This is also true for $\Tor\scrt(K\crt(A), K\crt(B))$. The details of the tensor product construction (and the $\Tor$ construction) for $CRT$-modules can be found in \cite{BoersemaKT},

To better understand these constructions, see the detailed calculation of 
$$KO_*(\mathcal{O}\pr_n \otimes \mathcal{O}\pr_m)$$ 
that is carried out in \cite{BoersemaKT} using the K\"unneth Theorem, where $\mathcal{O}\pr_n$ is the real Cuntz algebra described in Section~\ref{KKSection}.

\item The complex part of the $CRT$-module $K\crt(A \otimes B)$ is isomorphic to $K_*(A\sc \otimes B\sc)$ and the complex part of the tensor product
$K\crt(A) \otimes\scrt K\crt(B)$ is simply isomorphic to $K(A\sc) \otimes K(B\sc)$.
In fact, the complex part of the K\"unneth Theorem sequence in Part (3) of Theorem~\ref{RealKF} reduces to the 
complex K\"unneth Theorem sequence in Part (3) of Theorem~\ref{ComplexKF}. 
\end{enumerate}
\end{Rem}
\index{K\"unneth Theorem|)}

 \vspace{.1in}
{\bf A Spectral Sequence for real $K$-theory and united $K$-theory.}
\index{Atiyah-Hirzebruch spectral sequence|(}
\vspace{.1in}

 Finally, we present here a deep result in spectral sequences. 

\begin{Def} Suppose that $A$ is a real or a complex $C^*$-algebra. We say that $A$ is {\emph{filtered}} 
\index{filtered $C^*$-algebra} if there is an increasing sequence of closed ideals 
\[
A_0 \subset A_1 \subset A_2 \subset \dots
\]
with $\overline{\bigcup _j A_j } = A$.  
\end{Def}
Filtered $C^*$-algebras occur naturally in several situations, typically in the structure theory of certain types of $C^*$-algebras.  The following result was established 
in \cite{TopI} Theorem 1.2 for complex $C^*$-algebras and the proof there is equally valid in the real context of real $C^*$-algebras, and can also be extended to united $K$-theory as shown in Part (3).
\begin{Thm} (\cite{TopI}, Theorem 1.2 )  \label{Thm-spectralsequence}
\begin{enumerate}
\item
Let $A = \overline{\bigcup _j A_j }$ be a complex filtered $C^*$-algebra. Then
there is a spectral sequence $\{E^r, d^r\} $ with $d^r : E_{p,q}^r \to E_{p - r,q+r - 1}^r $ 
which converges to the graded object associated to the natural filtration on 
\[
\lim_{n \to \infty} K_*(A_n) \cong K_*(A),
\]
 with 
\[
E_{p,q}^1 (A)   = K_{p+q}(A_p/A_{p-1} ),
\]
where $K_{p+q}$ is graded mod 2.
\item
Let $A$ be a real filtered $C^*$-algebra. Then
there is a spectral sequence $\{E^r, d^r\} $ with $d^r : E_{p,q}^r \to E_{p - r,q+r - 1}^r $ 
which converges to the graded object associated to the natural filtration on 
\[
\lim_{n \to \infty} KO_*(A_n) \cong KO_*(A),
\]
 with 
\[
E_{p,q}^1 (A)   = KO_{p+q}(A_p/A_{p-1} ),
\]
where $KO_{p+q}$ is graded mod 8.
\item
Let $A$ be a real filtered $C^*$-algebra. Then
there is a spectral sequence of $CRT$-modules 
$\{E^r, d^r\} $ with $CRT$-module homomorphisms 
$$d^r : E_{p,q}^r \to E_{p - r,q+r - 1}^r $$
which converges to the $CRT$-module associated to the natural filtration on 
\[
\lim_{n \to \infty}   K\crt(A_n) \cong K\crt(A),
\]
 with 
\[
E_{p,q}^1 (A)   = K\crt_{p+q}(A_p/A_{p-1} ) \; .
\]

\end{enumerate}
\end{Thm}

A special case of the above result, established much earlier by Atiyah and Hirzebruch \cite {AH}, 
the very first paper on topological $K$-theory, is the Atiyah-Hirzebruch spectral sequence. \index{Atiyah-Hirzebruch spectral sequence}  We state it for finite 
CW-complexes $X$, though it holds for compact spaces and, with a little modification, for locally compact spaces\index{simplicial complex}\footnote{For CW-complexes, see Section~\ref{top spaces} and the citations therein. Every compact manifold is equivalent to some finite CW-complex.}. 

\begin{Cor}  Suppose that $X$ is a finite CW-complex.
Let $\{x_o\}$ denote the space with a single point.
Then there is a spectral sequence $E_r^{p,q}$  ($r \geq 1$, $-\infty < p,q < +\infty $), with 
\[
E_1^{p,q} \cong C^p(X ; KO^q(\{x_o\}) 
\]
($p$- cochains taking values in $KO^q(\{x_o\})$, with $d^1 $ the usual coboundary operator, and 
\[
E_2^{p,q} \cong H^p(X ; KO^q(\{x_o\}) 
\]
\[
E_\infty ^{p,q} \cong G_pKO^{p+q} (X).
\]
where $G_pKO^{p+q} (X)$ is the associated graded group built from the skeletal filtration of $KO^*(X)$ built from the filtration of $X$. 
The differential  
\[
d_r :  E_r^{p,q} \longrightarrow E_r^{p+r,q-r+1} 
\]
 has degree $(r, -r+1)$ as shown. 
\end{Cor}

\index{Atiyah-Hirzebruch spectral sequence|)}

\newpage   

\section{\bf{$KK$ theory and the Universal Coefficient Theorem}} \label{KKSection}
\index{KK-theory@$KK$-theory|(}
This is a very brief introduction to the powerful world of Gennadi  Kasparov's $KK$-theory.  There are several sources for this, starting with Kasparov's fundamental papers 
on the subject. For a newcomer, \cite{KK} is probably the best reference. One can also profitably consult Blackadar \cite{Blackadar}, Schr\"oder \cite{Schr}, and Boersema-Loring-Ruiz \cite{BLR}, for instance, and see Rosenberg \cite{R} for some important applications.  

We will not present here a detailed treatment of $KK(A,B)$ for real or complex $C \sp *$-algebras $A$ and $B$. 
In fact, there are multiple approaches to the definition of $KK(A,B)$, which lead to the same theory, but all are quite technical.
Instead we will focus on the significant properties, which will be stated without proof. Then the reader can consult the references listed above for a full development.
Three key points to remember for this section:
\begin{enumerate}
\item Kasparov treats both the real and the complex situations.  We will do that here and call it to the reader's attention in every important statement by referring to $C^*$-algebras over $\FF$ where $\FF$ is the field $\RR$ or $\CC$. So we will often simultaneously be working in the categorical setting of real $C^*$-algebras with $KO$-theory and the categorical setting of complex $C^*$-algebras and complex $K$-theory.  Similarly, $KK$-theory must be considered either in one category or the other, but the development is the same.

When it is necessary to specify one setting or the other, we will indicate that.

\item Kasparov also deals with the situation that there is a fixed topological group $G$ acting upon the $C^*$-algebras and, of course, this leads to equivariant $K$ and $KK$-theory. We will not discuss the equivariant theory.  He also insists all algebras are $\ZZ _2$-graded. We won't discuss this either.
\item In order for the general theory to work, we stipulate in this section that all $C \sp *$-algebras are separable\footnote{This is somewhat too stringent but will do for a first view.}.  
\end{enumerate}

\begin{Thm} (Kasparov)  \label{KK(A,B)}
For separable $C \sp *$-algebras $A$ and $B$ over $\FF$, there exists an abelian group $KK\pf(A,B)$, such that
\begin{enumerate}
\item $KK\pf(A, B)$ is contravariant in the first variable. In particular, a homomorphism $\alpha \colon A_1 \rightarrow A_2$ induces a map
$$\alpha^* \colon KK\pf(A_1, B) \leftarrow KK\pf(A_2, B)$$
\item $KK\pf(A, B)$ is covariant in the second variable. In particular,  a homomorphism $\beta \colon B_1 \rightarrow B_2$ induces a map
$$\beta_* \colon KK\pf(A, B_1) \leftarrow KK\pf(A, B_2)$$
\item $KK\pf(A,B)$ is homotopy invariant in each variable. That is, the maps $\alpha^*$ and $\beta_*$ in Parts (1) and (2) only depend on $\alpha$ and $\beta$ up to homotopy.
\item $KK\pf(A,B)$ is stable in each variable. That is, the natural maps $A \to A\otimes \KK\subf $ and $B \to B\otimes \KK\subf$ induce isomorphisms.
 \item
There are natural isomorphisms 
\[
KK\pr(\RR , B) \cong KO_0(B) \quad \text{for any real $C \sp *$-algebra $B$}
\]
and 
\[
KK\pc(\CC , B) \cong K_0(B) \quad \text{for any complex $C \sp *$-algebra $B$}
\]
showing that both real and complex $K$-theory can be expressed respectively in terms of real and complex Kasparov groups.
\item For every $C \sp *$-algebra homomorphism $\alpha \colon A \rightarrow B$, there exists an element $[\alpha] \in KK\pf(A,B)$. Furthermore, the map $\alpha \mapsto [\alpha]$ defines an additive homomorphism $\Hom(A,B) \rightarrow KK\pf(A,B)$.
\end{enumerate}
\end{Thm}

\begin{Rem}
Because of Item (6) above, it is often helpful to think of $KK\pf(A,B)$ as some sort of group of generalized homomorphisms from $A$ to $B$. Further justifying this perspective is the existence of a homomorphism
$$KK\pc(A,B) \rightarrow \Hom(K_*(A), K_*(B))$$
for complex $C \sp *$-algebras $A$ and $B$; and a homomorphism 
$$KK\pr(A,B) \rightarrow \Hom_{KO_*(\RR)}(KO_*(A), KO_*(B))$$
for for real $C \sp *$-algebras $A$ and $B$, which we describe below in the discussion leading up to Theorem~\ref{Thm-UCT}.
Furthermore, the compositions
$$\Hom(A,B) \rightarrow KK\pc(A,B) \rightarrow \Hom(K_*(A), K_*(B))$$
and 
$$\Hom(A,B) \rightarrow KK\pr(A,B) \rightarrow \Hom_{KO_*(\RR)}(KO_*(A), KO_*(B))$$
are exactly the maps $\alpha \mapsto \alpha_*$ induced on $K$-theory by any $C \sp *$-algebra homomorphism $\alpha$.

\end{Rem}

If we take the group $KK\pf(A,B)$ as given in Theorem~\ref{KK(A,B)}, then we can define the graded group $KK\pf_n(A,B)$ as follows for all $n \in \ZZ$.

\begin{Def}
For separable $C \sp *$-algebras $A$ and $B$ over $\FF$ and for all $n \in \ZZ$, define $KK\pf_n(A,B) = KK\pf(A, S^n B)$,
again with reference to Definition~\ref{Def-AllSus} for positive and negative suspensions.
\end{Def}

%

\begin{Thm}  \label{Thm-KandKK}
   For separable $C \sp *$-algebras $A$ and $B$ over $\FF$,
\begin{enumerate} 
 \item  For a fixed $C^*$-algebra $A$, the functors $KK\pf_n(A, - )$ form an additive homology theory. 
\item For a fixed $C^*$-algebra $B$, the functors $KK\pf_n(-, B)$ form an additive cohomology theory. 
\item 
 The groups $KK\pf_*(A,B)$ are periodic in each variable, so there are natural isomorphisms
 \[
KK\pr_n(A,B) \cong KK\pr_n(S^8 A, B) \cong KK\pr_n(A, S^8 B) \quad \text{for any real $C \sp *$-algebras $A$ and  $B$} 
\]
and 
\[
KK\pc_n(A,B) \cong KK\pc_n(S^2 A, B) \cong KK\pc_n(A, S^2 B) \quad \text{for any complex $C \sp *$-algebras $A$ and $B$}  \]
\item
There are natural isomorphisms 
\[
KK\pr_n(\RR , B) \cong KO_n(B) \quad \text{for any real $C \sp *$-algebra $B$}
\]
and 
\[
KK\pc_n(\CC , B) \cong K_n(B) \quad \text{for any complex $C \sp *$-algebra $B$}
\]

\end{enumerate}
\end{Thm}

Thus the Kasparov groups are a major source for homology and for cohomology theories.  They have an additional critical feature which enables much 
more power, namely the Kasparov pairing. Again, we see from this that $KK\pf(A,B)$ behaves formally a lot like $\Hom(A,B)$.

\begin{Thm} \label{Thm-pairing} 
For separable $C \sp *$-algebras over $\FF$,
there is a natural product map 
\[
KK\pf(A_1, B_1 \otimes D ) \times KK\pf(D \otimes A_2, B_2) \longrightarrow KK\pf(A_1 \otimes A_2, B_1 \otimes B_2 )
\] 
denoted 
\[
(\alpha , \beta ) \longrightarrow \alpha \otimes _D \beta
\]
which is natural in each variable, associative, and extends to a graded pairing. 
With this pairing, we obtain the following structures:
\begin{itemize}
\item $KK\pf_*(D, D)$ the structure of a unital graded ring -- the unit is the element $[\id_D] \in KK\pf_0(D,D)$, using the homomorphism of Theorem~\ref{KK(A,B)}, Part (6).
\item $KK\pf_*(D\otimes A, B)$ the structure of a graded left $KK\pf_*(D, D)$-module, and  
\item $KK\pf_*(A, B\otimes D)$ the structure of a graded right $KK\pf_*(D, D)$ module. 
\end{itemize}
\end{Thm}

%
%
%
%
%
%
 
 \begin{Rem} It is difficult to sort out and give precise credit for these theorems.  No doubt the main credit goes to Gennadi Kasparov, starting with \cite{KK1980} and \cite{KK}.  There were 
 many strengthenings of the theory and Skandalis' paper \cite{Skandalis} stands out. These papers and others like them fine-tuning the $KK$-groups are very technical and not really meant for our audience.
 \end{Rem}

\begin{Rem} Note to algebraic topologists.  The formal structure of the Kasparov pairing should not come as a surprise. If we let $[X, Y]$ denote stable homotopy classes of maps from finite complexes $X$ to $Y$ and 
let $E$ denote a ring spectrum then we may form  the groups $[X, Y \wedge E ]$. Fixing $X$ we may generate a homology theory in $Y$.  Fixing $Y$ we can generate a cohomology 
theory in $X$.   Now fix some space or spectrum $D$.   Then there is a natural pairing
\[
[X_1, Y_1 \wedge E \wedge D ] \times [D \wedge X_2, Y_2 \wedge E ] \longrightarrow [X_1 \wedge X_2, Y_1 \wedge Y_2 \wedge E ]
\]
constructed as follows. Given $\alpha \colon X_1 \rightarrow Y_1 \wedge E \wedge D$ and $\beta \colon D \wedge X_2 \rightarrow Y_2 \wedge E$, then the 
element $\alpha \otimes_D \beta \colon X_1 \wedge X_2 \rightarrow Y_1 \wedge Y_2 \wedge E$ is defined by the composition
\begin{multline}
	X_1 \wedge X_2 \xrightarrow{\alpha \wedge \id_{X_2}} Y_1 \wedge E \wedge D \wedge X_2
	\xrightarrow{\id_{Y_1 \wedge E} \wedge \beta} Y_1 \wedge E \wedge Y_2 \wedge E \\
	\xrightarrow{~~\cong~~}  Y_1 \wedge Y_2 \wedge E \wedge E
	\xrightarrow{\id_{Y_1 \wedge Y_2} \wedge m} Y_1 \wedge Y_2 \wedge E.   \; 
\end{multline}
where $m \colon E \wedge E \rightarrow E$ is the product structure on $E$.
The Kasparov pairing, when restricted to $C^*$-algebras of the form $C(X)$ (real or complex, of course), 
corresponds to the pairing obtained from this construction with $E $ the $K$-theory spectrum (real or complex respectively).
\end{Rem}

There have been many applications of the Kasparov theory. We shall mention a few of them: the work of Brown-Douglas-Fillmore (BDF) in \cite{BDF1} and \cite{BDF2}; and the Universal Coefficient Theorem (UCT) in \cite{RS} and \cite{BoersemaUCT}, which in turn led to the classification of real purely infinite simple C*-algebras in \cite{phillips-class},  \cite{Kirchberg94} and \cite{BRS}.
\index{KK-theory@$KK$-theory|)}

\vspace{.1in}

{\bf{ BDF: Classification of Essentially Normal Operators}}

\vspace{.1in}

Recall that a bounded operator T on Hilbert space is {\it Fredholm} if it has closed range and if $\ker T$ and $\ker T^*$ are both finite-dimensional.  The Fredholm index of such an operator is defined by
              $$  \ind(T) = \dim \ker T    - \dim \ker T^*$$
For example, if $T = N + K$ , the sum of a normal and compact operator then 
                            $$  \ind(T) = \ind(N + K) = \ind(N) = 0$$
since $\dim \ker(N) = \dim \ker (N^*)$.   On the other hand, if $S$ is the unilateral shift then $S$ is essentially normal but  
       $$\ind(S) = \dim \ker(S) - \dim \ker (S^*) = 0 - 1 = -1 \; ,$$
showing that $S$ cannot be written in the form $S = N + K$ where $N$ is normal and $K$ is compact.  
P. Halmos asked 
  whether the index condition was necessary and sufficient to determine whether or not an essentially normal operator can be so decomposed. BDF found the answer:
 
\begin{Thm}(BDF \cite{BDF1}).   Suppose that $T$ is an essentially normal operator on complex Hilbert space.   
Let $X$ denote the essential spectrum\footnote{That is, the spectrum of the image of $T$ in the Calkin algebra.  It is a closed subset of the plane. \index{essential spectrum} }   of $T$.  
Then the following are equivalent:
\begin{enumerate}

\item $T = N + K$ for some normal operator $N$ and some compact operator $K$.

\item For each $\lambda \in \CC - X$, 
\[
\ind(T - \lambda I) = 0 \; .
\]

\end{enumerate}
\end{Thm}

BDF's method was to consider all extensions of $C^*$-algebras of the form
\[
0 \longrightarrow \KK \longrightarrow E \longrightarrow \complex(X) \longrightarrow 0
\]
and take unitary equivalence classes of them to form a group which they called $Ext(X)$.  
The trivial element of that group is formed by taking $E$ to be all of the elements of the form $N + K$ as above.
Sitting in that group is the element $[T]$ formed by taking $E$ to be the $C^*$-algebra 
generated by the compact operators and by $T$.   Then $T$ has the form $N + K$ if and only if $[T] = 0$.    

BDF then showed that from the Fredholm index one could construct a homomorphism
\[
\Ext(X) \longrightarrow \Hom (K^1(X), \ZZ) .
\]
They proved that this map is an isomorphism for $X$ a closed subset of the plane,
and that in turn yields a proof of the theorem.

Later, they realized that, for finite complexes $X$ that 
\[
\Ext(X) \cong K_1(X).
\]
  They were able to extend this to prove the existence of a homology 
theory $\Ext_*(X) $ of period 2 and could be identified with $K_*(X)$.   Finally, it was established that
\[
\Ext_*(X) \cong KK\pc_*(C(X), \CC ). 
\]
Moving from $X$ a finite complex to it being a compact metric space was a serious problem, but this was resolved by Kaminker and Schochet \cite{KS}. 

\vspace{.1in}
{\bf{The Universal Coefficient Theorem for Complex and Real C*-Algebras}}
\vspace{.1in}

\index{universal coefficient theorem|(}

In this subsection, we present The Universal Coefficient Theorem (UCT) for both complex and real $C^*$-algebras and show their application to $KK$-equivalence.
The original sources for this material are \cite{RS} in the complex case and \cite{BoersemaUCT} in the real case.
The UCT, in both the complex and real settings, expresses the Kasparov $KK$-theory $KK\pf(A,B)$ in terms of the $K$-theory of the two algebras $A$ and $B$. To set the stage, note that the pairing of Theorem~\ref{Thm-pairing} gives a graded homomorphism
$$KK\pr_*(\RR, A) \otimes KK\pr_*(A, B) \xrightarrow{\otimes_A}  KK\pr_*(\RR, B)$$
which is equivalent to having a graded homomorphism
$$\gamma \colon KK\pr_*(A, B) \rightarrow \Hom(KK\pr_*(\RR, A), KK\pr_*(\RR,B))$$
(as the functors $\Hom(N, -)$ and $N \otimes -$ are adjoint to each other\footnote{See Definition~2.3.9 in \cite{W} for the definition of adjoint functors. The adjoint relationship also connects the homotopy constructions $\Omega$ and $S$ in Theorem~\ref{Thm-homotopy}\index{adjoint functors}.}).
Then from Theorem~\ref{Thm-KandKK}, this can be written as
$$\gamma\so \colon KK\pr_*(A, B) \rightarrow \Hom(KO_*(A), KO_*(B))  \; .$$


Now Theorem~\ref{Thm-pairing} also gives us a graded homomorphism
$$KK\pr_*(\RR, \CC \otimes A) \otimes KK\pr_*(A, B) \xrightarrow{\otimes_A} KK\pr_*(\RR, \CC \otimes B)$$
which is equivalent to
$$\gamma\su \colon KK\pr_*(A, B) \rightarrow \Hom(K_*(A\sc), K_*(B\sc))  \; .$$
Similarly there is a homomorphsm 
$$KK\pr_*(\RR, T \otimes_A) \otimes KK\pr_*(A, B) \xrightarrow{\otimes A} KK\pr_*(\RR, T \otimes B)$$
which can be written as 
$$\gamma\sst \colon KK\pr_*(A, B) \rightarrow \Hom(KT_*(A ), KT_*(B ))  \; .$$

For any $\xi \in KK\pr_*(A,B)$, the trio of homomorphisms $\Gamma(\xi) = (\gamma\so(\xi), \gamma\su(\xi), \gamma\sst(\xi))$ are compatible with the natural transformations $\Theta\crt$ and thus there exists a homomorphism of $CRT$-modules
$$\Gamma(\xi) \colon K\crt(A)  \rightarrow K\crt(B)  \; .$$

Therefore $\Gamma$ is a homomorphism
$$\Gamma \colon KK\pr_*(A, B) \rightarrow \Hom\scrt(K\crt(A), K\crt(B))$$
where $\Hom\scrt(K\crt(A), K\crt(B))$ representes the graded group of all graded $CRT$-module homomorphisms.

The UCT then states that $\Gamma$ is a surjection, under certain broad hypotheses, and furthermore it identifies the kernel of $\Gamma$.

\begin{Thm}[Universal Coefficient Theorem] \label{Thm-UCT}
\begin{enumerate}
\item
Let $A$ and $B$ be separable complex $C^*$-algebras
with $A$ in the bootstrap class.
Then there exists a natural short exact sequence,
$$0 \rightarrow \Ext^1(K_*(A), K_*(B)) \xrightarrow{\delta} KK\pc_*(A,B) \xrightarrow{\gamma} \Hom(K_*(A), K_*(B)) \rightarrow 0 \; $$
where $\delta$ has degree $-1$ and $\gamma$ has degree $0$. This sequence splits, though not naturally.
\item Let $A$ and $B$ be separable real $C^*$-algebras with $A$ in the pre-bootstrap class. Then there exists a natural short exact sequence,
\[
	0 \rightarrow \Ext\scrt(K\crt(A), K\crt(B)) \xrightarrow{\delta} KK\pr_*(A,B)  \xrightarrow{\Gamma} \Hom\scrt(K\crt(A), K\crt(B)) \rightarrow 0 \; \]
	where $\delta$ has degree $-1$ and $\Gamma$ has degree $0$.
	This sequence does not split in general.
\end{enumerate}
\end{Thm}

\begin{Rem}
In the real case, the object $\Hom\scrt(K\crt(A), K\crt(B))$ is understood as the group of all graded $CRT$-module homomorphisms 
$\alpha = (\alpha\po, \alpha\pu, \alpha\ppt)$ from $K\crt(A)$ to $K\crt(B)$, as in Definition~\ref{Def-CRTmodule}.
Then 
$$\Ext\scrt(K\crt(A), K\crt(B))$$ involves the derived functor as in Definition~2.5.2 of \cite{Weibel:HomAlg}.
\end{Rem}

One of the crucial consequences of the Universal Coefficient Theorem for both the real and complex cases, is that it provides a characterization of $KK$-equivalence.

\begin{Def} \label{Def-KKequiv}
Two $C^*$-algebras $A$ and $B$ over $\FF$ are said to be {\it $KK$-equivalent} \index{KK-equivalence@$KK$-equivalence} if there exists elements $\alpha \in KK\pf(A,B)$ and $\beta \in KK\pf(B,A)$ such that
$$\alpha\otimes \beta = 1_A  \in KK\pf(A,A) $$
and 
$$\beta \otimes \alpha = 1_B \in KK\pf(B,B) \; .$$
\end{Def}

We note that $KK$-equivalence over $\FF$ is an equivalence relation on the class of $C ^*$-algebras over $\FF$, and this relation is strictly weaker than isomorphism. The next result follows from Theorem~\ref{Thm-UCT}, see Corollary 7.5 of \cite{RS} and Corollary~4.11 of \cite{BoersemaUCT}.

\begin{Cor} \label{Cor-KKequiv}
\begin{enumerate}
\item Let $A,B$ be complex $C \sp *$-algebras in the bootstrap class. Then $A$ and $B$ are $KK$-equivalent, as complex $C^*$-algebras, if and only if $$K_*(A) \cong K_*(B) \; .$$
\item Let $A,B$ be real $C \sp *$-algebras in the pre-bootstrap class. Then the following are equivalent:
\begin{enumerate}
\item $A$ and $B$ are $KK$-equivalent, as real $C ^*$-algebras.
\item $K\crr(A) \cong K\crr(B) \; $
as $CR$-modules.
\item
$K\crt(A) \cong K\crt(B) \; $
as $CRT$-modules.
\end{enumerate}
\end{enumerate}
\end{Cor}
\index{universal coefficient theorem|)}

For example, note that the complex $C \sp *$-algebras $A = \complex(S^1)$ and $B = \complex(S^3)$ are $KK$-equivalent, as they satisfy $K_*(A) \cong K_*(B)$, but they are not isomorphic $C \sp *$-algebras.
Similarly, the real $C \sp *$-algebras $A = \real(S^1)$ and $B = \real(S^9)$ are $KK$-equivalent, but not isomorphic.

\vspace{.1in}

{\bf{The Classification of Purely Infinite Simple C*-Algebras}}

\vspace{.1in}

In this section, we briefly consider the class of simple, purely infinite $C \sp *$-algebras, to relate the powerful classification theorem for this class of $C \sp *$-algebras, in both the real and complex setting. Recall from Definition~\ref{Def-simpleC*-alg} that a $C \sp *$-algebra is simple if it has no non-trivial norm-closed, *-closed, two-sided idealsl. As we see in Proposition~\ref{Pro-c-simple}, if a real $C \sp *$-algebra $A$ is not simple, then $A\sc$ is also not simple; but the converse does not hold. This leads us to made the following definition, allowing us to focus on real simple $C \sp *$-algebras that remain simple upon complexification.


\begin{Def} \index{simple $C \sp *$-algebra!c-simple $C \sp *$-algebra} 
A real $C \sp *$-algebra $A$ is said to be {\it c-simple} if $A\sc$ is simple.
\end{Def}

Recall, that if $A$ is a complex $C \sp *$-algebra, then we can consider $A$ to be a real $C \sp *$-algebra, by restricting the scalar action to $\RR$. This is a sort of forgetful functor 
$$\mathcal{F} \colon \mathscr{C}^*_{\rm{complex}} \rightsquigarrow \mathscr{C}^*_{\rm{real}}$$
from the category of complex $C \sp *$-algebras to the category of real $C \sp *$-algebras.

\begin{Pro} \label{Pro-c-simple}
Let $A$ be a real $C \sp *$-algebra.
\begin{enumerate}
\item If $A$ is c-simple, then $A$ is simple.
\item If $A$ is simple and $A = \mathcal{F}(B)$ for some simple complex $C \sp*$-algebra $B$, then $A\sc \cong B \oplus B$ (thus $A\sc$ is not simple).
\item If $A$ is simple and $A\sc$ is not simple, then in fact $A \cong \mathcal{F}(B)$ for some complex $C \sp *$-algebra $B$.
\end{enumerate}
\end{Pro}

\begin{proof}
If $I$ is a non-trivial proper ideal of $A$, then $I\sc$ is a non-trivial proper ideal of $A\sc$. This proves Part (1).

For Part (2), suppose that $A = \mathcal{F}(B)$. For clarity in this proof, we write the elements of $A\sc$ in the form $a + j b$ for $a,b \in A$, and we reserve the symbol $i$ for the inherent complex scalar multiplication on $B$.

Then let
\begin{align*}
I &= \{ b + j(ib) \mid b \in B\} \subset A\sc \\
\text{and}~ J &= \{ b - j(ib) \mid b \in B \} \subset A\sc \;.
\end{align*}
There exists an isomorphism $B \cong I$ defined by
$$b \mapsto \tfrac{1}{2} (b + j(ib))$$
and similarly an isomorphism $B \cong J$.
Furthermore, $I \cap J = \{ 0\}$ and $I + J = A\sc$. This proves Part (2).

For Part (3), we refer the reader to Proposition~3.2 of \cite{BRS}.

\end{proof}

Because of this definition, we focus on real $C \sp *$-algebras $A$ that are not just simple, but are c-simple, as in Theorem~\ref{Thm-classification} below.

\begin{Def} \index{purely infinite $C \sp *$-algebra}
\begin{enumerate}
\item A complex simple $C \sp *$-algebra $A$ is {\it purely infinite} if $A$ is not isomorphic to $\CC$ and if for every pair of non-zero elements $a$ and $b$ in $A$ there are elements $x$ and $y$ such that $b = xay$.
\item A real c-simple $C \sp *$-algebra $A$ is {\it purely infinite} if $A$ is not isomorphic to $\RR$, $\CC$, or $\HH$ and if for every pair of non-zero elements $a$ and $b$ in $A$ there are elements $x$ and $y$ such that $b = xay$.
\end{enumerate}
\end{Def}

Significantly, we know that a real c-simple $C \sp *$-algebra $A$ is purely infinite if and only if $A\sc$ is purely infinite by Theorem 3.3 of \cite{Stacey2003-PI} and Theorem~3.9 of \cite{BRS}.

There are equivalent definitions of purely infinite in Lemma~3.8 of \cite{BRS} in the real case. Also, see Proposition~4.1.1 of \cite{RordamYellow} for alternative characterizations of this definition in the complex case, but all are quite technical. This class of $C \sp *$-algebras does not include any commutative $C \sp *$-algebras, nor any finite dimensional $C \sp *$-algebras, nor any direct limits of finite dimensional $C \sp *$-algebras. 
However, there are a rich set of examples that are simple and purely infinite, including the Cuntz algebras described in Example~\ref{NG-Models}.

The following proposition indicates that for a simple purely infinite $C \sp *$-algebra the group $K_0(A)$ can expressed without using matrix algebras, and without making use of the abstract Grothendieck construction. This gives a glimpse into why classification using $K$-theory is possible in this setting. 

\begin{Pro}[See Proposition~4.1.4 in \cite{RordamYellow} and Proposition~3.13 in \cite{BRS}]
\begin{enumerate}
\item If $A$ is a complex simple purely infinite $C \sp *$-algebra,
then any element of $K_0(A)$ can be represented as $[p]$ for some projection $p \in A$.
\item If $A$ is a real c-simple purely infinite $C \sp *$-algebra,
then any element of $KO_0(A)$ can be represented as $[p]$ for some projection $p \in A$.
\end{enumerate}
\end{Pro}

\begin{Pro}[Zhang's dichotomy. See Proposition~4.1.3 in \cite{RordamYellow} and Theorem~3.5 in \cite{BR}]
Suppose that $A$ is a complex or real simple purely infinite $C \sp *$-algebra. Then $A$ is either unital or stable.
\end{Pro}

The following theorems summarize the powerful classification theorems for simple purely infinite $C \sp *$-algebra. These results are due to Eberhard Kirchberg \cite{Kirchberg94} and N. Christopher Phillips \cite{phillips-class} in the complex case, and they are summarized in Mikael R\o rdam's monograph \cite{RordamYellow}. In the real case, the results are due to Boersema in \cite{BoersemaSurj} and Boersema, Ruiz, and Stacey in \cite{BRS}.

In the complex case, a simple separable nuclear purely infinite $C \sp *$-algebra is known as a {\it Kirchberg algebra}. 
\index{Kirchberg algebra} 
The proof of the classification Theorem~\ref{Thm-classification} below uses the Universal Coefficient Theorem in an essential way, so the hypothesis of the classification includes that the Kirchberg algebras must be in the bootstrap class. However, it is not known whether all Kirchberg algebras are in the bootstrap class.  In the real case, the classification applies to exactly those real $C \sp *$-algebras whose complexification is a Kirchberg algebra in the bootstrap class. \index{bootstrap category}

\begin{Thm} \label{Thm-classification} \index{classification theorem}
\begin{enumerate}
\item Let $A$ and $B$ be complex simple purely infinite stable $C \sp *$-algebras in the bootstrap class. Then $A \cong B$ if and only if $K_*(A) \cong K_*(B)$.
\item Let $A$ and $B$ be complex simple purely infinite unital $C \sp *$-algebras in the bootstrap class. Then $A \cong B$ if and only if $(K_*(A), [1]) \cong (K_*(B), [1])$.
\item Let $A$ and $B$ be real c-simple purely infinite stable $C \sp *$-algebras in the pre-bootstrap class. Then $A \cong B$ if and only if $K\crr(A) \cong K\crr(B)$.
\item Let $A$ and $B$ be real c-simple purely infinite unital $C \sp *$-algebras in the pre-bootstrap class. Then $A \cong B$ if and only if $(K\crr(A), [1]) \cong (K\crr(B), [1])$.
\end{enumerate}
\end{Thm}     

\begin{Rem}
In Part (2) above, the condition that $(K_*(A), [1]) \cong (K_*(B), [1])$ means that there is an isomorphism $\phi \colon K_*(A) \rightarrow K_*(B)$ and that this isomorphism carries the class of the identity in $K_0(A)$ to the class of the identity in $K_0(B)$, that is
$\phi([1_A]) = [1_B]$.
The condition in Part (4) is interpreted similarly.
\end{Rem}

\begin{Thm} \label{Thm-realization2}
\begin{enumerate}
\item For every ordered pair $(G_0, G_1)$ of countable abelian groups, there exists a complex simple purely infinite stable $C^*$-algebra $A$ in the bootstrap category such that 
$$K_*(A) \cong (G_0, G_1) \; .$$
\item For every ordered pair $(G_0, G_1)$ of countable abelian groups and element $g \in G_0$, there exists a complex simple purely infinite unital $C^*$-algebra $A$ in the bootstrap category such that 
$$(K_*(A), [1]) \cong (G_0, G_1, g) \; .$$
\item For every countable acyclic $CRT$-module $M$ there exists a real c-simple purely infinite stable $C^*$-algebra $A$ in the pre-bootstrap category such that 
$$K\crt_*(A) \cong M\; .$$
\item For every countable acyclic $CRT$-module $M$ and element $m \in M\po_0$, there exists a real c-simple purely infinite unital $C^*$-algebra $A$ in the pre-bootstrap category such that 
$$(K\crt_*(A), [1])  \cong (M, m) \; .$$
\end{enumerate}
\end{Thm}

Recall here that a $CRT$-module is acyclic if the sequences displayed in Theorem~\ref{Thm-realization} are exact. Note that Part (3) of Theorem~\ref{Thm-realization2} is the same as Theorem~\ref{Thm-realization}, but here we have extra properties stipulated on $A$.

\newpage   
  
 \section{\bf{$KO$-theory with Coefficients}}  \label{Section-KwithCoefficients}

 \index{K-theory@$K$-theory!$K$-theory with coefficients|(} 
  
  
    In \cite{TopIV}, the second author introduced mod $p$ coefficients into complex $K$-theory, obtained by tensoring with a $C \sp *$-algebra where the $K_0$ group is a prescribed group $G$. 
    In the present section we do something 
 similar with real $K$-theory. This situation is rather more complicated because the target category is graded modules over $KO_*(\RR)$, which has components isomorphic to $\ZZ$ and others isomorphic to $\ZZ_2$. 

  \begin{Def}
 Let $G$ be a countable abelian group. Then $N_G$ shall denote any complex $C^*$-algebra in the bootstrap class \index{bootstrap category} satisfying
 \[
  K_0(N_G) = G \quad \text{and} \quad K_1(N_G) = 0 \; .
   \]
\end{Def}

That such a $C^*$-algebra exists is a standard result (see for example the proof of Proposition~7.4 of \cite{RS} where such an $N_G$ is 
constructed and is also commutative). Furthermore, any two such $C^*$-algebras are $KK$-equivalent by Proposition~7.3 of \cite{RS} (this is a consequence of the Universal Coefficient Theorem \index{universal coefficient theorem} for complex $K$-theory). Furthermore, by the K\"unneth Theorem in \cite{TopII}, the following definition of $K$-theory with coefficients is well-defined (independent of the choice of $N_G$).

\begin{Def} In the context of complex $C \sp *$-algebras, $K$-theory with $G$-coefficients is defined as
 $$K_*(A; G) = K_*(A \otimes N_G) \; .$$ 
 \end{Def}
 
 \begin{Def}  \label{Def-NG}
\index[Notation]{NG@$N_G\pr$}
Let  $G$ be a countable abelian group. Then $N\pr_G$ shall denote any real $C^*$-algebra in the pre-bootstrap category such that
 \begin{itemize}
 \item $(N\pr_G)\sc = N_G$
 \item $KO_0(N\pr_G) = G$
  \item $c_0 \colon KO_0(N\pr_G) \rightarrow K_0((N\pr_G)\sc)$ is an isomorphism.
 \end{itemize}
 \end{Def}
 
 Note that we are making no restriction on $KO_i(N\pr_G)$ for $1 \leq i < 8$. Those groups are determined by the uniqueness part of Proposition~\ref{uniqueNG}.
 
 \begin{Pro} \label{uniqueNG}
 
 Let $G$ be a countable abelian group.  Then
 \begin{enumerate}
 \item There exists a real $C^*$-algebra  $N\pr_G$ in the pre-bootstrap category satisfying the three conditions listed above.
 \item The real $C^*$-algebra  $N\pr_G$ is unique in the sense that if there are two real $C^*$ algebras    $N\pr_G$
and  $\tilde N\pr_G$, each satisfying the three conditions listed above, then there is an isomorphism 
\[
 K\crt(N\pr_G) 
 \cong K\crt(\widetilde N\pr_G).
 \]
It follows that $N_G\pr$ and $\widetilde N_{G}\pr$ are $KK$-equivalent.
 \end{enumerate}
 \end{Pro}
 
 \begin{proof}
 In Section~8.4 of \cite{Hewitt}, Beatrice Hewitt presents a construction of an abstract $CRT$-module $M = (M\po, M\pu, M\ppt)$ having the structure of Definition~\ref{Def-NG}. In other words,
 \begin{itemize}
 \item $(M\pu)_0 = G$
  \item $(M\pu)_1 = 0$
 \item $(M\po)_0 = G$
 \item $c_0 \colon (M\po)_0 \rightarrow (M\pu)_0$ is an isomorphism.
 \end{itemize}
 Furthermore, $M$ is the unique $CRT$-module with these properties. 
 Then by Theorem~\ref{Thm-realization} or Theorem~\ref{Thm-realization2},
 there exists a concrete real $C \sp *$-algebra $N\pr_G$ with $K\crt(N\pr_G) \cong M$, and by 
Corollary~\ref{Cor-KKequiv}, this real $C \sp *$-algebra is unique up to $KK$-equivalence.
 \end{proof}
 
 Theorem~8.4.4 in \cite{Hewitt} describes the structure of $N_G\pr$ in general. For reference, we show the groups comprising $N\pr_\ZZ$ and $N\pr_\QQ$ below. 
 Note in particular that we can take $N\pr_\ZZ = \RR$. Later in this section we will also encounter $N\pr_{\Z_n}$.
 
  \begin{align*}
   &\begin{tabular} { c c c c c c c c c}
  $i$ & $0$ & $1$ & $2$ & $3$ & $4$   & $5$    &  $6$    & $7$ \\
  \hline
  $KO_i(N_{\ZZ}\pr)$ & $\ZZ$ & $\ZZ_2$  & $\ZZ_2$ & $0$ & $\ZZ$ & $0$ & $0$ & $0$ \\
  $K_i((N_{\ZZ}\pr)\sc)$ & $\ZZ$ & $0$   & $\ZZ$ & $0$ & $\ZZ$ & $0$ & $\ZZ$ & $0$ \\
    \end{tabular} \\~ \\
   &\begin{tabular} { c c c c c c c c c}
  $i$ & $0$ & $1$ & $2$ & $3$ & $4$   & $5$    &  $6$    & $7$ \\ 
  \hline
  $KO_i(N_{\QQ}\pr)$ & $\QQ$ & $0$  & $0$ & $0$ & $\QQ$ & $0$ & $0$ & $0$ \\
  $K_i((N_{\QQ}\pr)\sc)$ & $\QQ$ & $0$   & $\QQ$ & $0$   & $\QQ$ & $0$   & $\QQ$ & $0$ \\
    \end{tabular} \
   \end{align*}

 \begin{Def}
 \index[notation]{KOG@$KO_*(A;G)$}
 For a real $C^*$-algebra $A$, define
 $$KO_*(A; G) = K_*(A \otimes N_G\pr) \; .$$
 \end{Def} 
 
  
 \begin{Thm}
  $KO_*(A; G)$ is well-defined, independent of the choice of $N_G\pr$. The functor 
  $$A \rightsquigarrow KO_*(A; G)$$ 
 is a homology theory in sense of Definition~\ref{D:homology}.
 \end{Thm}
 
 \begin{proof}
 Suppose $N_G\pr$ and $\widetilde N_G\pr$ are both real $C^*$-algebras satisfying the necessary properties. 
 By Proposition~\ref{uniqueNG}, there is a $KK$-equivalence $\alpha \in KK\pr(N\pr_G, \widetilde N\pr_G)$. Then we get a map from the groups of this exact sequence to the groups of
the same sequence associated with $\widetilde{N}_G$.
Since $\alpha$ induces an isomorphism from $K\crt(N_G\pr)$ to $K\crt(\widetilde N_G\pr)$,
the five lemma implies that the homomorphism 
$$K\crt(A \otimes N_G\pr) \rightarrow K\crt(A \otimes \widetilde N_G\pr)$$
is an isomorphism.
In particular, there is an isomorphism
$$KO_*(A \otimes N_G\pr) \cong KO_*(A \otimes \widetilde N_G\pr) \; $$
and this proves that $KO_*(A; G)$ is well-defined.

That $KO_*(A; G)$ is a homology theory follows from Proposition~\ref{Proposition:homologytheory}.
\end{proof}

 \begin{Thm} \label{Prop-group-homomorphism}
 Let $\alpha \colon G \rightarrow H$ be a homomorphism of countable abelian groups. Then there exists a natural transformation 
 $$\theta_\alpha \colon KO_*(A;G) \rightarrow KO_*(A; H)$$
 such that $(\theta_\alpha)_0 = \alpha$.
 \end{Thm}
 
 \begin{proof}
 Suppose that $\alpha \colon G \rightarrow H$ is a homomorphism of countable abelian groups and let
$M_G$ and $M_H$ be the associated  $CRT$-modules as discussed in the first paragraph of the proof of Proposition~\ref{uniqueNG}. Then from the construction of $\mathcal{M}_G$ in Section 8.4 of \cite{Hewitt} it is clear that there is a homomorphism of $CRT$-modules $$\widetilde \alpha \colon {M}_G \rightarrow {M}_H$$ such that $(\widetilde \alpha)^O_0 = \alpha$.
 
 Now, let $N_G\pr$ and $N_H\pr$ be purely infinite simple stable real $C^*$-algebras in the bootstrap category, as in Theorem~\ref{uniqueNG}, such that
 $$K\crt(N_G\pr) \cong {M}_G \quad \text{and} \quad K\crt(N_H\pr) \cong {M}_H \; .$$ 
 The Universal Coefficient Theorem for real $C^*$-algebras from \cite{BoersemaUCT} tells us that
 $$\Gamma \colon KK\pr(N_G\pr, N_H\pr) \rightarrow \Hom\scrt ({M}_G, {M}_H)$$
 is surjective, so there exists a $KK$-element $\mathfrak{  a}$ such that $\Gamma(\mathfrak a ) = \widetilde \alpha$.
 Finally, by Theorem~\ref{Thm-UCT} there is a homomorphism 
 $$\psi_\alpha \colon N_G\pr \rightarrow N_H\pr$$
 such that $(\psi_\alpha)$ passes to $\mathfrak{a}$ under the natural homomorphism
 $$\Hom(N_G\pr , N_H\pr) \rightarrow KK\pr(N_G\pr , N_H\pr) \; .$$
 Then the desired natural transformation is given by
 $$\theta_\alpha = (\id \otimes \psi_\alpha)_* \colon KO_*(A \otimes N_G\pr) \rightarrow KO_*(A \otimes N_H\pr) \; .$$
   \end{proof}
  
\begin{Thm} \label{Thm:rationalKtheory} \index{K-theory@$K$-theory!rational $K$-theory}
For any real $C^*$-algebra $A$, there is an isomorphism
$$KO_i(A; \QQ) = KO_i(A) \otimes \QQ \; .$$
\end{Thm}

\begin{proof}
Let $\{\alpha\}_{n \in \NN}$ be a sequence of homomorphisms $\alpha_n \colon \ZZ \rightarrow \ZZ$ such that 
$$\lim_{n \to \infty} (\ZZ, \alpha_n) = \QQ \; .$$
Each homomorphism $\alpha_n$ is multiplication by some integer $k_n$.
Let 
$$\phi_n \colon N_\ZZ\pr \rightarrow N_\ZZ\pr$$
be a homomorphism of real $C^*$-algebras such that $\phi_n$ induces the homomorphism $\alpha_n$ on real $K$-theory in degree 0.
Then we calculate the $K$-theory of $\lim_{n \to \infty} (N_\ZZ\pr, \phi_n)$ using the direct limit structure and find that
$$ N_\QQ\pr = \lim_{n \to \infty} (N_\ZZ\pr, \phi_n) \; .$$
(Note we use the equal sign here, but we really mean that these $C \sp *$-algebras are $KK$-equivalent. After all, $N_\QQ\pr$ is only well-defined up to $KK$-equivalence.)

Then we also have
$$ \lim_{n \to \infty} (A \otimes N_\ZZ\pr, \id \otimes \phi_n) =A \otimes N_\QQ\pr \; .$$
Note that each map 
$$(\id \otimes \phi_n)_* \colon KO_*(A \otimes N_\ZZ\pr) \rightarrow KO_*(A \otimes N_\ZZ\pr)$$
is multiplication by $k_n$ on $KO_*(A \otimes N_\ZZ\pr) = KO_*(A)$.

Thus
\begin{align*}
	KO_*(A; \QQ) &= KO_*(A \otimes N_\QQ\pr) \\
		&=  \lim_{n \to \infty} (KO_*(A \otimes N_\ZZ\pr), (\id \otimes \phi_n)_*) \\
		&=  \lim_{n \to \infty} (KO_*(A), k_n) \\
		&= KO_*(A) \otimes \QQ
\end{align*}
\end{proof}

\begin{Cor} \label{Cor:RationalPeriodicity} 
Suppose that $G$ is a countable abelian torsion-free 2-divisible group. Then $K$-theory with coefficients in $G$ has period 4. In particular, there are isomorphisms
$$KO_i(A; \QQ) = KO_{i+4}(A; \QQ) \qquad \text{and} \qquad
 KO_i(A; \ZZ[\tfrac 12] ) = KO_{i+4} (A; \ZZ[\tfrac 12] ) \; .$$
\end{Cor}

\begin{proof}
The proof of this Corollary requires the K\"unneth Theorem and results about $CRT$-modules. Here we will only prove it in the case $G = \QQ$.
In fact, we will prove that there is an isomorphism
$$KO_i(A) \otimes G \cong KO_{i+4}(A) \otimes G$$ 
where $G$ is torsion-free and 2-divisible, which implies the result for $G = \Q$ using the isomorphism
$KO_i(A ; \QQ) \cong KO_i(A) \otimes \QQ$
from Theorem~\ref{Thm:rationalKtheory}.

Consider the natural transformation 
$$\xi \colon KO_i(A) \rightarrow KO_{i+4}(A)$$
given by multiplication by $\xi \in KO_4(\RR)$.
We know that $\xi^2 = 4 \beta_O$, by Theorem~\ref{thm:KO-relations}, where $\beta_O$ is the real Bott Periodicity isomorphism. This implies that 
$$(\xi \otimes \id)^2 \colon KO_i(A) \otimes G \rightarrow KO_{i+8}(A) \otimes G$$
is given by multiplication by 4. Then the hypotheses on $G$ imply that $(\xi \otimes \id)^2$ is an isomorphism, and this implies that
$\xi \otimes \id$ is an isomorphism.
\end{proof}

\index{K-theory@$K$-theory!$K$-theory with coefficients in $\Z_n$|(}

 In the remainder of this section, we investigate in more detail the properties of real $K$-theory with $\ZZ_n$-coefficients. Following the calculations of Section~5.1 of \cite{BoersemaKT}, we have the following structure for the $K$-theory of the real $C^*$-algebra $N_{\ZZ_n}\pr$:     
 
 \begin{align*}
 \text{For $n$ odd: \qquad }
   &\begin{tabular} { c c c c c c c c c}
  $i$ & $0$ & $1$ & $2$ & $3$ & $4$   & $5$    &  $6$    & $7$ \\
  \hline
  $KO_i(N_{\ZZ_n}\pr)$ & $\ZZ_n$ & $0$  & $0$ & $0$ & $\ZZ_n$ & $0$ & $0$ & $0$ \\
  $K_i((N_{\ZZ_n}\pr)\sc)$ & $\ZZ_n$ & $0$   & $\ZZ_n$ & $0$ & $\ZZ_n$ & $0$ & $\ZZ_n$ & $0$ \\
    \end{tabular} \\~ \\
    \text{For $n \equiv 2 \pmod 4$: \qquad }
   &\begin{tabular} { c c c c c c c c c}
  $i$ & $0$ & $1$ & $2$ & $3$ & $4$   & $5$    &  $6$    & $7$ \\ 
  \hline
  $KO_i(N_{\ZZ_n}\pr))$ & $\ZZ_n$ & $\ZZ_2$  & $\ZZ_4$ & $\ZZ_2$ & $\ZZ_n$ & $0$ & $0$ & $0$ \\
  $K_i((N_{\ZZ_n}\pr)\sc)$ & $\ZZ_n$ & $0$   & $\ZZ_n$ & $0$ & $\ZZ_n$ & $0$ & $\ZZ_n$ & $0$ \\
    \end{tabular} \\ ~ \\
  \text{For $n \equiv 0 \pmod 4$: \qquad }
   &\begin{tabular} { c c c c c c c c c}
  $i$ & $0$ & $1$ & $2$ & $3$ & $4$   & $5$    &  $6$    & $7$ \\ 
  \hline
  $KO_i(N_{\ZZ_n}\pr)$ & $\ZZ_n$ & $\ZZ_2$  & $\ZZ_2^2$ & $\ZZ_2$ & $\ZZ_n$ & $0$ & $0$ & $0$ \\
  $K_i((N_{\ZZ_n}\pr)\sc)$ & $\ZZ_n$ & $0$   & $\ZZ_n$ & $0$ & $\ZZ_n$ & $0$ & $\ZZ_n$ & $0$ \\
    \end{tabular}
   \end{align*}

\begin{Exam} \label{NG-Models}
We present here three natural ways to obtain models for $N_{\ZZ_n}\pr$.  

\begin{enumerate}

\vspace{.1in}
 \item {\bf Commutative Models.}
 
 Let $D^2$ be the unit disk and let $S^1$ be the unit circle in the complex plane.
 For $n \geq 2$, let $\zeta_n = \exp(2 \pi i /n) \in S^1$ be the principal $n$th root of unity and define an equivalence relation on $S^1 = \partial (D^2) \subset \CC$ generated by $z \sim \zeta_n \cdot z$. Then define the lens space by
 $$L_n = D^2/ \sim$$
and the reduced lens space by \index{lens space}
\index[notation]{ln@$\L_n$}
\index[notation]{lnz@$\Len_n$}
 $$\Len_n = (D^2/\sim) - [1] \; .$$
 Then we have
 $$\complex (\Len_n) \cong \{ f \colon D^2 \rightarrow \CC \mid \text{ $f(1) = 0$ and $f(z) = f(\zeta_n \cdot z)$ for all $z \in S^1$}  \} \; .$$
  There is a short exact sequence of complex $C^*$-algebras 
 $$0 \rightarrow \complex_0(\rm{int}(D^2)) \xrightarrow{~~\iota~~} \complex_0(\Len_n) \xrightarrow{\rm{~~ev~~}} \complex_0((S^1/\sim) - \{[1] \}) \rightarrow 0$$
 where $\rm{ev}$ is the evaluation map at the boundary. 
 We note that there are isomorphisms
 \[ \complex_0(\rm{int}(D^2)) \cong S^2 \CC \quad \text{and} \quad \complex_0((S^1/\sim) - \{[1] \}) \cong S \CC \]
 so in the resulting 6-term exact sequence on $K$-theory, the differential map is 
 $$\partial_1 \colon K_1(S \CC) \rightarrow K_0( S^2 \CC)$$ is multiplication by $n$ (from $\ZZ$ to $\ZZ$) and it follows that $$K_*(\complex_0(\Len_n)) = (0, \ZZ_n).$$ Thus $\complex_0(\Len_n)$ is (up to a suspension) an example of $N_{\Z_n}$. It is the complex model used by the second author in \cite{TopIV}.
 
 The same short exact sequence for real $C^*$-algebras
 \begin{equation} \label{Lens}  
0 \rightarrow \real_0(\rm{int}(D^2)) \xrightarrow{\iota} \real_0(\Len_n) \xrightarrow{\rm{~ev~}} \real((S^1/\sim) - \{[1] \})  \rightarrow 0
 \end{equation}
 yields a 24-term exact sequence on $KO$-theory. 
In the real case as in the complex case, we find that there are isomorphisms
 \[ \real_0(\rm{int}(D^2)) \cong S^2 \RR \quad \text{and} \quad \real_0((S^1/\sim) - \{[1] \}) \cong S \RR \; . \]
Then
 $$\partial_1 \colon KO_{-1}(S \RR) \rightarrow KO_{-2} ( S^2 \RR)$$
 is again multiplication by $n$ and from the 24-term long exact sequence we obtain
 \begin{equation*} 
 KO_{-2}( \real_0(\Len_n)) = \ZZ_n \quad \text{and} \quad KO_{-3}((\real_0(\Len_n)) = 0 \; .
 \end{equation*}
 This is a double shift of what is specified in the definition, so to correct for this we take
 $N_{\ZZ_n}\pr = S^{-2} \real_0(\Len_n)$.

  \vspace{.1in}
 \item {\bf Matrix Models.}
 
 \index{dimension-drop algebra} 
 \index[notation]{In@$\mathbb{I}_n\pr$}
 The real ``dimension-drop algebra" $\mathbb{I}_n\pr$ for any $n \geq 2$ is defined by
 $$\mathbb{I}_n\pr = \{ f \in C_0((0,1], M_n(\RR)) \mid f(1) \in \RR \cdot I_n \}$$
 as in Section~6 of \cite{Boersema2020}. The $K$-theory of these algebras is calculated in Theorem~6.1 of \cite{Boersema2020} and we find that $N_{\ZZ_n}\pr = S^{-1} \mathbb{I}_{n}\pr$.
 This is the real analog of the complex $C^*$-algebra used in \cite{DL96} for complex $K$-theory with coefficients.
  
  \vspace{.1in}
 \item {\bf Unital, Simple, Purely Infinite Models}. 
 
 \index{Cuntz algebra}   
 \index[notation]{On@$\mathcal{O}_n\pr$}
 The real Cuntz algebra $\mathcal{O}\pr_{n+1}$ is defined to be the unital real $C^*$-algebra generated by partial isometries $s_1, \dots, s_{n+1}$ that satisfy
 $$s_n^* s_n = 1 \text{~for all $k$} \quad \text{and} \quad \sum_{k=1}^{n+1} s_k s_k^* = 1 \; .$$
 Then by the real $K$-theory calculations for $\mathcal{O}\pr_{n+1}$ in
 Theorem~1.6.8 of \cite{Schr} and Section~5.1 of \cite{BoersemaKT}
 we have $N_{\ZZ_n}\pr =  \mathcal{O}\pr_{n+1} $.
Note that by Theorem~10.2 of \cite{BRS},
any unital simple purely infinite model for $N_{\ZZ_n}\pr$ is isomorphic to $\mathcal{O}_{n+1}\pr$ (this is stronger than what is stated in Proposition~\ref{uniqueNG}).
See Part (4) of Theorem~\ref{Thm-classification}.
 \end{enumerate} 
 \end{Exam}

The next two results generalize Proposition 1.6  in \cite{TopIV}  (the Bockstein exact sequence) and Part (2) of Proposition 2.6 in \cite{TopIV}  
from the complex case to the real case.

 \begin{Pro}
 For an integer $n \geq 2$, there exist natural transformations
 \begin{align*}
 \rho_n &\colon KO_*(A) \rightarrow KO_*(A; \ZZ_n) && \text{(of degree $0$)} \\
  \beta_n &\colon KO_*(A; \ZZ_n) \rightarrow KO_*(A) && \text{(of degree $-1$)}
 \end{align*}
  such that the following Bockstein sequence is exact
  $$\cdots \rightarrow KO_{k}(A) \xrightarrow{\rho_n} KO_k(A; \ZZ_n) \xrightarrow{\beta_n} 
  			KO_{k-1}(A) \xrightarrow{\times n} KO_{k-1}(A) \rightarrow \cdots$$
for $A$ in the pre-bootstrap category. \index{bootstrap category!pre-bootstrap category}
   \end{Pro}
 
 \begin{proof}
The existence of $\rho_n$ follows from Theorem~\ref{Prop-group-homomorphism} based on the surjective homomorphism $\ZZ \rightarrow \ZZ_n$, but to obtain $\beta_n$ and the exact sequence we take a different approach 
using the lens space model of $N_{\ZZ_n}\pr$. Start with the short exact sequence 
shown in displayed Equation (\ref{Lens}) above and tensor with a real $C^*$-algebra $A$. Since commutative $C^*$-algebras are exact, the resulting sequence 
\begin{equation*} 0 \rightarrow C_0(\rm{int}(D^2), A)  \xrightarrow{\iota} C_0(\Len_n, A) \xrightarrow{\rm{~ev~}} C_0(\topr, A) \rightarrow 0
 \end{equation*}
is still exact.
Using the isomorphisms 
$$C_0(\rm{int}(D^2), A) \cong S^2 A \quad \text{and} \quad C_0(\topr, A) \cong SA \; ,$$
the long exact sequence on $K$-theory can be formulated as
$$\cdots \rightarrow KO_{k}(S^2 A) \xrightarrow{\iota_*} KO_{k}(C_0(\Len_n, A)) \xrightarrow{\ev_*} KO_{k}(SA) 
	\xrightarrow{\partial} KO_{k-1}(S^2 A) \rightarrow \cdots \; .$$
Now we have $N_{\ZZ_n}\pr = S^{-2} \real_0(\Len_n)$, so
$$KO_i(A; \Z_n) = KO_i(A \otimes N_{\ZZ_n}\pr) 
	= KO_{i-2} ( C_0(\Len_n, A)) \; .$$
Using this, our long exact sequence can be written as
$$\cdots \rightarrow KO_{k+2}(A) \xrightarrow{\iota_*} KO_{k+2}(A;\ZZ_n)) \xrightarrow{\ev_*} KO_{k+1}(A) 
	\xrightarrow{\partial} KO_{k+1}(A) \rightarrow \cdots$$

This gives the natural transformations $\rho_n = i_*$ and $\beta_n = \ev_*$. It remains to show that the 
differential $\partial \colon KO_*(A) \rightarrow KO_*(A)$ is multiplication by $n$. 
The map arising from the same construction in the case of complex $K$-theory  is known to be multiplication by $n$ in for complex $K$-theory (see Lemma~1.4 of \cite{TopIV}). In particular, consider the $\partial$ maps when $A = \RR$ and when $A = \CC$ which commute with the complexification map as in the diagram below,  where each group is isomorphic to $\ZZ$ and the vertical maps are isomorphisms. Since the map $\partial$ is multiplication by $n$ on $K_0(\CC)$, the same is true for $\partial$ on $KO_0(\RR)$.

  \[
\xymatrix{ KO_0(\RR) \ar[r]^\partial \ar[d]_{c_0} & KO_0(\RR) \ar[d]_{c_0}  \\ 
K_0(\CC) \ar[r]^\partial & K_0(\CC) }
\]

Then it must also be multiplication by $n$ in the real case $KO_0(\RR)$ since there is an isomorphism 
$$c \colon KO_0(\RR) \rightarrow K_0(\CC) \; .$$ 

Now for any real $C \sp *$-algebra $A$ in the pre-bootstrap category, we have an isomorphism
$$KO_i(A) \otimes KO_0(\RR) \xrightarrow{} KO_i(A \otimes \RR)$$
by a special case of the K\"unneth Theorem in \cite{BoersemaKT};
and furthermore by naturality of the product map, the diagram below commutes
 \[
\xymatrix{  KO_i(A) \otimes KO_0(\RR) \ar[r] \ar[d]^{1 \otimes \partial} & KO_i(A \otimes \RR) \ar[d]^{\partial} & \\ 
KO_i(A) \otimes KO_0(\RR) \ar[r] & KO_i(A \otimes \RR) 
 }
\]
where the horizontal maps are isomorphisms of the K\"unneth Theorem.
Since the vertical map $1 \otimes \partial$ is multiplication by $n$, it then follows that $\partial$ is multiplication by $n$ for any $A$ in the pre-bootstrap category.  \end{proof}

\begin{Pro}
 For integers $m, n \geq 1$, there are  natural transformations
\begin{align*}
	 \kappa_{n,m} & \colon KO_i(A; \ZZ_m) \rightarrow KO_i(A; \ZZ_n)  \\
	 \beta_{n,m} & \colon KO_i(A; \ZZ_m) \rightarrow KO_{i-1}(A; \ZZ_n)  
\end{align*}
such that the following sequence is exact:
  $$\cdots \rightarrow KO_{k}(A; \ZZ_m) \xrightarrow{\kappa_{mn,m}} 
  			KO_k(A; \ZZ_{mn}) \xrightarrow{\kappa_{m,mn}} 
  			KO_{k}(A; \ZZ_n) \xrightarrow{\beta_{m,n} }
			KO_{k-1}(A; \ZZ_m) \rightarrow \cdots$$
 \end{Pro}

 \begin{proof}
 The existence of $\kappa_{n,m}$ follows from Theorem~\ref{Prop-group-homomorphism} based on the homomorphism
 $$ \ZZ_m \rightarrow \ZZ_n \; $$
 which is multiplication by $n/\gcd(m,n)$.
 Also, $\beta_{n,m}$ is defined to be the composition $\beta_{n,m} = \rho_m \circ \beta_n$.
 
 For the long exact sequence, we use the lens space model of $N_{\ZZ_n}\pr$. There is an obvious quotient map 
 $\Len_{n} \rightarrow \Len_{mn}$ which induces a injective homomorphism
 $$\alpha \colon N_{\ZZ_{mn}}\pr \rightarrow N_{\ZZ_{n}}\pr \; .$$
 It can be shown that the map induced by $\alpha$ on real $K$-theory  
 $$KO_0(N_{\ZZ_{mn}}\pr) \rightarrow KO_0( N_{\ZZ_{n}}\pr) \; ,$$ 
 is the natural surjective map from $\ZZ_{mn}$ to $\ZZ_n$. (This is also true of the map on complex $K$-theory.)
 Now we apply the mapping cone construction (Definition 4.9) to get the short exact sequence
 $$0 \rightarrow SN_{\ZZ_{mn}}\pr \rightarrow C \alpha \rightarrow N_{\ZZ_n}\pr \rightarrow 0 \; $$
 and the corresponding long exact sequence is
 $$\cdots \rightarrow KO_{i-1}(N_{\ZZ_{n}}\pr) \xrightarrow{\alpha_*} 
 	KO_{i-1} (N_{\ZZ_{mn}}\pr) \rightarrow
 	KO_i (C \alpha) \rightarrow
	KO_i( N_{\ZZ_n}\pr) \rightarrow \cdots  \; .$$
Note that the partial map arising from this long exact sequence coincides with the map $\alpha_*$ by Theorem 3.5 of \cite{TopIII}, which applies in the case of real $K$-theory just as for complex $K$-theory.
Using the fact that $\alpha_*$ is surjective on
$$KO_{0}(N_{\ZZ_{mn}}\pr) \rightarrow 	KO_{0} (N_{\ZZ_{n}}\pr)
	\quad \text{(which is $ \ZZ_{mn} \rightarrow \ZZ_{n}$)}$$
we can calculate that $KO_{-2}(C \alpha) = 0$ and $KO_{-1}(C \alpha) = K_{-1}((C \alpha)\sc) = \ZZ_m$. It follows that $C\alpha = S N_{\ZZ_m}\pr$ and that the map
$$KO_{i} (N_{\ZZ_{mn}}\pr) \rightarrow
 	KO_{i-1} (C \alpha)$$
is equivalent to the map
$$KO_i(N_{\ZZ_{mn}}\pr) \rightarrow KO_i(N_{\ZZ_{m}}\pr)$$
induced by multiplication by the quotient map $\ZZ_{mn} \rightarrow \ZZ_m$.
 \end{proof}
 
 \index{K-theory@$K$-theory!$K$-theory with coefficients in $\Z_n$|)}

We finish the presentation with a short exact sequence which can be used to compute $K$-theory with general coefficients, which is a simple consequence of the K\"unneth Theorem in Section~\ref{Section-Kunneth}.

\begin{Thm}
Let $A$ be any complex $C^*$-algebra and let $G$ be any countable abelian group. Then there is a split exact sequence of $K_*(\CC)$-modules,
$$0 \rightarrow K_i(A) \otimes G \xrightarrow{\alpha} K_i(A; G) \xrightarrow{\beta} \rm{Tor}(K_{i+1}(A), G) \rightarrow 0 \; $$
for $i = 0,1$, where $\alpha$ has degree $0$ and $\beta$ has degree $-1$.
\end{Thm}

\begin{proof}
We know that $N_G$ is in the bootstrap cateogory, so the K\"unneth Theorem (Theorem~\ref{ComplexKF}, Part (3)) applies, giving a short exact sequence with $K_*(A; G) = K_*(A \otimes N_G)$ as the middle term.
Since $K_0(N_G) = G$ and $K_1(N_G) = 0$, the outer terms are as shown.
\end{proof}

To be more explicit about this sequence in the complex case, there are in fact two split exact sequences:
\begin{align*}
0 \rightarrow K_0(A) \otimes G \xrightarrow{\alpha} &K_0(A; G) \xrightarrow{\beta} \rm{Tor}(K_{1}(A), G) \rightarrow 0 \;  \\
0 \rightarrow K_1(A) \otimes G \xrightarrow{\alpha} &K_1(A; G) \xrightarrow{\beta} \rm{Tor}(K_{0}(A), G) \rightarrow 0 \; 
\end{align*}
so each group of $K_i(A; G)$ can be expressed as a trivial extension of groups.

In the real case, the K\"unneth Theorem for real $C^*$-algebras must be expressed in terms of united $K$-theory, giving the result below. Again it follows from the 
Theorem~\ref{RealKF}, Part (3) giving us an expression of $K\crt(A; G)$ as an (not necessarily trivial) extension of two $CRT$-modules.
 
 \begin{Thm}
 Let $A$ be any real $C^*$-algebra and let $G$ be any countable abelian group. Then there is a short exact sequence of $CRT$-modules,
 $$ 	0 \rightarrow K\crt(A) \otimes\scrt K\crt(N_G) \xrightarrow{\alpha} K\crt(A; G) 
 		\xrightarrow{\beta}  
			\Tor\scrt(K\crt(A), K\crt(N_G)) \rightarrow 0 \; 
$$
where $\alpha$ has degree $0$ and $\beta$ has degree $-1$.
 \end{Thm}
 \index{K\"unneth Theorem}
 \index{K-theory@$K$-theory!$K$-theory with coefficients|)} 
 
 \newpage    
      \section{\bf{$K$-theory and Unitaries  }}
      \label{Section-KwithUnitaries}
      \index{K-theory@$K$-theory!unitary picture of $K$-theory|(}

In this section, we will present a description of real $K$-theory in which each group $KO_i(A)$ is expressed entirely in terms of certain unitaries. The picture of $K$-theory that we describe here is
well-known in the complex case for $K_0(A)$ and $K_1(A)$, as well as in the real case for $KO_0(A)$ and $KO_1(A)$, but is a recent development for $KO_i(A)$ in general.
It was developed in \cite{BL} to incorporate all cases of the Ten-Fold Way, and can be especially useful to understand the elements of $KO_i(A)$ in a concrete way.
 
The real $K$-theory groups $KO_i(A)$ of a real $C^*$-algebra $A$ were defined in Definitions~\ref{defn:realKOi} and \ref{defn:realKOi.v2} in terms of homotopy groups, or in terms of suspensions. The picture of 
$K_i(A\sc)$ and $KO_i(A)$ that we describe here is equivalent to those definitions, though it is not obvious that they are equivalent in most cases. It has the advantage of giving us a way to represent all elements of $KO_i(A)$ in terms of a specific (homotopy class of a) unitary in a matrix algebra over $A\sc$. This concrete representation also allows physicists to develop formulas for topological markers.

This unitary picture of $K$-theory is summarized in Table~\ref{table:summary1} below, which is reproduced from \cite{BL}. The table needs to be read very carefully, and we provide a detailed ``Guide to Table~\ref{table:summary1}" to explain how to read it. Also, the content of each row of the table is stated explicitly in the Propositions~\ref{prop:K0} through~\ref{prop:Klast} which follow in the exposition below.
For example, Propositions~\ref{prop:K0} and \ref{prop:KO0} below are summarized in the lines of Table~\ref{table:summary1} corresponding to $K_0(A\sc)$ and $KO_0(A)$, respectively. 
The other lines in Table~\ref{table:summary1} similarly show how to represent $K_i(A\sc)$ and $KO_i(A)$ for $i \neq 0$ in this unitary picture. At the end of the section, we provide some examples --- we give specific formulas for the unitaries representing non-trivial elements of $KO_i(A)$, for $A = \RR, A = S\RR$, and $A = S^{-1} \RR$.

\begin{Table}  
\centerline{\bf {Unitary Picture of $K$-theory:
Summary of Propositions~\ref{prop:K0} through \ref{prop:Klast} } \label{table:summary1} }
\vglue .5in
\begin{center}

\begin{tabular}{c|c|c|c|c|} 
\hhline{~----}
& K-group \TTT \BBB & $\ell$ & $\mathcal{S}$ & $I$ \\ \hhline{=====}
\multicolumn{1}{|c|}{\multirow{2}{*}{$A$ complex}}  & $K_0(A)$ \TTT \BBB & 2 & $u = u^*$ & $\smat{1}{0}{0}{-1}$ \\ \hhline{~----}
\multicolumn{1}{|c|}{} & $K_1(A)$ \TTT \BBB & 1 & -- & $1$ \\ \hhline{=====}
\multicolumn{1}{|c|}{\multirow{2}{*}{$A$ real}} & $KO_0(A)$ \TTT \BBB & 2 &  $u = u^*$, $u^{\tau} = u^* $ & $\smat{1}{0}{0}{-1}$  \\ \hhline{~----}
\multicolumn{1}{|c|}{} & $KO_1(A)$ \TTT \BBB & 1 &  $u^{\tau} = u^* $ & $1$ \\ \hhline{~----}
\multicolumn{1}{|c|}{} & $KO_2(A)$ \TTT \BBB & 2 &  $u = u^*$, $u^\tau = -u$ & $\smat{0}{i \cdot 1}{-i \cdot 1}{0}$ \\\hhline{~----}
\multicolumn{1}{|c|}{} & $KO_{3}(A)$ \TTT \BBB & 2 & $u^{\sharp \otimes \tau} = u$ & $ 1_2 $ \\ \hhline{~----}
\multicolumn{1}{|c|}{} & $KO_4(A)$ \TTT \BBB & 4 &  $u = u^*$, $u^{\sharp \otimes \tau} = u^*$ & ${\diag}(1_2,-1_2)$ \\ \hhline{~----}
\multicolumn{1}{|c|}{} & $KO_5(A)$ \TTT \BBB & 2 &  $u^{\sharp \otimes \tau} = u^*$ & $1_2$ \\ \hhline{~----} 
\multicolumn{1}{|c|}{} & $KO_6(A)$ \TTT \BBB & 2 &  $u = u^*$, $u^{\sharp \otimes \tau} = -u$ & $\smat{0}{i \cdot 1}{-i \cdot 1}{0}$ \\  \hhline{~----}
\multicolumn{1}{|c|}{}  & $KO_{7}(A)$ \TTT \BBB & 1 & $u^\tau = u $ & $1$ \\ \hhline{=====}
\end{tabular}
\end{center}
\end{Table} 

\vspace{1cm}

\newpage

\centerline{\bf{ Guide to Table~\ref{table:summary1}}}
\begin{itemize}
\item For a real $C \sp *$-algebra $A$, we work in the complexification $A\sc$ and we also make use of the antimultiplicative involution $\tau$ on $A\sc$ given by 
\[
\tau (a + ib) \equiv (a + ib)^\tau   = a^* + ib^*.
\] 
Also, let $\tau$ be the extended antimultiplicative involution defined by $\tau = \rm{Tr} \otimes \tau$
on $M_n(A\sc) = M_n(\RR) \otimes A\sc$.
Note that we are using the condensed notation $a^\tau = \tau(a)$ for $a \in A$ or $a \in M_n(A\sc)$ and recall that an element $a$ in $M_n(A\sc)$ is in the underlying real $C^*$-algebra $M_n(A)$ if and only if it satisfies $a^\tau = a^*$. 

\item
For a real $C \sp *$-algebra $A$, define
 $\mathcal{U}_n(A\sc)$ to be the set of unitary elements in $M_n(A\sc)$ and let 
\begin{align*} 
\mathcal{U}_{1, \infty}(A\sc) &= \bigvee_{n = 1}^\infty \mathcal{U}_{n}(A\sc) \\
\mathcal{U}_{2,\infty}(A\sc) &= \bigvee_{n = 1}^\infty \mathcal{U}_{2n}(A\sc) \\
\mathcal{U}_{4,\infty}(A\sc) &= \bigvee_{n = 1}^\infty \mathcal{U}_{4n}(A\sc) \; .
\end{align*}
\index[notation]{Mkzz@$M_{k, \infty}(A\sc)$}
\index[notation]{Ukzz@$\mathcal{U}_{k, \infty}(A\sc)$}


\item The value of $\ell $ is $1,2$ or $4$ in the second column of Table~\ref{table:summary1}. This indicates that $KO_i(A)$ is represented using equivalence classes of unitaries in 
$\mathcal{U}_{1, \infty}(A\sc)$, $\mathcal{U}_{2, \infty}(A\sc)$, or $\mathcal{U}_{4, \infty}(A\sc)$.

\item The entry in Column $\mathcal{S}$ gives us the critical relation that unitaries must satisfy to represent an element of $KO_i(A)$. For example, 
$KO_1(A)$ is represented using equivalence classes of unitaries in $\mathcal{U}_{n}(A\sc)$ that satisfy the relation $u^\tau = u^*$
and $KO_2(A)$ is represented using equivalence classes of unitaries in $\mathcal{U}_{2n}(A\sc)$ that satisfy both relations $u = u^*$ and $u^\tau = -u$. 


\item Two unitaries $u$ and $v$ in $\mathcal{U}_{\ell n}(A\sc)$ that satisfy the relation $\mathcal{S}$ for $KO_i(A)$ are equivalent if there is a continuous path of unitaries $U_t$ from $u$ to $v$ such that $U_t$ satisfies the relation $\mathcal{S}$ for all $t \in [0,1]$. We then write $[u] = [v]$ to indicate that $u$ and $v$ represent the same equivalence class in $\mathcal{U}_{n}(A\sc)$.

\item We identify unitaries in $\mathcal{U}_{\ell n}(A\sc)$ with unitaries in $\mathcal{U}_{\ell n+\ell}(A\sc)$ using the formula
\[ [a] = \left[ \left( \begin{smallmatrix} a & 0 \\ 0 & I \end{smallmatrix} \right) \right] \quad \text{in $KO_i(A)$} \]
where $I$ is the $\ell \times \ell$ matrix shown in the last column of Table~\ref{table:summary1}.

\item In the unitary picture of $K$-theory, 
\[ KO_i(A)  = \{ [u] \mid u \in  \mathcal{U}_{\ell, \infty}  \}  \]
where the equivalence class $[u]$ is determined by the equivalences described in the two previous points.

\item The group operation for $KO_i(A)$ is given by
$$[u] \oplus [v] = \left[ \left( \begin{smallmatrix} u & 0 \\ 0 & v \end{smallmatrix} \right) \right] \;  . $$
In particular, the group operation is {\it not} performed by simply multipliying representative unitaries (except for $KO_1(A)$ and $K_1(A\sc)$).



%
%
%

\item Finally, the unitary $I \in \mathcal{U}_{\ell, \infty}(A\sc)$ in the last column of the Table represents the trivial group element in $K_i(A\sc)$ or $KO_i(A)$. That is, 
\[
	[I] = 0  \; .
\]
\end{itemize}

In what remains of this section, we state the propositions that form the basis of this picture, and we present some examples.
For full proofs of these propositions see Section~5 of Boersema and Loring, \cite{BL} for the even case and 
Section 6 for the odd case. We first discuss the cases for $K_i(A\sc)$ and $KO_i(A)$ when $i$ is even. As we know, generally speaking, $KO_0(A)$ is built from projections in $A$ and $K_0(A\sc)$ is built from projections in $A\sc$. Note that if $p$ is a projection in a real or complex unital $C \sp *$-algebra $A$, then $u = 2(p-1)$ is a  self-adjoint unitary element. Conversely, if $u$ is a  self-adjoint unitary, then $p = \tfrac{1}{2}(u+1)$ is a projection. Using this relationship, we can express $K_0(A\sc)$ and $KO_0(A)$ in terms of unitaries instead of projections, as stated in the following propositions.

\begin{Pro} \label{prop:K0}
Let $A$ be a real unital $C \sp *$-algebra. Then $K_0(A\sc)$ is isomorphic to the group of homotopy classes of self-adjoint elements in 
$\mathcal{U}_{2,\infty}(A\sc)$.
\end{Pro}

\begin{Pro} \label{prop:KO0}
Let $A$ be a real unital $C \sp *$-algebra. Then $KO_0(A)$ is isomorphic to the group of homotopy classes of self-adjoint elements in 
$\mathcal{U}_{2,\infty}(A)$.
Equivalently, $KO_0(A)$ is isomorphic to the group of homotopy classes in $\mathcal{U}_{2,\infty}(A\sc)$ of self-adjoint elements $u$ that satisfy 
$u = u^{\rm{Tr} \otimes \tau}$.
\end{Pro}

%

Recall that for $KO_0(A)$ the group operation is given by $[u] + [v] = \left[ \left( \begin{smallmatrix} u & 0 \\ 0 & v \end{smallmatrix} \right) \right]$, and the inverse rule is given by $[u]^{-1}  = [-u]$. The group operation will be the same in all cases going forward, but the inverse rule is not always the same.

Remarkably, by making a small sign change to the relation $\mathcal{S}$ characterizing $KO_0(A)$, we obtain the unitary picture of $KO_2(A)$, as in Proposition~\ref{Prop:KO2} and as shown in the line of Table~\ref{table:summary1} for $KO_2(A)$.

\begin{Pro} \label{Prop:KO2}
Let $A$ be a real unital $C \sp *$-algebra. Then $KO_2(A)$ is isomorphic to the group of homotopy classes in $\mathcal{U}_{2,\infty}(A\sc)$ of self-adjoint elements $u$ that satisfy $u = -u^\tau$.
\end{Pro}

In the group of unitaries for $KO_2(A)$ described above, the identity element is given by $[I]$ where $I = \left( \begin{smallmatrix} 0 & i \cdot 1 \\ -i\cdot 1 & 0 \end{smallmatrix} \right)$. 

For $n = 4$ and $n = 6$ we let $\sharp$ be the involution 
\[\smat abcd \mapsto
		\smat{d}{-b}{-c}{a}  \]
on $M_2(A)$, and let $\sharp \otimes \tau$ be the corresponding involution on $M_{2n}(A) = M_2(\RR) \otimes M_n(A)$. Then we have the following propositions. For these propositions, the unitary $I$ such that $[I]$ represents the identity element of $KO_i(A)$ is listed in Table~\ref{table:summary1}.

\begin{Pro}
Let $A$ be a real unital $C \sp *$-algebra. Then $KO_4(A)$ is isomorphic to the group of homotopy classes in $\mathcal{U}_{4,\infty}(A\sc)$ of self-adjoint unitaries $u$ that satisfy $u = u^{\sharp \otimes \tau}$.
\end{Pro}

\begin{Pro}
Let $A$ be a real unital $C \sp *$-algebra. Then $KO_6(A)$ is isomorphic to the group of homotopy classes in $\mathcal{U}_{2,\infty}(A\sc)$ of self-adjoint unitaries $u$ that satisfy $u = -u^{\sharp \otimes \tau}$.
\end{Pro}

This completes the discussion of the even cases of Table~\ref{table:summary1}.
In the odd cases, we first discuss the well-known exponential map in the complex case. The following is Theorem 10.1.3 \cite{RLL}.

\begin{Pro} \label{complexK1}
Let $A$ be a real unital $C \sp *$-algebra. Then $K_1(A\sc)$ is isomorphic to the group of homotopy classes in $\mathcal{U}_{1,\infty}(A\sc)$.\end{Pro}


\begin{proof}[Sketch of Proof]
We have $K_1(A\sc)= K_0(S A\sc)$ (as presented in Theorem~\ref{T:K_*})
 and for the moment we 
let $K_1^{u}(A\sc)$ denote that group of homotopy classes in 
$\mathcal{U}_{1,\infty}(A\sc)$.
Recall that a unitary $u \in M_{n}(A\sc)$ is identified here with 
$\smat u001 $ in $M_{n+1}(A\sc)$. The identity is the class of the unit $[1]$ and the group inverse of an element $[u]$ in $K_1^u(A\sc)$ is given by $[u^*]$, since the unitary 
 $\smat u00{u^*} $ can be connected by a path of unitaries to 
 $\smat 1001$.
 
 We will sketch the definition of the map $\Theta$ that gives an isomorphism between $K_1^{u}(A\sc)$ and $K_0(S A\sc)$, following that in \cite{RLL}.

Let $u$ be a unitary in $M_n(A\sc)$. Let $x_t$ be a path of unitaries in $\mathcal{U}_{2n}(\CC)$ from 
$\smat{0}{1_n}{1_n}{0} $ to $1_{2n}$, and let $v_t$ be the path
$$v_t =  \smat u00{1_n}   x_t \smat{1_n}00{u^*} x_t^* $$
from $1_{2n}$ to $\smat u00{u^*} $.

Then $p_t = v_t  \smat{1_n}000  v_t^*$ is a projection in $C([0,1], M_{2n}(A\sc))$ with 
$p_0 = p_1 =\smat{1_n}000$.
The element 
$$ [p_t] - \left[ \smat{1_n}000\right]
	\in K_0( C([0,1], M_{2n}(A\sc))$$
is in the kernel of $(\ev_0)_* \colon K_0( C([0,1], A\sc ) ) \rightarrow K_0(A\sc)$, which therefore defines an element $\chi_u \in K_0(SA\sc)$.
Then we define $\Theta([u]) = \chi_u$. 

We refer to the proof of Theorem~10.1.3 for the details that $\Theta$ is well-defined and gives an isomorphism.
\end{proof}

\begin{Pro}
Let $A$ be a real unital $C \sp *$-algebra. Then $KO_1(A)$ is isomorphic to the group of homotopy classes in $\mathcal{U}_{\infty}(A\sr)$, which is isomorphic to the group of homotopy classes in $\mathcal{U}_{\infty}(A\sc)$ of unitaries $u$ that satisfy $u^* = u^\tau$. 
\end{Pro}

\begin{proof}[Sketch of Proof]
The second statement is true simply because the condition $u^* = u^\tau$ (for an element in $M_n(A\sc)$)  is equivalent to the condition that $u$ is in $M_n(A)$.

For the first statement, the proof is the same as in the complex case, with one extra point. 
Note that $\left( \begin{smallmatrix} 0 & 1_n \\ 1_n & 0 \end{smallmatrix} \right) $ is not in the same homotopy class of $1_{2n}$ in $M_{2n}(\RR)$ when $n$ is odd. However, we solve this by moving to a larger matrix algebra.
Specifically, for $n$ odd
$$\left[ \begin{smallmatrix} 0 & 1_n & 0  \\ 1_n & 0 & 0 \\ 0 & 0 & -1 \end{smallmatrix} \right] $$
is in the same homotopy class of unitaries as $1_{2n+1}$. 
Using this, we can produce a path $v_t$ of unitaries from $1_{2n+1}$ to $\diag(u,u^*,1)$ in $M_{2n+1}(\RR)$. The rest of the proof is the same as for Proposition~\ref{complexK1}.
\end{proof}

\begin{Pro}
Let $A$ be a real unital $C \sp *$-algebra. 
Then $KO_{-1}(A)$ is isomorphic to the group of homotopy classes in $\mathcal{U}_{\infty}(A\sc)$ of elements $u$ that satisfy $u = u^\tau$. 
\end{Pro}

\begin{Pro}
Let $A$ be a real unital $C \sp *$-algebra. Then $KO_3(A)$ is isomorphic to the group of homotopy classes in $\mathcal{U}_{2, \infty}(A\sc)$ of elements $u$ that satisfy 
$u= u^{\sharp \otimes \tau}$. 
\end{Pro}

\begin{Pro} \label{prop:Klast}
Let $A$ be a real unital $C \sp *$-algebra. Then $KO_5(A)$ is isomorphic to the group of homotopy classes in $\mathcal{U}_{2, \infty}(A\sc)$ of elements $u$ that satisfy 
$u^*= u^{\sharp \otimes \tau}$.
\end{Pro}

To understand better how this picture of real $K$-theory works, we will discuss some examples of real $C \sp *$-algebras and see how some of the specific non-trivial $KO$-classes can be represented by unitaries satisfying their appropriate relations. These examples and more can be found with more detailed exposition in \cite{BL} and \cite{Boersema2020}.

\begin{Exam}
$KO_i(\RR)$.
\end{Exam}

The table below shows the unitaries representing generators of all the non-trivial $KO_i(\RR)$ groups. Recall that $\RR\sc \cong \CC$, and the corresponding antimultiplicative involution $\tau$ on $\CC$ is the identity.

$$
\begin{array}{|c|c|c|}  \hline
  & \text{isomorphism class} & \text{unitary representing a generator} \\ \hline \hline
K_0(\RR\sc) & \ZZ & y_0 = 1_2   \\ \hline 
K_1(\RR\sc) & 0 &   \\ \hline  \hline
KO_0(\RR) & \ZZ & x_0 = 1_2   \\ \hline 
KO_1(\RR) & \ZZ_2 & x_1 = -1  \\ \hline 
KO_2(\RR) & \ZZ_2 & x_2  = \smat{0}{-i}{i}{0}   \\ \hline 
KO_3(\RR) & 0 & \\ \hline 
KO_4(\RR) & \ZZ & x_4 = 1_4   \\ \hline 
KO_5(\RR) & 0 & \\ \hline 
KO_6(\RR) & 0 & \\ \hline 
KO_7(\RR) & 0 & \\   \hline \hline 
\end{array}$$

As indicated by the table, the self-adjoint untiary $x_0 = 1_2$ represents a generator of $KO_0(\RR) \cong \ZZ$ 
(and recall that according to Table~\ref{table:summary1}, the self-adjoint unitary $I = \diag(1, -1)$ represents the trivial class of $KO_0(\RR)$).
More generally, if $x$ is a self-adjoint unitary in $M_{2n}(\RR)$, then $\tfrac{1}{2} \rm{tr}(x)$ is an integer, which determines the $KO_0(\RR)$ class represented by $x$. That is,  $[x] \mapsto \rm{tr}(x)$ is well-defined and gives an isomorphism from the unitary picture of $KO_0(\RR)$ to $\ZZ$. The same is true for $K_*(\RR\sc) = \ZZ$.

The unitary $x_1 = -1$ represents the non-trivial element of $KO_1(\RR) \cong \ZZ_2$. More generally, if $x$ is a unitary in $M_{2n}(\RR)$, then $\det(x)$ is either $1$ or $-1$, and this map characterizes the class represented by $x$ in $KO_1(\RR)$.

In the same way, the $KO_2(\RR)$ class of a self-adjoint, skew-symmetric matrix is determined by its pfaffian.

\vspace{1cm}

\begin{Exam} \index{suspension!the suspension algebra}
$KO_i(S\RR)$ 
\end{Exam}

$$\begin{array}{|c|c|c|c|}  \hline
  & \text{isomorphism class} & 
 	\multicolumn{1}{|c|}{ \text{unitary representing a generator}}  \\ \hline \hline
 K_0((S \RR)\sc) & 0 &   \\ \hline
 K_1((S \RR)\sc) & \ZZ & y_1 = x + iy  \\ \hline \hline
KO_0(S \RR) & \ZZ_2 & x_0 = \smat{x}{y}{y}{-x}  \\ \hline 
KO_1(S \RR) & \ZZ_2 & x_1 =  \smat{x}{-y}{y}{x}  \\ \hline
KO_2(S \RR) & 0 & & \\ \hline 
KO_3(S \RR) & \ZZ & x_{3}  \smat{x+iy}{0}{0}{x+iy}  \\ \hline
KO_4(S \RR) & 0 & \\ \hline 
KO_5(S \RR) & 0 & \\ \hline 
KO_6(S \RR) & 0 & \\ \hline 
KO_{7}(S \RR) & \ZZ & x_{-1} = x + iy   \\ \hline 
\end{array}$$

Since $S \RR = \realo(\RR)$ is not unital, we must consider these elements in the unitization 
$$(S \RR)^+ \cong  \real(S^1)$$ and represent $KO_*(S \RR)$ as the kernel of the evaluation map
$$\pi_* \colon KO_*( \real ( S^1))  \rightarrow KO_*(\RR)$$
where $\pi$ is evaluation at $z = 1 \in S^1$. The unitaries in the third column of the table show the generators of the various non-trivial $KO_i$ groups in this picture. 
Furthermore, to make sense of this table, we need to think about the complex $C \sp *$-algebra $\complex (S^1)$ and the corresponding involution $\tau$ on $\complex (S^1)$. In this case we have $\tau = \id$.

For $KO_0(S \RR)$, we observe that $x_0 = \smat{x}{y}{y}{-x}$ is a self-adjoint unitary matrix in $M_2( \real (S^1))$. Since $\pi(x_0) = \diag(1,-1)$ represents the trivial element of $KO_0(\RR)$, we find that $x_0$ is a well-defined element of $KO_0(S \RR)$. According to the table, $[x_0]$ represents the non-trivial element of $KO_0(\RR) \cong \ZZ_2$ (whereas the trivial element of $KO_0(S \RR)$ would be represented by the constant matrix $\diag(1,-1)$). We also point out that one can confirm the fact that $KO_0(S \RR) \cong \ZZ_2$ by verifying that $\diag(x_0, x_0)$ can be connected to the identity by a homotopy through self-adjoint unitaries. 

Now $KO_{-1}(S \RR) = \ZZ$ is represented by unitaries $u$ that satisfy the relation 
$u = u^\tau$, which in this case collapses to $u = u^{\rm{Tr}}$.
This relation is clearly satisfied by the unitary $x_{-1} = x + iy \in M_1(\complex(S^1))$ 
and we also have that $\pi[x_{-1}] = [1]$ which is trivial in $KO_{-1}(\RR) = 0$.
Hence, we have $[x_{-1}] \in KO_{-1}(S \RR)$, and in fact $x_{-1}$ represents one of the generators of $KO_{-1}(S \RR)$.

Finally, the reader can check that the unitary $x_3$ satisfies $x_3^{\sharp \otimes \tau} = x_3$. This element turns out to be a non-trivial element of $KO_3(S \RR) \cong \ZZ$ and represents a free generator of that group.

\vspace{1cm}

\begin{Exam} \index{suspension!the desuspension algebra}
$KO_i(S^{-1} \RR)$.
\end{Exam}

$$\begin{array}{|c|c|c|c|}  \hline
  & \text{isomorphism class} & 
 	\multicolumn{1}{|c|}{ \text{unitary representing a generator}}  \\ \hline \hline
KU_0(S^{-1} \RR) & 0 &  \\ \hline
KU_1(S^{-1} \RR) & \ZZ & y_{1} = x+iy  \\ \hline  \hline
KO_0(S^{-1} \RR) & 0 &  \\ \hline
KO_1(VS^{-1} \RR) & \ZZ & x_{1} = x+iy  \\ \hline 
KO_2(S^{-1} \RR) & \ZZ_2 & x_{2} =  \smat{y}{ix}{-ix}{y} \\ \hline 
KO_3(S^{-1} \RR) & \ZZ_2 & x_{3} = \smat{x+iy}{0}{0}{x-iy} \\ \hline 
KO_4(S^{-1} \RR) & 0 &   \\ \hline 
KO_5(S^{-1} \RR) & \ZZ & x_{5} = \smat{x+iy}{0}{0}{x+iy}  \\ \hline 
KO_6(S^{-1} \RR) & 0 & \\ \hline 
KO_7(S^{-1} \RR) & 0 & \\  \hline 
\end{array}$$

\vspace{1cm}
Recall that 
$$S^{-1}  \RR \cong \{f \in C_0(\RR, \CC) \mid \text{$f(-x) = \overline{f(x)}$ for all $x \in \RR$} \}$$ and that 
$$KO_*(S^{-1} \RR) = \Sigma^{-1} KO_*(\RR) \; .$$ 
The algebra $S^{-1} \RR$ and the one $S \RR$ in the previous example both have the same complexification (and the same unitization of the complexification $\complex(S^1)$). The difference is the associated involution on the complex $C \sp *$-algebra $\complex(S^1)$ is now 
$f^\tau(z) = f(\overline{z})$ for $f \in \complex(S^1)$ and $z \in S^1 \subset \CC$. 
The element $x_1 = x+iy$ satisfies $x_1^{\tau} = x_1^*$, so it represents an element in $KO_1(S^{-1} \RR)$ (whereas in the previous example the same unitary represented an element in $KO_{-1}(S \RR)$).

We leave it to the reader to check that the given elements $x_2$, $x_3$, and $x_5$ are also unitaries and that they satisfy the appropriate relations to represent the appropriate group $KO_i( S^{-1} \RR)$. See \cite{Boersema2020} for further details.

We find it elucidating to note that $KO_0(S^{-1} \RR)$ is represented by unitaries in $M_{2n}(\complex(S^1))$ satisfying 
$u = u^* = u^\tau$. Furthermore, the condition that $\pi(1)$ represents the trivial element of $KO_0(\RR)$ implies that the trace of such a unitary must be 0. Using the functional calculus, then one can easily show that $u$ is equivalent to a constant unitary with trace 0. Thus we see in a very concrete way that $KO_0(S^{-1} \RR) = 0$.

 \index{K-theory@$K$-theory!unitary picture of $K$-theory|)}

\newpage
  \part{$KO^*$-theory for Spaces}
  \vglue .3in
   
  Part 3 deals with spaces, and there is a basic question: which spaces? There are four choices in the literature:
  \begin{itemize}
\item   Finite CW complexes (which includes smooth compact manifolds)
\item        Infinite CW complexes (including the finite ones but also spaces like BO) 
\item       compact spaces (includes the finite CW complexes but not spaces like BO)
\item       compactly generated spaces (includes all of the above, to be explained)
\end{itemize}

Each of these leads to $KO^*(X)$.  The second option cannot be described directly using vector bundles and so we shall not consider it.  Most $K$-theory books concentrate 
upon $KU$-theory for finite CW complexes. Nevertheless, we take the third option in this monograph as this corresponds to all unital commutative complex $C^*$-algebras and to all of the unital commutative real $C^*$-algebras (when equipped with an involution). The motivation gets very concrete as well. We want to use $K$-theory to understand essentially normal operators. The 
essential spectrum of such an operator may be any compact subset of the plane! So the only choice for us is the third option.

When we get to $K$-theory in Part 3 we focus almost entirely upon the real case but we will also comment on the complex case for context. 
Atiyah's book \cite{Atiyah} is devoted to complex $K$-theory for compact spaces.
Karoubi's book \cite{K} covers both real and complex $K$-theory and does it for compact spaces. 

There are many books on
classical algebraic topology; we particularly recommend 
G. W. Whitehead \cite{W} for basics, 
J.F. Adams \cite{A} for spectra, cohomology theories, and related matters,
and J.P. May \cite{May} for a modern treatment.  
   \section{\bf{The Basics: Topological Spaces } } \label{top spaces}

 In this section we review the information we need about topological spaces,
  loop spaces and classifying spaces of $\oo $ and $\uu $,    in order to present $KO$-theory for spaces.   
    
   \index[notation]{uuz@$\uu$} 
    
 \begin{Def} A {\emph{based topological space}} \index{based topological space} $(X, x_o)$ is a topological space $X$ together with a choice of basepoint $x_o \in X$.  A {\emph{based map}}   
 \[
 f: (X,x_o) \to (Y, y_o) 
 \]
 of based topological spaces (sometimes simply written $f: X \to Y$ when the context is clear) is a continuous map  from $X$ to 
   $Y$ with the property that $f(x_0) = y_0$.  We denote the set of based maps from $X$ to $Y$ by $F(X, Y)$ and
   based homotopy classes of based maps by $[X,Y]$.    In particular, 
   we may define the k'th {\emph{homotopy group}} \index{homotopy group} of $X$ by 
   \[
   \pi_k(X) = [S^k, X].
   \]
   
   \end{Def}
   Note that on occasion we will want to deal with spaces and maps without specified basepoints, in which case we will denote the set of maps by $F(X,Y)^{ub}$, where $ub$ refers to 
  ``unbased." \footnote{We considered using the term ``baseless" but that seemed too derogatory. }
 As noted previously, we follow Bourbaki and  understand compact spaces to be compact and Hausdorff.

          In order to proceed, and as noted in the Introduction, we must put a small point-set topological restriction on our topological spaces: we follow G. Whitehead \cite {W} and May \cite{May}.
 
 \begin{Def} A space $X$ is {\emph{compactly generated}}  \index{compactly generated space} if $X$ is Hausdorff   and if it has the following property:  if $Y$ is a subset of $X$ and $Y \cap C$ is closed in $X$  for each compact 
 subset $C$ of $X$, then $Y$ itself is closed in $X$. 
 \end{Def}
 
 Every locally compact space is compactly generated. So is every metrizable space, so this includes smooth manifolds and in particular all of the topological groups that we will need as well as the homogeneous spaces \index{homogeneous space} such as the space $O_{2k}/(O_k \times O_k)$. 
 All   CW complexes, including spaces like $\oo $  are compactly generated.     If you confront a space $X$ that is not compactly generated then it is easy to retopologize it canonically to get a new compactly generated  space 
 $kX$ and a map $kX \to X$ which is the identity function and is continuous.\footnote{See May (\cite{May} p. 37) for details.} 
 
 {\bf{We assume henceforth  that all spaces are compactly generated.}}

 Next we need some basic homotopy constructions.  Some of the constructions need to be retopologized (as just described) to 
 make them stay in the category of compactly generated spaces, and we will assume that that  automatically happens when needed.  For example, if we take the function space $F(X,Y)$ of all based maps 
 from $X$ to $Y$ with the compact-open topology then that requires a slight adjustment to make it compactly generated.  The basepoint of $F(X,Y)$ will be the constant function $f(x) = y_0$ for all $x \in X$.
 
 \begin{Def}  \label{Def-homotopy}
 Suppose that $X$ and $Y$ are based compactly generated spaces. Then:
  \begin{enumerate}
  \item
  The {\emph{product}} \index{product of spaces} $X \times Y$ is given by  
\[
X \times Y = \{(x,y) \mid x \in X, y \in Y\}
\index[notation]{XtimesY@$X \times Y$}
\]
 with basepoint $(x_0, y_0)$ and the product topology.
  \item
The {\emph{join}} (or {\emph{one-point union}}) \index{join of spaces}  of $X$ and $Y$, denoted $X \vee Y$, is the disjoint union of $X$ and $Y$ with $x_0 $ identified 
 with $y_0$.  Alternatively, think of this as a subset of $X \times Y$, namely 
 \[  X \vee Y = \{ (x, y) \in X \times Y \mid \text{$x = x_0$ or $y = y_0$ (or both)}  \}  \]
 \index[notation]{xveey@$X \vee Y$}
 topologized as a closed subset of $X \times Y$. 
  For example, the join of two circles is a figure 8 and the point that the two circles have in common is the basepoint.
   \item
The {\emph{smash product}} \index{smash product of spaces}  $X \wedge Y$ is given by  
\[
X \wedge Y = (X \times Y)/(X \vee Y),
\index[notation]{xwedge@$X \wedge Y$}
\]
 the product $X \times Y$ modulo the collapsing of the join 
$X \vee Y$ to a point, which becomes the basepoint. This is topologized as a quotient of $X \times Y$.    
For example, $S^1 \times S^1$ is the 2-torus, whereas $S^1 \wedge S^1$ is the 2-sphere  $S^2$.
 
  \item
 \index[notation]{Sx@$SX$}
 The {\emph{suspension}} \index{suspension!of a space} of a topological space $X$ is the space 
 \[
 SX = S^1 \wedge X.
 \]
   We let 
 \[
 S^nX = S^n \wedge X
 \]
  and note that 
  \[
  S(S^nX) \cong S^{n+1}X.
  \] 
  Sometimes in the literature this is called the {\emph{reduced suspension}} of $X$, with the term ``suspension" reserved for $S^1 \times X$.   
 
  \item   The {\emph{loop space}} \index{loop space} of a topological space $X$ is defined to be  
  \index[notation]{Ozz@$\Omega X$}
  \[
  \Omega X = F(S^1, X)
  \]
   and higher loop spaces are given by  
   \[
   \Omega^n  X \cong F(S^n, X).
   \]
        Note that 
        \[
        \Omega (\Omega ^n X ) \cong \Omega ^{n+1} X.
        \]
 
  \item The {\emph{pullback}} \index{pullback} construction is defined exactly as for $C^*$-algebras. That is, given a diagram 
 
 \[
\CD
 @.      X \\
@.    @VVfV \\
Y   @>g>>    Z
\endCD
\]
 then the {\emph{pullback}} $P$ is defined to be 
 \[
 P = \{ (x,y) \in X \times Y \mid f(x) = g(y)  \}
 \]
 with the subspace topology and with associated commuting diagram
  \[
\CD
 P @>>>      X \\
@VVV    @VVfV \\
Y   @>g>>    Z
\endCD
\]
 
  \item The {\emph{pushout}} \index{pushout} construction is dual. That is, given a diagram
  \[
\CD
 W @>f>>      X \\
@VVgV    @. \\
Y   @. 
\endCD
\]
then the {\emph{pushout}} $Q$ is defined by taking the disjoint union of $X$ and $Y$ modulo the equivalence relation generated by   $f(w)  \sim g(w)$ for each $w \in W$.  
 \[
\CD
 W @>f>>      X \\
@VVgV    @VVV\\
Y   @>>>   Q
\endCD
\]

There is a very important special case . Suppose 
that $Y = D^n$, the $n$-dimensional disk with boundary sphere $S^{n-1}$, $\io $ is the  inclusion map 
\[
\io: S^{n-1} \to D^n
\]
and we are given a map $f: S^{n-1} \to X$.    Then we have a diagram 
 \[
\CD
S^{n-1} @>f>>      X \\
@VViV    @ VVV \\
D^n   @>>>  X \bigcup _f D^n 
\endCD
\]
and the associated pushout is denoted $X \bigcup _f D^n$.
This construction is called {\emph{adding an $n$-cell to $X$}}.   Here's a specific example: if $f: S^3 \to S^2$ is the Hopf map, then 
\[
S^2 \bigcup _f D^4 \cong \CC P^2
\]
the complex projective plane.
 \end{enumerate}
 \end{Def}
  
 Both $S$ and $\Omega $ are covariant functors.  If $f : X \to Y$ is a based map then there are natural induced maps $Sf : SX \to SY$ and $\Omega f : \Omega X \to \Omega Y$ with the requisite 
 properties.  Better yet, we have the following classical result, which is basic in homotopy theory and of critical importance in understanding the physics invariants to come.

 \begin{Thm}  \label{Thm-homotopy}  Suppose that $X$ and $Y$ are based and compactly generated. Then   the functors $S$ and $\Omega $ are adjoint functors \index{adjoint functors}.  That is, there is a natural homeomorphism 
 \[
 F(SX, Y) \overset{\cong}\longrightarrow F(X, \Omega Y).
 \]
 and, by iterating, 
 \[
 F(S^nX, Y) \overset{\cong}\longrightarrow F(X, \Omega ^n Y).
 \]
 \end{Thm} 
   
  \begin{proof}  We will prove something stronger, namely that if $W$ is also based and compactly generated, then 
  there is a natural homeomorphism 
  \[
 F(W \wedge X, Y) \overset{\cong}\longrightarrow F(X, F(W, Y)).
 \] 
  Suppose that $f : W\wedge X \to Y$.  Define $\phi (f) :  X \to F(W,Y) $ by 
  \[
  \phi (f) (x) (w) = f(w\wedge x).
  \]
     One checks that this map is well-defined, continuous, preserves basepoints, 
  and then write down its inverse as follows.   If $g : X \to F(W,Y)$ then define 
  $\psi (g) : W \wedge X \to Y$ 
  by 
  \[
  \psi (g)(w\wedge x) = g(x)(w).
  \]
    Then check that $\psi $ is well-defined, continuous, preserves basepoints, and that $\Psi$ is the inverse for $\phi $. 
  \end{proof}

Next, we must briefly discuss  CW-complexes. \index{CW-complex|(}This is a big classical subject and is very nicely developed in many texts, such as G.W. Whitehead \cite{W} Ch. 2, J. P. May \cite{May} Ch. 10, and John M. Lee \cite{LeeTopMan} Ch. 5.   We will give a small taste with no details, which will suffice for our needs. 
  
 We have explained what is meant by {\emph{adding an $n$-cell}} above.
  There are easy examples. If $X$ is a point then $X \bigcup _f D^n \cong S^n$.  If $X = S^2$ and $f: S^3 \to S^2$ is the Hopf map then $S^2 \bigcup _f D^4 \cong CP^2$.
  
  \begin{Def}  \label{CWcomplex}  A {\emph{CW-complex}} $X$ is a topological space $X$ which is the union of spaces $X^0 \subseteq X^1  \subseteq X^2  \subseteq X^3 \subseteq \dots $    where $X^0$ is 
  a discrete set of points, called vertices,  and inductively, $X^n$ (called the {\emph{n-skeleton of $X$}}) is obtained from $X^{n-1}$ by adding on some\footnote{Perhaps none and perhaps an infinite number}  n-cells.   A map $f: X \to Y$
  of CW-complexes is said to be {\emph{cellular}} if $f(X^n) \subseteq Y^n$ for each n.  If only a finite number of cells are added in total then $X$ is called a 
  {\emph{finite CW-complex}} \index{CW-complex!finite CW-complex}
  and, of course, $X^k = X$ for $k$ sufficiently large.
  \end{Def}
  
  It is easy to give examples of CW-complexes, though the structure may not be obvious at first.  For instance, any compact connected manifold $M$ is homotopy equivalent to a finite CW-complex with one vertex.  Thus spaces like $U_n$ and $SO_n$ fit this description, as does $O_n$ if we allow two vertices. On a much deeper level, we may, for example, regard the sequence 
  \[
  O_1 \subset O_2 \subset O_3  \subset \dots \subset O_n \subset \dots \oo
  \]
  as a sequence of  skeleta with $X^n = O_{n+1}$, showing that $\oo $ is a CW-complex.     
  Similarly, the sequence 
  \[
  RP^1 \subset RP^2 \subset RP^3 \subset  \dots \subset RP^\infty
  \]
  tells us how real projective spaces are constructed.

  How does all of this help us?    Here's a nice consequence.
  
  \begin{Thm}\label{T:homotopylimit}     Suppose that $X$ is a CW-complex with skeleta $X^n$ and associated maps $\rho_n : X^n \to X$. 
  Then  
  \begin{enumerate}
  \item
  The map 
  \[
  (\rho_n)_* : \pi_k(X^n) \longrightarrow \pi _k(X)
  \]
  is an isomorphism for $k < n$ and is surjective for $k = n$.   
  \item The natural induced map $\rho$
  \[
\lim_{n \to \infty}  \pi _k(X^n) \overset{\rho}\longrightarrow  \pi_k(X)
  \]
  is an isomorphism.
  \end{enumerate}
  \end{Thm}
  
  For example, it is not hard to show that $\pi_1(RP^2) = \ZZ _2$ and hence 
  $\pi_1(RP^n) = \ZZ_2$ for all $n \geq 2$ and   $\pi_1(RP^\infty ) = \ZZ_2$

 \begin{proof}  This is an easy consequence of the Cellular Approximation Theorem \index{cellular approximation theorem} (cf. G. W. Whitehead \cite{W} II. 4.5) which says, in part, that if $f: Y \to X$ is 
 a map of CW-complexes, then $f$ is homotopic to a cellular map $f'$,  a map that has the property that $f'(Y^j) \subseteq X^j$ for each $j$.    Taking $Y = S^j$, 
 a sphere, this tells us that a map $f: S^j \to X$ is homotopic to a map $S^j \to X^j$.  So, for example, if $n > 2$ then  any map $S^2 \to S^n$ is null-homotopic.  
 \end{proof}
 \index{CW-complex|)}

          \newpage
    \section{\bf{The Basics: Fibrations, Fibre Bundles, Principal Bundles, and $\bbo$}}
       \label{Section:Fibrations}


  
  Fibrations and fibre bundles 
 underlie large portions of mathematics and physics, albeit sometimes in disguise.  For example, real and complex vector bundles are 
  going to be prominent for us later in this work.  For now we sketch out the basics. 
  
  There are several fine sources for this material, and we would particularly mention 
  Husemoller \cite{H}, Milnor \cite{M1} and \cite{M3}, Steenrod \cite{St}, G. W. Whitehead \cite{W}, 
  and, more recently,  J. P. May \cite{May}.  
  We present here merely a brief outline of the theory of fibrations, based upon May \cite{May} Chapter 7,  and urge the interested reader to explore further in those sources. 
      
   Lie groups come up in many chapters of this monograph, either as examples (such as $O_n$ or $U_n$) or in general. We recommend Chevallay~\cite{chev} or Chapter 6 of Lee \cite{Lee} as a basic reference.

       \begin{Def} A surjection $p: E \to B$ is a {\emph{(Hurewicz) fibration}} \index{fibration!of spaces}if it satisfies the covering homotopy property (CHP). 
       \index{covering homotopy property}  Given a commuting diagram 
       
         \[
\xymatrix{Y \ar[d]^{\iota} \ar[r]^f & E\ar[d]^p   \\ 
Y \times I \ar@{.>}[ur]^{\tilde{h}} \ar[r]^{h} & B}
\]
     with $\iota (y) = (y, 0)$, then there exists a map $\tilde{h} $ making the diagram commute.  In other words, if one end of the homotopy $h$ may be covered (by $f$)  to $E$ 
     then the whole homotopy may be    covered by $\tilde h$. 
       
       For each $b \in B$,   the space  $F_b = p^{-1}(b)$  is called  the {\emph{fibre}}\footnote{Some people spell the word ``fiber" but we follow Norman Steenrod, the second author's mathematical 
 ancestor who, quite literally, wrote the (first!) book on ``Fibre Bundles" \cite{St} in 1951.} over the point $b$.       
       \end{Def}
       
     We shall see many examples in future sections;  we mention only a few here.
     
     \begin{enumerate}
     
     \item {\underline{The trivial fibration}}
     
      If $E = B \times F$ with $p = \pi _1: B \times F \to B $ the projection on the first factor, then 
     $E \overset{p}\longrightarrow B$ is a fibration, and then $F_b \cong F$ for each $b \in B$.
     
     \item {\underline{The M\"{o}bius band}}   \index{M\"{o}bius band}
     
     Form the rectangle $[0,17] \times [0,1]$. If we then identify $(0, t)$ with $(17,t)$ we obtain a cylinder of width $1$ and circumference  $17$ and this is 
     a trivial fibration over the circle of circumference $17$.  However, if we identify $(0,t)$ with $(17, 1-t)$ then we obtain a non-trivial fibration 
     over the circle, with total space the M\"{o}bius band.
     
      \item {\underline{Coset spaces}} \index{coset spaces}
      
      If $G$ is a Lie group and $H$ a closed subgroup with coset space $G/H$ then the sequence
      \[
      H \longrightarrow G \longrightarrow G/H
      \]
      is a fibration with fibre $H$. 
      
      \end{enumerate}

     \begin{Pro} If $p: E \to B$ is a fibration and $f: B' \to B$ is some map then the pullback   $p': p^{-1}(E) \to B'$ in the associated pullback diagram
     \[
     \begin{CD}
  p^{-1}(E) @>>>    E \\
     @VV{p'}V       @VVpV \\
     B' @>f>> B
     \end{CD}
     \]
     is a fibration. For any $x \in B'$, the fibre $(p')^{-1}(x)$ is isomorphic to the fibre $p^{-1}(f(x))$.
     \end{Pro}

     \begin{Def} \index{numerable open cover} An open cover $\mathcal {U} = \{U_\alpha \}$ of a space $B$ is said to be {\emph{numerable}} if for each open set $U$ in the cover there is a continuous map $\lambda _U : B \to [0,1] $ such that 
     \[
     \lambda _U^{-1}(0, 1] = U
     \]
     and the cover is locally finite; that is, each $b \in B$ appears in only finitely many open sets of the cover.
 \end{Def}
     
  Any open cover of a paracompact\footnote{A space is {\emph{paracompact}} \index{paracompact} if every open cover of the space has a locally finite refinement. Examples include 
  compact spaces and metrizable spaces. }     space has a numerable refinement.
     
     \begin{Thm} (cf. \cite{May}, p. 49.)  Suppose that $\mathcal{U}$ is a numerable open cover of $B$. Then the following are equivalent:
     \begin{enumerate}
     \item The map $p: E \longrightarrow B$ is a fibration. 
     \item For every open set $U \in \mathcal{U}$, the pullback  
     \[
     \begin{CD}
     p^{-1}(U) @>>> E \\
     @VV{p'}V     @VVpV  \\
     U @>\subset >>    B
     \end{CD}
     \]

     is a fibration.
     \end{enumerate}
     \end{Thm}

     \begin{Thm} (cf. \cite{May} p. 52) If the base space $B$ of a fibration $p: E \to B$ is 
     path-connected\footnote{A space $X$  is {\emph{path-connected}} \index{path connected} if given any two points $x_0,  \, x_1 \in X$ there is a continuous path $x_t$ in $X$ connecting   them. }
      then any two fibres $F_x$ and $F_y$ are homotopy 
     equivalent.
     \end{Thm}

   Next, we record what will be for us one of the most useful properties of fibrations.  
     
     \begin{Thm} (cf. May \cite{May} p. 64).  Suppose that $p: E \to B$ is a fibration and $B$ is path-connected with basepoint $b_0$. Let $F = p^{-1}(b_0) $ and fix some basepoint in $F$.  Then 
     there is a natural long exact sequence of homotopy groups
     \[
     \to \pi_n(F)   \to \pi_n(E)   \overset{p_*}\longrightarrow \pi_n(B)  \overset{\delta }\longrightarrow \pi_{n-1}(F)   \to  \dots  \to \pi_0(B ) \to 0 .
     \]
     
     \end{Thm}

\begin{Pro}\label{loops}  Suppose that 
\[
F \longrightarrow E \longrightarrow B 
\]
is a fibration of CW-complexes and 
$E$ is contractible. 
Then $F$ is homotopy equivalent 
 to $\Omega B$.
\end{Pro}

\begin{proof}
The long exact homotopy sequence for fibrations degenerates, since $\pi _j(E) = 0$ for all $j$, and we are left with 
\[
0 \longrightarrow \pi_{j+1}(B) \overset{\cong}\longrightarrow \pi_{j}(F) \longrightarrow 0
\]
and so 
\[
  \pi _{j}(F) \cong \pi _{j+1}(B) \cong \pi_{j}(\Omega B) \qquad \forall j. 
\]
  If $F$ and $B$ are CW-complexes then so is $\Omega B$,  and this implies by Whitehead's theorem that they are homotopy equivalent.
  \end{proof}

     \begin{Def}  A map $p: E \to B$ is a {\emph{fibre bundle}} \index{fibre bundle}
     if  $B$ possesses a numerable open cover $\mathcal {U}$ and there is a fixed topological space $F$ and 
     homeomorphisms $\phi _U $ for each open set $U$ in the cover such that the diagram
     \[
     \begin{CD}
     U \times F @>{\phi _U}>\cong> p^{-1}(U) \\
     @VV{\pi _1}V       @VVpV \\
     U @>\subset>> B 
     \end{CD}
     \]
     commutes.   Every fibre bundle is a fibration.
     \end{Def}
     
     Perhaps the most well-known examples of fibre bundles are real and complex vector bundles, \index{vector bundle} 
     such as the tangent bundle of a smooth manifold. We will be 
     discussing those in future sections. For now we concentrate upon a different class of examples, namely principal $G$-bundles.
      This will be very sketchy and we recommend Husemoller \cite{H} starting on page 42 for a detailed treatment.

     \begin{Def} Suppose that $G$ is a topological group and $p: E \to B$ is a fibre bundle together 
     with  a continuous right action $E \times G \to G$. Then  $p: E \to B$ is a {\emph{principal $G$-bundle }} 
     \index{principal $G$-bundle}
     and is homeomorphic to $E \to E/G$ 
     if the following properties hold:
     \begin{enumerate}
     \item  For each $b \in B$, the action of $G$ restricts to an action on the fibre 
     \[
     F_b \times G \to F_b.
     \]
     \item  $G$ acts freely on $F_b$ ; that is, $xg = x$ for all $x \in X$  implies that  $g $ is the identity element of $G$.
      \item   $G$ acts transitively on $F_b$; that is, if $x, y \in F_b$ then there exists some $g \in G$ with $xg = y$.
      \item For each $b \in B$ and $x \in F_b$, the map  $\phi _{b,x}: G \to F_b$ given by 
      \[
      \phi _{b,x} (g) = xg
      \]
      is a homeomorphism.   
      \end{enumerate} 
      Hence $F_b \cong G$ for each $b \in B$ and the orbit space $E/G$ is homeomorphic to $B$.  So, up to isomorphism, we have
      \[
      G \longrightarrow E \longrightarrow E/G .
      \]
      If $G$ is a Lie group, the action on $E$ is smooth,  and $E \to B$ is a map of smooth manifolds then the construction may be carried out smoothly.

      \end{Def}

   Examples of  principal $G$ bundles that will occur many times in this work are easy to state.  Take $G$ to be a compact   Lie group,  $H$ a closed 
      subgroup, and $G/H$ the homogeneous space. Then $$H \rightarrow G \rightarrow G/H$$ is a principal $H$-bundle. 
      For example, we have the coset fibrations
      \[
      O_n \longrightarrow O_{n+1} \longrightarrow O_{n+1}/O_n \cong S^n
\]
and
\[
O_n \times O_n \longrightarrow O_{2n} \longrightarrow  O_{2n}  /  (O_n \times O_n)  .
\]
  These are both principal bundles.
  
  Here's another example.  Take $G = \ZZ_2$ and let $\ZZ_2$ act on the sphere $S^n$ by the antipodal action. Then there is a natural 
  principal $\ZZ_2$ bundle, namely 
  \[
  \ZZ_2 \longrightarrow S^n \longrightarrow RP^n
\]
with base space real projective space.  There is an obvious map of principal $\ZZ _2$- bundles
\[
\begin{CD}
  \ZZ_2 @>>> S^n @>>>  RP^n \\
@VVV    @VVV  @VVV \\
  \ZZ_2 @>>> S^{n+1} @>>> RP^{n+1}
\end{CD}
\]
  and taking direct limits we obtain
  \[
   \ZZ_2 \longrightarrow S^\infty \longrightarrow RP^\infty
\]
  The space $S^\infty$ is contractible, and so this presents us with an example of a {\emph{universal}} 
   principal $G$-bundle:
      
      \begin{Def} (cf. May \cite{May} 196-198)     A principal $G$-bundle 
      \[
      G \longrightarrow E \longrightarrow B
      \]
      is a {\emph{universal principal $G$-bundle}}  \index{principal $G$-bundle!universal} if the total space $E$ is contractible.

      Any two universal principal $G$-bundles are homotopy equivalent.  The usual notation for the universal principal 
      $G$-bundle is 
      \[
      G \longrightarrow EG \longrightarrow {{BG}}
      \]
      The space $BG$ is called the {\emph{classifying space}} \index{classifying space} of $G$-bundles for reasons that will 
      become apparent.
      \end{Def} 
      
      So, for example, 
       \[
   \ZZ_2 \longrightarrow S^\infty \longrightarrow RP^\infty
\]
  is the universal principal $\ZZ_2$-bundle and $RP^\infty \cong {{B\ZZ_2}}$.
  
  \begin{Rem}\label{classify}  The universal principal $G$-bundles play a critical role in the classification of $G$-bundles.  The basic fact is that the pullback 
  of a principal $G$-bundle is itself a principal $G$-bundle. Then for some fixed space $B$ (e.g. paracompact to simplify matters) and 
  some continuous function $f: B \to BG$, the pullback diagram
  \[
  \begin{CD}
  E @>>> EG \\
  @VVV   @VVV \\
  B @>f>> BG
  \end{CD}
   \]
   produces a principal $G$-bundle $E \to B$.  This gives us a map 
   \[
   [B, BG ]^{ub} \longrightarrow \Prin_G(B)
  \]
  from homotopy classes of maps $B \to BG$ to isomorphism classes of principal $G$-bundles over $B$.  This map is an isomorphism (cf. 
  Husemoller \cite{H} 12.2 ).  So this justifies calling $BG$ a ``classifying space" since homotopy classes of maps into it classify principal $G$-bundles.  
  
  The cases that are of greatest interest to us are when $G = O_n$ or $U_n$, because with some additional work one can show that maps into $BO_n$ 
  classify real vector bundles of real dimension $n$ and similarly maps into $BU_n$ classify complex vector bundles of complex dimension $n$.   We will 
  come to this when we discuss real and complex vector bundles.
  
  \end{Rem}

  There is one other construction that will be very useful to us.

  \begin{Def}
  Suppose that $E \longrightarrow B$ is a principal $G$-bundle and $F$ is a $G$-space with an effective\footnote{This means 
  that if $gf = f$ for all $f$ then $g$ is the identity element of the group.}  left $G$-action. Then $G$ acts on $E \times F$ 
  by 
 \[
  g(e,f) = (eg^{-1}, gf).
  \]
    Define the {\emph{associated bundle}} \index{associated bundle} by
  \[
  E \times _G F =  (E \times F)/G \longrightarrow B.
  \]
  The associated bundle is a fibre bundle 
  \[
  F \longrightarrow E \times _G F \longrightarrow B
  \]
 with  group $G$.   Every fibre bundle with group $G$ arises up to isomorphism in this manner. 
  
  \end{Def}
  
Here are a  very important couple of examples.  If we take $G = O_n$ and $F = \RR ^n$  with the usual action of $O_n$ on $\RR ^n $ then for any $G$-space $Y$, 
    \[
  Y \times _G \RR ^n  \longrightarrow B \; 
  \]
will be a real vector bundle of dimension $n$.   Similarly, if $G = U_n$ acts on $\CC ^n$ then we obtain the associated complex vector bundle 
   \[
  Y \times _G \CC ^n  \longrightarrow B \; .
  \]

 We need to carefully construct the  topological space $\bbo$. 
 \index[notation]{BO@$\bbo$}
  \index[notation]{BU@$\bbu$}
   \index[notation]{Bsp@$\bbsp$}
Hatcher \cite{Hatcher} and Husemoller \cite{H} are excellent references for this section and May \cite{May}  gives a very useful presentation.  This is a very brief summary.     The construction of the spaces $\bbu $ and $\bbsp $ in the complex and symplectic settings
 is done in the same way and we omit 
discussion of them.      Along the way we will discover some very important fibre bundles. 

Form the Grassmann manifold $G_k(\topr ^n)$, which consists of the space of all real $k$-dimensional vector subspaces through the origin in $\RR ^n$. It is a based compact Hausdorff space
and a finite CW-complex. There are natural inclusion maps 
\[
G_k(\RR ^n)       \longrightarrow G_k(\RR ^{n+ 1})
   \]
   and we define $BO_k$     to be their union. 
   The name is not accidental:
   
   \begin{Pro} (May \cite{May} pp. 184-5)     
   \begin{enumerate}
   \item
   There is a contractible $O_n$-space $EO_n$ and a universal principal $O_n$-bundle
   \[
   O_n \longrightarrow EO_n \longrightarrow BO_n
   \]
   and similarly for $BU_n$ and $ BSp_n$.\footnote{For example, we have shown that there is a fibre bundle 
   \[
   \ZZ_2 \longrightarrow S^\infty \longrightarrow RP^\infty
   \]
   with $S^\infty $ contractible, and thus $B\ZZ_2 \cong RP^\infty$.  There is an analogous construction that shows that $BS^1 \cong CP^\infty$. Thus
   since $S^1 = U_1$, we have $BU_1 \cong CP^\infty$. } 
   \item 
    There are natural inclusion maps
 \[
G_k(\RR ^n)       \longrightarrow G_{k+1}(\RR ^{n+ 1})
   \]  
   which induce inclusion maps 
   \[
   BO_k \longrightarrow BO_{k+1} .
   \]
Define $\bbo $ to be their union.  It is  a CW-complex and there is a universal  principal bundle with $G = \oo = \lim_{n \to \infty}  O_n$ 
   of the form
   \[
   \oo \longrightarrow E\oo \longrightarrow \bbo .
   \]
   \item
 
   The analogous constructions for the unitary   groups yields a universal principal $\uu $- bundle 
     \[
   \uu \longrightarrow E\uu \longrightarrow \bbu .
 \]
 
 \item 
  The analogous constructions for the symplectic   groups yields a universal principal $\pp$-bundle 
     \[
   \pp \longrightarrow E\pp \longrightarrow \bbsp .
 \]

\end{enumerate}   
   \end{Pro}

\begin{Cor} The natural maps 
   \[
\lim_{n \to \infty}  \pi _*( G_k(\RR ^n)  ) \longrightarrow \pi _*(BO_k )
   \]
   and 
   \[
\lim_{k \to \infty}  \pi _*( BO_k) \longrightarrow \pi _*(\bbo )
   \]
   are isomorphisms. Similar statements hold for $\bbu$ and $\bbsp$.
   \end{Cor}

   This implies by Proposition~\ref{loops} that $\Omega BO_n \simeq O_n$  and hence 
    \[
    \pi _j(BO_n) \cong \pi _{j-1}(O_n) \qquad {\text{for each}}\,\, j. 
    \]
  Proposition~\ref{loops}  tells us that 
       \[
   \Omega\bbo \simeq \oo \qquad and \qquad \pi _n(\bbo ) \cong     \pi _{n-1} (\oo ).
   \]
     
     It turns out that $\bbo$ has natural addition and multiplication maps which turn any $[X, \bbo]^{ub}$ into a commutative semiring.  Here is
      a sketch of the construction.

The addition operation ultimately arises from the maps 
\[
O_j \times O_k \longrightarrow O_{j+k} 
\]
given by 
\[
(A, B) \longrightarrow  \left[ \begin{matrix} 
A & 0 \\
0 & B
\end{matrix}\right]
\]
which leads to a map $\oo \times \oo \to \oo$ and then to 
\[
\bbo \times \bbo \overset{\alpha}\longrightarrow \bbo. 
 \]
  
 The multiplication map is more complicated. It arises from the Kronecker \index{Kronecker product} 
 product of matrices $M \times N \mapsto M\otimes N$
which gives rise to the product
 \[
 O_j \times O_k \longrightarrow O_{jk} 
 \; . \]
 and leads to a map 
 \[
 \bbo \times \bbo \overset{\mu}\longrightarrow \bbo .
 \]


   \newpage
    \section{\bf{Bott Periodicity Revisited: 
     First View of the Homotopy Road Map}}
     \label{Section:BottRevisited}

    \index{homotopy group|(}
    
    One of the principal goals of this monograph is to plot a path that takes us from classical Bott Periodicity to physics, stopping along the way at Clifford algebras, \index{Clifford algebra}
    symmetric spaces and Cartan classification, and the TCS physics symmetries.  Our guide will be the {\bf{Homotopy Road Map}} \index{Homotopy Road Map} 
    which we will construct along the way.  Step One takes us from Bott Periodicity to the Big Ten spaces.

    We have glossed over the proofs of real and complex Bott Periodicity. There are many proofs and we have referenced several.  In the previous 
    sections we have seen that real Bott Periodicity can be stated in terms of 
    \[
    KO_{j+8}(A) \cong KO_j(A)
    \]
    for all real $C^*$-algebras, or equivalently as  
    \[
    \Omega ^8 \oo \simeq \oo .   \index{loop space}
    \]
    Similarly, complex Bott Periodicity can be stated as 
       \[
    K_{j+2}(A) \cong KO_j(A)
    \]
    for all complex $C^*$-algebras, or equivalently as  
    \[
    \Omega ^2 \uu  \simeq \uu .
    \]
    We need to take a much more careful look at the individual spaces $\Omega ^k \oo$ and $\Omega ^k \uu $.  We shall simply state the result that 
    we need, and we refer the reader to Theorem~\ref{Wood}  and to Part IV of Milnor  \cite{M3} for details.

   \index{classical group|(}
 
 \begin{Def}
 \begin{enumerate}
 \item
  Define
 \[
 \cR_0 = \bbo \times \ZZ \qquad\qquad     \cR_1 = \oo 
 \]
 and for $j = 2,3, \dots 7 $  define  $ \cR_j = \Omega ^{j-1}\oo $.
 \item
 Similarly, define
  \[
 \cR_0\pc = \bbu \times \ZZ \qquad\qquad     \cR_1\pc  = \uu 
 \]
 \end{enumerate}
 \end{Def} 
 
 \index[notation]{rj@$\cR_j$}
 
 Here is the Homotopy Road Map, followed by an explanation of its entries.   
 
  \begin{Table} \label{road1}
 \centerline{\bf{The Homotopy Road Map: Step I }} \index{Homotopy Road Map}  
 
 \begin{center}
     \begin{tabular}{   r  |r |   r|  r | r }  \hline
finite approx   &   Big Ten    &    Bott &      $\cR_k$ &  $\pi _0(\cR_k)$ \\  \hline \hline
  $U_n $ &  $\uu $ &    $\Omega ^0 \uu$   &   $\cR _1\pc $      &       0  \\  \hline
  
    $U_{2n}/(U_n \times U_n)  $ &  $\bbu \times \ZZ$ &     $\Omega ^1\uu $  &     $\cR _0\pc  $      &       $ \ZZ $    \\  \hline \hline
  $O_n $ &  $\oo $  &  $\Omega ^0 \oo$  &     $\cR _1$      &       $ \ZZ_2 $    \\  \hline
    $O_{2n}/U_n $ &  $\oo / \uu $  &  $\Omega ^1\oo$  &     $\cR _2$      &       $ \ZZ_2 $    \\  \hline
    $U_n/Sp_n $ &  $\uu / \pp $  &  $\Omega ^2 \oo$  &     $\cR _3$      &       $ 0 $    \\  \hline
    $Sp_{2n}/(Sp_n \times Sp_n) $ &  $\bbsp \times \ZZ $  &   $\Omega ^3 \oo$  &     $\cR _4$      &       $ \ZZ  $    \\  \hline
    $Sp_n $ &  $\pp $  &   $\Omega ^4\oo $  &     $\cR _5$      &       $ 0 $    \\  \hline
    $Sp_{n}/U_n $ &  $\pp/\uu $  &   $\Omega ^5\oo$  &     $\cR _6$      &       $ 0 $    \\  \hline
    $U_n /O_n $ &  $\uu /\oo $  &    $\Omega ^6\oo$  &     $\cR _7$      &       $ 0$    \\  \hline
    $O_{2n}/(O_n \times O_n) $ &  $\bbo  \times \ZZ $  &   $\Omega ^7\oo$  &     $\cR _0$      &       $ \ZZ  $    \\  \hline
\end{tabular}
   \end{center}
 \end{Table}
    \vglue .5in

   \begin{enumerate}
   \item The first column indicates the family of Lie groups \index{Lie group} or homogeneous spaces \index{homogeneous space} associated with the row. These are all compact spaces.
  
   \item The second column is the direct limit of the first column. It indicates description of the limiting Big Ten\footnote{The second author grew up
    following Big Ten football and couldn't resist the 
   opportunity to call the ten spaces by their true name.}        homogeneous space in terms of the iterated loop space of the orthogonal group (or the unitary 
   group, for the first two rows).   Note that in the three cases $\bbu \times \ZZ$,  $\bbsp \times \ZZ$, and $\bbo \times \ZZ$,
   the space mentioned is the limit of the unstable space times the integers. 
   So what is true, for example in the bottom row, is that 
   \begin{enumerate}
   \item 
   \[
   O_{2n}/(O_n \times O_n)  \longrightarrow \oo /(\oo \times \oo ) \simeq \bbo   
   \]
   and
   \item 
   \[
   \bbo  \times \ZZ = \cR_0  \simeq \Omega ^7 \oo
   \]
   
   \end{enumerate}
   
  \item The third column gives the identification of the Big Ten spaces in terms of iterated loop spaces of $\uu $ and of $\oo $ and illustrates their essential 
  tie with Bott Periodicity.  This will be explained in future sections -- for the moment we ask that it be taken on faith.
   \item The fourth column indicates the name of the space in the physics literature. We will use this terminology with the hope that it will make the connection with 
   physics clearer.
     \item The fifth column indicates  $\pi _0(\cR_j)$.  Remember that saying that  $\pi _0 = 0$ means that the space is path-connected. Saying that  
   $\pi _0 =  \ZZ _2$ says that the space has two connected components. If $\pi _0 = \ZZ$ then there are an infinite number of path components.

 \end{enumerate}

  Guideposts along the way:
   \begin{enumerate}
   \item Bott Periodicity:
\[
\Omega ^2 \uu \simeq \uu \qquad\qquad    \Omega ^8 \oo  \simeq \oo  
\]
\item Notation:
 \[
 \Omega \cR _j  \simeq  \cR_{j+1}  \qquad \pmod 8 \qquad\qquad    \Omega \cR _j\pc  \simeq  \cR_{j+1}\pc   \qquad \pmod 2
 \]
 \item Restatement of Bott Periodicity:
   \[
 { \pi_j(\cR _k) =   \pi _0(\cR _{j + k}) }        \cong \pi _{j+k}(\cR _0)  \qquad \pmod 8 \,\, \footnote{since
 $
 { \pi_j(\cR _k) = [S^j, \cR _k ] = [S^0, \Omega ^j \cR_k ] = \pi _0(\cR _{j + k}) }.$
 } 
 \] 
   and
    \[
 { \pi_j(\cR _k\pc  ) =   \pi _0(\cR _{j + k}\pc  ) }        \cong \pi _{j+k}(\cR _0\pc  )  \qquad  \pmod 2 \,\, \footnote{since
 $
 { \pi_j(\cR _k\pc ) = [S^j, \cR _k \pc  ] = [S^0, \Omega ^j \cR_k\pc  ] = \pi _0(\cR _{j + k}\pc ) }.
  $
 } 
 \] 
 \item There are natural maps
    \[
  \cR_i \wedge \cR_j \longrightarrow \cR_{i+j} 
  \]
  in the real case, (with indexing $\pmod 8$), and 
  \[
  \cR_i\pc  \wedge \cR_j\pc   \longrightarrow \cR_{i+j} \pc  
  \]
  in the complex case (with indexing $\pmod 2$).     These maps correspond to tensor product of vector bundles.

  \end{enumerate}

Take for example the fourth row.  The group $U_n$ may be included naturally in the group $O_{2n}$ with homogeneous space \index{homogeneous space} 
$O_{2n}/U_n$. 
As $n \to \infty $ we obtain the homogeneous space $\oo / \uu$.  
It turns out that $\oo / \uu $ is homotopy equivalent to $\Omega ^1 \oo = \cR_2$
and so 
\[
\pi _0(\oo/\uu ) = \pi _0 (\cR_2) \cong  \ZZ _2.
\]
  \vglue .3in

 The identification of the spaces $\cR_j$ with the homogeneous spaces \index{homogeneous space} mentioned follows from some of the proofs of Bott Periodicity. We shall explore this in much more detail when we reach Clifford algebras \index{Clifford algebra} (see Theorem \ref{Wood} in particular).  
So for the moment consider this theorem as a taste of the Homotopy Road Map, \index{Homotopy Road Map}  which will stretch from homotopy land and Bott Periodicity to the Ten-Fold classification of symmetric spaces and to the Ten-Fold Way in physics.

\begin{Rem}
 Each of the spaces $\cR_j $ is a CW-complex \index{CW-complex} and that means, for example, that the map
\[
\lim_{n \to \infty}  \pi _j (O_n) \longrightarrow \pi_j (\oo ) \cong \pi _j(\cR _1)  
\]
is an isomorphism.  However, that does NOT mean that $\pi _j (O_n) \to \pi_j (\oo ) $ is an isomorphism for all $n$.   Here's an example:
we calculate $\pi _5(O_k) $ for all $k$. 
     \[
    \begin{CD}
    \pi_5(O_2) @>>> \pi_5(O_3) @>>>  \pi_5(O_4) @>>> \pi_5(O_5) @>>>  \pi_5(O_6) @>>> \pi_5(O_k)  \\
   @VV{\cong}V     @VV{\cong}V     @VV{\cong}V     @VV{\cong}V     @VV{\cong}V    @V{k \geq 7}V{\cong}V  \\
   0 @>>>   \ZZ_2 @>>>  \ZZ_2\oplus \ZZ_2 @>>>  \ZZ_2 @>>> \ZZ  @>>>  0      
    \end{CD}
\]
 Thus
 \[
 \pi_5(O_k) = \pi _0 (\Omega ^5 O_k) = 0 \qquad k \geq 7
 \]
 and hence the eighth row lists $\pi_0(\Omega ^5 (\oo ) ) = 0$
 even though 
 \[
 \pi _5(O_3) \neq 0 \qquad   \pi _5(O_4) \neq 0 \qquad   \pi _5(O_5) \neq 0 \qquad   \pi _5(O_6) \neq 0  .
 \]
 In Section~\ref{Section:LowHomotopyGroups} we will see many examples of this behavior.
 \end{Rem}
 
  \index{classical group|)}
  \index{homotopy group|)}


 \newpage
      \section{\bf{Next Contact with Physics: Low Homotopy Groups of the Ten-Fold Way}}
      \label{Section:LowHomotopyGroups}
 
  \index{homotopy group|(}
   \index{classical group|(}
In this section we will compute $\pi_0$, $\pi_1$      and  $\pi_2$ of the classical Lie groups \index{Lie group}  and homogeneous 
spaces that arise 
  in the Ten-Fold Way.   
  This section contains material that is  very well-known in the math world but not in the physics world (yet).\footnote{Lest 
  we give the wrong impression, a great deal more is known about the higher homotopy groups of compact Lie groups 
  and their associated homogeneous spaces \index{homogeneous space} . We are simply looking at the lowest dimensional calculations.} 
  It is useful there because some physicists (cf. Orion and Akkermans \cite{OA} )     need not just the ``stable" results provided by the Bott Periodicity theorem 
but also the low-dimensional results to which this section is devoted. It is fascinating to see how soon the Bott results take shape here.   

Suppose that $G$ is a compact Lie group \index{Lie group} and $H$ is a closed subgroup. Then
 
 \[
H \longrightarrow G \longrightarrow G/H
\] 
  is a fibre bundle and there is an associated long exact sequence in homotopy that in low dimensions takes 
  the form

\begin{multline} \label{A}
0 = \pi_2(G) \to  \pi _2(G/H) \to \pi _1(H) \to \pi _1(G) \to \pi _1(G/H) \to  \\
\to \pi _0(H) \to \pi _0(G) \to \pi _0(G/H) \to 0 
\end{multline}

\noindent The group  $\pi _2(G)  = 0 $  for all compact Lie groups \index{Lie group}   as previously mentioned.
If the group $H$ is path-connected
so that $\pi_0(H) = 0$,
then the sequence breaks apart into two simpler sequences-  

\begin{equation} \label{B}
\begin{split} 
0 \to  \pi _2(G/H) &  \to \pi _1(H) \to \pi _1(G) \to  \pi _1(G/H) \to 0 \; , \\
&\text{and}\\
 & 0 \to \pi _0(G) \to \pi _0(G/H) \to 0 \; .
\end{split} 
\end{equation}
In particular, there is an isomorphism $\pi _0(G) \cong \pi _0(G/H)$. 

 These sequences are natural.  That is, given a map (i.e. a continuous homomorphism)  of compact groups $f: G \to G'$, closed subgroups $H \subset G$ and $H' \subset G'$ and 
 a commuting diagram
 \[
 \begin{CD}
 H @>f>>    G \\
 @VVV   @VVV  \\
 H' @>f>>    G'
 \end{CD}
 \]
then there is a commuting diagram of long exact sequences.

 It's easy to show that if $X$ and $Y$ are topological spaces then 
 \[
 \pi _j (X \times Y) = \pi _j(X) \oplus \pi _j(Y) \qquad \forall j.
 \]
 We will use this fact many times in what comes below.  Another fact to keep in mind is that homotopy groups 
 revolve about basepoint-preserving maps.  This implies that for $j > 0$ that $\pi _j(G)$ depends only on the 
 path-component of $G$ that contains the basepoint, which we always take to be the identity when dealing with 
 topological groups.  In particular,  the inclusion $SO_n \to O_n$ induces isomorphisms
 \[
 \pi _j(SO_n) \overset{\cong}\longrightarrow \pi _j(O_n)  \qquad\qquad     j > 0.
 \]
 
  Here is an outline of this section. Note that there are ten groups of spaces and these correspond 
  to the   spaces $\cR_j$ discussed in Section~\ref{Section:BottRevisited}
 
 \index{classical group} \index{homogeneous space|(} 
   
{\bf{

\begin{enumerate}
\item    The Groups $U_n$ and $SU_n$
\item  The Homogenous Spaces $U_{n}/(U_k \times U_{n-k}) $ and 
 $SU_{n}/(SU_k \times SU_{n-k}) $
\item   The Groups $O_n$ and $SO_n$
\item   The Homogeneous Spaces  $O_{2n}/U_n$,   $SO_{2n}/SU_n$,  and    $SO_{2n}/U_n$
\item     The Homogeneous Spaces $U_n/Sp_n $ and $SU_n/Sp_n$
\item    The Homogeneous Spaces $ Sp_{n}/(Sp_k \times Sp_{n-k}) $
\item     The Groups $ Sp_n $
\item     The Homogeneous Spaces $ Sp_n/U_n $ and $ Sp_n/SU_n $
\item     The Homogeneous Spaces $U_n/O_n $, $SU_n/SO_n$ and $U_n/SO_n$
\item     The Homogeneous Spaces $O_{n}/(O_k \times O_{n-k}) $ and 
    $SO_{n}/(SO_k \times SO_{n-k}) $
 \end{enumerate}
  }}

\vspace{.1in} \paragraph{{\bf 1. \underline{The Groups $U_n$ and $SU_n$}}}
 \indent

   The group $U_n$ is path-connected for all $n$.  Note that $U_1 \cong S^1$ and 
maps from $S^1$ into $U_n$ factor up to homotopy through the $2$-skeleton, which is $S^1$, by cellular approximation.    (For example, $U_2$ 
looks like $U_1$ with an additional $S^3$ attached).  Hence  
 \[
  \pi _1(U_n) = \pi _1(U_1) = \ZZ \qquad \forall n \geq 1. 
  \]
   The group $SU_n$ is simply connected,\footnote{A space $X$ is {\emph{simply connected}} if $\pi _1(X) = 0$.   In the case of $SU_n$, even more is true. We have $SU_1 = 0$, $SU_2 \cong S^3 $ and, for larger $n$, $SU_n$ is obtained from $SU_{n-1}$ by adjoining even higher-dimensional cells.   }  and hence for all $n$ we have

 \begin{center}
     \begin{tabular}{    r|   c| c| c  }  
 $X$ & $\pi _0(X)$    &  $\pi _1(X)$   &   $\pi _2(X)$       \\  \hline\hline
  $U_n $ & $0$    & $\ZZ $ & $0$ \\ \hline
    $SU_n $ & $0$    & $ 0$ & $0$ \\ \hline
  
\end{tabular}
   \end{center}

   These results connect up with the stable group $\cR_1\pc = \uu$ since $ \lim_{n \to \infty}  U_n = \uu.$

  \vspace{.1in}
   \paragraph{{\bf 2. \underline{The Homogenous Spaces $U_{n}/(U_k \times U_{n-k}) $ and 
 $SU_{n}/(SU_k \times SU_{n-k}) $}}}
  \indent

      To eliminate a trivial case and get rid of duplications, let us assume throughout that $0 < k \leq n - k  $. 
         
       First we consider the path-connected homogeneous spaces   $U_{n}/(U_k \times U_{n-k}) $.  
      The map 
 \[
 \phi = \diag(A, B)  : U_k \times U_{n-k} \to U_{n} 
 \]
 takes a pair $(A,B)$ of   matrices and places them in the upper left and lower right parts of an $n \times n $ matrix (with zeros in the remaining two quarters). 
      Sequence~(\ref{B}) applies, and so we have 
 \[
0 \to  \pi _2(U_{n}/(U_k \times U_{n-k}) ) \to \pi _1(U_k \times U_{n-k})  \overset{\phi _*}\longrightarrow \pi _1(U_{n}) \to  \pi _1(U_{n}/(U_k \times U_{n-k}) ) \to 0  
\]
Plugging in what we know gives 
 
\[
0 \to  \pi _2(U_{n}/(U_k \times U_{n-k}) ) \to \ZZ ^2  \overset{\phi _*}\longrightarrow  \ZZ \to  \pi _1(U_{n}/(U_k \times U_{n-k}) ) \to 0  \]

   Let $f: S^1 \to U_k$ represent the generator 
 of   $\pi _1(U_k)$ and let 
 \[
 \delta : U_k \to U_k \times U_{n-k} 
 \]
  be given by $\delta (x) = (x, 1)$.     Then the map $\phi _* \delta _*( f)$ generates 
 $\pi _1(U_{n}) \cong \ZZ$ and  hence $\phi _*$ is surjective.
 This implies that  
 \[
  \pi _1(U_{n}/(U_k \times U_{n-k}) ) = 0        \qquad and  \qquad
  \pi _2(U_{n}/(U_k \times U_{n-k}) ) = \ZZ .
  \]
   
   Determining   $\pi _j(SU_{n}/(SU_k \times SU_{n-k}))  $ is easy: it is zero for $j  \leq 2$ because 
   $\pi _j(SU_k) = 0 $ for $j \leq 2$.

 We pause to note that what we are seeing here is a small manifestation of Bott Periodicity.  The space 
 $U_n$ stabilizes to the infinite unitary group $\uu $ and the homogeneous spaces  $U_{2n}/(U_n\times U_{n})$ stabilize to 
 the classifying space $\bbu $ and $\bbu \times \ZZ    \cong \cR_0\pc$ as $n \to \infty$.  Periodicity implies that $\pi _1(\uu ) \cong \pi _0(\bbu )   $
 and that is what we are seeing already at the earliest stage:
 
     \begin{center}
     \begin{tabular}{    r|   c| c| c| c  } 
 $X$ & $\pi _0(X)$    &  $\pi _1(X)$   &   $\pi _2(X)$       \\  \hline\hline
  $U_n  $ & $0$    & $\ZZ $ & $0$ \\ \hline
   $U_{n}/(U_k \times U_{n-k}) \times \ZZ $ & $\ZZ $    & $0 $ & $\ZZ $    &  $0 <  k < n $\\ \hline    
  $SU_n  $ & $0$    & $0$  & $0$ \\ \hline
$SU_{n}/(SU_k \times SU_{n-k})  \times \ZZ $ & $\ZZ $    & $0 $ & $0 $ &      $0 <  k < n $        \\ \hline

\end{tabular}
   \end{center}

 Bott Periodicity gives us 
\[
\pi _n(\uu) \cong \pi_{n+2}(\uu) 
 \]
 and putting that together with the easy statement that $\Omega (\bbu ) \simeq \uu $ gets us the full complex periodicity statement.    To see this when looking at   $SU_{n}/(SU_k \times SU_{n-k})$ one has to go to higher homotopy groups.

    \vspace{.1in} \paragraph{{\bf 3. \underline{The Groups $O_n$ and $SO_n$}}}
    \indent

   The groups  $O_n $ have two path components, rotations and reflections, and
    $SO_n$ is the closed, index $2$ subgroup of rotations. Thus 
    \[
     \pi _0(O_n) = \ZZ_2  \qquad\qquad      \pi_0(SO_n) = 0 \qquad \forall n \geq 1,
     \]
     and
     \[
     \pi_j(O_n) = 
\pi _j(SO_n ) \qquad\qquad \forall j \geq 1.
\]     
    Now $SO_2 \cong S^1$ as is easy to see.  Less obvious\footnote{Steenrod  (\cite{St} \S 22)  has 
    a detailed exposition of this (22.3) and related matters.}  is the fact that $SO_3 \cong RP^3$.  This implies that  
      
    \[
     \pi_1(O_1) \cong  \pi _1 (SO_1) = 0          \qquad  \pi _1 (O_2) \cong  \pi _1 (SO_2) = \ZZ    \qquad   
        \pi _1(O_3) \cong  \pi _1(SO_3) = \ZZ_2 
           \]
    and, by cellular approximation,   
      \[
    \pi _1 (O_n) \cong     \pi _1(SO_n) = \ZZ_2 \quad \forall n \geq 3 .
    \]
 
To summarize: 
    
 \begin{center}
     \begin{tabular}{    r|    c|  c|  c |  cl } 
 $X$ & $\pi _0(X)$    &  $\pi _1(X)$   &   $\pi _2(X)$       \\  \hline\hline
   $O_1 $ & $\ZZ_2 $    &$0$  & $0$ \\ \hline
  $O_2 $ & $\ZZ_2 $    & $\ZZ $ & $0$ \\ \hline
      $O_n $ & $\ZZ_2 $    & $ \ZZ_2$ & $0$ & $n  \geq 3$\\ \hline\hline
       $SO_1 $ & $0$    & $0 $ & $0$ \\ \hline
      $SO_2 $ & $0$    & $\ZZ $ & $0$ \\ \hline
     $SO_n $ & $0$    & $\ZZ_2 $ & $0$ & $ n \geq 3$\\ \hline

    \end{tabular}
   \end{center}

These results tie to the group $\cR_1 = \oo$, of course, since $ \lim_{n \to \infty}  O_n \cong \oo$.

As mentioned earlier, much more is known about these groups.  For example, Steenrod 
shows \cite{St} that we may extend the chart as follows.  (Note $\pi_j(O_n) = \pi_j(SO_n)$ for $j > 0$ since $SO_n$ is the path-component of the identity of $O_n$. )
 \begin{center}
     \begin{tabular}{    r|    c| c| c| c  } 
 $X$ & $\pi _3(X)$    &  $\pi _4(X)$   &   $\pi _5(X)$       \\  \hline\hline
          $O_3 $ & $\ZZ $    & $\ZZ_2$ & $\ZZ_2$ \\ \hline
      $O_4 $ & $\ZZ^2$    & $\ZZ_2 \oplus \ZZ_2 $ & $\ZZ_2 \oplus \ZZ_2$ \\ \hline
     $O_5 $ & $\ZZ $    & $\ZZ_2 $ & $\ZZ_2$ & \\ \hline
     $O_k $ & $\ZZ$    & $0 $ & $\ZZ $ & $ k  > 5$\\ \hline
    \end{tabular}
   \end{center}

    
  \vspace{.1in} \paragraph{{\bf 4. \underline{The Homogeneous Spaces  $O_{2n}/U_n$,   $SO_{2n}/SU_n$,  and    $SO_{2n}/U_n$}}}
  \indent


    Let 
    \[\phi _n : U_n \to O_{2n} \]
     denote the standard inclusion.
    The groups $U_n$ are path-connected for all $n$. So we have from Sequence~(\ref{B}),
\begin{equation*} \begin{split}
0 \to  \pi _2(O_{2n}/U_n) & \to \pi _1(U_n) \overset{(\phi_n)_* }\to \pi _1(O_{2n}) \to 
 \pi _1(O_{2n}/U_n) \to 0  \\
\text{and~} &0  \to   \pi_0(O_{2n} ) \to \pi _0( O_{2n}/U_n ) \to 0.
\end{split} \end{equation*}

First, consider the case $n = 1$.  Then the sequence becomes

\[
0 \to  \pi _2(O_{2}/U_1) \to\ZZ \overset{\phi _*}\longrightarrow \ZZ \to  \pi _1(O_{2}/U_1) 
\to 0 \to \ZZ_2 \to  \pi _0(O_{2}/U_1) \to 0
 \]
The map $\phi _*$ is an isomorphism since $U_1 = SO_2 \subset O_2$  and this implies that
\[
 \pi _0 ( O_{2}/U_1) \cong \ZZ_2 \qquad            \pi _1(O_{2}/U_1)    =    \pi _2(O_{2}/U_1) = 0.
\]
Now consider the case $n > 1$.    
Then the sequence becomes
 
\[
0 \to  \pi _2(O_{2n}/U_n) \to\ZZ \overset{\phi _*}\longrightarrow \ZZ_2 \to  \pi _1(O_{2n}/U_n) 
\to 0 \to \ZZ_2 \to  \pi _0(O_{2n}/U_n) \to 0
 \]
with the map $(\phi _n)_*$ surjective \footnote{The diagram 
\[
\CD
U_1 @>\phi _1>>     O_2 \\
@VVV     @VVV \\
U_n    @>\phi _n >> O_{2n}
\endCD
\]
commutes; then apply $\pi _1$. } for all $n$,
 implying that 
\[
 \pi _1(O_{2n}/U_n) = 0       \qquad \text{and} \qquad  \pi _2(O_{2n}/U_n) = \ZZ  
\]
for all $n$.

For $SO_{2n}/SU_n$ we apply Sequence~(\ref{B}) since $SO_{2n}$ is path-connected. We get 
 
\[
0 \to  \pi _2(SO_{2n}/SU_n) \to\pi _1(SU_n)   \longrightarrow  \pi _1(SO_{2n})  \to  \pi _1(SO_{2n}/SU_n) 
\to 0  
 \]
which simplifies   to 
\[
0 \to  \pi _2(SO_{2n}/SU_n) \to    0   \longrightarrow  \pi _1(SO_{2n})  \to  \pi _1(SO_{2n}/SU_n) 
\to 0  
 \]
Thus $\pi _2(SO_{2n}/SU_n) = 0$ and  
\[
 \pi _1(SO_{2n}/SU_n) \cong  \pi _1(SO_{2n}) 
 \]
  which is $\ZZ$ 
for $n = 1$ and $\ZZ _2 $ for $n > 1$. 

Finally, we consider $SO_{2n}/U_n$.   The case $n = 1$ is easy, since $U_1 = SO_2$ and so all of the 
homotopy groups vanish.  So assume that $n > 1$. As $SO_{2n}$ is path-connected we may use 
Sequence~(\ref{B}).   It reads 
\[
0 \to \pi _2 (SO_{2n}/U_n) \longrightarrow  \pi _1(U_n) \overset{\phi _*}\longrightarrow \pi _1(SO_{2n} ) \longrightarrow
 \pi _1 (SO_{2n}/U_n) \to 0
\]
which simplifies to 
\[
0 \to \pi _2 (SO_{2n}/U_n) \longrightarrow \ZZ  \overset{\phi_*}\longrightarrow \ZZ_2  \longrightarrow
 \pi _1 (SO_{2n}/U_n) \to 0
\]
To determine the map $\phi_*$, we can consider the commutative diagram 
\[ \xymatrix{
\pi _1(U_1) \ar[r]^\cong \ar[d]^\cong
& \pi_1(SO_2) \ar[r]^\cong \ar@{->>}[d]
& \ZZ \ar@{->>}[d] \\
 \pi _1(U_n) \ar[r]^{\phi_*} 
& \pi_1(SO_{2n}) \ar[r]^\cong
& \ZZ_2
} \]

which implies that the map $\phi _*$ is surjective.  Thus 
\[
 \pi _1 (SO_{2n}/U_n) = 0      \qquad \text{and} \qquad  \pi _2 (SO_{2n}/U_n) \cong \ZZ .   
\]
In summary,

 \begin{center}
     \begin{tabular}{    r|   c| c| c| c} 
 $X$ & $\pi _0(X)$    &  $\pi _1(X)$   &   $\pi _2(X)$       \\  \hline\hline
        $O_2/U_1 $ & $\ZZ_2 $    & $ 0$ & $0$ & \\ \hline
            $O_{2n}/U_n$ & $\ZZ_2 $    & $ 0$ & $\ZZ $ & $n  \geq 2$\\ \hline  
        $SO_2/SU_1 $ & $\ZZ $    & $ \ZZ $ & $0$ & \\ \hline
            $SO_{2n}/SU_n$ & $0 $    & $ \ZZ _2 $ & $0 $ & $n  \geq 2$\\ \hline  
          $SO_{2n}/U_n$ & $0 $    & $ 0 $ & $ \ZZ $ &  \\ \hline   
              \end{tabular}
   \end{center}

These results connect up with $\cR_2 =  \lim_{n \to \infty}   O_{2n}/U_n $.


   \vspace{.1in} \paragraph{{\bf 5. \underline{The Homogeneous Spaces $U_n/Sp_n $ and $SU_n/Sp_n$ }}} 
\indent

These are path-connected spaces. We have  $\pi _j(Sp_n)     = 0$ for  $j\leq 2$ (see Item 7 below), so Sequence~(\ref{B}) simplifies to
    
\[
                 0 \to \pi _2(U_n/Sp_n)    \to                      0   \to \ZZ \to  \pi _1(U_n/Sp_n ) \to 0 \; . 
                 \]
 
 Hence
 \[
  \pi _0(U_n/Sp_n )    =  0        \qquad                       \pi _1(U_n/Sp_n ) \cong \ZZ \qquad 
     \pi _2(U_n/Sp_n ) \cong 0 \qquad \forall n
 \]
 
 For $SU_n/Sp_n$, the  sequence simplifies to

\[
                 0 \to \pi _2(SU_n/Sp_n)    \to                      0   \to 0 \to  \pi _1(SU_n/Sp_n ) \to 0 
\]
 \vglue .in
and so 
\[
\pi _j(SU_n/Sp_n) = 0 \quad j \leq 2 \qquad \forall \,\, n
 \]

 \begin{center}
     \begin{tabular}{    r|   c | c | c } 
 $X$ & $\pi _0(X)$    &  $\pi _1(X)$   &   $\pi _2(X)$       \\  \hline\hline
        $U_n/Sp_n $ & $0 $    & $ \ZZ $ & $0$  \\ \hline
         $SU_n/Sp_n $ & $0 $    & $ 0 $ & $0$  \\ \hline
    
    \end{tabular}
   \end{center}

These results tie in with $\cR_3 =  \lim_{n \to \infty}  U_n/Sp_n$.

  \vspace{.1in}  \paragraph{{\bf 6. \underline{The Homogeneous Spaces $ Sp_{n}/(Sp_k \times Sp_{n-k}) $ }}}
 \indent

 We will see in the next section that $\pi _j(Sp_n) = 0$  for $j \leq 2$.  This implies that in Sequence~(\ref{A}) 
 two out of every three groups is zero, and hence 
 \[
\pi _j( Sp_{n}/(Sp_k \times Sp_{n-k})) = 0 \qquad j = 0, 1, 2   \qquad \forall \,\, n
 \]

 \begin{center}
     \begin{tabular}{    r|    c | c | c}
     $X$ & $\pi _0(X)$    &  $\pi _1(X)$   &   $\pi _2(X)$       \\  \hline\hline
        $Sp_{n}/(Sp_k \times Sp_{n-k}) $ & $0 $    & $ 0 $ & $0$  \\ \hline

    \end{tabular}
   \end{center}

 These results connect with $\cR_4  = {\bf{BSp}} \times \ZZ $ since $  Sp_{2n}/(Sp_n \times Sp_{n}) $ converges to ${\bf{BSp}}$.

       \vspace{.1in}  \paragraph{{\bf 7. \underline{The Groups $ Sp_n $}}}
  \indent

    This is the easiest section to write. The groups $Sp_n$ have their lowest-dimensional cell in dimension 3, which 
   implies that 
   \[
   \pi _j(Sp_n) = 0  \quad \forall \,\,n  \quad j \leq 2.  
    \]

 \begin{center}
     \begin{tabular}{    r|    c| c | c } 
 $X$ & $\pi _0(X)$    &  $\pi _1(X)$   &   $\pi _2(X)$       \\  \hline\hline
        $Sp_{n}  $ & $0 $    & $ 0 $ & $0$  \\ \hline

    \end{tabular}
   \end{center}

These groups connect up with $\cR_5 =  \lim_{n \to \infty}  Sp_n $.

  \vspace{.1in}  \paragraph{{\bf 8. \underline{The Homogeneous Spaces $ Sp_n/U_n $ and $ Sp_n/SU_n $ }}}
 \indent

 The spaces $ Sp_n/U_n $ are path-connected.  Sequence~(\ref{B}) gives an exact sequence 
 \[
 \to \pi_2(Sp_n) \to \pi _2(Sp_n/U_n) \to \pi_1(U_n) \to \pi _1(Sp_n) \to \pi_1(Sp_n/U_n) \to 0
 \]
 which, after evaluating, becomes the sequence
 \[
 \to 0 \to \pi _2(Sp_n/U_n) \to\ZZ \to 0 \to \pi_1(Sp_n/U_n) \to 0
 \]
 from which we deduce that 
 \[
  \pi_1(Sp_n/U_n)  = 0   \qquad\qquad          \pi_2(Sp_n/U_n) \cong \ZZ.
 \]
  
  For $ Sp_n/SU_n$, Sequence~(\ref{B}) gives an exact sequence
   \[
0 \to \pi _2(Sp_n/SU_n) \to 0 \to 0 \to \pi_1(Sp_n/SU_n) \to 0
 \]
  and so   $\pi _j(Sp_n/SU_n ) = 0$  for all $n$ and for $j \leq 2$.

 \begin{center}
     \begin{tabular}{    r|    c| c | c}
 $X$ & $\pi _0(X)$    &  $\pi _1(X)$   &   $\pi _2(X)$       \\  \hline\hline
        $Sp_{n}/U_n  $ & $0 $    & $ 0 $ & $ \ZZ $   \\ \hline
              $Sp_{n}/SU_n  $ & $0 $    & $ 0 $ & $ 0 $   \\ \hline

          \end{tabular}
   \end{center}
   
   These groups relate to $\cR_6 =   \lim_{n \to \infty}   Sp_n/U_n $.

  \vspace{.1in}  \paragraph{{\bf 9. \underline{The Homogeneous Spaces $U_n/O_n $, $SU_n/SO_n$ and $U_n/SO_n$ }}}
 \indent

First consider $U_n/O_n $.      As $O_n$ is not path-connected we use Sequence~(\ref{A})  and obtain
 
\[
  0 \to \pi_2(U_n/O_n) \to \pi _1(O_n) \overset\alpha\to  \pi _1(U_n) \to \pi _1(U_n/O_n) \to \pi _0(O_n) \to \pi _0(U_n) \to
   \pi _0(U_n/O_n) \to 0 
   \]
The fact that $U_n$ is path-connected implies that  $ \pi _0(U_n/O_n) = 0$ for all $n$
and so the sequence simplifies to
\[
  0 \to \pi_2(U_n/O_n) \to \pi _1(O_n) \overset\alpha\to  \pi _1(U_n) \to \pi _1(U_n/O_n) \to \ZZ_2  \to 0 
    \]

Suppose first that  $n = 1$.  Then $O_1 \to U_1 \to U_1/O_1$ is  really just the double cover  
\[
\ZZ_2 \longrightarrow S^1 \overset{\pi}\longrightarrow U_1/O_1 \cong S^1
\]
    and so 
\[
\pi_0(U_1/O_1) = 0 \qquad         \pi_1(U_1/O_1) = \ZZ  \qquad   \pi_2(U_1/O_1) = 0.  
\]

   Suppose next that $n = 2$.  The commuting diagram
  \[
 \CD
 \pi _1(O_2) @> \cong >>   \pi_1(U_1) \\
  @VV\cong V      @VV\cong V \\
 \pi _1(O_2) @>\alpha  >>   \pi_1(U_2) 
 \endCD
  \]
 shows that the map $ \alpha :    \pi _1(O_2) \to   \pi_1(U_2) $ is an isomorphism. Thus Sequence~(\ref{A})  reduces down to 
 
 \[
  0 \to  \pi_2(U_n/O_n) \longrightarrow\ZZ \overset\cong\longrightarrow  \ZZ  \longrightarrow \pi _1(U_n/O_n) \longrightarrow \ZZ_2 \to 0 
      \]
which implies that
\[
\pi _1(U_n/O_n) \cong  \ZZ_2 \qquad \text{and} \qquad \pi_2(U_n/O_n) = 0.
 \]

Suppose next that $n > 2$.  Then we have 
 \[
    0 \to \pi_2(U_n/O_n)  \to \ZZ_2   \overset\alpha\to \ZZ \overset{\beta}\longrightarrow \pi _1(U_n/O_n) \to \ZZ_2 \to 0 
    \]
  The map $\alpha = 0$ and so it follows that
$   \pi_2(U_n/O_n) = \ZZ_2 $
and that $\pi _1(U_n/O_n)$ is either isomorphic to $\ZZ$ or to $\ZZ \oplus \ZZ_2$. We will settle 
this point in a moment.

We turn next to $SU_n/SO_n$.   The space  $SU_n$ is path-connected and so  
\[
\pi _0(SU_n/SO_n) = 0.
\]
Turning to Sequence~(\ref{B}) we obtain 
\[
  0 \to \pi_2(SU_n/SO_n) \to \pi _1(SO_n) \overset\alpha\to  \pi _1(SU_n) \to \pi _1(SU_n/SO_n)  
  \to 0 
   \]
which becomes
\[
  0 \to \pi_2(SU_n/SO_n) \to \pi _1(SO_n) \overset\alpha\to 0 \to \pi _1(SU_n/SO_n)  
  \to 0 
   \]
Hence $ \pi _1(SU_n/SO_n) = 0$ and $\pi_2(SU_n/SO_n) \cong \pi _1(SO_n)$. This group is $\ZZ$ when $n = 2$ and 
$\ZZ_2$ for $n > 2$. 

We return to the identification of  $\pi_1( U_n/O_n)$ above.   There is a commutative diagram 
\[
\CD
SO_n @>>> SU_n @>>> SU_n/SO_n   \\
@VVV @VVV     @VVV    \\
O_n @>>> U_n  @>\beta >>   U_n/O_n  \\
@VVV      @VVV       @VVV     \\
O_1 = \{\pm 1\}    @>>> U_1 = S^1     @>>>   U_1/O_1 = S^1
\endCD
\]

\noindent with each vertical column and horizontal row a fibration.  Taking homotopy groups we obtain 
the following diagram with exact rows and columns

\[
\CD
@VVV     @VVV       @VVV    \\
\pi_1(SO_n) @>>> \pi_1(SU_n) = 0 @>>> \pi_1(SU_n/SO_n) = 0  @>>> 0  @>>> 0\\
@VVV @VVV     @VVV    @VVV   \\
\pi _1(O_n) @>>> \pi_1(U_n) = \ZZ  @>\beta >>  \pi_1( U_n/O_n) @>>>  \ZZ_2 @>>> 0 \\
@VVV      @V\cong VV       @V\zeta VV     @V\cong VV  \\
\pi_1(O_1) = 0     @>>> \pi _1(U_1) = \ZZ      @>{\bar\beta }>>   \pi_1(U_1/O_1) = \ZZ  @>>> \ZZ_2 @>>> 0
\endCD
\]

The map $\bar\beta $ is multiplication by $2$ which implies that the map $\zeta$ must be an inclusion 
of  $\pi_1( U_n/O_n)$ into $\ZZ$ and so $ \pi_1( U_n/O_n) \cong \ZZ$.   This settles the extension question mentioned above.

Finally, consider $U_n/SO_n$.  The space $U_n$ is path-connected.  If $n = 2$ then  
Sequence~(\ref{B}) gives us 
\[
0 \to \pi _2(U_2/SO_2) \longrightarrow \pi_1(SO_2) \overset{\cong}\longrightarrow \pi _1(U_2) 
\longrightarrow   \pi _1(U_2/SO_2) \to 0
\]
which simplifies\footnote{Exercise \# 1 for the reader: If 
\[
0 \to H \to G_1 \overset{\cong}\longrightarrow G_2 
\]
is an exact sequence of abelian groups, then $H = 0$.      Exercise \# 2: If 
\[
  G_1 \overset{\cong}\longrightarrow G_2  \to H \to 0
\]
is an exact sequence of abelian groups, then $H = 0$.  } to 
\[
0 \to \pi _2(U_2/SO_2) \longrightarrow \ZZ \overset{\cong}\longrightarrow \ZZ
\longrightarrow   \pi _1(U_2/SO_2) \to 0
\]
and hence $\pi _j(U_2/SO_2) = 0$ for $j \leq 2$. 

If $n > 2$ then the sequence
 \[
0 \to \pi _2(U_n/SO_n) \longrightarrow \pi_1(SO_n) \longrightarrow \pi _1(U_n) 
\longrightarrow   \pi _1(U_n/SO_n) \to 0
\]
simplifies to 
 \[
0 \to \pi _2(U_n/SO_n) \longrightarrow \ZZ_2  \longrightarrow \ZZ
\longrightarrow   \pi _1(U_n/SO_n) \to 0.
\]
The map $Z_2 \to Z$ is the trivial map, of course, and so 
\[
 \pi _1(U_n/SO_n) \cong \ZZ   \qquad \text{and}  \qquad         \pi _2(U_n/SO_n) \cong \ZZ_2
\]

To summarize:
    
 \begin{center}
     \begin{tabular}{    r|   c| c | c | c}
 $X$ & $\pi _0(X)$    &  $\pi _1(X)$   &   $\pi _2(X)$       \\  \hline\hline
 $U_1/O_1 $ & $0  $    & $ \ZZ_2$ & $0$ & \\ \hline
    $U_2/O_2 $ & $0 $    & $ \ZZ_2$ & $0$ & \\ \hline
            $U_n/O_n$ & $0 $    & $ \ZZ $ & $\ZZ_2 $ & $n  \geq 3$\\ \hline 
     $SU_2/SO_2 $ & $0 $    & $ 0$ & $\ZZ $ & \\ \hline
            $SU_n/SO_n$ & $0 $    & $ 0$ & $\ZZ_2 $ & $n  \geq 3$\\ \hline
       $U_2/SO_2 $ & $0 $    & $ 0$ & $0 $ & \\ \hline
            $U_n/SO_n$ & $0 $    & $\ZZ $ & $\ZZ_2 $ & $n  \geq 3$\\ \hline
    \end{tabular}
   \end{center}

These groups relate to $\cR_7 =  \lim_{n \to \infty}  U_n/O_n $.

  \vspace{.1in} \paragraph{{\bf 10. \underline{The Homogeneous Spaces $O_{n}/(O_k \times O_{n-k}) $ and 
    $SO_{n}/(SO_k \times SO_{n-k}) $ }}}
 \indent
 
 We will need the following information frequently this section, so we copy it from Item 3 above.
  
 
 \begin{center}
     \begin{tabular}{    r|    c |   c |  c |  c l }
 $X$ & $\pi _0(X)$    &  $\pi _1(X)$   &   $\pi _2(X)$       \\  \hline\hline
   $O_1 $ & $\ZZ_2 $    &$0$  & $0$ \\ \hline
  $O_2 $ & $\ZZ_2 $    & $\ZZ $ & $0$ \\ \hline
      $O_n $ & $\ZZ_2 $    & $ \ZZ_2$ & $0$ & $n  \geq 3$\\ \hline 
       $SO_1 $ & $0$    & $0 $ & $0$ \\ \hline
      $SO_2 $ & $0$    & $\ZZ $ & $0$ \\ \hline
     $SO_n $ & $0$    & $\ZZ_2 $ & $0$ & $ n \geq 3$\\ \hline
\end{tabular}
\end{center}

In what follows we will always assume that $0 < k \leq  n-k $.

   Consider 
 path-components first. We have 
 \[
 \to \pi_0(O_k \times O_{n-k})\overset{\gamma}\longrightarrow    \pi_0(O_{n})     \to    
    \pi_0(O_{n}/(O_k \times O_{n-k})) \to 0
 \]
 which simplifies  to
 \[
  \to \ZZ_2 \oplus \ZZ_2 \overset{\gamma}\longrightarrow \ZZ_2 \to \pi_0(O_{n}/(O_k \times O_{n-k})) \to 0.
   \]
 The map $\gamma $ is surjective and hence  $ \pi_0(O_{n}/(O_k \times O_{n-k})) = 0 $. 
 A similar and easier argument works for $SO$ and hence
\[
\pi _0(O_{n}/(O_k \times O_{n-k}) ) = 0 \qquad \text{and}  \qquad   \pi _0(SO_{n}/(SO_k \times SO_{n-k}) ) = 0
  \]
  for $0 < k \leq n-k < n$.
  
  Next we calculate  $\pi _1 $ and $\pi _2 $ for seven special cases.


\vspace{.1in} \paragraph{  
\bf 1. \underline{The spaces $O_2/(O_1 \times O_1)$ and $SO_2/(SO_1 \times SO_1)$}}
   \indent
   
  There is a natural identification 
  \[
  O_2/O_1  \cong SO_2/SO_1 = SO_2 = S^1
  \]
  with associated fibration
  \[
  O_1 \longrightarrow O_2/O_1 \longrightarrow O_2/(O_1 \times O_1)
  \]
  which we may identify with the double cover
  \[
  \ZZ_2 \longrightarrow S^1 \longrightarrow O_2/(O_1 \times O_1) 
  \]
  Thus  $O_2/(O_1 \times O_1) \cong S^1 $ and hence 
  \[
\pi_1(  O_2/(O_1 \times O_1)) \cong \ZZ \qquad \text{and}  \qquad            \pi_2 (  O_2/(O_1 \times O_1) = 0.
  \]
The fibration
  \[
  SO_1 \longrightarrow SO_2/SO_1 \longrightarrow SO_2/(SO_1 \times SO_1)
  \]
  trivializes since $SO_1$ is one point, so $SO_2/(SO_1 \times SO_1) \cong S^1$ and hence
  
   \[
\pi_1(  SO_2/(SO_1 \times SO_1)) \cong \ZZ \qquad \text{and}  \qquad            \pi_2 (  SO_2/(SO_1 \times SO_1) = 0.
  \]
  \vspace{.1in} \paragraph{  
\bf 2. \underline{The space $O_3/(O_1 \times O_2)$ }}
   \indent
   
  We need some preliminary calculations. First of all, 
  \[
  O_3/O_2 \cong SO_3/SO_2 \cong RP^3/S^1.
  \]
  The long exact sequence of the associated fibration  has the form
  \[
  \to \pi_2(RP^3) \to \pi _2(RP^3/S^1) \to \pi _1(S^1) \overset{p}\longrightarrow \pi _1(RP^3) \to \pi _1(RP^3/S^1) \to 0
\]
  which simplifies to 
   \[
   0 \to \pi _2(RP^3/S^1) \to \ZZ  \overset{p} \longrightarrow \ZZ_2 \to \pi _1(RP^3/S^1) \to 0
\]
with $p $ surjective, and this implies that 
\[
  \pi _1(RP^3/S^1)    = 0         \qquad \text{and}  \qquad        \pi _2(RP^3/S^1) = \ZZ
  \]
  and, translating back, we have
  \[
  \pi _1(O_3/O_2)    = 0         \qquad \text{and}  \qquad        \pi _2(O_3/O_2) = \ZZ
  \]
  and
  \[
  \pi _1(SO_3/SO_2)    = 0         \qquad \text{and}  \qquad        \pi _2(SO_3/SO_2) = \ZZ
  \]
  
  Now, apply these results to the long exact sequence associated to the fibration 
  \[
  O_1 \longrightarrow O_3/O_2 \longrightarrow O_3/(O_1 \times O_2)
  \]
  and one obtains 
  \[
  0 \to \pi_2(O_3/O_2) \cong \ZZ \longrightarrow \pi _2 (O_3/(O_1 \times O_2)) \to 0
  \]
  and 
  \[
   0 \longrightarrow \pi _1 (O_3/(O_1 \times O_2) \to \pi _0(O_1)\cong \ZZ_2   \to 0
  \]
  so that 
  \[
  \pi _1 (O_3/(O_1 \times O_2)  = \ZZ _2 \qquad \text{and}  \qquad \pi _2 (O_3/(O_1 \times O_2) \cong \ZZ .
  \]
  
\vspace{.1in} \paragraph{  
\bf 3. \underline{The space $SO_3/(SO_1 \times SO_2)$}}
   \indent
  
  Similarly, we have 
   \[
  0 \to \pi_2(SO_3/SO_2) \cong \ZZ \longrightarrow \pi _2 (SO_3/(SO_1 \times O_2)) \to 0
  \]
  and 
  \[
   0 \longrightarrow \pi _1 (SO_3/(SO_1 \times SO_2) \to \pi _0(SO_1)\cong 0  \to 0
  \]
so that 
  \[
  \pi _1 (SO_3/(SO_1 \times SO_2)  =  0 \qquad \text{and}  \qquad \pi _2 (SO_3/(SO_1 \times SO_2) \cong \ZZ .
  \]
 \vspace{.1in} \paragraph{  
\bf 4. \underline{The space $O_4/(O_2 \times O_2)$  }}
     \indent
     
  Applying Sequence~(\ref{A}) to $O_4/(O_2 \times O_2)$ gives us 
  \[
  0 \to \pi _2 (O_4/(O_2 \times O_2)) \to \pi _1(O_2 \times O_2) \to \pi_1(O_4) \to \pi _1 (O_4/(O_2 \times O_2)) \to
  \]
  \[
  \to \pi _0(O_2 \times O_2) \to \pi_0(O_4) \to \pi _ 0(O_4/(O_2 \times O_2)) \to 0
  \]
  which, after applying what we already know, becomes 
 \[
  0 \to \pi _2 (O_4/(O_2 \times O_2)) \to\ZZ^2 \twoheadrightarrow  \ZZ _2
  \overset{p}\longrightarrow   \pi _1 (O_4/(O_2 \times O_2))
   \to \ZZ_2 \oplus \ZZ_2  \to\ZZ_2 \to  0
  \]
  This implies that $p$ must be the zero map, and then we have 
  \[
  0 \to \pi _2 (O_4/(O_2 \times O_2)) \to\ZZ^2  \twoheadrightarrow \ZZ _2
 \to 0
 \]
 and
 \[
 0  \to  \pi _1 (O_4/(O_2 \times O_2))
   \to \ZZ_2 \oplus \ZZ_2  \to\ZZ_2 \to  0
  \]
  so that
   
  \[
      \pi _1 (O_4/(O_2 \times O_2))  \cong \ZZ_2     \qquad \text{and}  \qquad    \pi _2 (O_4/(O_2 \times O_2)) \cong \ZZ^2
  \]
  
 \vspace{.1in} \paragraph{  
\bf 5. \underline{The space  $SO_4/(SO_2 \times SO_2)$}}
   \indent

  The sequence for $SO_4 / (SO_2 \times SO_2) $ is simpler.  We wind up with 
  \[
  0 \to    \pi _2 (SO_4/(SO_2 \times SO_2)) \longrightarrow \pi _1(SO_2 \times SO_2)  \twoheadrightarrow
  \pi_1(SO_4) \longrightarrow   \pi _1 (SO_4/(SO_2 \times SO_2)) \to 0
  \]
  which simplifies to
    \[
  0 \to    \pi _2 (SO_4/(SO_2 \times SO_2)) \longrightarrow \ZZ^2  \twoheadrightarrow
 \ZZ_2 \overset{p}\longrightarrow   \pi _1 (SO_4/(SO_2 \times SO_2)) \to 0
  \]
The map $p$ must be the zero map, and so 
\[
 \pi _1 (SO_4/(SO_2 \times SO_2)) \cong 0 \qquad \text{and}  \qquad     \pi _2 (SO_4/(SO_2 \times SO_2)) \cong \ZZ^2
\]

 \vspace{.1in} \paragraph{  
\bf 6. \underline{The space  $O_{n}/(O_k \times O_{n-k})$}}
  \indent
  
Suppose as usual  that $2 < k \leq n-k$. Then the sequence 
  \[
  0 \to \pi_2(O_{n}/(O_k \times O_{n-k})) \to \pi _1(O_k \times O_{n-k} ) \overset\gamma\to  \pi _1(O_{n}) \to \pi _1(O_{n}/(O_k \times O_{n-k})) \to 
  \]
  \[
 \to  \pi _0(O_k \times O_{n-k})  \to \pi _0(O_n) \to    \pi_0(O_{n}/(O_k \times O_{n-k}))
   \]
simplifies to
  \[
  0 \to \pi_2(O_{n}/(O_k \times O_{n-k})) \to \ZZ_2 \oplus \ZZ_2 \overset{\bar\gamma}\longrightarrow  \ZZ_2 
 \overset{\zeta}\longrightarrow \pi _1(O_{n}/(O_k \times O_{n-k})) \to \ZZ_4 \overset{\gamma}\longrightarrow  \ZZ_2\to 0       
   \]
The map $\bar\gamma $ is surjective, and hence   $\zeta $ is the zero map.  The map $\gamma $ is also surjective. Thus 
\[
\pi _1(O_n/(O_k \times O_{n-k})) \cong  \ZZ_2 \qquad \text{and}   \qquad \pi_2(O_{n}/(O_k \times O_{n-k})) \cong \ZZ_2 .
\]
 
 \vspace{.1in} \paragraph{  
\bf 7. \underline{The space $ SO_{n}/(SO_k \times SO_{n-k})$}}
  \indent
  
 The case of 
$ SO_{n}/(SO_k \times SO_{n-k})$ is much simpler. The space $SO(n)$ is path-connected so Sequence~(\ref{B}) gives us 
\[
0 \to \pi_2  (SO_{n}/(SO_k \times SO_{n-k})) \longrightarrow 
  \pi_1  (SO_k \times SO_{n-k}) \longrightarrow \pi _1(SO_{n} ) \to  \pi_1  (SO_{n}/(SO_k \times SO_{n-k})) \to 0.
 \]
  Since  $k > 2$, the sequence simplifies to 
 \[
0 \to \pi_2  (SO_{n}/(SO_k \times SO_{n-k})) \longrightarrow 
  \ZZ_2 \oplus \ZZ_2 \twoheadrightarrow \ZZ_2  \to  \pi_1  (SO_{n}/(SO_k \times SO_{n-k})) \to 0.
 \]
and so 
\[
 \pi_1  (SO_{n}/(SO_k \times SO_{n-k}))  = 0 \qquad \text{and}  \qquad         \pi_2  (SO_{n}/(SO_k \times SO_{n-k})) = \ZZ_2.
\]

 To summarize:

 \begin{center}
     \begin{tabular}{    r|    c| c | c | c }
 $X$ & $\pi _0(X)$    &  $\pi _1(X)$   &   $\pi _2(X)$       \\  \hline\hline
        $O_{2}/(O_1 \times O_1) $ &     $0 $    &   $ \ZZ $           & $0$ & \\ \hline
          $O_{3}/(O_1 \times O_2) $ &     $ 0 $    &   $ \ZZ_2 $           & $\ZZ $ & \\ \hline
            $O_{4}/(O_2 \times O_2) $ &   $0 $    &  $ \ZZ_2 $         & $\ZZ^2 $ &  \\ \hline  
               $O_{n}/(O_k \times O_{n-k}) $ & $0 $    & $ \ZZ_2$         & $\ZZ_2 $ & $2 < k \leq n-k$\\ \hline  
          $SO_{2}/(SO_1 \times SO_1) $ &     $0 $    &   $ \ZZ  $           & $0$ & \\ \hline
          $SO_{3}/(SO_1 \times SO_2) $ &     $ 0 $    &   $ 0 $           & $\ZZ $ & \\ \hline       
            $SO_{4}/(SO_2 \times SO_2) $ &   $0 $    &  $0 $         & $\ZZ^2 $ &  \\ \hline  
               $SO_{n}/(SO_k \times SO_{n-k}) $ & $0 $    & $ 0$         & $\ZZ_2 $ & $2 < k \leq n-k$\\ \hline
    \end{tabular}
   \end{center}

  Here the connection is to $\cR_8 = {\bf{BO}} \times \ZZ  $ and $ O_{2n}/(O_n \times O_{n})$ converges to $ {\bf{BO}}.$

\vspace{.1in}

  \paragraph{{\bf\underline{Summary: The Eight Real Groups and Homogeneous Spaces}}}

The chart below summarizes the information we have about the stable behavior of the first three
homotopy groups of the eight real groups and homogeneous spaces.  As before, we insert the ``$\times \ZZ $"
as needed and adjust the results accordingly.    

Note that the bottom four rows, marked with **, are simply a repetition of the top four rows, inserted so the diagonal pattern is more evident.
 \vglue .5in
 
 \begin{Table}
 \label{Table:HomogeneousSpaces}
 \centerline{\bf{{The Eight Real Groups and Homogeneous Spaces} }}

  \begin{center}
     \begin{tabular}{    r|   c| c | c | c }
 $X$ & $\pi _0(X)$    &  $\pi _1(X)$   &   $\pi _2(X)$       \\  \hline\hline
      $O_n $ & $\ZZ_2 $    & $ \ZZ_2$ & $0$ & $n  \geq 3$\\ \hline
      $O_{2n}/U_n$ & $\ZZ_2 $    & $ 0$ & $\ZZ $ & $n  \geq 2$\\ \hline  
        $U_n/Sp_n $ & $0 $    & $ \ZZ $ & $0$ & \\ \hline
        $Sp_{n}/(Sp_k \times Sp_{n-k}) \times \ZZ$ & $\ZZ $    & $ 0 $ & $0$ & \\ \hline
         $Sp_{n}  $ & $0 $    & $ 0 $ & $0$ & \\ \hline
    $Sp_{n}/U_n  $ & $0 $    & $ 0 $ & $ \ZZ $   \\ \hline
$U_n/O_n$ & $0 $    & $ \ZZ $ & $\ZZ_2 $ & $n  \geq 3$\\ \hline  
 $O_{n}/(O_k \times O_{n-k}) \times \ZZ $ & $\ZZ $    & $ \ZZ_2$         & $\ZZ_2 $ & $2 < k \leq n-k$\\ \hline 
  **    $O_n $ & $\ZZ_2 $    & $ \ZZ_2$ & $0$ & $n  \geq 3$\\ \hline
   **   $O_{2n}/U_n$ & $\ZZ_2 $    & $ 0$ & $\ZZ $ & $n  \geq 2$\\ \hline  
      **  $U_n/Sp_n $ & $0 $    & $ \ZZ $ & $0$ & \\ \hline
         **   $Sp_{n}/(Sp_k \times Sp_{n-k}) \times \ZZ$ & $\ZZ $    & $ 0 $ & $0$ & \\ \hline
    \end{tabular}
   \end{center}
   \end{Table}

   Look carefully and you will see that parts of Bott Periodicity are revealed, noting the entries on the diagonal lines moving up from left to right. 
   
 \index{homotopy group|)}
  \index{classical group|)}
  \index{homogeneous space|)}

\newpage
  \section{\bf{The Basics: Real Vector Bundles}}
  \label{Section:RealVectorBundles}
 
 In this basic section we provide the definitions and results needed to define $KO^*(X)$ for compact spaces in the following section. 
   
 The elementary theory of {\emph{complex}} vector bundles and $K^*(X)$ has been explicated beautifully, particularly by Atiyah \cite{Atiyah}, Karoubi \cite{K}, Husemoller \cite{H}, Blackadar \cite{Blackadar}, 
 and more recently by Hatcher \cite{Hatcher}, Park \cite{Park}.\footnote{We shall refer to  these luminaries as the ``complex authorities".}   So one is tempted to say ``read these sources but everywhere you see {\emph{complex}}  insert {\emph{real.}}"   This 
 approach falls apart when one gets to graded $K$-theory and more sophisticated matters.  So we will briefly sketch out the {\emph{real}} facts of life and urge the reader 
 to look at the complex sources as well. At several points in the next few sections we will state a result in detail for the real setting 
and then very briefly indicate that the complex analog holds as well. 
 
 A separate question arises to the practical physicist: Why isn't this Chapter 1?  Why do we need all this operator algebra stuff?  Our response to this is mostly practical. 
 First, we need to use algebras such as $C(X, M_n(\RR))$. For a deeper answer, consider that 
 if we stick to real vector bundles then we can suspend but we cannot desuspend. We need the real commutative algebras such as $\realo (i\topr ) $ and $\real (iS^1)$ to perform 
 this trickery, and they are NOT of the form $\realo (X) $ or $\real (X)$.

    Deeper yet, we know from studying complex theory that deep results such as the K\"unneth Theorem and the Universal 
 Coefficient Theorem are quite useful in many places in physics, and we will see that analogous real results will help us as well. 
 Thus  we need $KO$ for real $C^*$-algebras and not just for real vector bundles.  
 
 Having said that, we recognize that $KO^*(X)$ is important in its own 
 right, as its use of vector bundles brings us much closer to topological and geometric properties that we need.

   \begin{Def} A {\emph{real vector bundle}} 
   \index{vector bundle}
   $(E, p, X)$ (frequently denoted by just $E$) consists of a continuous map $p: E \to X$ of topological spaces such that 
 \begin{enumerate}
 \item The space $E_x = p^{-1}(x)$ is a real vector space for each $x \in X$ called the fibre over $x$,  and 
 \item The following {\emph{local triviality}} condition holds: for each $x \in X$ there is an open set $U_x \subset X$ with  $x \in U_x$  and a homeomorphism $h $ such that 
 the diagram 
 \[
 \begin {CD}  
 p^{-1}(U_x) @>h>>  U_x \times \RR  ^k  \\
 @VVpV                  @VVp_1 V  \\
 U_x @>\id>>          U_x      
 \end{CD}
 \]
 commutes, where $p_1$ is projection on the first coordinate.  In addition, for each $x \in X$ the induced map $h : E_x \to  \RR ^k  $ must be linear.
 \end{enumerate}
 \end{Def}

 Note that we could avoid some  work here by using principal $O_n$-bundles, as their associated $\RR ^n$-bundles are essentially real vector bundles. That is, a real vector bundle $E \longrightarrow X$ of dimension $n$ over some compact space is essentially equivalent to an associated bundle of the 
 form 
 \[
 Y \times _{O_n} \RR ^n \longrightarrow X .
 \]
 Complex vector bundles are defined similarly, substituting ``complex" for ``real" and ``$\CC$" for ``$\RR$" everywhere.
 
 The simplest examples of vector bundles over $X$ are the product bundles 
 \[
 X \times \RR^k \to X
 \]
  denoted $\theta _k$ when $X$ is understood.\footnote{There is a bit 
 of debate about whether $\theta _0$ exists; if you believe that there is a vector space of dimension $0$ then your answer is ``yes".}

 \begin{Def} Suppose that $(E, p,  X) $ and $(E', p',  X')$
 are real vector bundles. 
  A {\emph{morphism}} or {\emph{map}} of real vector bundles    is a commuting diagram 
  \[
  \begin{CD} E' @>{\tilde{f}}>>     E \\
  @VVp' V       @VV pV      \\
  X'    @>f>>        X
  \end{CD}
  \]
  of continuous maps, 
  where 
  \[
  {\tilde{f}} : E_y' \longrightarrow E_{f(y)}
  \]
    is a linear transformation for each $y \in X'$. 
   If $\tilde{f}$  is a homeomorphism then it  is called an {\emph{isomorphism}} of vector bundles. \index{vector bundle!isomorphism of}
 Any bundle isomorphic to a product bundle is called a {\emph{trivial}} bundle. \index{vector bundle!trivial}
 \end{Def} 
 
 If $X$ is a path-connected space then the local triviality condition insures that  fibres of a vector bundle must all be of the same dimension: this is called the {\emph{rank}} 
 or {\emph{dimension}} of the vector bundle.   If $X$ has several path components then a vector bundle over $X$ may well have a different rank on each path component. 
 
 Perhaps the most important real vector bundles  arise as   {\emph{tangent bundles}} \index{vector bundle!tangent bundle} of   compact smooth manifolds   $M$. Say that  $M$ has dimension $k$.  
 Each point $m \in M$ has a $k$-dimensional real vector space $(TM)_m$ associated to it, namely the vectors tangent to the manifold $M$ at the point $m$.  These coalesce 
 to form the $k$-dimensional tangent bundle $TM $ with projection map $p: TM \to M$ with $p^{-1}(m) = (TM)_m$. 
 We denote the tangent bundle as $TM \to M$ or simply $TM$.     Local triviality says that each point $m$ of the manifold has an open set $U$ surrounding it so that if we restrict the bundle to $U$ then, up to isomorphism,  the bundle takes the form 
 \[
 U \times \RR ^k \longrightarrow U.
 \]
 We may form the dual  {\emph{cotangent bundle}} \index{vector bundle!cotangent bundle} $T^*M \to M$ similarly, with 
 \[
 (T^*M)_m \cong \Hom ((TM)_m, \RR ) \qquad \forall m \in M. 
 \]
 
 A further example. Think of the sphere $S^k \subset \RR ^{k+1}$.  Then each point $m \in S^k $ has vectors based at the point $m$ and perpendicular  to the sphere. This forms a one-dimensional 
 bundle $N(S^k)$ called the {\emph{normal bundle}} \index{vector bundle!normal bundle} of $S^k$. 
 
 \begin{Def} Suppose that $f: X \to Y$ and $E \to Y$ is a real vector bundle over $Y$.  Then we define the
  real vector bundle $f^*E \to X$ to be the {\emph{pullback}} \index{pullback!vector bundle} \index{vector bundle!pullback}  in the commuting diagram
 \[
 \begin{CD}
 f^*(E) @>{\widetilde f}>>     E   \\
 @VVV      @VVpV   \\
 X  @>f>>      Y.
 \end{CD}
 \]
 That is, 
 \[
 f^*(E) = \{ (x, e) \in X \times E \mid f(x) = p(e) \}.
 \]
 \end{Def}
 
 It is a nice exercise to prove that $f^*(E) \to X$ is really a real vector bundle, of the same dimension as $E$, and each 
 $\widetilde f : f^*(E)_x \to E_{f(x)} $ is an isomorphism on fibres.  If $X$ is 
 a subset of $Y$ and $f: X  \hookrightarrow Y$ is the inclusion map,    then $f^*(E)$ is called the {\emph{restriction}} 
 \index{vector bundle!restriction} of the bundle to $X$ and sometimes it is  written $E|_X$. 
 
 \begin{Def} Suppose that $p: E \to X$ is a real vector bundle of $X$. A {\emph{cross-section}} or  {\emph{section}} 
\index{vector bundle!section} of the bundle is a continuous 
 function $s: X \to E $ with the property that $ps : X \to X$ is the identity.    It is {\emph{nowhere zero}} if $s(x) \neq 0 \in E_x$ for every $x \in X$. 
 The set of all sections of the real vector bundle $E$ is denoted $\Gamma(E)$.   
  \end{Def}
 For example, the trivial bundle $\theta _n \to X$ has $n$ linearly independent nowhere zero sections $s_1,\dots s_n$.  
 That is, for each $x \in X$, the set of vectors
 \[
 \{ s_1(x), \dots s_n(x) \}
 \]
 is linearly independent.

 On the other hand, the tangent bundle of $S^2$ has the property that 
 every section $s: S^2 \to T(S^2)$ must be zero at some point.\footnote{If this were not the case then  you could comb the hair of a coconut. This was establshed by H. Poincar\'{e} in 1885. }

 \begin{Def} 
 \index[notation]{vectrx@$\rm{Vect}\pr (X)$}
 Suppose that $E$ and $F$ are real vector bundles over the compact space $X$.  We define their {\emph{sum}} \index{vector bundle!direct sum} to be  the real vector bundle $E \oplus F$ where the fibre over 
 the point $x$ is $E_x \oplus F_x$.   Similarly, the {\emph{tensor product}} \index{vector bundle!tensor product} of the bundles $E \otimes F$ is the real vector bundle where the fibre 
 over the point $x$ is $E_x \otimes _{\RR} F_x$.\footnote{This requires more careful explanation but is identical to the complex situation, so please refer to the complex authorities.}  
 For example, if $E$ is a real vector bundle over $X$ and  $\theta _1$ is the one-dimensional trivial bundle, then $E \otimes \theta _1 \cong E$,  generalizing the 
 fact that $\RR ^n \otimes \sr \RR \cong \RR ^n$. 
 Finally let $\rm{Vect}_k\pr (X)$ denote the set of isomorphism classes of real vector bundles of dimension $k$ and $\rm{Vect}\pr (X) $ denote the union of the $\rm{Vect}_k\pr (X)$. 

 \end{Def}
 
 It is easy to see that the operations $E \oplus F$ and $E\otimes F$ give $\rm{Vect}\pr (X) $ the structure of an abelian semiring with unit $\theta _1$. 
 It turns out that we can calculate it as well in homotopy terms.

 It is best to restrict at this point to bundles over compact spaces.  One very important consequence of this restriction is the following fact: If $E$ is any finite-dimensional real vector bundle 
 over a compact space  $X$ then there exists another finite-dimensional real vector bundle $F$ over $X$  such that $E \oplus F$ is a trivial bundle.\footnote{cf. Corollary 4.14 of  Atiyah \cite{Atiyah}.} 
  
  Over a compact space, isomorphism classes of real vector bundles correspond to homotopy classes of certain maps. Here is the result:

 \begin{Thm} (cf. Husemoller \cite{H}, page 34, Theorem 7.2 ) 
 \begin{enumerate}
 \item
 For each $k$, there is a canonical universal $k$-dimensional real vector bundle $\gamma _k: EO_k  \to BO_k$ such that for any compact  space $X$ there 
 is an isomorphism of commutative semirings 
 \[
 [ X, BO_k ]^{ub} \overset\cong\longrightarrow \rm{Vect}_k\pr (X )
 \]
 induced by sending a function $f: X \to BO_k$ to $f^*\gamma _k$. 
   \item
 For each $k$, there is a canonical universal $k$-dimensional complex vector bundle $\gamma _k: EU_k  \to BU_k$ such that for any compact  space $X$ there 
 is an isomorphism of commutative semirings 
 \[
 [ X, BU_k ]^{ub} \overset\cong\longrightarrow \rm{Vect}_k\pr (X )
 \]
 induced by sending a function $f: X \to BU_k$ to $f^*\gamma _k$. 
 \end{enumerate}
 \end{Thm}
 
 \begin{Rem} One should compare this result to Remark~\ref{classify} where the classification of principal $G$-bundles is discussed. \index{principal $G$-bundle} The vector bundle $\gamma _k$
 is strongly related to the universal principal $O_n$-bundle.  
  
 \end{Rem}
 
 \begin{Def} Suppose that $\pi : E \to X$ is a real vector bundle over the compact space $X$.        Then 
 the set of sections $\Gamma (E) $ is a module \index{projective module} over the commutative ring $\real (X)$ with the following operations:
 \begin{enumerate}
 \item Addition defined by 
 \[
 (s_1 + s_2)(x) = s_1(x) +  s_2(x) \; .
 \]

 \item Scalar multiplication defined by 
 \[
 (rs)(x) = r(s(x)) \qquad r \in \RR , s \in \Gamma(E) \; .
 \]
 
 \item The action of $C\pr (X)$ on $\Gamma (E)$ is defined by 
 \[
 (f \circ s) (x) = f(x)s(x) \qquad f \in  \real (X),   s \in \Gamma (E) \; .
 \]
 \end{enumerate}
 \end{Def}

It is a basic result that a real bundle $E$ of dimension $n$ has $n$  linearly independent non-vanishing 
sections if and only if $\Gamma (E) \cong \real (X)^n$ and hence $E$ is isomorphic to the trivial real bundle $\theta _n$ of dimension $n$.
Similarly, for a complex vector bundle $\pi \colon E \to X$, the set of sections $\Gamma(E)$ is a module over $\complex(X)$, and $\Gamma(E) \cong \complex(X)^n$ if and only if $E$ is isomorphic to the trivial complex vector bundle $\theta_n$.

As remarked above, every real vector bundle $E \to X$ over a compact space is a direct summand of a trivial bundle:
 \[
 E \oplus F \cong \theta _n \;  .  
 \] 
 Take sections on both sides and note that 
 \[
 \Gamma (E) \oplus \Gamma (F) \cong   \Gamma (E\oplus F) \cong  \Gamma (\theta _n) \cong \real (X)\oplus  \overset{n} \cdots  \oplus \real (X) \; .
  \] 
  Thus $\Gamma (E)$ is a direct summand of a finitely generated free $\real (X)$-module and hence it is a  finitely generated projective $\real (X)$-module.
 This implies that the map 
 \[
 p: \real (X)^n \cong  \Gamma (E \oplus F) \longrightarrow \Gamma (E) \hookrightarrow \real (X)^n 
 \]
  is an orthogonal projection:  \index{projection}
  \[
  p = p^2 = p^*  \in M_n(\real (X)).
  \]
      Thus to each real vector bundle $E$ over a compact 
  space $X$ we have associated a projection  $ p = p(E) \in M_n(\real (X)) $. We denote this map as 
     \[
 \mathcal{P} : \rm{Vect}\pr (X) \longrightarrow \mathcal{P}_\infty (\real(X))
 \]
from the category $\rm{Vect}\pr  (X)$ of isomorphism classes of vector bundles over $X$ to the category $\mathcal{P}_\infty (\real (X) )$ of  projections in matrix rings over $C(X)$.

\begin{Rem} 
As noted earlier, all of the basic properties of real vector bundles on compact spaces noted in this section hold as well, with 
obvious modifications, for complex vector bundles. 

\begin{enumerate}

\item Locally, a complex bundle has the structure $U \times \CC ^k \longrightarrow U$.
\item
A complex bundle of complex dimension $n$  has an associated $U_n$-bundle 
\[
E \times _{U_n} \CC^n \to B. 
\]
\item The space $\Gamma (E)$ of complex sections of a complex bundle $E \to B$ over a compact space $X$ is a module over the ring $\complex (X)$.

\item
\[
[X, BU_k]^{ub} \cong \rm{Vect}_k\pc (X).
\]
\end{enumerate}

\end{Rem}

  Here's a similar construction that we will need in the next section.
  
  \begin{Def} Suppose that $p: E \to X$ is a real vector bundle of real dimension $n$. Then we may form the associated {\emph{complexified vector bundle}} 
  \index{vector bundle!complexification}
  \[
  E \otimes \CC = E \times _{\RR^n} \CC^n
  \] 
  of complex dimension $n$.
  On each fibre $E_x \cong \RR ^n$ we are simply forming 
  \[
  E_x \otimes _{\RR^n} \CC \cong \CC^n.
  \]
      This construction induces the {\emph{complexification  map}}  \index{complexification}
  \[
  c:    \rm{Vect}\pr _n (X) \longrightarrow  \rm{Vect}\pc _n (X)
  \]
  which is a map of semirings.
  \end{Def}
  
  Principal bundles \index{principal $G$-bundle} may be pulled back exactly as is done with vector bundles, and in fact, for $X$ compact, isomorphism classes of principal $G$-bundles correspond 
  exactly to elements of $[X, BG]$. 
  
  \index{vector bundle|)}  
 
\newpage
    \section{\bf{The  Basics: $KO$-theory for Spaces, Fredholm Operators, and Characteristic Classes}}
  \label{Section:KOforSpaces}
  
 In this section we present the basic properties of real $K$-theory $KO^*(X)$ for $X$ compact.  These are completely analogous to complex $K$-theory 
 $K^*(X)$\footnote{with the obvious adjustments caused by the fact that $2 \neq 8$} and we will briefly put these on the record as well. 
 
 \begin{Def} For a compact space $X$, define 
 \[
 KO^0(X) = \rm{Groth} ( \rm{Vect} \pr (X))
 \]
 \index[notation]{KO0X@$KO^0(X)$}
 the Grothendieck group of $\rm{Vect}\pr (X)$.  Elements of $KO^0(X)$ may be represented as formal 
 differences of vector bundles $[E] - [F]$  over $X$. Then $KO^0(X)$ is a commutative ring with identity $1 = [\theta _1]$, the class of the one-dimensional trivial bundle over $X$. 
 
 If $f: X \to Y$ then there is an induced map 
 \[
 f^* : KO^0(Y) \longrightarrow KO^0(X)
 \]
 defined on generators by 
 \[
 f^*[E] = [f^*(E) ]
 \]
 and then $KO^0$ is a contravariant\footnote{A functor $F$ is {\emph{contravariant}} if a map $f: A \to B$ induces a map 
 $f^* : FB \to FA$  with the usual naturality conditions.}  functor from compact spaces to commutative unital rings. 
 \end{Def}

 Here are a few  examples. 
 \begin{enumerate}
 \item If $X$ is contractible, such as the closed unit ball, then $X$ has only trivial bundles, and hence $KO^0(X) \cong \ZZ $.   
 In particular, let $ \pt$ denote the space with a single point, and then we have $KO^0(\pt) \cong \ZZ$ as commutative unital rings. Thus if $\pt$ is the basepoint of some compact space $X$ then 
 there is an inclusion map $\pt  \rightarrow X$ and an
 associated ``dimension" map 
 \[
 KO^0(X) \overset{\dim}\longrightarrow KO^0(\pt) \cong \ZZ    \; .
 \]
  It's called the dimension map because if $X$ is compact connected and $E \to X$ is a vector bundle over $X$ of dimension $n$, then $dim([E]) = n$.

 \item If $X$ is the circle  $S^1$ then the tangent bundle is trivial, as there is an obvious nowhere zero section. However, it turns out that there is a non-trivial real bundle $M$   of dimension $1$. It is essentially the Mobius strip!   Note that if you go around the 
 Mobius strip twice then you get a trivial bundle $\theta _2 = 2\theta _1$.  So $KO^0(S^1)$ is generated by $[\theta _1]$ and by $[M]$ with 
 $2([M]  - \theta _1) = 0$. Then we have 
 the split exact sequence 
 \[
 0 \longrightarrow  \ZZ_2 \longrightarrow KO^0(S^1)\overset{~\dim~~}\longrightarrow \ZZ \longrightarrow 0
 \]
 and hence $KO^0(S^1) \cong \ZZ \oplus \ZZ _2$ with the summands being generated by $\theta _1$ and $ [M] - \theta _1$ respectively and 
 \[
 \dim ( [M] - \theta _1 ) = 0.
 \]

 \item (Kaplansky) Suppose that $X = S^n$ some sphere. Then we have the tangent bundle $TS^n$ 
  \index{vector bundle!tangent bundle}
  and also the normal bundle 
 \index{vector bundle!normal bundle}
 $N( S^n)  \cong \theta _1$, consisting of all vectors at a point of the sphere that are perpendicular to the tangent plane of the sphere. It  is easy to see that 
 \[
 TS^n \oplus N(S^n) \cong \theta _{n+1}
 \]
  so, in particular, $TS^n$ plus a trivial bundle  is isomorphic to a  trivial bundle, whether $TS^n$  is trivial or not.  Such a bundle is called {\emph{stably trivial}}.    So the tangent bundle of $S^n$ is stably 
  trivial for all $n$.  This implies $[TS^n] = [\theta _n] $ in $KO^0(S^n)$ for all $n$.

  It is a rather deep result that $TS^n $ is isomorphic to a {\emph{trivial}}  bundle if and only if $n = 1,\,3,$ or $7$.  (This is related to the existence of the complex numbers, the quaternions, 
  and the Cayley numbers in dimensions $2,\,4,$ and $8$ and the fact that there are no other such algebras in other dimensions.) 
     \end{enumerate}
 
 We are particularly interested in relating $KO$ for spaces to $KO$ for $C^*$-algebras.  Here is the main result.

 \begin{Thm} (Swan \cite{Sw}) \label{Thm:Swan}  \index{Swan's theorem} Let $X$ be a compact space. 
 \begin{enumerate}
 \item There is a natural 
   isomorphism of semigroups
   \[
   \rm{Vect}\pr (X)  \simeq  \mathcal{P}_\infty(\real (X))
 \]
 which induces a natural isomorphism of groups
 \[
  KO^0(X) 
\cong    KO_0(\real(X))  .
 \]
 \item In the complex case, there is a natural 
   isomorphism \[
   \rm{Vect} (X)  \simeq  \mathcal{P}_\infty(\complex (X))
 \]
 which induces a natural isomorphism  
 \[
  K^0(X) 
\cong    K_0(\complex(X))  .
 \]
\end{enumerate}
   \end{Thm}  
 
 We will put off a discussion of the proof of Swan's Theorem and related results until the next section. See Theorems~\ref{T:swanK} and \ref{T: swan}.

Here's an immediate consequence.

\begin{Def}  
\index[notation]{KOnx@$KO^{n}(X)$} Let $X$ be a compact space. Then the  groups $KO^{-n}(X) $ are defined by 
\[
KO^{-n}(X) = KO_n(\real (X) ) = KO_0(S^n \real(X)).
\]
If $f: X \to Y$ then this definition yields a map 
\[
f^*: KO^{-n}(Y)  \longrightarrow    KO^{-n}(X)
\]
and  each $KO^{-n}$ is a contravariant, homotopy-invariant functor from compact spaces to abelian groups.
\end{Def}

\begin{Def}
\begin{enumerate}

\item
For a locally compact space $X$ with one-point compactification $X^+$,  define
\[
KO^0(X) = \ker \big[ KO^0(X^+) \longrightarrow KO^0(+) \cong \ZZ  \big]
\]
so there is a natural short exact sequence
\[
0 \longrightarrow KO^0(X) \longrightarrow KO^0(X^+) \longrightarrow \ZZ \longrightarrow 0
\]
and an (unnatural) isomorphism 
\[
KO^0(X^+) \cong KO^0(X) \oplus \ZZ
\]
A proper\footnote{Proper maps of locally compact spaces are those that extend uniquely to one-point compactifications}   map $X \to Y$ of locally compact spaces induces a homomorphism $KO^0(Y) \to KO^0(X)$, and $KO^0$ is a contravariant 
functor from the category of locally compact spaces and proper maps and homomorphisms to abelian groups. 

\item 
 For $n > 0$ and  $X$ a locally compact space,  define
 \[
 KO^{-n}(X) = KO_n(\realo(X))
\]

\item
As $KO_n(A)$ is periodic of period 8 by Bott Periodicity for all real $C^*$-algebras, we conclude that the same holds for the special case of $KO^*(X) $ and hence we may extend $KO^*(X)$ to a $\ZZ$-graded 
collection of contravariant functors which satisfy 
\[
KO^n(X) \cong KO^{n-8}(X) \qquad \forall n \in \ZZ .
\]

\item We similarly define complex $K$-theory $K^{-n}(X)$ for $X$ compact, and locally compact; and we note that it too extends to a $\ZZ$-graded collection of contravariant functors which satisfy 
\[
K^n(X) \cong K^{n-2}(X) \qquad \forall n \in \ZZ .
\]
 \end{enumerate}
 \end{Def}
 
 \begin{Def} 
 \index[notation]{KOnxred@$\widetilde{KO}^{n}(X)$}
 For $X$ a compact space with basepoint $x_o$, we define the {\emph{real reduced $KO$-theory groups}} by
 \[
 \widetilde{KO}^n (X) = \ker \big[ KO^n(X) \longrightarrow KO^n(x_o) \big]      .
 \]
 and note at once that there is an unnatural isomorphism
                                   \[  KO^n(X)     \cong     \widetilde{KO}^n(X)  \oplus   \ZZ  \; . \]

\end{Def}

 \begin{Thm}  Suppose that $X$ is a based compact space. Then 
 \begin{enumerate}
 \item
 \[
 \widetilde {KO}^0(X) \cong [X, \bbo \times \ZZ ] \cong [X, \cR _0 ]. 
 \]
 
 \item 
 For $0 \leq n \leq 7$, 
 \[
  \widetilde {KO}^{-n}(X) \cong   \widetilde {KO}^0(S^nX) \cong [S^nX, \bbo \times \ZZ ]   \cong 
 \index{suspension!of a space} \index{loop space}
 \]
 \[
  \cong [S^nX, \cR _0 ] 
 \cong [X , \Omega ^n \cR_0 ]  \cong [X, \cR_n]  .
 \]

 \end{enumerate}
 \end{Thm}
 	 
 Here is a comprehensive statement of the main properties of the reduced cohomology theory $ \widetilde{KO }^*(X)$. All of these facts are found in Karoubi \cite{K} and Hatcher \cite{Hatcher}
 and hold just as well for the real and complex settings.
 We could also deduce these properties for $KO^*(X)$ from the corresponding properties of $KO_*(A)$, using Swan's Theorem --  Theorem~\ref{Thm:Swan}.
   
   \begin{Thm}  \index{cohomology theory!reduced cohomology theory}    
  $ \widetilde{KO }^*(X)$ is a reduced cohomology theory on compact spaces.
   That is, for each $j \geq 0$, 
  \begin{enumerate}
  \item 
  $\widetilde{KO }^{-j}(-)$ is a contravariant functor 
   from based compact spaces to  abelian groups.
   
    \item
    
   $\widetilde{KO }^{*}(-)$ is a contravariant functor from based compact spaces to 
   $\ZZ$-graded commutative rings.
      
  \item
  
  If $f \simeq g$ are based homotopic maps $X \to Y$ then 
  \[
  f^* = g^* : \widetilde{KO }^{-j}(Y) \to \widetilde{KO }^{-j}(X).
  \]
     
  \item

  There is a suspension isomorphism \index{suspension!isomorphism} 
  \[
  \widetilde{KO }^{-j}(SX) \cong \widetilde{KO }^{-j-1}(X) .
  \]
  \item If $\iota : Y \to X$ is an inclusion of a closed based subset of $X$ with associated sequence
 \[
 Y \overset{\iota}\longrightarrow X \overset{\pi }\longrightarrow  X/Y
  \]
     then there is a long exact sequence
  \[
\dots \longrightarrow 
    \widetilde{KO}^{-j-1}(X) \overset{\iota ^*}\longrightarrow 
    \widetilde{KO }^{-j-1}(Y) \longrightarrow 
     \widetilde{KO }^{-j}(X/Y) \overset{\pi ^*}\longrightarrow 
     \widetilde{KO }^{-j}(X) \overset{\iota ^*}\longrightarrow  \dots
 \]
    
  \item For $X$ and $Y$ based compact spaces, there is a natural short exact sequence 
 
  \[
  0 \longrightarrow 
    \widetilde{KO }^{-j} (X \wedge Y) 
  \longrightarrow
   \widetilde{KO }^{-j} (X \times Y)
   \longrightarrow
    \widetilde{KO }^{-j} (X)  \oplus \widetilde{KO }^{-j} ( Y) 
   \longrightarrow
   0
   \]
   which splits unnaturally.

    \item 
  If   
 $X \cong \Invlim X_j $, then the natural projection maps $X \to X_j $ induce maps $KO ^*(X_j) \longrightarrow KO ^*(X) $ and these maps coalesce to yield an 
 isomorphism 
 \[
 \lim_{j \to \infty}  KO ^*(X_j) \overset{\cong}\longrightarrow KO ^*(X) .
 \]
 and similarly for $\widetilde {KO^*}(X)$.
 
 \item
 Bott Periodicity:  There is a natural isomorphism \index{Bott Periodicity}
 \[
 \widetilde {KO}^*(X) \cong \widetilde {KO}^{*-8}(X).
    \]
       
\item K\"unneth Theorem (special case) \index{K\"unneth Theorem} \cite{BoersemaKT}.  If $X$ and $Y$ are compact spaces and $Y$ has the homotopy type of a sphere or a product 
of spheres, then the internal product induces  a natural isomorphism
\[
\widetilde{KO}^*(X) \otimes _{KO^*(\pt)} \widetilde{KO}^*(Y) \overset{\cong}\longrightarrow \widetilde{KO}^*(X \wedge Y).
\]

    \end{enumerate}
    \end{Thm}

 \begin{Rem}\label{algebras} If $X$ is an infinite CW-complex \index{CW-complex} (see Definition~\ref{CWcomplex}) such as $RP^\infty $, $BO_n$,  or $\bbo $,  then there are several ways to generalize $KO $-theory to include them. However, 
 we generally must give up the idea of describing elements of the resulting group via  real  vector bundles. 
 \end{Rem}

 \vspace{.1in}
{\bf{Fredholm operators}}
\vspace{.1in}

\index{Fredholm operator}

 There are strong analytic ties between $KO$- and $K$-theories for spaces and bounded operators on real/complex Hilbert space.
In this section we start on the most elementary level, that of Fredholm operators.  We state the definitions for the real case, as the complex
case is easy to look up, such as in \cite{Atiyah} or \cite{Mu}

\begin{Def} Suppose that $T \in \Lin(\hil)$ is a bounded operator on real Hilbert space.  We say that $T$ is a {\it{Fredholm operator}} if there
exists a bounded operator $S$ such that $TS - I$ and $ST - I$ are compact operators. Equivalently, we may require instead that
$T$ has closed range and
$\ker(T)$ and $\ker(T^*)$ are both finite-dimensional.  Let $\mathcal{F}\sr \subset  \Lin(\hil)$ denote the space of all Fredholm operators. We define the {\it{Fredholm index}} of such an operator $T$ by
\[
\ind(T) = \dim \ker(T)  - \dim \ker(T^*)  .
\]
so that  $\ind : \mathcal{F}\sr  \longrightarrow \ZZ $.
\end{Def}

There are easy examples. Obviously, every invertible operator $T$ is Fredholm, and has $\ind(T) = 0$.   If we take $e_0, e_1, e_2 \dots $ to be the basis of the Hilbert space and define the unilateral shift by $T(e_n) = e_{n+1} $ then $\ker(T) = 0$ and $\ker(T^*)$ is one-dimensional with basis $e_0$. and hence
\[
\ind(T) = 0 - 1 = -1.
\]

The following theorem provides a secure link between Fredholm operators and topological K-theory. Atiyah \cite{Atiyah} proved it 
in the complex setting. K. Janich \cite{janich} independently proved it in the real,
complex, and quaternionic settings. (Much more sophisticated results in the real setting were proved years later by Mingo \cite{mingo} and 
Kasimov \cite{Kasimov};
see Section 2.2 of Schr\"oder \cite{Schr}.)

%

\begin{Thm} For all compact spaces $X$ there are natural isomorphisms
\[
 [ X , \mathcal{F}\sc ]^{ub}   \cong     K^0(X)
 \]
and
\[
[ X , \mathcal{F}\sr ]^{ub}  \cong   KO^0(X)
\]
induced by the index map.
\end{Thm}

Indeed, if we take $X$ to be a single point, then a map from $X$ to $\mathcal{F}$ simply chooses a path component of $\mathcal{F}$
and these path components are determined by the index of the operators contained in each path component.

Using some serious homotopy theory, which we omit, we may conclude that the spaces of Fredholm operators actually classify:

\begin{Cor} There are homotopy equivalences
\[
\mathcal{F}\sc  \simeq \bbu \times \ZZ \qquad \text{and} \qquad  \mathcal{F}\sr  \simeq \bbo \times \ZZ.
\]
and so the spaces  $\mathcal{F}\sc$      and   $\mathcal{F}\sr $ are classifying spaces  \index{classifying space} for complex and real $K$-theory respectively.
\end{Cor}
       
\begin{Rem} $KO^*$ is thus a cohomology theory in the classical Eilenberg-Steenrod sense, 
and so it has a spectrum\footnote{The word ``spectrum" has many diverse meanings within science.  Here it 
is important to note that this use of the word has nothing at all to do with the spectrum of an operator on Hilbert space.  The word spectrum here means a sequence of spaces with related properties. 
For example, it might be the sequence $\{ \Omega ^n X \} $ for some space $X$ but with no assumption of 
periodicity. For much more information, see J.F. Adams \cite{A} Part III (2), p. 131 and following.}  
and an associated homology theory (cf. May \cite{May}). Unfortunately, this homology theory
does not have a concrete, simple characterization and one cannot, for instance, easily write down elements of the $KO$-homology of nice spaces. The Kasparov 
framework, using 
\[
KO_*(X) \cong KK_*\pr (\real (X), \RR ) 
\]
doesn't seem easy to use either.
Perhaps it will be of some use in the future to 
physicists, but that is out of our vision for now. 
     \end{Rem}

  {\bf{Characteristic Classes}}. \index{characteristic class}

  Both complex and real $K$-theory for spaces are related to the cohomology of the space.  This connection is useful but complicated. The following is a minimal sketch of the idea, first in the complex case and then in the real case. 
  
  Each complex vector bundle $E \to X$ has associated to it certain cohomology classes 
  \[
  c_i(E) \in H^{2i}(X ; \ZZ)
  \]
  called {\emph{Chern classes}}. \index{characteristic class!Chern class} (cf.  Milnor-Stasheff \cite{MS}, Chapter 14). These may 
  be defined using differential 
  geometry or purely topologically.  Here's a quick sketch. Suppose $E\to X$ 
  is a complex bundle of 
  dimension $n$. We use the fact that the integral cohomology of the space $BU_n$ is a polynomial algebra over the integers with
  chosen  generators $\{c_1, \dots , c_n\}$ with $\dim(c_j) = 2j$:
  \[
  H^*(BU_n ; \ZZ ) \cong \ZZ [c_1, \dots , c_n ] .    
\]
Assuming $X$ is compact, then there is a continuous map $f: X \to BU_n $   
such that $E \cong f^*(\gamma _n) $ is isomorphic to the pullback of the universal  bundle $\gamma _n$ over $BU_n$
and $f$ is unique up to based homotopy.
Define the Chern class $c_j(E)$ of the vector bundle $E$  by 
\[
c_j(E) = f^*(c_j) \in H^{2j}(X ; \ZZ) .
\]

 Next, there is a certain rational combination of the classes $c_j(E)$ called the {\emph{Chern 
 character}} \index{characteristic class!Chern character|(} and written
 \[
 \crn (E) \in \bigoplus _k H^{2k}(X ; \QQ  ).
 \]
  Here are its principal properties on the category of compact spaces.
  
  \begin{Pro}  (cf. Hatcher \cite{Hatcher}, Section 4.1.)
  \begin{enumerate}
  \item The Chern character is a natural transformation
  \[
  \crn: K^0(X) \longrightarrow \bigoplus _k H^{2k}(X ; \QQ  )
  \]
  and using suspension this may be extended to 
   \[
  \crn : K^{-1}(X) \longrightarrow \bigoplus _k H^{2k-1}(X ; \QQ  ).
  \]
  \item The Chern character is a homomorphism of graded commutative rings. Specifically, 
  \[
  \crn(E \oplus F) = \crn(E) + \crn(F)
  \]
  and
  \[
    \crn(E \otimes F) = \crn(E) \crn(F).
  \]
  \item The Chern character induces an isomorphism   of graded commutative rings
   \[
  \crn : K^*(X)\otimes \QQ  \overset{\cong}\longrightarrow \bigoplus _k H^k (X ; \QQ  ).
  \]
\end{enumerate}
\end{Pro}
   
   It is very important to note that in general the Chern character 
   takes values generally in rational cohomology $\oplus _k H^{2k}(X ; \QQ )$.
Only rarely does  it take values in integer cohomology  $\oplus _k H^{2k}(X ; \ZZ )$.
   
   Here's a good example. Take $X = S^{2n}$, the sphere of even dimension. Then its cohomology vanishes except in dimensions $0$ and $2n$, so the only Chern class of a complex bundle 
   $E \to S^{2n}$ that has a chance of being non-zero is $c_n(E)$.  Here is what happens.

  For the even spheres, the Chern character simplifies considerably. If $E$ is any complex bundle
  over $S^{2n}$, then 
  \[
  \crn (E) = \dim (E) + \frac{(-1)^{n+1}c_n(E)}{(n-1)! } .
  \]
  \index{characteristic class!Chern character|)}
Using this information, we have the following result.

   \begin{Pro} (cf. Hatcher \cite{Hatcher}, Corollary 4.4) 
  A class in $H^{2n}(S^{2n} ; \ZZ )\cong \ZZ  $ occurs as a Chern class $c_n(E)$ if and only 
   if it is divisible by $(n-1)!$.
 
   \end{Pro} 
   
      Moving on, we come to the case of $KO^*(X)$, and for this we refer to 
   Milnor-Stasheff \cite{MS} Chapter 15 and to Hatcher \cite{Hatcher}. 
   Suppose that  $E \to X$ is a real vector bundle    with complexification  $E\otimes \CC$. \index{complexification} It turns out that 
   \[
  2 c_{2j + 1}(E \otimes \CC) = 0 \quad \forall j      
   \]         
   so the only elements of the form
   $c_{k}(E \otimes \CC)$
    that may be nonzero after tensoring with $\QQ$  have $k$ even,  so that the dimension of the element  is 
   some multiple of $4$.

  Define the $i$'th {\emph{Pontrjagin}} \index{characteristic class!Pontrjagin class}
  class $p_i(E)$ of the real bundle $E \to X$ by 
  \[
  p_i(E) = (-1)^i c_{2i}(E \otimes \CC ) \in H^{4i}(X ; \ZZ ) .
   \]
   Define the  {\emph{Pontrjagin character}} \index{characteristic class!Pontrjagin character} $Ph$ to be the composite
   \[
   \pont : KO^{0}(X) \overset{c}\longrightarrow K^0(X) \overset{\crn}\longrightarrow 
    \bigoplus _k H^{4k}(X ; \QQ  ) .
   \]
   
   Here are its principal properties on the category of compact spaces. For much more, see Grady and Sati \cite{GS} Section 2.1.

  \begin{Pro}  (cf. Hatcher \cite{Hatcher}, Section 4.1.)
  \begin{enumerate}
  \item The Pontrjagin character is a natural transformation
  \[
  \pont : KO^0(X) \longrightarrow \bigoplus _k H^{4k}(X ; \QQ  )
  \]
  and using suspension this may be extended to 
   \[
\pont : KO^*(X) \longrightarrow \bigoplus _k H^k (X ; \QQ  ).
  \]
  \item The Pontrjagin character is a homomorphism of graded commutative rings. Specifically, 
  \[
  \pont(E \oplus F) = Ph(E) + Ph(F)
  \]
  and
  \[
    \pont (E \otimes F) = Ph(E)Ph(F).
  \]
  \item The Pontrjagin character induces an isomorphism of graded commutative rings
   \[
 \pont: K^0(X)\otimes \QQ  \overset{\cong}\longrightarrow \bigoplus _k H^{4k}(X ; \QQ  ).
  \]
\end{enumerate}
\end{Pro}

 
   \begin{Rem}
    (Grady and Sati \cite{GS}, \S 2.1.)
 It's interesting to look at the simplest case, when $X = \pt$ (a one-point space).  Then $KO^{-k}(\pt) = KO_k(\RR)  $  which is torsion or zero except when $k = 0$ or $4$ mod $8$.
 So $KO^*(\pt)\otimes \QQ = \QQ$ when $* = 0, 4$ and is $0$ otherwise. The generator of $KO^0(\pt)$ is the identity of the ring, 
 and it is represented by a one-dimensional trivial bundle over the point. It goes to the generator of $H^0(\pt ; \ZZ) \cong \ZZ \subset H^0(\pt ; \QQ ) $ 
 since $Ph $ is a ring map. The group $KO^{-4}(\pt)$ is more interesting. We have a natural isomorphism 
 \[
 KO^{-4}(\pt) \cong KO^0(S^4)
 \]
 and then 
 \[
 \ZZ \cong \pi _4(S^4) \overset{\cong}\longrightarrow   KO^{0}(S^4 ) \overset{Ph}\longrightarrow  H^4(S^4 ; \QQ)
 \]
 and the image of the identity is twice the  generator of $H^4(S^4 ; \ZZ)$, as may be seen by pairing it with the fundamental class of $S^4$. 
 The upshot of this is that the map 
 \[
 KO^*(\pt) \otimes \QQ \overset{Ph}\longrightarrow \bigoplus _k H^{4k}(\pt ; \QQ ) 
 \]
 actually is integral, taking the form\footnote{One may reasonably ask why don't we see  $H^{-4}(\pt ; \ZZ )$ rather than 
 $ H^4(\pt ; \ZZ )$.  The answer is that we should! The problem is that ordinary 
 homology and cohomology were invented first- H. Poincar\'e's 1895 paper {\emph{Analysis Situs}}
 is a reasonable starting point, but he used Betti numbers rather than homology groups per se, and then in the 
 1920's and 1930's several people started to use our present notation. There is no reason to think that they were 
 troubled by pairing two classes of degree $n$ and winding up with a class of degree $0$.    } 
  \[
 KO^0(\pt) \overset{\cong}\longrightarrow H^0(\pt ; \ZZ ) \cong \ZZ 
 \qquad \text{and}\qquad
 KO^{-4}(\pt) \overset{\times 2}\longrightarrow H^4(\pt ; \ZZ ) \cong \ZZ   .
 \]
 
\end{Rem}
 
 Finally, we touch on Chern numbers.  \index{Chern number} We do this gingerly, since the name {\emph{Chern number}} is used in certain areas of physics for 
 pretty much any number that can be calculated via some sort of topological construct.  Here is the classical definition.  
 
 \begin{Def} Suppose that $M^n$ is an oriented connected manifold and the tangent bundle $TM$ can be given the structure of a complex vector bundle  with associated
  Chern classes $  c_j(TM) \in H^{2j}(M^n ; \ZZ)$.    Choose some polynomial 
  \[
  p = p(c_1(TM), \dots  c_n(TM) )
  \]
  where each term has total degree $n$.\footnote{For example, we could have 
  \[
  p = 17 c_1(TM)^4 + 37 c_2(TM)^2
  \]
  as each term has total degree $8$. }
     Let $[M] \in H_n(M ; \ZZ)$ 
  denote the fundamental class of the manifold.   
 Then evaluating $p \in H^n(M ; \ZZ) $ against the fundamental class 
 $[M] \in H_n(M ; \ZZ )$ (i.e. integrating $p$ over the manifold)     gives 
 a {\emph{Chern number}}  
 \[
 < p, [M] > \,\in\, \ZZ .
 \]
 More generally, we can replace the tangent bundle by any complex bundle over $M$ and define Chern numbers associated to that new bundle.
 \end{Def} 
 
 \newpage
 
      \section{\bf{Relating $KO_*$-theory for $C^*$-Algebras and $KO^*$-theory for Spaces}}
      \label{Section:RelatingKO-theory}
       
In this section, we prove theorems that connect the topological $KO$-theory of a compact or locally compact space $X$ to the $C \sp *$-algebraic $KO$-theory of the $C \sp *$-algebra $A = \real(X)$ or $A = \realo(X)$ respectively, and similarly in the complex setting. Some of these theorems are well-known among operator algebraists, but references are hard to find.  

\begin{Thm}\label{T:swanK}
Let $X$ be a compact topological space. Then:
\begin{enumerate}
\item
 There is a natural isomorphism
\[
       KO_0(\real(X))   \cong           KO^0(X) .
\]
Similarly, if $X$ is a locally compact topological space, then there is a natural isomorphism
\[
 KO_0(\realo(X))    \cong     KO^0(X) .
\]

\item
 There is a natural isomorphism
\[
       K_0(\complex(X))   \cong           K^0(X) .
\]
Similarly, if $X$ is a locally compact topological space, then there is a natural isomorphism
\[
  K_0(\complexo(X))    \cong     K^0(X) .
\]
\end{enumerate}

\end{Thm}

 This theorem is an immediate consequence of a theorem of R.G. Swan\footnote{CS: Swan
 was on faculty at the U. of Chicago where I studied. He was very shy - would time his lectures 
 so that he was standing by the exit door at the end.  One day in my homological algebra class the instructor, Saunders MacLane, got stuck on a proof. The class ended and he told us he would finish the proof next time.  At the next class MacLane said that, being still stuck, he went for help to the person that everyone in the Department went to when they needed help - to Dick Swan. Swan gave him the answer. }  that we mentioned previously:

 \begin{Thm} \label{T: swan} (R. G. Swan \cite{Sw})   Let $X$ be a compact space and let $\FF $ denote either the real or the complex numbers.\footnote{Swan also includes 
 the quaternions.}
 \index[notation]{theta@$\Theta_X$}
  Then there is  a natural isomorphism  
 \[
 \Theta = \Theta_X: \mathcal{P}_\infty (C^\FF (X))\simeq    \rm{Vect}^\FF (X) 
 \]
 between   $\rm{Vect}^\FF  (X)$, isomorphism classes of $\FF$ - vector bundles over $X$, and     $\mathcal{P}_\infty (C^\FF (X) )$, projections in matrix rings over $C^\FF (X)$. 
 
   \end{Thm} 
 
 Swan's introduction is worth quoting:
 
 ``{\emph{Serre \dots  has shown that there is a one-to-one correspondence between algebraic vector bundles over an affine variety and 
 finitely generated projective modules \index{projective module} over its coordinate ring. For some time, it has been assumed that a similar correspondence exists between topological 
 bundles over a compact Hausdorff space $X$ and finitely generated projective modules over the ring of continuous real-valued functions on $X$... However, 
 no rigorous treatment of the correspondence seems to have been given. I will give such a treatment here... }}"

 Theorem \ref{T: swan} implies Theorem \ref{T:swanK} by simply applying the Grothendieck construction.
 If $\FF = \RR$ then 
\[
\rm{Groth}( \mathcal{P}_\infty(A) ) = KO_0(A).  
\]
and
\[
\rm{Groth}( \rm{Vect}\pr (X)) = KO^0(X).
\]
Then use the fact that the Grothendieck construction is natural with respect to homomorphisms of abelian semigroups.

 There are several nice proofs of the complex case of these theorems, but proofs of the real case are scarce.  We shall supply a new proof of the 
 real case of Theorem \ref{T: swan} which is somewhat more analytic than Swan's proof. 
 
 We first focus on the case with $X$  compact.
 For any projection 
 \[
 p \in \mathcal{P}_n(A) \subseteq C(X, M_n(\RR)),
 \]
  let
\begin{align*} E_p = {\rm range~} p &= \{ (x, v) \in X \times \RR^n \mid v \in p(x)  \RR^n \} \\
&= \{ (x, v) \in X \times \RR^n \mid p(x) v = v \} \; .
\end{align*}
Then $E_p$ is a real vector bundle over $X$ with associated projection map $\pi \colon E_p \rightarrow X$ given by $\pi(x,v) = x$.
It is clear that $E_p \cong E_{p \oplus 0}$ where
\[
p \oplus 0 = \diag(p,0) \in \mathcal{P}_{n+1}(A).
\]

 \begin{Lem}\label{T:lemma123}
Let $X$ be a compact space and let $A = \real(X)$. If $p,q \in \mathcal{P}_n(A)$ are unitarily equivalent, then there is an isomorphism of real bundles $E_p \cong E_q$.
 \end{Lem}
 
 \begin{proof}
    Let 
\begin{align*} E_p &= \{ (x, v) \in X \times \RR^n \mid p(x) v = v \} \\
E_q &=  \{ (x, v) \in X \times \RR^n \mid q(x) v = v\} \\
\end{align*}
and let $u \in C(X, M_n(\RR))$ be a unitary satisfying $u p u^* = q$.
There is an homomorphism  $\phi : E_p \to E_q$ defined by $\phi  (x, v) = (x, u(x) v)$.
Indeed, we check that if $p(x) v = v$ then 
\[
q(x) \cdot u(x) v = u(x) p(x) v = u(x)v.
\]
 The inverse is defined using $u^*(x)$.
Furthermore, $\phi$ clearly commutes with the projection map $\pi$.
Thus $E_p$ and $E_q$ are isomorphic objects in the category of real vector bundles. 
\end{proof}

Let $A = \real(X)$. It follows that the construction $p \mapsto E_p$ gives a well-defined map
\[
\Theta \colon \mathcal{P}_\infty(A) \rightarrow \rm{Vect}\pr (X)
\] 
and it is easy to show that this is a semigroup homomorphism, that is
\[
\Theta (p \oplus q) \cong \Theta(p) \oplus \Theta(q).
\]

\begin{Lem}\label{T:Lemma124} Suppose that $X$ is a compact space. Then the map $\Theta $ is an isomorphism.
\end{Lem}

\begin{proof}
To prove that $\Theta$ is injective, assume that $p,q$ in $M_n(A)$ are projections that satisfy $E_p \cong E_q$. So there exists $a \in M_n(A) = C(X, M_n(\RR))$ where $a_x$ is an isomorphism from $p(x) \RR^n$ to $q(x) \RR^n$ for each $x$
(and $a_x$ vanishes on $p^\perp$).
Thus we have $ap = a$ and $qa = a$. Using  continuous functional calculus, let $s = f(a^* a) a^*$, where
\[ f(t) = \begin{cases} 0 &  t = 0 \\ t^{-1/2} & t > 0 . \end{cases} \]
Since $X$ is compact, the non-zero part of the spectrum of $a$ is bounded below. So $f$ is continuous on $\rm{sp}(a)$ and $s \in M_{n}(A)$. Then it can be checked that $s$ is a partial isometry that satisfies $s s^* = q$ and $s^* s = p$.

Now let
\[ u = \begin{bmatrix} s & 1-q \\ 1-p & s^* \end{bmatrix}  \;  \]
in $M_{2n}(A) = C(X, M_{2n}(\RR))$.
Check that $u$ is a unitary and that it satisfies
\[ u \begin{bmatrix} p & 0 \\ 0 & 0 \end{bmatrix} u^* = \begin{bmatrix} q & 0 \\ 0 & 0 \end{bmatrix} \; . \]
Thus $\diag(p, 0)$ is unitarily equivalent to $\diag(q, 0)$ in $M_{2n}(A)$, so $[p] = [q]$ in $\rm{Proj}(A)$. Thus $\Theta$ is injective.

To show that $\Theta$ is surjective, let $E$ be any locally trivial vector bundle with projection map $\pi \colon E \rightarrow X$.
Since $E$ is locally trivial, there is a finite open cover $\{U_i\}_{i = 1}^k$ of $X$ such that for each $i  $ the restriction $E_{U_i}$ is isomorphic to $U_i \times \RR^n$. Let $\phi_i \colon E_{U_i} \rightarrow U_i \times \RR^n$ be this isomorphism.
Now, let $\{\psi_i\}_{i = 1}^k$ be a partition of unity subordinate to $\{ U_i\}$, a collection of continuous functions from $X$ to $[0,1]$ such that
$\sum_{i=1}^k \psi_i = 1$ and $\rm{supp}(\psi_i) \subseteq U_i$.
Then the function 
\[ f(v) = (\pi(v),  \psi_1(v) \phi_1(v) , \dots, \psi_k(v)  \phi_k(v) ) \]
defines an injective map of vector bundles from $E$ to $X \times \RR^{nk}$.
Using the standard inner product on $\RR^{nk}$, for each $x \in X$, let $p_x$ be the orthogonal projection from $\RR^{nk}$ to the range of $f$ over $x$. Then $p_x$ is a continuous function of $x$ and is a projection in the algebraic sense that $p^2 = p$ and $p^* = p$. 
Thus $[p] \in \mathcal{P}_\infty(A)$ and by construction we have $\Theta([p]) \cong [E]$.
 \end{proof}
 
 \begin{proof}[Proof of Theorem~\ref{T:swanK}.]  
 When $X$ is compact and $\FF = \RR$, this is Lemma~\ref{T:Lemma124}. Now, suppose that $X$ is locally compact. Let $X^+ = X \bigcup \{ +\}$ be the one-point compactification. As in Definition~11.3 we have
\begin{equation} \label{ker1} KO^0(X) = \ker \left[ KO^0(X^+) \xrightarrow{\iota^*} KO^0(+) \right]  \;  \end{equation}
where the map is induced by the inclusion map $\iota \colon \{+\} \hookrightarrow X^+$.
On the other hand, if we set $A = \realo(X)$, then $A$ is a non-unital real $C \sp *$-algebra 
and from Definition~6.5 we have
\begin{equation} \label{ker2}  KO_0(A) = \ker \left[ KO_0(A^+) \xrightarrow{\pi_*} KO_0(\RR) \right] \; \end{equation}
where the map is induced by the canonical map $\pi \colon A^+ \rightarrow  A^+/A \cong \RR$ associated with a unitization.

In this situation, we have an isomorphism $A^+ \cong \real(X^+)$ and under this isomorphism the map $\pi$ corresponds to the evaluation map
\[ \rm{ev} \colon \real(X^+) \rightarrow \RR \; , \]
specifically evaluation at the point $+$. This evaluation map is dual to the inclusion map $\iota$. This shows that the maps shown in Equations~\ref{ker1} and \ref{ker2}  are identical and hence there is an isomorphism $KO^0(X) \cong KO_0(A)$.
 
Finally, the proof in the complex case (for $X$ compact and locally compact) is exactly the same as the proof above, including the proofs of Lemmas~\ref{T:lemma123} and \ref{T:Lemma124}, replacing projections in $M_n(\real(X))$ with projections in $M_n(\complex(X))$, replacing real vector bundles with complex vector bundles, and replacing the semigroup $\rm{Vect}\pr(X)$ with $\rm{Vect}\pc(X)$.
\end{proof}

 \begin{Cor} Suppose that $X$ is a connected compact space.  If $i > 0$, then  \begin{enumerate}
 \item  
 \[
  KO^{-i}(X)  \cong  KO^0( \topr ^i \times X) \;
  \]
 understood as a subgroup of $K$-theory of the compactification
 \[
 KO^0( ( \topr ^i \times X)^+).
 \] 
 \item 
 In particular, $KO^0( \topr ^i \times X)$ is the subgroup of $KO^0(( \topr ^i \times X)^+)$ generated by formal differences of vector bundles $[E] - [F]$ over
$ ( \topr ^i \times X)^+$ where $E$ and $F$ have the same dimension.
\item
$KO_0( \realo( \topr ^i \times X))$ is the subgroup of $KO_0( \real(( \topr ^i \times X)^+))$ generated by formal differences of projection classes $[p] - [q]$ where $\ev_*[p] \cong \ev_*[q]$ where $\ev \colon \real( (\topr^i \times X)^+) \rightarrow \RR$ is point evaluation.
 
 \item
We can deal with negative degrees by writing (again for $i > 0$),
\[ 
KO^{-i}(X)    
	= KO^0( \topr ^{i} \times X, \hat \tau) \; 
	\]
  where $\hat{\tau}$ is the involution on $\topr ^i \times X$ given by $(t_1, \dots, t_i, x)  \mapsto (-t_1, \dots, -t_i, x)$.
  \item
In a similar way to the previous paragraph, this $K$-theory group can be expressed as a subgroup of 
$KR^0( ( \topr ^i \times X, \tau)^+)$ consisting of differences $[(E, \sigma_1)] - [(F, \sigma_2)]$ where $E$ and $F$ have the same dimension; or as the kernel of $\ev_*$ in $KO_0(A)$ where 
$$A = \{ f \in \complex(( \topr ^i \times X)^+) \mid f(t_1, \dots, f_i, x) = \overline{ f(-t_1, \dots, -t_i, x) } \} \; . $$
\end{enumerate}
\end{Cor}
 
 Note that in Part (5) above, ``$KR^0(~)$" refers to Atiyah's ``Real $K$-theory" which is described in Section~\ref{Section:KRTheory}.
 
 \begin{proof}
 The first statement follows from Swan's Theorem:
 \[
  KO^{-i}(X) = KO_i( \real(X)) = KO_0( S^i \real(X))
	= KO_0( \realo(\topr^i \times X)) = KO^0( \topr^i \times X). \; 
	\]
 In the negative degree setting we have, again from Swan's Theorem,
 \[
  KO^{-i}(X) = KO_{i}( \real(X)) = KO_0( S^{-i} \real(X))
	= KR^0( \topr ^{i} \times X, \hat \tau) \;
	\]
 and the other results follow.
 \end{proof}
   
\begin{Rem} We have discussed topological $K$-theory for real and complex vector bundles.  There are many other families of compact \index{Lie group}   
Lie groups, such as $Spin(n)$, $Spin ^c(n)$, $Sp(n)$, etc. with natural actions on $\RR ^n$,  $\CC ^n$, or $\HH ^n$.   Let $G$ be one of them, and assume 
it acts on $\RR ^n$.  Then we may form 
$G$-vector-bundles by starting with principal $G$-bundles $P \to X$ (with $X$ compact for safety) and forming $V = P \times _G \RR ^n$. 
The set of isomorphism classes of such $G$-bundles is closed under direct sum, and so we may form a Grothendieck group to obtain 
$KG^0(X)$.   Some of these groups turn out to be very useful indeed, but they take us too far away from our main object of interest.  To get 
a taste of what is possible, check out Stong \cite{Stong} Chapter XI, which has the great title {\emph{``{$Spin$, $Spin ^c$, and Similar Nonsense."}}} 
Stong is writing about cobordism, but the two subjects are strongly related in some cases.\footnote{Stong was an unusual person, to put it 
mildly.  He would get to his office around 3 AM, work diligently until 8 or 9 AM, and then lounge in the coffee room for the rest of the day, 
so it looked like he was doing nothing.  He is well-known, among other things, for the Hattori-Stong theorem.  J.F. Adams once commented that 
it should have been called the Stong-Hattori theorem, since Stong proved it first, but people called it Hattori-Stong because they liked the emphasis 
implied by saying ``Stong" second.  }
\end{Rem}

\begin{Rem} Suppose that $G$ is a topological group (compact for simplicity).  Then we may focus attention on a  compact space $X$ for 
which there is a natural action $G \times X \to X$. Call $X$ a {\emph{$G$-space}}.  \index{G-space@$G$-space} Suppose then given  (real or complex) vector bundles $V \to X$ where $G$ acts on $X$ and the projection map commutes with the $G$-action.  Then we get $G$-vector bundles.  Taking Grothendieck 
groups again we obtain
 {\emph{equivariant}} $K$-theory \index{equivariant $K$-theory} \index[notation]{KGX@$K^G(X)$}
  $K^G(X)$.\footnote{We can do this more generally for operator algebras but this is 
abstract enough for a remark.}     These groups turn out to be computable (sometimes) and quite interesting (always).  The cases $G$ finite and $G$ 
compact connected turn out to be quite different, and we refer the interested reader to N. C. Phillips \cite{Phillips} and to Rosenberg-Schochet \cite{RS} 
for examples of their use. 
\end{Rem}

\newpage
      \section{\bf{Atiyah's $KR$-theory}}
      \label{Section:KRTheory}

Now we consider $K$-theory in the category of topological spaces with involution, a cohomology theory invented by M. Atiyah \cite{AtiyahR}.  We give a new proof 
of his result that it can be 
realized as $KO$ of an appropriate real $C^*$-algebra. This should be considered as an extension of Swan's Theorem to the case of topological spaces with involution.

  Atiyah described the birth of $KR$-theory 
as follows:

{\emph{
In attempting to understand reality questions of elliptic operators Singer and I were, for a long time, held up by the fact that real $K$-theory behaves differently from complex
$K$-theory. Eventually, Singer pointed out that a different notion of reality was needed. The essential point is that the Fourier transform of a real-valued function is not 
real but instead satisfies the relation $f(-x) = {\overline {f(x)} } $.  It was therefore necessary to develop a new version of real $K$-theory for spaces with involution. This was carried 
out in \cite{AtiyahR}. In fact $KR$-theory, as I called it, turned out to provide a better approach even for purely topological purposes. (\cite{AtiyahColVol2}, p. 4) }}

An object in this category is a pair $(X, \tau)$ where $X$ is a topological space (compact or locally compact) with involution  $\tau$. 
(This is what Atiyah refers to  as a ``Real space." We will not use that terminology.) 
Given such a pair $(X, \tau)$ when $X$ is compact, we will define a  {\emph{complex vector bundle with involution}} 
\index{vector bundle!with involution}
to be a pair $(E, \sigma)$ where $E$ is a complex vector bundle over $X$ and $\sigma$ is a conjugate-linear involution on $E$ that satisfies $\tau( \pi(e)) = \pi (\sigma(e))$ for all $e \in E$.

Then, again following Atiyah, we define
$\rm{Vect}(X, \tau) $ 
\index[notation]{vectXt@$\rm{Vect}(X, \tau)$}
to be the abelian semigroup of isomorphism classes of 
complex vector bundles with involution over $(X, \tau)$ and
\[ KR^0(X, \tau) = \rm{Groth}( \rm{Vect}(X, \tau)). \]
\index[notation]{KRX@$KR^0(X, \tau)$}
This is a contravariant functor on the category of topological spaces with involution.
We assert that $KR^*$-theory for spaces with involution can be expressed in terms of $KO_*$- theory for an associated 
real $C^*$-algebra. Here is the precise statement.

\begin{Thm}\label{T:AtiyahR} (Atiyah)
\begin{enumerate}
\item
Let $(X, \tau)$ be a compact topological space with involution.  Define a real $C^*$-algebra 
\[
A = A(X, \tau ) =  \{ f \in \complex(X) \mid f(\tau(x)) = \overline{f(x)} \} . \; 
\]
Then there is a natural isomorphism
\[
\Theta = \Theta_{X,\tau} :     KO_0(A)     \overset{\cong}\longrightarrow     KR ^0(X, \tau) .
\index[notation]{thetaXt@$\Theta_{X,\tau}$}
\]
\item  Let $(X, \tau)$ be a locally compact topological space with involution. Define
\[
A = A(X, \tau) =  \{ f \in \complex_0(X) \mid f(\tau(x)) = \overline{f(x)} \} \; .
\]
Then there is a natural isomorphism
\[
\Theta:     KO_0(A)     \overset{\cong}\longrightarrow     KR ^0(X, \tau) .
\]
\end{enumerate}
  \end{Thm}

One should compare Theorem 4.1 of \cite{Knumber102} (Karoubi) for a related result.
Before beginning the proof, we establish a basic property of complex involutive bundles over spaces with involution. 

Note that the trivial bundle over $(X, \tau)$ of dimension $n$ is given by $(X \times \CC^n, \tau \times \psi_n)$ where $\psi_n$ is pointwise conjugation on $\CC^n$.

\begin{Lem} \label{sub-trivial}
Let $(X, \tau)$ be a compact space with involution and let $(E, \sigma)$ be a complex vector bundle over $X$ with a conjugate linear involution satisfying $\sigma(E_x) = E_{\tau(x)}$. Then $(E, \sigma)$ is isomorphic (as a complex vector bundle with involution) to a sub-bundle $(E', \sigma)$ of the trivial bundle $(X \times \CC^n, \tau \times \psi_n)$ for some $n$.  Further, there is a subbundle 
$(F, \rho )$ of the trivial bundle such that 
\[
(E', \sigma ') \oplus (F , \rho) \cong (X \times \CC^n, \tau \times \psi_n)
\]
as bundles with involution.
\end{Lem}

\begin{proof}
First we work locally with the claim that for all $x_0 \in X$ there exists a neighborhood $\mathcal{U}$ of $x_0$ in $X$ such that $\tau(\mathcal{U}) = \mathcal{U}$ and $(E_{\mathcal{U}}, \sigma)$ is isomorphic to $(\mathcal{U} \times \CC^n, \tau \times \psi_n)$

Let $x_0 \in X$ be given and let $\mathcal{V}$ be a neighborhood of $x_0$ on which $E$ is trivial. 
If $\tau(x_0) \neq x_0$ we can arrange without loss of generality that $\mathcal{V} \cap \tau( \mathcal{V} ) = \emptyset$. Let $\phi \colon E_{\mathcal{V}} \rightarrow \mathcal{V} \times \CC^n$ be an isomorphism; we use the notation $\phi_x \colon E_x \rightarrow \CC^n$ for the map on the fibres. Let $\mathcal{U} = \mathcal{V} \bigcup \tau(\mathcal{V})$ and extend $\phi$ to 
$$\phi \colon E_\mathcal{U} \rightarrow \mathcal{U} \times \CC^n $$ 
by defining
$$\phi(x, v) = (x,  \psi (\phi_{\tau(x)}( \sigma_{x} (v)))) \; \quad \text{for $x \in \tau(\mathcal{V})$.}$$
Now $\phi \colon E_\mathcal{U} \rightarrow \mathcal{U} \times \CC^n$ is defined and satisfies
$(\tau \times \psi_n) \circ \phi = \phi \circ \sigma$.

The other case to consider is if $\tau(x_0) = x_0$. Now, we can arrange so that $\mathcal{U}$ is a neighborhood of $x_0$ on which $E$ is trivial and such that $\tau(\mathcal{U}) = \mathcal{U}$ (but it may not be the case that $\tau = \id$ on $\mathcal{U}$). Let $\phi \colon E_\mathcal{U} \rightarrow \mathcal{U} \times \CC^n$ be a vector bundle isomorphism, but which may not carry the involution on $E$ to the standard involution (conjugation) on $\mathcal{U} \times \CC^n$.
Thus, let $\widetilde{\sigma}_x$ be the involution on $\CC^n$ induced from $\sigma_x$ by $\phi$. In other words, $\phi$ is a bundle-with-involution isomorphism between
$(E_\mathcal{U}, \sigma)$ and $(\mathcal{U} \times \CC^n, \tau \times \widetilde{\sigma} )$. It remains to show that there is an isomorphism $\theta$ from $(\mathcal{U} \times \CC^n, \tau \times \widetilde{\sigma} )$ to $(\mathcal{U} \times \CC^n, \tau \times \psi_n)$.

Since $\widetilde{\sigma}_x$ is conjugate-linear and depends continuously on $x$, then $\psi_n \circ \widetilde{\sigma}$ is linear and can be represented an $n \times n$ matrix $A_x$ depending continuously on $x$. Thus we can write 
\[
\widetilde{\sigma}_x = \psi_n \circ A_x \; ,
\]
conflating the matrix $A_x$ with the linear transformation on $\CC^n$ that it represents.
The condition that $\widetilde{\sigma}$ is an involution implies that ${A_x } \overline{A_{\tau(x)}} = I_n$ for all $x \in \mathcal{U}$. 
Indeed, $(\widetilde{\sigma})^2 = \id$ can be rewritten as 
\[
\psi_n A_{\tau (x)} \psi_n A_{x} v = v
\qquad \text{for} \qquad  v \in \CC^n
\]
and then again as $\overline{A_{\tau(x)}} A_x v = v$, which is equivalent to 
\[
\overline{A_{\tau(x)}} A_x = I_n \; . 
\]

Furthermore, we can assume without loss of generality that $A_{x_0} = I_n$, by use of the isomorphism 
$$\id \times A_{x_0}^{-1} \colon (\mathcal{U} \times \CC^n) \rightarrow (\mathcal{U} \times \CC^n) \; .$$

The isomorphism $\theta$ that we are looking for is from $\mathcal{U} \times \CC^n$ to itself, so it can be represented by a continuous matrix-valued function $B_x$ for $x \in \mathcal{U}$. We need $\theta$ to carry the involution $\widetilde{\sigma}$ to the involution $\psi_n$ on $\CC^n$. In other words, we need the diagram
\[ \xymatrix{
\mathcal{U} \times \CC^n \ar[r]^ \theta \ar[d]^{ \widetilde{\sigma}}
& \mathcal{U} \times \CC^n  \ar[d]^{\tau \times \psi_n}\\
\mathcal{U} \times \CC^n \ar[r]^\theta
& \mathcal{U} \times \CC^n 	} \]
to commute, which can be expressed as 
\[
(B_{\tau(x)} \circ \psi_n \circ A_x)(v)  = (\psi_n \circ B_x)(v).
\]
This condition is equivalent to matrix equations $B_{\tau(x)} \overline{A_x} = \overline{B_x}$ for all $x \in \mathcal{U}$.

There is a continuous square root defined on square matrices in a neighborhood of the identity. We replace $\mathcal{U}$ with a smaller neighborhood to insure that square roots are defined on $A_x$ for $x \in \mathcal{U}$, and let $B_x = (A_x)^{1/2}$. Then using the formula
$A_{\tau(x)} = (\overline{A_x})^{-1}$ we have
$B_{\tau(x)} =  (\overline{A_x})^{-1/2}$
so
$$ {B_{\tau(x)}} \overline{A_x} = (\overline{A_x})^{-1/2} \overline{A_x} = (\overline{A_x})^{1/2} = \overline{B_x}$$
as desired.

Now, using a partition of unity, these local maps can be patched to form an injective homomorphism of vector bundles-with-involution  
\[
(E, \sigma) \longrightarrow  (X \times \CC^m, \tau \times \psi_n)
\]
 for some large enough integer $m$, as in the proof of Lemma~\ref{T:Lemma124}.

 Finally, take  $(E, \sigma) = (E, \tau \times  \psi^n)$  to be   a subbundle of the trivial bundle. Then take $(F, \sigma) = (F, \tau \times \psi^n)$  to be the point-wise orthogonal complement to $E$. Check easily that since vector conjugation on $\CC ^n$  respects the inner product, we have $(\tau \times  \psi^n)(F) = F$.
Then we clearly have 
   \[
   (E, \tau \times \psi^n) \oplus (F, \tau \times  \psi^n) = (X \times \CC^n, \tau \times \psi^n) 
\]
completing the proof.
\end{proof}
 
The rest of this section is devoted to a new proof of Theorem~\ref{T:AtiyahR}, using Lemma~\ref{sub-trivial}.  Let $(X, \tau)$ be a compact space with involution and let 
$A = \real(X, \tau )$.
We shall show that there is a natural  isomorphism of abelian semigroups 
 \[
 \Theta \colon \mathcal{P}_\infty (A) \rightarrow \rm{Vect}(X, \tau)
 \]
which implies that there is an induced natural isomorphism of groups  
\[
\Theta \colon KO_0(A) \rightarrow KR^0(X, \tau).   
\]

\begin{proof}[Proof of Theorem  \ref{T:AtiyahR}]

%
%

First let $p \in M_n(A)$ be a projection. Then $p \in C(X, M_n(\CC))$ and $p$ satisfies
$$p(\tau(x)) = \overline{p(x)}  \; $$ 
where the conjugation is entry-wise conjugation on matrices.
Let
\begin{align*} E = {\rm ran~} p &= \{ (x, v) \in X \times \CC^n \mid v \in p(x) \cdot \CC^n \} \\
&= \{ (x, v) \in X \times \CC^n \mid p(x) v = v \} 
\end{align*}
with the associated projection map $\pi \colon E \rightarrow X$ by $\pi(x,v) = x$.
Now $E$ is a complex vector bundle and has a conjugate-linear involution $\sigma = \tau \times \psi_n$ given by $(x,v) \mapsto (\tau(x), \overline{v})$.
This involution $\sigma = \tau \times \psi_n$ certainly commutes with the respective involutions on $X$ as required. 
To show that $(\tau(x), \overline{v}) \in E$ we
check that 
\[ p(\tau(x)) \overline{v} = \overline{p(x)} \, \overline{v} 
	= \overline{p(x) v} 
	= \overline{v}  \;   \; . \]
Hence $(E, \sigma) \in \rm{Vect}(X, \tau)$.    
We define $\Theta (p) = (E, \sigma)$.

We must show that $\Theta $ is well-defined on $\mathcal{P}_\infty$.
Suppose that $p$ and $q$ are projections in $M_n(A)$ and
that $u \in M_n(A)$ is a unitary that satisfies $upu^* = q$. 
Let 
\begin{align*} E_p &= \{ (x, v) \in X \times \CC^n \mid v \in p(x) v = v \} \\
E_q &=  \{ (x, v) \in X \times \CC^n \mid v \in q(x) v = v\} \\
\end{align*}
As in the proof of Lemma \ref{T:lemma123}, there is a vector bundle homomorphism  
\[
\phi: E_p \to E_q
\]
 defined by 
 \[
 \phi  (x, v) = (x, u(x) v).
 \]
The bundle homomorphism $\phi$ commutes with the involution $\sigma$:
$$ \phi( \sigma(x,v)) = \phi( \tau(x), \overline{v}) = 
	( \tau(x), u( \tau(x))  \overline{v}) = (\tau(x), \overline{u} \;  \overline{v} ) = \sigma( x, u(x) v) = \sigma( \phi(x,v)) \; .$$
It follows that $E_p$ and $E_q$ are isomorphic objects in the category of 
complex vector bundles with involution over $(X, \tau)$.

Furthermore, it is straightforward to see that if $p$ is a projection in $M_n(A)$, and if $p' = p \oplus 0 \in M_{n+m}(A)$, then $E_p \cong E_{p'}$.
More generally, if $p \in M_m(A)$ and $q \in M_n(A)$ are projections, then $p \oplus q$ is a projection in $M_{m+n}(A)$ 
and there is an isomorphism
$$E_{p \oplus q} \cong E_p \oplus E_q \; .$$
Thus the construction $p \mapsto (E, \sigma)$ is well-defined and gives a homomorphism of abelian semigroups 
\[
\Theta \colon \mathcal{P}_\infty (A) \rightarrow \rm{Vect}(X, \tau) .
\]

Next, we claim that $\Theta$ is injective. Assume that $p,q$ are projections in $M_n(A)$ that satisfy $E_p \cong E_q$ (again as complex vector bundles with involution over $(X, \tau)$). Then there exists 
an element $a \in  C(X, M_n(\CC))$ that implements this isomorphism (as in the first paragraph of the proof of Lemma \ref{T:Lemma124}). So $a_x$ is an isomorphism from $p(x) \CC^n$ to $q(x) \CC^n$ for each $x \in X$.
Furthermore, since the isomorphism commutes with the involutions on $E_p$ and $E_q$, we have
$a_{\tau(x)}( \overline{v}) = \overline{ a_{x}(v)}$ for all $x \in X$ and all $v \in \pi^{-1}(x)$.
This implies that $a_{\tau(x)}( \overline{v}) = \overline{ a_{x}} \,\overline{v}$ and hence
$a_{\tau(x)} = \overline{a_x}$ for all $x \in X$. Therefore $a \in M_n(A)$.

As in the proof of Lemma \ref{T:Lemma124}, we let $s = f(a^* a) a^*$ where
\[ f(t) = \begin{cases} 0 &  t = 0 \\ t^{-1/2} & t > 0  \end{cases} \]
noting in particular that because $a \in M_n(A)$, we also have $s \in M_n(A)$.
As before, $s$ is a partial isometry that satisfies $s s^* = p$ and $s^* s = q$.
Since \[ u = \begin{bmatrix} s & 1-q \\ 1-p & s^* \end{bmatrix}  \;  \]
is in $M_{2n}(A)$
and since
\[ u \begin{bmatrix} p & 0 \\ 0 & 0 \end{bmatrix} u^* = \begin{bmatrix} q & 0 \\ 0 & 0 \end{bmatrix} \; , \]
it follows that $\diag(p, 0)$ is unitarily equivalent to $\diag(q, 0)$ in $M_{2n}(A)$, so $[p] = [q]$ in $\rm{Proj}(A)$.
Thus $\Theta $ is injective.

Next, we show that $\Theta$ is surjective.
Let $(E, \sigma)$ be any complex vector bundle with involution over $(X, \tau)$.
Using Lemma~\ref{sub-trivial},  identify $E$ with a sub-bundle of  
$X \times \CC^n$ in such a way that the involution $\sigma$ on $E$ corresponds to conjugation on $\CC^n$. In other words, $\sigma$ can be expressed by the formula
$$\sigma(x,v) = (\tau(x), \overline{v})$$
where $\overline{v} = \psi_n(v)$ is component-wise conjugation on $\CC^n$.
Define $p = p(x) \in C(X, M_n(\CC))$ by letting $p(x)$ be the orthogonal projection of $\CC^n$ to $E_x$. Then $p$ is continuous in $x$, and is a 
projection: $p^2 = p^* = p$. 

Now, let $v \in E_x$. This implies that $\overline{v} \in E_{\tau(x)}$. Therefore we have
$p(x) \cdot v = v$ and $p({\tau(x)}) \cdot \overline{v} = \overline{v}$. This implies also that
$$\overline{p({\tau(x)})} \cdot v = \overline{p{(\tau(x))} \cdot \overline{v}}  = \overline {\overline{v}} = v \; .$$
As $\overline{p(\tau(x))}$ and ${p(x)}$ are both orthogonal projections with the same rank which agree on $E_x$, they must be equal, implying that $\overline{p(\tau(x))} = {p(x)}$
and hence $p(\tau(x)) = \overline{p(x)} \; $.
From this it follows that $p$ is an element in $M_n(A)$. 
By construction we have $(E, \sigma) \cong E_p$, which shows that $\Theta$ is surjective.
 
Thus we have established that  if $X$ is compact then there is an isomorphism   
 of abelian semigroups 
 \[
 \Theta \colon \mathcal{P}_\infty (A) \rightarrow \rm{Vect}(X, \tau)
 \]
which implies that there is an induced group isomorphism 
\[
\Theta \colon KO_0(A) \rightarrow KR^0(X, \tau).   
\]

It remains to establish that $\Theta $ is natural with respect to maps of spaces with involution. 
Suppose that
\[
f \colon (X_1, \tau _1 )  \rightarrow (X_2, \tau _2 ) 
\]
is a continuous homomorphism of compact spaces with involution. Then there is a $C \sp *$-algebra homomorphism 
 \[
 f^* : A_2 = A(X_2, \tau _2 ) \longrightarrow A_1 = A(X_1, \tau _1) 
 \]
 with associated map 
 \[
 f^*: \mathcal{P}_\infty (A_2)  \longrightarrow  \mathcal{P}_\infty (A_1) .
 \]
 Similarly, the pullback construction induces a natural map 
 \[
 f^*: \rm{Vect}(X_2, \tau _2)    \longrightarrow        \rm{Vect}(X_1, \tau _1)
 \]
 and it is easy to check that the diagram
 \[
 \begin{CD}
  \mathcal{P}_\infty (A_2)       @>{\Theta}>>        \rm{Vect}(X_2, \tau _2)  \\
 @VV{f^*}V                  @VV{f^*}V      \\
  \mathcal{P}_\infty (A_1)         @>{\Theta}>>          \rm{Vect}(X_1, \tau _1)
 \end{CD}
 \]
 commutes.  This implies immediately that the diagram 
 \[
 \begin{CD}
 KO_0(A_2)      @>{\Theta}>>       KR^0(X_2, \tau _2)  \\
 @VV{f^*}V                  @VV{f^*}V      \\
 KO_0(A_1)      @>{\Theta}>>       KR^0(X_1, \tau _1)  
 \end{CD}
 \]
 commutes, establishing naturality.  This completes the proof of Theorem \ref{T:AtiyahR} in the case where $X$ is compact.

The proof in the locally compact case follows the same strategy as that in the proof of Theorem \ref{T: swan}. Note that if $X$ is locally compact then any involution $\tau$ on $X$ extends naturally to an involution $\tau^+$ on $X^+$, which satisfies $\tau^+( +) = +$. Thus
$KR^0(X, \tau)$ is defined to be 
\[ 
KR^0(X, \tau) = \ker \left[ KR^0(X^+, \tau^+) \xrightarrow{\iota_*} KR^0(+, \id) \right]  \;  .
 \]
 We leave the details to the reader.
 \end{proof}

\begin{Rem}
Since $(\real(X, \tau))$ is commutative, there is a product structure on 
 $KO_*( \real (X, \tau ))$. As far as we know this product structure has not been investigated, except in the classical case with $\tau = \id$. If we identify $KO_*(\real (X, \tau ) ) $ with 
 Atiyah's Real $K$-theory $KR^*(X )$, then there is a natural coproduct
 \[
 (X, \tau ) \longrightarrow (X, \tau ) \times (X , \tau ) 
 \]
 that induces a product structure on $KO_*(\real (X, \tau )) \cong KR^*(X )$ exactly analogous to the product structures on $KO^*(X)$ and 
 $K^*(X)$.   This needs more attention than we can give it at this time.
 \end{Rem}
 


 \newpage
\part{Clifford Algebras and Symmetric Spaces} \index{Clifford algebra}

 
 At first glance, there is no particular reason to suspect that a space such as $\Omega ^j(\oo )$ has any special relationship to physics.  It turns out that  
 these spaces have a close relationship with certain symmetric spaces, and in fact the eight spaces in the real Bott Periodicity Theorem and the two in the Complex Bott Periodicity Theorem correspond  to the renowned Ten-Fold Way in physics.  We can see this connection to symmetric spaces best via Clifford algebras\index{Clifford algebra}. So in Part 4 we will introduce Clifford algebras and show their connection to Bott Periodicity and to symmetric spaces.   Then
  we will relate symmetric spaces to physics directly.
 
 \vglue .5in
        \section{\bf{Clifford Algebras and Their Associated Lie Groups}} \index{Clifford algebra|(}
 \label{Section:CliffordAlgebras}
 
 \index{Lie group}   
  The classic reference for Clifford algebras and their connection to $K$-theory  is Atiyah, Bott, and Shapiro \cite{ABS}.  We refer to this paper and these authors by ABS.  Husemoller \cite{H} has a fine treatment. For connections with physics we rely on Baez \cite{Baez}, Stone et.al. \cite{Stone} and G. Moore \cite{Moore}.  Notation for Clifford algebras varies a lot. When in doubt, please consult the chart below, which gives our notation.
  
\begin{Def} The {\emph{real Clifford  algebra}} 
\index[notation]{cliff@$\CCC_k$}
$\CCC_k$ is defined to be the universal $\RR$-algebra generated by a unit $1$ and symbols  $e_1, \dots , e_k$  
 modulo the relations
 \[
 e_j^2 = -1 \qquad\qquad e_ie_j + e_je_i = 0 \qquad i \neq j
  \]
  The {\emph{complex Clifford algebra}} $\CCC _k\pc $ is defined by 
  \[
  \CCC _k\pc = \CCC_k \otimes _{\RR} \CC.
  \]
  There are natural inclusions
  \[
  \CCC _k \hookrightarrow \CCC _{k+1} \qquad\qquad and \qquad\qquad  \CCC _k\pc  \hookrightarrow \CCC _{k+1}\pc .
  \]
  \end{Def}
  
 Let's look at the first few of these. For $ \CCC _0 $ we are handed $\RR $ itself.  The algebra $\CCC _1 $ is formally given by
 \[
 \CCC_1 \cong (\RR \cdot 1) \oplus (\RR \cdot e_1) 
 \]
which is easily seen to be isomorphic to the complex numbers $\CC $. Similarly, 
 $\CCC _2 $ is formally given by
  \[
  \CCC_2 \cong (\RR \cdot 1) \oplus (\RR \cdot e_1) \oplus (\RR \cdot e_2) \oplus (\RR \cdot e_1 e_2) 
  \]
 which is isomorphic to the quaternions $\HH $.    

Complexifying is easy as well. 
\[
\CCC _0\pc \,\cong\, \RR \otimes \sr \CC \,\cong\, \CC \qquad\qquad\qquad \CCC _1\pc \,\cong\, \CC \otimes \sr \CC \,\cong\, \CC \oplus \CC .
\]
 Here is a list of the first two complex Clifford algebras, the first   eight  real Clifford algebras, and a statement of their periodicity.
  
 \begin{Thm} \cite{ABS} \label{T:york}
 \begin{enumerate}
 
 \item In low dimensions, the Clifford algebras are given as follows:
 
  \begin{center}
     \begin{tabular}{    r  |  r   } 
  {\rm{name}} & {\rm{Clifford algebra}}                          \\  \hline\hline
  $\CCC _0\pc $ & $\CC $     \\ \hline
  $\CCC _1\pc $ & $\CC \oplus \CC $            \\ \hline\hline  
  $\CCC _0 $ & $\RR $           \\ \hline
  $\CCC _1 $ & $\CC $                     \\ \hline
   $\CCC _2 $ & $\HH $                 \\ \hline
  $\CCC _3 $ & $\HH \oplus \HH  $ \\ \hline
  $\CCC _4 $ & $M_2(\HH ) $                                  \\ \hline
  $\CCC _5 $ & $M_4 (\CC) $                               \\ \hline
  $\CCC _6 $ & $M_8(\RR)  $                                \\ \hline
   $\CCC _7 $ & $M_8(\RR)\oplus M_8(\RR) $                \\ \hline
\end{tabular}
   \end{center}

\noindent where $M_k(\FF )$ denotes $k \times k$-matrices with entries in $\FF $.

\item All other real and complex Clifford algebras are determined by the periodicity formulas:
\[
\CCC _{k + 8} \,\cong\, \CCC _k \otimes_{\RR} M_{16}(\RR)  \qquad\qquad \CCC \pc _{k + 2}  \,\cong\, \CCC \pc _k \otimes_{\CC} M_{2}(\CC) 
\]

\item
\[
\dim\sr (\CCC _k ) = 2^k \qquad\qquad    \dim\sc (\CCC _k \pc ) = 2^k
\]		
\end{enumerate}
\end{Thm}

 Each Clifford algebra $A$  has a Lie \index{Lie group}    group $Lie (A) \subset A$  associated to it, namely the elements $a \in A$  that are ``unitary" in the sense that 
  $aa^* = a^*a = 1$.  
  We write 
  \[
  \LL _j = Lie (\CCC _j ) \qquad\qquad         \LL _j\pc = Lie (\CCC _j\pc ). 
  \]
  
  We can easily see what these are in low dimensions. We can also see what happens as $n \to \infty$ and 
  that the result is periodicity of period 2 (in the complex case) and 8 (in the real case.)
  
   \begin{center}
     \begin{tabular}{    r|    r|  r|  r   } 
  name & Clifford algebra    &   Lie group \,\,\,\, $\LL _j$    &     stable Lie          \\  \hline\hline
  $\CCC _0\pc $ & $\CC $    & ${\mathcal {L}}_0\pc = U_1$ & $\uu $ \\ \hline
  $\CCC _1\pc $ & $\CC \oplus \CC $   &   $\LL_1\pc = U_ 1 \times U_1 $ & $\uu \times \uu $   \\ \hline\hline

  $\CCC _0 $ & $\RR $     &  $ \LL_0 = O_1 $       &$ \oo $    \\ \hline

  $\CCC _1 $ & $\CC $           &$ \LL_1 = U_1$     &$\uu $      \\ \hline

   $\CCC _2 $ & $\HH $        &$\LL_2 = Sp_1$       & $ \pp $  \\ \hline
 \index[notation]{spk@$Sp_k$}
 \index[notation]{spz@$\pp$}
  $\CCC _3 $ & $\HH \oplus \HH  $ & $ \LL_3 =Sp_1 \times Sp_1   $  & $\pp \times \pp $ \\ \hline

  $\CCC _4 $ & $M_2(\HH ) $       &  $\LL_4 =Sp_2 $        &$ \pp $              \\ \hline

  $\CCC _5 $ & $M_4 (\CC) $            & $\LL_5 =U_4 $       & $\uu $           \\ \hline

  $\CCC _6 $ & $M_8(\RR)  $                & $\LL_6 =O_8 $       & $\oo $         \\ \hline

   $\CCC _7 $ & $M_8(\RR)\oplus M_8(\RR) $  & $\LL_7 = O_8 \times O_8 $  &$\oo \times \oo $           \\ \hline

   $\CCC _8 $ & $M_{16}(\RR)  $  & $\LL_8 = O_{16}  $  &$\oo  $           \\ \hline
   
\end{tabular}
   \end{center}

 Note that the {\emph{order}} of the Clifford algebras
 
  \[
 \CCC_0\pc \hookrightarrow \CCC _1\pc \hookrightarrow \CCC _2\pc  \hookrightarrow \CCC _3\pc \hookrightarrow \dots
 \]
 and

 \[
 \CCC_0 \hookrightarrow \CCC _1 \hookrightarrow \CCC _2 \hookrightarrow \CCC _3 \hookrightarrow \dots
 \]
 corresponds nicely to the order of the inclusions of the stable groups. The group $\uu = \bigcup U_n(\CC)$ 
 repeats in two steps
  \[
   \uu \,\hookrightarrow\, \uu /(\uu \times \uu )   \,\hookrightarrow\, \uu 
  \]
  and similarly in the real and symplectic setting there are eight steps, 
 namely  
 \[
  \oo \hookrightarrow  \uu   \hookrightarrow \pp   \hookrightarrow \pp \times \pp   \hookrightarrow \pp    \hookrightarrow \uu  \hookrightarrow \oo   \,\hookrightarrow  \oo \times \oo   \hookrightarrow \oo
\]

 \begin{Rem}\label{R:chart1}
  We are particularly interested in the relationship between the Lie groups, \index{Lie group}    and so we reconstruct the chart as follows, showing the various
  natural inclusions. The complex situation is fairly simple:
  
\[ 
\begin{CD} 
j \,\,\,\,\, @. {\text{Clifford algebra $\CCC_j$}}   @. \LL_j  @.         \LL _j / \LL_{j-1}  @.@.  {\text{ limit space}}    \\ \\
 0@. \CCC\pc_0 = \CC   @.   U_1                     \\
 @.      @VVV      @VVV  (U_1 \times U_1)/U_1 @. \quad @. (\uu \times \uu)/\uu  \simeq \uu\\
 1@.  ~~~~ \CCC\pc_1 =\CC \oplus \CC ~~ @. ~~U_1   \times U_1  ~~               \\
 @.      @VVV      @VVV   ~~U_2/(U_1 \times U_1)~~  @.@. \uu / (\uu \times \uu ) \simeq \bbu \\
 2@. \CCC\pc_2 =M_2(\CC )   @.     U_2          \\
 \end{CD}
 \]
\newpage 
 The real case is analogous, but of course more complicated. Note that {\bf{every vertical arrow below is an inclusion.}}

  \[
 \begin{CD}
j \,\,\,\,\, @. {\text{Clifford algebra $\CCC_j$}}   @. \LL_j  @.         \LL _j / \LL_{j-1}  @.@.  {\text{ limit space}}    \\ \\
0@. \CCC_0 = \RR  @.    O_1        @.   @.              \\
@. @VVV @VVV  U_1/O_1           @.@.  \uu /\oo  \\
1 @.   \CCC_1 =\CC  @.       U_1                      \\
@. @VVV @VVV   Sp_1/U_1   @.  @. \pp /\uu   \\
  2 @.\CCC_2 =\HH  @.     Sp_1          \\
@. @VVV @VVV (Sp_1 \times Sp_1)/Sp_1   @.  @. (\pp  \times \pp )/\pp = \pp  \\
  3 @.\CCC_3 = \HH \oplus \HH  @.     Sp_1 \times Sp_1                                  \\
@. @VVV @VVV Sp_2/(Sp_1 \times Sp_1)@.  @. \pp /(\pp \times \pp) \simeq \bbsp \\
  4 @.\CCC_4 = M_2(\HH )  @.       Sp_2      @.                        \\
@. @VVV@VVVU_4/Sp_2    @.   @. \uu/\pp \\
  5 @.\CCC_5 =    M_4(\CC )  @.     U_4       @.                               \\
@. @VVV @VVV                  O_8/U_4    @. @. \oo /\uu\\
   6 @.\CCC_6 =   M_8(\RR ) @.      O_8       @.                           \\
@. @VVV @VVV            (O_8 \times O_8)/O_8    @.   @. (\oo \times \oo )/\oo = \oo \\
 7@.\CCC_7 =M_8(\RR ) \oplus M_8(\RR )  @.       O_8 \times O_8     @.                        \\
@. @VVV @VVV  O_{16}/(O_8 \times O_8)  @.  @. \oo /(\oo \times \oo) \simeq \bbo \\
 8@. \CCC_8 =M_{16}(\RR )  @.      O_{16}       @. \\
   \end{CD}
  \]
 Note that the ``limit space" column is exactly the list of the Big Ten spaces in the order previously presented in 
 the Homotopy Road Map.  Here it is again:
 \newpage
 
 \begin{Table}\label{road1redux}
   \centerline{\bf{The Homotopy Road Map: Step I (Table~\ref{road1} } redux)} 
 \index{Homotopy Road Map}

    \begin{center}
     \begin{tabular}{   r  |r |   r|  r | r  } 
finite approx   &   Big Ten    &    Bott &      $\cR_k$ &  $\pi _0(\cR_k)$ \\  \hline\hline
  $U_n $ &  $\uu $ &    $\Omega ^0 \uu$   &   $\cR _1\pc $      &       0  \\  \hline
  
    $U_{2n}/(U_n \times U_n)  $ &  $\bbu \times \ZZ$ &     $\Omega ^1\uu $  &     $\cR _0\pc  $      &       $ \ZZ $    \\  \hline \hline
  $O_n $ &  $\oo $  &  $\Omega ^0 \oo$  &     $\cR _1$      &       $ \ZZ_2 $    \\  \hline
    $O_{2n}/U_n $ &  $\oo / \uu $  &  $\Omega ^1\oo$  &     $\cR _2$      &       $ \ZZ_2 $    \\  \hline
    $U_n/Sp_n $ &  $\uu / \pp $  &  $\Omega ^2 \oo$  &     $\cR _3$      &       $ 0 $    \\  \hline
    $Sp_{2n}/(Sp_n \times Sp_n) $ &  $\bbsp \times \ZZ $  &   $\Omega ^3 \oo$  &     $\cR _4$      &       $ \ZZ  $    \\  \hline
    $Sp_n $ &  $\pp $  &   $\Omega ^4\oo $  &     $\cR _5$      &       $ 0 $    \\  \hline
    $Sp_{n}/U_n $ &  $\pp/\uu $  &   $\Omega ^5\oo$  &     $\cR _6$      &       $ 0 $    \\  \hline
    $U_n /O_n $ &  $\uu /\oo $  &    $\Omega ^6\oo$  &     $\cR _7$      &       $ 0$    \\  \hline
    $O_{2n}/(O_n \times O_n) $ &  $\bbo  \times \ZZ $  &   $\Omega ^7\oo$  &     $\cR _0$      &       $ \ZZ  $    \\  \hline
\end{tabular}
   \end{center}
 \end{Table}
    \vglue .5in

This  shows a striking connection between the structure of Clifford algebras  (a clearly purely algebraic matter, 
at least at the start) to Bott Periodicity.  This very deep connection was first observed by  Atiyah, Bott, Shapiro  \cite{ABS} (1964),  who say in their introduction:

{\emph{``It is this theorem which shows the significance of Clifford algebras in $K$-theory and strongly suggests that one should look for a  proof of the 
periodicity theorem using Clifford algebras.  Since this paper was written a proof on these lines has in fact been found by R. Wood. \cite{Wood} (1966)" }}

Wood did more than this - he generalized $K$-theory from compact spaces to Banach algebras. It is his proof of periodicity in that context, 
together with Karoubi's monumental work (we list many papers of his in the bibliography, some with Villamayor)
which made 
possible the use of $K$-theory   in functional analysis and the subsequent explosive developments due to Brown-Douglas-Fillmore, Kasparov, Atiyah,  Karoubi, Connes, and many others in the 1970's, 1980's and to this day.  

\end{Rem}

\index[notation]{Lj@$\LL _k$}
 We may realize  the Lie groups $\LL _k $ \index{Lie group}    of the Clifford algebras in a different way.  Fix some large index $t$ and integer $k$.  Fix  a faithful  $C^*$-representation
 \[
\phi :  \CCC _k \longrightarrow M_{16t}(\RR). 
\]
Following the custom in physics, we let $J_i = \phi (e_i) \in O_{16t} $.  Then $J_rJ_s = - J_sJ_r $ for $r \neq s$ and $J_r^2 = -1$.  Let 
\[
\langle  O_{16t}; \, J_1, \dots J_k \rangle = \{ a \in O_{16t} \mid [a, J_i] = 0 \text{ for all $i$} \} 
\]
denote the subgroup of $O_{16t}$ consisting of all  of the matrices that commute with $J_1, \dots J_k$. 
Then we have the following result.

\begin{Thm}\label{T:chart2}
\[
\LL _k =\langle  O_{16t}; \, J_1, \dots J_k \rangle  \notag
\]
More precisely,
\begin{align*}
U_{8t} &\cong  \langle  O_{16t}; \,  J_1  \rangle \\
Sp_{4t}&\cong   \langle O_{16t} ; \,  J_1, J_2  \rangle \\
Sp_{2t} \times Sp_{2t}  &\cong   \langle  O_{16t}, \,\,\,J_1, J_2, J_3  \rangle \\
Sp_{2t} &\cong   \langle  O_{16t} ; \,  J_1, J_2, J_3, J_4  \rangle \\
U_{2t}&\cong   \langle O_{16t} ; \, J_1, J_2, J_3, J_4, J_5 \rangle \\
O_{2t} &\cong   \langle O_{16t} ; \,  J_1, J_2, J_3, J_4, J_5, J_6 \rangle \\
O_t \times O_r &\cong   \langle  O_{16t}; \,  J_1, J_2, J_3, J_4, J_5, J_6, J_7   \rangle \\
O_t &\cong   \langle  O_{16t} ; \,  J_1, J_2, J_3, J_4, J_5, J_6, J_7, J_8 \rangle \
\end{align*} 
The sequence may be extended in each direction when $t$ is a sufficiently large power of $2$. 
The analogous result holds in the complex setting as well. 

\end{Thm}

The proof of this result in \cite{Stone} is quite concrete. It references Milnor \cite{M3} for the construction.
We shall return to this construction when we introduce symmetric spaces.

 \begin{Rem}  We want to highlight something that has confused us no end when looking casually at the situation.  
Remark \ref{R:chart1} gives us the following   list for the Lie groups \index{Lie group}   , with inclusions {\bf{ left to right}} corresponding to the increasing sequence of Clifford algebras:
\begin{equation*} \label{increasing}
 O_1  \hookrightarrow   U_1   \hookrightarrow  Sp_1   \hookrightarrow  Sp_1 \times Sp_1   \hookrightarrow Sp_2   \hookrightarrow U_4  \hookrightarrow O_8   \hookrightarrow O_8 \times O_8   \hookrightarrow O_{16}     
\end{equation*}
whereas Theorem  \ref{T:chart2}, following Stone et al \cite{Stone}  gives the following list, corresponding to the increasing number of $J_i$ being used and the corresponding 
decreasing sequence of Lie groups \index{Lie group}  and Lie algebras \index{Lie algebra} and  with inclusions {\bf{right to left}}:
\begin{equation*} \label{decreasing} 
O_{16}  \hookleftarrow   U_8  \hookleftarrow   Sp_4  \hookleftarrow   Sp_2 \times Sp_2  \hookleftarrow  Sp_2  \hookleftarrow  U_2 \hookleftarrow  O_2
  \hookleftarrow  O_1 \times O_1  \hookleftarrow  O_{1}   
\end{equation*} 
If we take larger and larger groups, to reach stability, we wind up with two contrasting sequences, namely
\[
  \oo \,\overset{1}\hookrightarrow\, \uu 
  \,\overset{2}\hookrightarrow\, \pp
  \,\overset{3}\hookrightarrow\, \pp \times \pp
  \,\overset{4}\hookrightarrow\, \pp
  \,\overset{5}\hookrightarrow\, \uu
  \,\overset{6}\hookrightarrow\, \oo
  \,\overset{7}\hookrightarrow\, \oo \times \oo 
  \,\overset{8}\hookrightarrow\, \oo. 
\]
from Sequence (\ref{increasing}) and
\[
  \oo \, \overset{6}\hookleftarrow\, \uu 
  \, \overset{5}\hookleftarrow\, \pp
  \, \overset{4}\hookleftarrow\, \pp \times \pp
  \, \overset{3}\hookleftarrow\, \pp
  \, \overset{2}\hookleftarrow\, \uu
  \, \overset{1}\hookleftarrow\, \oo
  \, \overset{8}\hookleftarrow\, \oo \times \oo 
  \, \overset{7}\hookleftarrow\, \oo. 
\]
from Sequence (\ref{decreasing}).
So the groups appear in opposite order. What is really neat is that the maps are the same - just re-ordered! We have inserted the numbering of the maps 
simply to facilitate comparison. 
The same thing happens in the complex setting, but there are only two groups involved, and so it passes without notice. 
We want to refer to these two different points of view.

\begin{Def} \begin{enumerate}
\item The {\emph {\york perspective}} will entail regarding the ordering 
\[
\overset{8}\hookrightarrow\,
  \oo \,\overset{1}\hookrightarrow\, \uu 
  \,\overset{2}\hookrightarrow\, \pp
  \,\overset{3}\hookrightarrow\, \pp \times \pp
  \,\overset{4}\hookrightarrow\, \pp
  \,\overset{5}\hookrightarrow\, \uu
  \,\overset{6}\hookrightarrow\, \oo
  \,\overset{7}\hookrightarrow\, \oo \times \oo 
  \,\overset{8}\hookrightarrow\, \oo  \overset{1}\hookrightarrow\
\]
with associated homogeneous spaces ordered correspondingly:
\[
  \uu /   \oo
  \,\quad\, \pp/\uu
  \,\quad\, (\pp \times \pp) / \pp
  \,\quad\, \pp /(\pp \times \pp) 
  \,\quad\, \uu /\pp
  \,\quad\, \oo /\uu
  \,\quad\, ( \oo \times \oo) /\oo
  \,\quad\, \oo / (\oo \times \oo )
\]
 arising from the inclusions of Clifford algebras $\CCC _0 \subset \CCC_1 \cdots $ as fundamental.

\item  The {\emph{\lancaster perspective}} will entail regarding the ordering
 \[
  \overset{7}\hookleftarrow\,
  \oo \, \overset{6}\hookleftarrow\, \uu 
  \, \overset{5}\hookleftarrow\, \pp
  \, \overset{4}\hookleftarrow\, \pp \times \pp
  \, \overset{3}\hookleftarrow\, \pp
  \, \overset{2}\hookleftarrow\, \uu
  \, \overset{1}\hookleftarrow\, \oo
  \, \overset{8}\hookleftarrow\, \oo \times \oo 
  \, \overset{7}\hookleftarrow\, \oo
   \overset{6}\hookleftarrow\, 
\]
with associated homogeneous spaces
\[
  \,\quad\, \oo /\uu
  \,\quad\, \uu /\pp
  \,\quad\, \pp /(\pp \times \pp) 
  \,\quad\, (\pp \times \pp) / \pp
    \,\quad\, \pp/\uu
  \,\quad  \uu /   \oo 
   \,\quad\, \oo / (\oo \times \oo )
\,\quad\, ( \oo \times \oo) /\oo
  \]
arising from using the $J_i$ to shrink the Clifford algebras as fundamental.  
Note that the maps in the \lancaster perspective match exactly to the maps in the \york perspective
and hence exactly the same homogeneous spaces appear, albeit in reverse order!  \index{Clifford algebra|)}
\end{enumerate}
\end{Def}
 \end{Rem}
 
\newpage
  \section{\bf{Clifford Algebras  and Bott Periodicity: the Real Case}}
  \label{Section:CliffordAlgebrasBottPeriodicity}
  
 In this section and the next we outline the Wood-Karoubi proofs of Bott Periodicity \index{Bott Periodicity} via Clifford algebras. \index{Clifford algebra|(}  We try to balance the detail, 
and we cite but don't present the one hard homotopy step involved in the proof.   The real case is significantly longer in the sense that 8 
is greater than 2, but the mathematical steps are essentially the same. In this section we present the real case. In the next 
section we present the complex case.

 \begin{Def}
 \index[notation]{LLA@$\LL A$}
  \index[notation]{GGA@$\GL A$}
Suppose that $A$ is a unital real Banach algebra. Let $\LL A$ denote the Lie group \index{Lie group}  associated to $A$.  Let   $\GL A$ 
denote the group of invertible 
elements of $A$.   Then $\LL A \subset \GL A $ and this map is a deformation retraction (essentially, by Gram-Schmidt orthogonalization) and 
hence a homotopy equivalence. 
\end{Def}

Using the larger 
groups has some technical conveniences in homotopy theory, and they are used by Wood, and so we shall use them as well. 

Here's some notation that we will use:
\begin{enumerate}

\item   $A\otimes B = A \otimes \sr B$  
 
\item $\GL_k(A) = \GL(A \otimes M_k(\RR ) )$  

\item $\GL_\infty (A) = \lim_{k \to \infty}  \GL_k(A) $, the direct limit with the weak topology. 
\index[notation]{GLKA@$\GL_k(A)$}
\index[notation]{GLzA@$\GL_\infty(A)$}

\item If $X$ is a based topological space, then $X_o$ denotes the path component of the basepoint of $X$.  (If $X$ is a topological group 
then we take the identity of the group as the basepoint unless there is a good reason not to.)

\end{enumerate}

There are natural inclusions of real  Clifford algebras $\CCC _r \longrightarrow \CCC_{r+1} $  which induce a natural  commutative diagram with \york perspective

\[
\begin{CD}
\GL(A \otimes \CCC _r )        @>>>                   \GL_k(A \otimes \CCC _r )               @>>>          \GL_\infty(A \otimes \CCC _r ) \\ 
@VVV         @VVV            @VVV    \\
\GL(A \otimes \CCC _{r+1 } )       @>>>                   \GL_k(A \otimes \CCC _{r+ 1} )               @>>>          \GL_\infty(A \otimes \CCC _{r+ 1} ) \\ 
@VVV         @VVV            @VVV    \\
\GL(A \otimes \CCC _{r+1 })/ \GL(A \otimes \CCC _{r })   @>>>                   \GL_k(A \otimes \CCC _{r+ 1} )/ \GL_k(A \otimes \CCC _{r} )
                @>>>          \GL_\infty(A \otimes \CCC _{r+ 1} )  /   \GL_\infty(A \otimes \CCC _{r} )     \\ 
\end{CD}
\]

We have seen instances of these spaces before. For example, taking $A = \RR$, we obtain (after a deformation retraction)
\[
\GL_\infty( \CCC _{7} )  /   \GL_\infty( \CCC _{6} )   \simeq   (\oo \times \oo ) / \oo  = \oo  = \cR_1
\]
and 
\[
\GL_\infty( \CCC _{6} )  /   \GL_\infty( \CCC _{5} )   \simeq     \oo / \uu    =   \cR_2.
\]

Here is an example of the contrast/confusion generated between the \york perspective and the \lancaster perspective in the numbering.  We go ``down" from $7\dots 6$ 
to $6\dots 5$ in the \york indexing of the quotients associated to the Clifford algebras and the number goes up from $1$ to $2$ in the \lancaster indexing of the $\cR_j$.

 Next, we place ABS Clifford periodicity on record in the form that we need.

\begin{Pro} (Atiyah-Bott-Shapiro  \cite {ABS} ) \label{ABS}  Let $A$ be a unital real Banach algebra. Then 
there is a natural homotopy equivalence
\[
GL_\infty(A \otimes  {\CCC}_{r+ 8})   \simeq     GL_\infty(A \otimes \CCC _{r} ) .
\]
 \end{Pro}

\begin{proof}
ABS prove that 
\[
\CCC _{r+8} \cong \CCC _r \otimes M_{16}(\RR)
\]
which implies immediately that 
\[
A \otimes \CCC _{r+8} \cong  A \otimes \CCC _r \otimes M_{16}(\RR) .
\]
Thus for each $k$,  
\[
\GL _k(A \otimes \CCC _{r+8} ) \cong \GL _k( A \otimes \CCC _r \otimes M_{16}(\RR) ) \cong \GL _{k+ 16}(A \otimes \CCC _{r} ) 
\]
and letting $k \to \infty $  implies the result.  

\end{proof}

Next, we invoke Wood's result.  It is essentially topological, using the internal structure of the topological groups involved and their symmetries. 
   (Karoubi (cf. his excellent book \cite{K} ) packages the matter differently but ultimately proves 
the analogous result.)    For yet a different packaging, see Milnor \cite{M1}, \S 24.    For specifics, recall that the maps 
\[
 \LL _r \hookrightarrow  \GL_\infty(  \CCC _{r} )
\]
are homotopy equivalences, and then look at Remark \ref{R:chart1} for a list of the various $\LL _r$ and, most importantly, $\LL_r / \LL _{r-1}$.

\begin{Thm} Wood \cite{Wood}\label{Wood} ( 4.7, 4.9)   Suppose that $A$ is a real unital Banach algebra.  Then:
\begin{enumerate}
\item
There is a weak homotopy equivalence\footnote{A based map $f: X \to Y$ is a {\emph{weak homotopy equivalence}} if the induced map 
\[
f_* : \pi _*(X) \to \pi_*(Y)
\]
 is an isomorphism.  Every homotopy equivalence is a weak homotopy equivalence but not conversely.}
\[
 \big( \GL_\infty(A \otimes \CCC _{r} ) / \GL_\infty(A \otimes \CCC _{r-1} )\big)_o    \overset{\simeq}\longrightarrow    \Omega ^1\big(  \GL_\infty(A \otimes \CCC _{r+1} )/  \GL_\infty(A \otimes \CCC _{r} ) \big)_o.
\]
\item  There is a weak homotopy equivalence 
\[
\big( \GL_\infty(A \otimes \CCC _{r} ) / \GL_\infty(A \otimes \CCC _{r-1} )\big)_o    \overset{\simeq}\longrightarrow    \Omega ^8\big(  \GL_\infty(A \otimes \CCC _{r+8} )/  \GL_\infty(A \otimes \CCC _{r+ 7} ) \big)_o.
\]
\item There is a weak homotopy equivalence
\[
\big( \GL_\infty(A \otimes \CCC _{r} ) / \GL_\infty(A \otimes \CCC _{r-1} )\big)_o    \overset{\simeq}\longrightarrow
\Omega ^8 \big(  \GL_\infty(A \otimes \CCC _{r} ) / \GL_\infty(A \otimes \CCC _{r-1}) \big)_o   
\]

\item Taking  $A = \RR$ and   $r = 7$ in Part (3) above, there is a weak homotopy equivalence
\[
\oo _o  \overset{\simeq}\longrightarrow \Omega ^8 \oo _o .
\]
\item Rephrasing in the \lancaster perspective and recalling that $\cR_1 = \oo$, there is a natural homotopy equivalence 
\[
\cR_s \simeq \Omega ^s \cR_0
\]
which implies that 
\[
\pi _a(\cR_b) \cong \pi _0(\cR_{a+b} ) \cong \pi_{a+b}(\cR_0)
\]
and, in particular,
\[
\pi_k(\oo) = \pi_k(\cR_1) = \pi _0(\cR_{k+1}) 
\]
with the indexing understood to be modulo 8.
So we know all of the homotopy groups of $\oo$.\footnote{This cannot be said for very many spaces. Topologists build Eilenberg-MacLane spaces, for instance, that have exactly one non-zero homotopy group. However, in situations (such as the 2-sphere $S^2$, to take a really simple example) where the higher homotopy groups do not all vanish, very little is 
actually known. Keep in mind that if somebody says that they know the first million homotopy groups of $S^2$ then they are talking about 0\% of all of the homotopy groups of $S^2$.  
} 
\end{enumerate}
\end{Thm}

Note that part (2) follows from part (1): you simply iterate (1) eight times.  For part (3)  you use the ABS result, Theorem \ref{ABS}.  
Part (4) is simply the recognition that 
\[
 \GL_\infty( \CCC _{7} ) / \GL_\infty(\CCC _{6} )_o  \cong \oo _o   
\]
The hard part, then, is Part 1, and this is where Wood makes a major contribution.

Please note that this strongly ties Clifford algebras, or, more precisely, their Lie groups, \index{Lie group} to the Homotopy Road Map, \index{Homotopy Road Map} 
as shown in Table~\ref{road1redux}.  
Rather than regarding 
the first column as coming to us from Bott Periodicity, we may think of the column as simply presenting yet a different look at the \york and 
\lancaster perspectives.  
  
 \index{Clifford algebra|)}

\newpage
  \section{\bf{Clifford Algebras  and Bott Periodicity: the Complex Case}} \index{Bott Periodicity}
   \label{Section:CliffordAlgebrasComplex}
 
 Here we vaguely sketch Wood's proof of Bott Periodicity in the complex setting.
  \index{Clifford algebra|(}
 
 \begin{Def}
Suppose that $A$ is a unital complex Banach algebra.  
Here's some notation that we will use in this section:
\begin{enumerate}

\item   $A\otimes B = A \otimes \sc B$  
 
\item $\GL_k(A) = \GL(A \otimes M_k(\CC ) )$  

\item $\GL_\infty (A) = \lim_{k \to \infty} \GL_k(A) $, the direct limit with the weak topology. 
\end{enumerate}
 \end{Def}
  
There are natural inclusions of complex  Clifford algebras $\CCC  \pc_r  \longrightarrow \CCC_{r+1}\pc  $  which induce a natural commutative diagram

\[
\begin{CD}
\GL(A \otimes \CCC  \pc_r)        @>>>                   \GL_k(A \otimes \CCC  \pc_r )               @>>>          \GL_\infty(A \otimes \CCC _r \pc) \\ 
@VVV         @VVV            @VVV    \\
\GL(A \otimes \CCC _{r+1 }\pc )       @>>>                   \GL_k(A \otimes \CCC _{r+ 1}\pc )               @>>>          \GL_\infty(A \otimes \CCC _{r+ 1}\pc) \\ 
@VVV         @VVV            @VVV    \\
\GL(A \otimes \CCC _{r+1 }\pc)/ \GL(A \otimes \CCC _{r }\pc) )  @>>>                   \GL_k(A \otimes \CCC _{r+ 1} \pc)/ \GL_k(A \otimes \CCC _{r} \pc)
                @>>>          \GL_\infty(A \otimes \CCC _{r+ 1}\pc )  /   \GL_\infty(A \otimes \CCC _{r}\pc )     \\ 
\end{CD}
\]

Taking $A = \CC$, we obtain (after a deformation retraction)
\[
\GL_\infty( \CCC _{1}\pc  )  /   \GL_\infty( \CCC _{0}\pc  )   \simeq   (\uu \times \uu ) / \uu  = \uu .
\]

 Next, we place ABS Clifford complex periodicity on record in the form that we need.

\begin{Pro} (Atiyah-Bott-Shapiro  \cite {ABS}) \label{ABS}  Let $A$ be a unital complex Banach algebra. Then 
there is a natural homotopy equivalence
\[
GL_\infty(A \otimes \CCC _{r+ 2}\pc  )   \simeq     GL_\infty(A \otimes \CCC _{r}\pc ) .
\]
 \end{Pro}

\begin{proof}
ABS prove that 
\[
\CCC _{r+2}\pc \cong \CCC _r \pc \otimes M_{2}(\CC)
\]
which implies immediately that 
\[
A \otimes \CCC _{r+2}\pc  \cong  A \otimes \CCC  \pc_r  \otimes M_{2}(\CC) .
\]
Thus for each $k$,  
\[
\GL _k(A \otimes \CCC _{r+2}\pc  ) \cong \GL _k( A \otimes \CCC  \pc_r  \otimes M_{2}(\CC) ) \cong \GL _{k+ 2}(A \otimes \CCC _{r} \pc ) 
\]
and taking limits, this implies the result.  

\end{proof}

Next, we invoke Wood's result.   

\begin{Thm} Wood \cite{Wood} (Proposition~4.7 and Remark~4.9).   Suppose that $A$ is a complex unital Banach algebra.  Then:
\begin{enumerate}
\item
There is a weak homotopy equivalence
\[
 \big( \GL_\infty(A \otimes \CCC _{r}\pc  ) / \GL_\infty(A \otimes \CCC _{r-1}\pc  )\big)_o    \overset{\simeq}\longrightarrow    
 \Omega ^1\big(  \GL_\infty(A \otimes \CCC _{r+1}\pc  )/  \GL_\infty(A \otimes \CCC _{r}\pc  ) \big)_o.
\]
\item  There is a weak homotopy equivalence 
\[
\big( \GL_\infty(A \otimes \CCC _{r} \pc) / \GL_\infty(A \otimes \CCC _{r-1}\pc )\big)_o    \overset{\simeq}\longrightarrow    
\Omega ^2\big(  \GL_\infty(A \otimes \CCC _{r+2}\pc )/  \GL_\infty(A \otimes \CCC _{1}\pc ) \big)_o.
\]
\item There is a weak homotopy equivalence
\[
\big( \GL_\infty(A \otimes \CCC _{r}\pc ) / \GL_\infty(A \otimes \CCC _{r-1}\pc )\big)_o    \overset{\simeq}\longrightarrow
\Omega ^2 \big(  \GL_\infty(A \otimes \CCC _{r} \pc) / \GL_\infty(A \otimes \CCC _{r-1}\pc) \big)_o   
\]

\item Taking  $A = \CC$,  there is a weak homotopy equivalence
\[
\uu _o  \overset{\simeq}\longrightarrow \Omega ^2 \uu _o .
\]

\end{enumerate}
\end{Thm}

Note that part (2) follows from part (1): you simply iterate (1) twice.  For part (3)  you use the ABS result,  Theorem \ref{ABS}.  
Part (4) is simply the recognition that 
\[
 \GL_\infty( \CCC _{1}\pc  ) / \GL_\infty(\CCC _{0}\pc  )_o  \cong \uu _o .  
\]
  
 Our comment at the end of the previous section regarding the Homotopy Road Map applies here as well. 
    \index{Clifford algebra|)} \index{Homotopy Road Map}
   
  \newpage 
  \section{\bf{The Basics: Symmetric Spaces and a Timely Example}} \index{symmetric space|(}
       \label{Section:SynmetricSpaces}

 We rely on Helgason \cite{Helgason} and Gorodski \cite{Gorodski}  for our basic references for symmetric spaces.  This is actually a vast topic and we are simply skimming the surface.

   \begin{Def} Suppose that $(M, g)$ and $(N, h)$ are Riemann manifolds\footnote{A {\emph{Riemann metric}} on a smooth manifold is a smooth choice of a positive definite inner product for each tangent 
   space of the manifold. A {\emph{Riemann manifold}} is a smooth manifold together with a Riemann metric.}
   \index{Riemannian manifold}
    and $f: M \to N$ is a diffeomorphism (that is, $f$ is a homeomorphism which is smooth and has a smooth inverse.) 
 Then $f $ is an {\emph{isometry}} if $f$ preserves distance.  In shorthand, we may say $g = f^*h$, and more explicitly, 
 \[
 g_p(u, v) = f_{f(p)} (df_p(u), df_p(v)) 
 \]
 for all $p \in M$ and $u, v \in T_pM$. 
 \end{Def}

 \begin{Def} Suppose that $M$ is a connected Riemann manifold and let $p \in M$.  There is an exponential map 
 \[
 \exp _p : T_pM \to M. 
 \]
   We may choose a small open ball $B_p$ about the origin of 
 $T_pM$ and label its image as 
 \[
 V_p  \equiv  \exp_p(B_p) \subset M.  
 \]
 If $B_p$ is small enough then the map 
 \[
 \exp_p : B_p \longrightarrow   V_p  
 \]
 is a diffeomorphism.    Define the {\emph{geodesic symmetry at $p$}},  \index{geodesic symmetry} 
 written   $s_p : V_p \to V_p$ to be the composite
  \[
V_p \overset{\exp _p^{-1}}\longrightarrow B_p \overset{~~u \mapsto -u~~}\longrightarrow  B_p  \overset{\exp _p}\longrightarrow   V_p.
 \]
 \end{Def}
 The key point is that $s_p$ reverses the direction of geodesics.
   
 \begin{Def} A Riemann manifold $M$ is a {\emph{symmetric space}} if for each $p \in M$ the map $s_p$ is a globally defined isometry of $M$.   
 (Sometimes these are referred to as {\emph{Riemannian symmetric spaces}} but we will leave off the first word.)
 \end{Def}
 
 Here are some examples:
 \begin{enumerate}
 \item  Euclidean space $\RR ^n$, with $s_p(y) = 2p - y$ 
 \item  spheres $S^n$
 \item projective spaces 
 \item compact Lie groups \index{Lie group}  
 \end{enumerate}
 
 In order to simplify matters and to concentrate immediately upon the cases that concern us,  we will stick to compact symmetric spaces.  We will also not  
 discuss the exceptional \index{Lie group}  Lie groups.   These seem to not play a role in physics, as yet.\footnote{ For example,   the lowest dimensional irreducible representation 
  of the group $E_8$ is in dimension 248. }  
 
 \index{Lie group|(} 
 
  \begin{Pro} (cf. Helgason \cite{Helgason}, \S 3-4.)
 Every symmetric space arises as a homogeneous space of Lie groups. 
  \end{Pro}

 For example, the symmetric space $S^k$ has isometry group $O_{k+1} $, its isometry group at the point $(1,0,0,\dots , 0) $ is $O_k$ and hence we may write 
 \[
 S^k \cong O_{k+1} / O_k .
 \]
 Note that had we used the connected component of the identity, then we would get the same outcome, since $SO_{k+1}/ SO_k \cong S^k$. 

 We follow, more or less, the 1926 approach of \' Elie Cartan. 
  
  First, any irreducible symmetric space $M$ fits into one of three categories:
  \begin{enumerate}
  \item $M$ is isometric to some Euclidean space.
  \item $M$ is compact. (It will then have nonnegative (but not identically zero) sectional curvature.)
  \item $M$ is not compact and has nonpositive (but not identically zero) sectional curvature.  This class is dual to the compact class.
  \end{enumerate}
  
  As mentioned, we are restricting to the second category.
  The classification of \' E. Cartan uses heavily the Lie algebra   associated to the Lie group \index{Lie group}  . Here's a quick introduction. We will  focus completely on  Lie algebras over either $\RR $ or $\CC $, and we denote 
  the field by $\FF $. 
  
  \begin{Def} A {\emph{Lie algebra}} \index{Lie algebra|(}   $\mathfrak{g}$ over a field $\FF $ is an $\FF$-vector space $\mathfrak{g}$  together with a map  
  \[
[\,\, , \,\,  ] :   \mathfrak{g} \times \mathfrak{g}   \longrightarrow  \mathfrak{g}
  \]
  such that 
  \begin{enumerate}
  \item The map is bilinear.
  \item $[x, x ] = 0$ for all $x \in \mathfrak{g} $.
  \item The pairing satisfies the Jacobi identity. That is, for all $x,y,z \in \mathfrak{g}$, 
  \[
  [x, [y,z ]] + [y, [z,x]] + [z, [x,y]] = 0.
  \]
  \end{enumerate}
  \end{Def}
  
  For example, take $ M_k(\FF)$ with $[A,B] = AB - BA$.\footnote{Everyone should do the following exercise once in their lives. Suppose
  that $\mathfrak {g} = M_k(\FF)$ with $[A,B] = AB - BA$.  Then check explicitly that $\mathfrak{g} $ satisfies the Jacobi identity by calculating 
   $[A, [B,C ]] + [B, [C,A]] + [C, [A,B]] $ for $A, B, C  \in M_k(\FF)$.}
 
  The most important examples for us are the Lie algebras  associated to particular  Lie groups.  These ideas date back to Sophus Lie (in the period 1869-1873), to his collaborators, particularly to  Felix Klein, and to his students, particularly {\' E}lie Cartan.
  
  \begin{Def}
   Fix a Lie group $G$. \index{Lie group}  Then it is a smooth manifold, and so at each point $h \in G$ there is a 
  tangent space $T_hG$.  Suppose that $L_g : G \to G$ is a vector field on $G$ induced by $L_g(x) = gx$.  We say that $L$ is  a {\emph{left invariant vector field}} if   
  for all $g, h$, 
  \[
  (dL_g)_h : T_hG \longrightarrow T_{gh}G .
  \]
  We denote by $Lie (G)$ the real vector space of all left-invariant vector fields on $G$.  It is a Lie algebra with Lie bracket defined by
  \[
  [L_g, L_h] = L_gL_h - L_hL_g
  \]
      There is special notation for the Lie algebras associated to the classical Lie groups. The Lie algebra \index{Lie group}  \index{Lie algebra}  of the 
  group $SO_n$ is ${\mathfrak{so}}_n $ and similarly for $O_n$, $U_n$, $SU_n$ and $Sp _n$.
 
     \end{Def}
     
     Now, what is the relationship between Lie groups and Lie algebras?  Here are the facts:

     \begin{Thm} (S. Lie)
     \index[notation]{LieG@$Lie(G)$}
     \begin{enumerate}
     \item If $G$ is a Lie group then $Lie (G)$ is a   Lie algebra, with $\dim (G)$ (as a compact smooth manifold) equal to $\dim (Lie(G))$ (as a real vector space.)
    \item  If $ \mathfrak g$    is a finite-dimensional Lie algebra,  then there exists a connected Lie group $G$ with $Lie (G) \cong \mathfrak g$.      Thus 
    the map 
    \[
    \{ \text{Lie Groups  } \} \overset{\text{Lie}(-)}\longrightarrow \{ \text{Finite-Dimensional Lie Algebras} \footnote{where both sides mean ``isomorphism classes".}  \}
    \]
         is surjective.
     \item Any two Lie groups with the same Lie algebra are locally isomorphic and, in particular, have the same universal cover. 
          \item  If we restrict attention to connected,  {\emph{simply connected\footnote{i.e. with $\pi _1(G) = 0$.}}} Lie groups, then  there is a bijective correspondence
      \[
    \{ \text{Simply Connected Lie Groups} \} \xrightarrow{~Lie~}  \{ \text{Finite-Dimensional Lie Algebras} \}.
     \]
     \item {\bf Cartan Decomposition.} If $K$ is a closed Lie subgroup of the Lie group $G$  with associated Lie algebras $\mathfrak k$ and $\mathfrak g$, then $\mathfrak g$ may be written as a direct sum 
     \[
     \mathfrak g \cong \mathfrak k \oplus \mathfrak m 
     \]
     as vector spaces, with the commutator Lie bracket attending to this decomposition as follows: if $k, k_1, k_2  \in \mathfrak k$ and $m, m_1, m_2 \in \mathfrak m$, then 
     \[
     [k_1, k_2 ] \in  \mathfrak k     \qquad    [ k, m ] \in  \mathfrak m     \qquad  [  m_1,  m_2 ] \in  \mathfrak k  .
     \] \index{Cartan decomposition}
  \end{enumerate}
 \end{Thm}    
    \vglue .3in

  \begin{Thm} (E. Cartan)   Ignoring the exceptional Lie groups, there are ten different families of  compact symmetric spaces. They fall into three groups:
    \begin{enumerate}
  \item     There are three compact connected  Lie groups that appear in the classification as is: they are $SO_n$,  $ SU_n$, and $Sp_n$. 
  \item     There are four symmetric spaces built in the $G/K$ format that are built out of two different groups. They are
  \[
  SU_n / SO_n \qquad SU_{2n}/Sp_n  \qquad 
   SO_{2n}/U_n \qquad\qquad Sp_n / U_n .
 \]
\item   There are three symmetric spaces, also in the $G/K$ format, that have the form $G_{p+q}/(G_p \times G_q)$
\[
 SO_{p+q}/(SO_p \times SO_q).  
 \qquad SU_{p+q}/(SU_p \times SU_q)
 \qquad Sp_{p+q}/(Sp_p \times Sp_q)  
 \]
\end{enumerate}
\end{Thm}

These ten different categories of spaces correspond exactly to the ten categories of spaces we have 
seen, namely the eight spaces of the type $GL_\infty (\CCC _r)/GL_\infty ( \CCC _{r-1}) $ associated to the real numbers and the two 
space types associated with the complex numbers! 

\index{symmetric space|)}

Here is an extension of the Homotopy Road Map demonstrating the correspondence.


\begin{Table}\label{road2}
\centerline{\bf{The Homotopy Road Map: Step 2 }} \index{Homotopy Road Map} 

    \begin{center}
     \begin{tabular}{ r |r  |r  |r |r |r |r  } 
label & name & finite approx   &   Big Ten    &    Bott &      $\cR_k$ &  $\pi _0(\cR_k)$ \\  \hline\hline
 A &   unitary & $U_n $ &  $\uu $ &    $\Omega ^0 \uu$   &   $\cR _1\pc $      &       0  \\  \hline
AIII & chiral unitary & $U_{2n}/(U_n \times U_n)  $ &  $\bbu \times \ZZ$ &     $\Omega ^1\uu $  &     $\cR _0\pc  $      &       $ \ZZ $    \\  \hline \hline
  D & BdG &     $O_n $ &  $\oo $  &  $\Omega ^0 \oo$  &     $\cR _1$      &       $ \ZZ_2 $    \\  \hline
   DIII &  BdG & $O_{2n}/U_n $ &  $\oo / \uu $  &  $\Omega ^1\oo$  &     $\cR _2$      &       $ \ZZ_2 $    \\  \hline
    AII & symplectic & $U_n/Sp_n $ &  $\uu / \pp $  &  $\Omega ^2 \oo$  &     $\cR _3$      &       $ 0 $    \\  \hline
    CII & chiral  symplectic &  $Sp_{2n}/(Sp_n \times Sp_n) $ &  $\bbsp \times \ZZ $  &   $\Omega ^3 \oo$  &     $\cR _4$      &       $ \ZZ  $    \\  \hline
    C & BdG & $Sp_n $ &  $\pp $  &   $\Omega ^4\oo $  &     $\cR _5$      &       $ 0 $    \\  \hline
    CI & BdG &  $Sp_{n}/U_n $ &  $\pp/\uu $  &   $\Omega ^5\oo$  &     $\cR _6$      &       $ 0 $    \\  \hline
    AI & orthogonal &$U_n /O_n $ &  $\uu /\oo $  &    $\Omega ^6\oo$  &     $\cR _7$      &       $ 0$    \\  \hline
    BDI & chiral orthogonal & $O_{2n}/(O_n \times O_n) $ &  $\bbo  \times \ZZ $  &   $\Omega ^7\oo$  &     $\cR _0$      &       $ \ZZ  $    \\  \hline
\end{tabular}
   \end{center}
\end{Table}

 \vglue .3in
 {\bf{A Timely Example}}
       \vglue .3in
          
For the rest of this section we are going to consider a timely example, a case that was used very recently by Orion and Akkermans \cite{OA} to understand the entanglement of two qubits.   \index{qubit}
They consider the case where the Lie group is $U_4$ and the closed subgroup is $O_4$ and then restrict to $SU_4$.   So let's look at the  chart, 
featuring as a special case 
 $U_{2t}/O_{2t}$.   Here are 
some basic topological facts:

\begin{Pro}
\begin{enumerate}
\item $U_n$ and $SU_n$ are compact smooth manifolds of dimension $n^2$ and $(n^2 - 1)$ respectively.  Hence their Lie algebras 
$\mathfrak{u}_n$ and $\mathfrak{su}_n$ 
have dimension $n^2$ and $(n^2 - 1)$ respectively as 
real vector spaces.
\item  $O_n$ and $SO_n$ are compact smooth manifolds, both of dimension $\frac{n(n-1)}{2}$.  Their Lie algebras are isomorphic:
\[
\mathfrak{o}_n \cong \mathfrak{so}_n 
\] 
 and have dimension $\frac{n(n-1)}{2}$ as real vector spaces.
 \item In particular, 
 \[
\dim \,\mathfrak{su}_4 = 15 \qquad \text{and} \qquad \dim\, \mathfrak{o}_4 = 6.
 \]
 \item Using  the Cartan decomposition, we may decompose
 \[
 \mathfrak{su}_4      \cong   \mathfrak{o}_4 \oplus \mathfrak{m}
 \]
as vector spaces, where  $\dim \, \mathfrak{m} = 9$. 
  Under this decomposition, the Lie bracket has the following property:
 \[
 [h_1, h_2] \in   \mathfrak{o}_4 \qquad\forall h_j \in \mathfrak{o}_4  
 \]
 \[
  [m_1, m_2] \in   \mathfrak{o}_4 \qquad \forall m_j \in \mathfrak{m}
 \]
 \[
 [h, m] \in \mathfrak{m} \qquad \forall h \in \mathfrak{o}_4 \quad m \in \mathfrak{m}
 \]
\end{enumerate}
   \end{Pro}
  
   In order to move forward we need some definitions known to every physicist but not, perhaps, to every mathematician. {\bf{We will 
   use physics notation for the remainder of this section and for the next. Beware!}}
   
   \index{Lie group|)} \index{Lie algebra|)} 
   
   \begin{Def} 
   \begin{enumerate}
   \item We denote the complex conjugate of a complex number $z$ by $z^*$ rather than $\bar{z}$. 
   \item If $A$ is a complex matrix then the {\emph{Hermitian adjoint}} $A^\dag $ of $A$ is defined by 
   \[
   [a_{ij}]^\dag = [a_{ji}^* ]
   \]
   \item $A$ is said to be {\emph{Hermitian}} (or self-adjoint)  \index{Hermitian} \index{self-adjoint}
   if $A^\dag = A$ and $A$ is {\emph{skew-Hermitian}}  \index{skew-Hermitian} if $A^\dag = -A$. 
   \item
   If $f: V \to W$ where $V$ and $W$ are complex vector spaces, then $f$ is said to be {\emph{antilinear}} or {\emph{conjugate linear}} 
   if for all $x, y \in V$ and $z \in \CC$, we have 
   \[
   f(x + y) = f(x) + f(y)  \qquad                {\text{and}}           \qquad f(zx) = z^* f(x) .
   \]
   \item If $U : \mathcal{H} \to \mathcal{H}$ where $\mathcal{H}$ is a complex Hilbert space, then $U$ is said to be {\emph{antiunitary}} \index{antiunitary} if it is antilinear and, in addition, if 
   for all $x, y \in \mathcal{H}$, 
   \[
   \langle Ux , Uy \rangle  =  \langle  x,y \rangle ^* .
   \]
   \item Given a basis $\beta = \left( v_1, \ldots, v_n \right)$ for $V$, we define the complex conjugation operator $K \colon  V \rightarrow V$ by
   \begin{align*} Kv_j &= v_j \\ K(iv_j) &= -iv_j  \;  \end{align*}  
  (and extended to all of $V$ by linearity over $\RR$). Then one checks that 
   $
   K^\dag K = K^2 = I.
   $
   \end{enumerate}
   \end{Def}
   
   We need a basis for the vector space $\mathfrak{su}_4$ that respects the Cartan decomposition
\[
 \mathfrak{su}_4      \cong   \mathfrak{o}_4 \oplus \mathfrak{m}.
 \]
It is constructed as follows.
Let 
\begin{equation*}
\sigma _1  = \begin{bmatrix} 0&1\\1 &0 \end{bmatrix}   \qquad
   \sigma _2  =\begin{bmatrix} 0&-i\\i &0 \end{bmatrix}  \qquad
   \sigma _3  = \begin{bmatrix} 1 &0\\0&-1 \end{bmatrix}   \qquad
 \sigma _0 = \begin{bmatrix} 1&0\\0&1 \end{bmatrix}
 \end{equation*}
   denote the $2 \times 2$  Pauli matrices in $M_2(\CC )$.       
   
   \begin{Pro} The Pauli matrices \index{Pauli matrices} satisfy the following properties: 
   \begin{enumerate}
   \item Each matrix is Hermitian.\footnote{ The matrices $\{ \sigma _0, \sigma _1,\sigma _2,\sigma _3 \}$
   form a basis for 
   the Hermitian elements of $M_2(\CC )$. }
      
   \item The Cayley table, their multiplication table,  shows the values of the row times the column:
     \begin{center}
     \begin{tabular}{   r  |r |   r|   r | r| }
 $\times$  & $\sigma _1$   & $\sigma _2$ & $\sigma _3$\\       \hline
  $\sigma _1 $ & $\mathbb{I}$ & $ i\sigma _3 $& $-i\sigma _2$      \\  \hline
    $\sigma _2 $ & $-i\sigma _3 $ & $ \mathbb{I}$& $i\sigma _1 $      \\  \hline
      $\sigma _3 $ & $i\sigma _2 $ & $ -i\sigma _1  $& $\mathbb{I}$      \\  \hline
       \end{tabular}
   \end{center}
(To read the table correctly, note that $\sigma _1 \sigma _2 =   i\sigma _3 $.) 
      
   \item In particular, for $i, j = 1,2,3$ and $i \neq j$, then 
   \[
   \sigma _i \sigma _j = - \sigma _j \sigma _i .
   \]
      
   \item Complex conjugation is given by 
   \[
      \sigma _1^* = \sigma _1 \qquad    \sigma _2^* = - \sigma _2 \qquad    \sigma _3^* = \sigma _3
   \]
   and
   \[
   (i \sigma _1)^* = -i \sigma _1 \qquad    (i \sigma _2)^* = i  \sigma _2 \qquad    (i \sigma _3)^* = -i \sigma _3
   \]

   \item The matrices $i \sigma _1, i \sigma _2, i\sigma _3 $ form a basis for the real Lie algebra $\mathfrak{su}_2 $. 
      
 \item 
 \[
 \det \,\sigma _j = -1  \qquad\qquad \rm{tr}\, \sigma _j = 0 \qquad \forall j
   \]
   and hence each matrix $\sigma _j $ has eigenvalues $+1$ and $-1$. 
   
   \end{enumerate}
   \end{Pro}

  Here is the  basis for $\mathfrak{su}_4 $ that respects the Cartan decomposition, as promised.
    
    \begin{Pro} \label{Prop-mathfrak-m}
    Regard 
       \[
   M_4(\CC) \cong M_2(\CC ) \otimes M_2(\CC ) .
   \]
   
    Then the Lie algebra  $\mathfrak{o}_4 $  is the span of the six basis 
   elements 
   \[
   i\sigma _j \otimes \mathbb{I}, \quad i\mathbb{I}\otimes \sigma _k \qquad\qquad  j, k \in \{1,2,3\} .
   \]
   The remaining space $\mathfrak{m}$ is the span of the nine basis elements
   \[
   i\sigma _j \otimes \sigma _k   \qquad\qquad  j, k \in \{1,2,3\} .
   \]
       
    \end{Pro}
  
  We now can define the operators that we need. 
  
  \begin{Def}
  \begin{enumerate}

  \item Let $K$ denote complex conjugation with respect to the basis we have constructed for $\mathfrak{su}_4 $ as defined above. We have noted
  that $K^2 = I$.
  
       \item   Define the antiunitary operator \index{antiunitary}
   \[
  \Theta :  \CC^4  \longrightarrow \CC^4
    \]
  by 
  \[
  \Theta  =  (\sigma _2 \otimes \sigma _2)K.
  \]
   Check that $\Theta ^2 = I$.
     
   \item Define 
   \[
   Q = K(\sigma _1 \otimes i\sigma _2 ).
   \]
   Check once again that $Q^2 = -I$.

  \item Define $T$
     by 
  \[
  T(A) = \Theta A \Theta ^\dagger .
  \]
  Then $T$ is a map
  
  \[
T :  \mathfrak{su}_4 \longrightarrow \mathfrak{su}_4 
    \]
  and $T^2 = I$. 
  
  \item If instead of using $\Theta $ we used $Q$,  defining 
  \[
  T'(A) = QAQ^\dagger \; ,
  \]
  then 
    \[
T' :  \mathfrak{su}_4 \longrightarrow \mathfrak{su}_4 
    \]
  and $(T')^2 = -I$.

  \end{enumerate}
  \end{Def}

  \begin{Thm}  Let $T : \mathfrak{su}_4 \longrightarrow \mathfrak{su}_4 $ as described above.  Then the decomposition $\mathfrak {su}_4 \cong \mathfrak{o}_4 \oplus \mathfrak{m} $ 
  coincides with the eigenspace decomposition of $T$.     That is:
  \begin{enumerate}
  \item \qquad
  $T(V) = V  \quad\quad {\text{if}} \quad V \in \mathfrak{o}_4    $.
  
  \item\qquad
$   T(V) = -V  \quad {\text{if}} \quad V    \in  \mathfrak{m} .$
  \end{enumerate}
 \end{Thm}
 
  We omit the proof and refer the reader to \cite{OA}
   
   \newpage

   \subsection{\bf{Topology and Physics, by Nadav Orion}}
~

In quantum mechanics, the physical description of a system consists
of two notions: 

\begin{enumerate}
\item {\bf{The state of the system}}. All information about the system, as it is in a given moment in time, is represented by a vector in a complex inner product space (specifically a Hilbert space), the state space, which in turn contains all possible system states. These vector spaces may be of any complex dimension, from two (such a system is called ``qubit") to infinity.


\item {\bf{The Hamiltonian $H$}}. It is a Hermitian operator on that space of
physical states that determines how a system evolves in time. $H$
generates the unitary evolution operator 
\[
\exp{(-iHt/\hbar)}
\]
where $t$
is the time and $\hbar$ is the reduced Planck constant. For a system described by the vector $v$ at $t=0$, the
state of the system at time $t$ is given by
\[
\exp{(-iHt/\hbar)}v.
\]
 For a finite interval $t\in\left[0,T\right]$, the equation above 
 describes a path in the state space. 
\end{enumerate}
In a lab setting, physicists control certain aspects of the Hamiltonian
(e.g. by applying electromagnetic fields), while other aspects depend
on the specific fabrication and noise and thus change from one experiment
to another. For these reasons a set of Hamiltonians, rather than a
single one, is usually considered. In many cases this set corresponds
to a Lie algebra \index{Lie algebra}, Clifford algebra \index{Clifford algebra}, or a well-defined subset of one of these.

Different elements in the set of Hamiltonians may produce different
paths in the state space. These paths are at the basis of many physical phenomena such as
\begin{enumerate}
\item Aharonov-Bohm effect, 
\item Berry phase, 
\item quantum Hall effect, and 
\item entanglement
\end{enumerate}
to name a few. These are captured as topological properties of Hamiltonians
as maps (or a set of maps) from the space of physical parameters,
both controlled and uncontrolled, to the set of Hamiltonians. The {\emph{space}} of Hamiltonians is defined by applying homotopic equivalence and additional equivalence relations, depending upon the situation.
The topological properties of the evolution operator space (defined
using the set of Hamiltonians) can also be examined, and they are
related to those of the space of Hamiltonians. This relation is given by Bott Periodicity when the state space dimension $\rightarrow \infty$, approximately true also for large enough dimensions. 

Let us consider the case of two qubits, discussed purely mathematically in ``A Timely Example'' in the previous section. For each qubit
the state space is $\mathbb{C}^{2}$, giving a state space of 
\[
\mathbb{C}^{2}\otimes\mathbb{C}^{2}\cong\mathbb{C}^{4}
\]
for the joint system (both spaces with the usual inner product). A
typical Hamiltonian for such a system contains both terms that affect
a single qubit ($\sigma_{j}\otimes I$ and $I\otimes\sigma_{k}$)
and terms that describe interactions between the two ($\sigma_{j}\otimes\sigma_{k}$),
such that the set of Hamiltonians is all of the Hermitian $4\times4$ matrices,
which is the Lie algebra \index{Lie algebra} $\mathfrak{u}_4$. The evolution operator space
is thus $U_4$, a well-known 
compact complex symmetric space.

In some physical cases not all terms are allowed, for example when
the system possesses time reversal symmetry, i.e. 
\[
T\left(H\right)\equiv\Theta H\Theta^{\dagger}=H.
\]
For two qubit systems the natural time reversal symmetry is 
\[
\Theta=\left(\sigma_{2}\otimes\sigma_{2}\right)K,
\]
thus 
\[
H\in\text{span}\left\{ \sigma_{j}\otimes\sigma_{k}\right\} _{j,k\in\left\{ 1,2,3\right\} } 
\]
which is (up to a factor of $i$) the space $\mathfrak{m}$ as
in Proposition~\ref{Prop-mathfrak-m}. The disallowed terms are the Lie algebra $\mathfrak{o}_4$.
To get the space of evolution operators one must then identify all
operators in the original $U_4$ space that differ by
the disallowed terms $O_4$, giving the homogeneous space
\[
 U_4/O_4,
 \]
 a known real symmetric space.

Entanglement is a valuable resource in physics, and two qubits
is the simplest case in which it can be observed. The Cartan decomposition
distinguishes between Hamiltonians that may create and destroy entanglement ($\mathfrak{m}$)
and those that cannot interact with it ($\mathfrak{o}_4$), and as a consequence also
distinguishes between paths where entanglement changes and paths where
it does not. The evolution operators in $U_4/O_4$
can change the entanglement and those in $O_4$ cannot.
Furthermore, the different Hamiltonians in $\mathfrak{m}$ generate
paths in $U_4/O_4$ that differ by the amount
of entanglement they can create. 

As mentioned before, the space of Hamiltonians is obtained by homotopic equivalence, among others, following that the topological properties of the Hamiltonian do not depend on the size of its eigenvalues. One may then ``flatten" $H$ by
deforming its eigenvalues to be $+1$  or  $-1$  (or $0$ in extremely special cases). 
This leads to the associated evolution operator paths mentioned 
being loops.\footnote {The  loop starts at $t = 0$. To define the number for the Hamiltonian, we stop the first time the loop returns to its starting place.}  These loops then are maps 
\[
S^1 \longrightarrow  U_4/O_4
\]
and their homotopy classes lie in the group $\pi_1(U_4/O_4)$, which is isomorphic to $\mathbb{Z}$. 
Thus to each time evolution is associated an integer, which is a topological invariant of that path. By homotopy theory, if a single loop of a given Hamiltonian is given the invariant 1, performing the loop again will give 2. We can thus define the topological number of the Hamiltonian by the number given by a single loop.

 \newpage
     \section{\bf{Quantum Symmetry}}
     \label{Section:QuantumSymmetry}
     
    
    Suppose that $G$ is a compact Lie group \index{Lie group} and $K$ is a Lie subgroup of $G$.  Then the Lie algebra \index{Lie algebra}  of $G$ 
      decomposes as real vector spaces 
           \[   
       Lie(G) \cong Lie (K) \oplus \mathfrak m  
       \]
as noted previously, and the decomposition preserves the Lie algebra structure in $Lie(K)$.  However, $\mathfrak m$ is not a Lie algebra as we saw in the Timely Example. 
       
     Recall that $G$ is in fact a compact smooth manifold with a closed submanifold $K$ and homogeneous space $G/K$. 
       Topologically, $Lie(G)$ corresponds to the tangent space $T_oG$ of $G$ at the origin, $Lie(K)$ corresponds to the subspace $T_oK$
         that represents the tangent space of $K$ at the origin.  Taking quotients, we see that $\mathfrak m$ corresponds to the tangent 
       space of $G/K$ at the point on the smooth manifold $G/K$ that is the image of the origin of $G$.  Thus symmetries of $\mathfrak m$ correspond to 
       symmetries of the tangent space of $G/K$ at the origin.             
    
   Next we provide a very simplified introduction to the $T,\, C,\, S$ symmetries so beloved by followers of the Ten-Fold Way.
   The assignment of these symmetries to the ten Cartan classes dates to 1996 and is attributed to Altland and Zirnbauer (cf. \cite{AZ}, \cite{Z}).\footnote{We 
   note that the whole TCS assignment is considered somewhat tentative by some physicists- Emil Prodan alerted us to this matter and has 
   an alternate approach which is too involved to be introduced in this work.}
   
  
  We start with a quantum system with some real Hamiltonian $H$, as in the previous section, and seek topological invariants. Generally the Hamiltonian is a map from some compact space X (usually low dimensional $d\leq 4$) to a set of matrices defined as in the previous section. Topological invariants are, roughly, quantities that pay no attention to gentle continuous perturbations of the Hamiltonian. For this to be meaningful, the Hamiltonian must have energy gaps, since if not then every Hamiltonian can be deformed to any other. An energy gap is closed when two or more eigenvalues of a matrix corresponding to the same point in X coincide. The topological invariants of a system are given by the homotopy class $[X,R_H]$, where $R_H$ is the space of Hamiltonians, usually a symmetric space. In many physical cases the matrices turn out to be very large and $R_H$ is accordingly is taken to be the stable space.
  
  We seek the simplest examples. If the Hamiltonian is invariant under a unitary operator then it may be block-diagonalized with 
  respect to the operator into smaller blocks. We do this until there are no further unitary operators, so that it is irreducible. 
  
  There are, however, other possibilities, the first of which being {\emph{time-reversal symmetry}} $T$.   Possessing this symmetry means that the system would look the same, even at the atomic 
  level, if time were reversed. \index{time reversal}
  So in general, as we saw in the previous section,  we can write time-reversal symmetry as 
  $T$ with $T^2 = I$ or $T^2 = -I $ and the custom is to write $0$ if it is not present.       The Timely Example was indeed an example of time-reversal 
  symmetry with $T^2 = I$. 
  
  The second sort of symmetry that might be present is 
 {\emph{charge conjugation symmetry}}, also known as {\emph{particle-hole symmetry}}, 
 \index{particle-hole symmetry}
 where we, at least symbolically, switch particles and holes. This is 
 modeled by an antiunitary \index{antiunitary} operator $C$ with $C^2 = 1$ or $C^2 = -1$ and we write  $0$ if it is not present.

A system with both time reversal and charge-conjugation symmetry will also be symmetric with respect to $S=CT$. It is also possible to have a system that has neither time reversal nor charge conjugation symmetry, but is nonetheless symmetrical under the combination $S=CT$. 

Finally, we complete the Homotopy Road Map. 
   \newpage

   \begin{Table}\label{roadcomplete}
 \centerline{\bf{The Complete Homotopy Road Map }} \index{Homotopy Road Map} 

  \begin{center}
   \begin{tabular}{   r | r  |   | c  |  p{.15in} |  p{.15in} | l  |  c   | c  | c |} \hline\hline
  Label  & Name &{\bf{T}} & {\bf{C}} & {\bf{ S}}\,\,      &$\hspace{.5cm} \text{Big Ten} $ \hspace{.5cm} & $\cR _j $  & $\pi _0(\cR_j)$ & $\pi _1(\cR_j)$  \\    \hline \hline
   A & unitary & 0 & 0 & 0 &      $\Omega^0 \uu \, \, \, \cong \, \, \,  \uu $  & $\cR_1\pc  $ & $0$  & $\ZZ$ \\ \hline
  AIII & ch. unitary & 0  &  0 &1 & $\Omega^0 \uu \, \, \, \cong \, \, \, \bbu \times \ZZ $ & $\cR_0\pc $ &$ \ZZ$ & $0$  \\ \hline \hline 
    D  & BdG  & 0  &  +1 & 0   & $\Omega^0 \oo \, \, \,  \cong \, \, \,   \oo  $ &         $\cR _1$ & $ \ZZ _2$ & $\ZZ_2$ \\ \hline 
  DIII & BdG & -1 &  +1 & 1    & $\Omega^1  \oo\, \, \,  \cong \, \, \,   \oo / \uu $      &     $\cR _2$      &$ \ZZ_2 $ & $0$ \\ \hline
   AII & symplectic & -1  &  0 & 0  & $\Omega^2 \oo  \, \, \,  \cong \, \, \,   \uu / \pp $    &     $\cR _3$        &$ 0$ & $\ZZ $ \\ \hline
    CII & ch. symplectic & -1  &  -1 & 1  & $\Omega^3 \oo \, \, \,  \cong \, \, \,   B\pp \times \ZZ $          &    $ \cR _4$ & $ \ZZ$ & $0$ \\ \hline
    C & BdG & 0  &  -1 & 0   & $\Omega^4 \oo \, \, \,  \cong \, \, \,   \pp $   &     $\cR _5 $&$0$ & $ 0$ \\ \hline
  CI  & BdG & +1  &  -1 & 1  & $\Omega^5 \oo \, \, \,  \cong \, \, \,   \pp / \uu  $   &    $ \cR _6$ &$0$ & $0$ \\ \hline
    AI  & orthogonal  & +1  &  0 & 0   & $\Omega^6 \oo \, \, \,  \cong \, \, \,   \uu / \oo  $   &     $\cR _7 $&$0 $ & $\ZZ$ \\ \hline
   BD1 &chiral orth & +1  &  +1 & 1     & $\Omega^7 \oo  \, \, \,  \cong \, \, \,   \bbo \times \ZZ $   &    $ \cR _0$ & $ \ZZ$ & $\ZZ_2$ \\ \hline \hline
\end{tabular}
\end{center}
 \end{Table}

Each of the ten rows of this chart corresponds to one of the ten types of compact symmetric spaces \index{symmetric space}
and hence to the Ten-Fold Way. 
 We have discussed most of the columns previously; here is a reminder.

\begin{enumerate} \index{symmetric space}
\item {\bf{Label}} gives the \'E. Cartan indicator for this class of symmetric spaces.
\item {\bf{Name}} gives the Cartan name for this class of symmetric spaces. 
\item{$\bf{T, C, S}$} indicate the types of symmetry present on this class of symmetric spaces. For example, the ``orthogonal" row,  second row from the bottom, 
includes the Timely Example, in which time-reversal symmetry $T$ was present with $T^2 = +1$ and $C$ and $S$ were not present.
\item {\bf{Big Ten}} This shows the stable version of this class of spaces.  For example, look at the ``orthogonal" row, second from the bottom. The space $\uu/\oo$
is the union of symmetric spaces $U_n/O_n$.  The space $\uu/\oo$ came to our attention first via the isomorphism 
\[
\uu/\oo \cong \Omega ^6\oo 
\]
as part of our exposition of Bott Periodicity.  We named this space $\cR _7$ following physics custom.   It arose a second time via the \lancaster perspective 
in the context of Clifford algebras  \index{Clifford algebra|} .  It comes up a third time via $U_4/O_4$ in the Timely Example of the previous section. 
 \item ${\bf{ \pi_0(\cR _j)}}$ shows the lowest homotopy group of one of the Big Ten spaces.
 In the case of  $\uu/\oo$
  the group is $0$, indicating that the space is path-connected. 
 Also recall that $ \Omega \cR_j \cong \cR_{j+1}$ (where the index is considered modulo 8) so that 
         \[  \pi_1(\cR_j) = [S^1, \cR_j] = [S^0, \Omega \cR_j] = [S^0, \cR_{j+1} ] = \pi_0(\cR_{j+1})  \; . \]
Similarly, in the complex case,
 	\[  \pi_1(\cR\pc _j) = \pi_0(\cR\pc_{j+1})  \; .\]
where the index is considered modulo 2.
Compare to Table~\ref{road1}. 
\end{enumerate}

       This is laid out in Stone  et al \cite{Stone}, \S 3 in very concrete detail.  They demonstrate there that each of the  symmetric spaces $\cR_j$  listed above has certain symmetries associated to it, or, more 
   precisely, associated to $\mathfrak{m}_j$.  
    Their Table 4 lists which symmetry occurs in each of 
   the eight real cases.        These symmetries are vital for the applications in physics.

Since
\[
\widetilde{KO}^{-n}(X) \cong [X, \cR _n]
\]
and 
\[
\widetilde{K}^{-n}(X) \cong [X, \cR\pc _n] ,
\]
we see at once that $KO$-theory and $K$-theory are intimately related to this story.  The homotopy properties of the $\cR _n$ are detected by $KO$ and hence connected to 
real vector bundles.  This allows topological properties of the various symmetries to be detected and utilized. Much current physics research in topological insulators 
and related matters uses $KO^*$, and this was one of the main motivations for writing this {\BS}.

It is important to remember that the $\cR_j$ are unions of compact symmetric spaces, as this structure helps us very much to understand the homotopy side of the story. For example, 
\[
\cR_2 = \oo / \uu =  \lim_{t \to \infty} O_{16t}/U_{8t}
\]
and the map 
\[
\pi _j( O_{16t}/U_{8t} ) \to  \pi _j(\oo /\uu ) = \pi _j(\cR_2 )
\]
is an isomorphism for $j << t $ by cellular approximation.  Thus, for example,  if $X$ is some compact space then 
\[
\widetilde{KO}^{-2}(X) \cong [X, \cR_2] = [X, O_{16t}/U_{8t} ]
\]
for $\dim (X) << t$.
 
For many of the applications we need $\pi_j$ of the finite approximations (e.g. $\pi _j(O_4/U_4)$ for $j = 0,1,2$). 
It is important to remember that
  these may not coincide with 
$\pi_j (\oo / \uu )$ so one has to be careful.  See Section~\ref{Section:LowHomotopyGroups} for a full discussion of this matter.
 
\newpage

\bibliographystyle{amsplain}
\bibliography{Kthy}
    
 \printindex
 \printindex[notation] 
 
  \end{document}